\documentclass[12pt]{book}

 \usepackage{caption}

   \usepackage{cancel}
   \usepackage{color}
\usepackage{xcolor}
\newcommand{\zProof}{{\bf\underbar{Proof}.}\ } 
 
\newcommand{\ZREF}[1]{
\smallskip\noindent\fbox{\fbox{
\parbox{5in}{#1}
}}
\medskip
}

\newcommand{\ztc}{\ \mbox{such that}\ }
 \newcommand{\zdia}{~~\rule{1mm}{2mm}\par\medskip}  
  
\newcommand{\zdiaform}{\mbox{~~\zdia}}  

\newtheorem{Theorem}{Theorem}  
\newtheorem{Corollary}[Theorem]{Corollary}  
\newtheorem{Lemma}[Theorem]{Lemma} 

\newtheorem{Remark}[Theorem]{Remark}
\newtheorem{Definition}[Theorem]{Definition} 
\newtheorem{Example}[Theorem]{Example} 
\newtheorem{Assumption}[Theorem]{Assumption}

 \newcommand{\zb}[1]{{\bf #1}}

  \def\zHG0{H_{\Gamma_0}^1(\ZOMq)}

 \newcommand{\Dom}{{\rm dom}}

 \newcommand{\carfunz}{\mathbbm{1}}
\newcommand{\heaviside}{\mathbbm{1}}

\newcommand{\charfun}{\mathbbm{1}}

\newcommand{\defin}[1]{%
{\sc{#1}\/}
}

\newcommand{\ZSUno}{{\sum_{n=1}^{+\ZIN}}}

 \renewcommand{\epsilon}{\varepsilon}
 \renewcommand{\rho}{\varrho}

\newcommand{\zzn}{\mathbb{N}}
\newcommand{\ZEP}{\varepsilon}

\newcommand{\zg}{\gamma}

\newcommand{\ZOM}{\omega}
\newcommand{\zaa}{\alpha}

\newcommand{\ZDE}{\delta}

\newcommand{\zthe}{\theta}

\newcommand{\ZLA}{\label}
\newcommand{\ZIN}{\infty}

\newcommand{\zzr}{\mathbb{R}}
\newcommand{\zzq}{\mathbb{Q}}

 \newcommand{\ZD}{\;\mbox{\rm d}}
\newcommand{\zl}{\lambda}
\newcommand{\ZSI}{\sigma}
\newcommand{\ZOMq}{\Omega}
  \newcommand{\Rpalla}{\overset{\circ}{ R}}

   \usepackage{url}
 \usepackage{hyperref}
  \usepackage{bm}
  \usepackage{bbm}
\newcommand{\limn}{\lim _{n\to+\ZIN}}
\newcommand{\limm}{\lim _{m\to+\ZIN}}

\newcommand{\Q}{\pmb Q}

\newcommand{\QP}{{{\pmb Q}_{\pmb P}}  }
 
\newcommand{\QN}{{{\pmb Q}_{\pmb N}}  }

  \usepackage{amssymb,amsmath,exscale}
  \DeclareMathOperator{\sinc}{sinc}
\usepackage{type1cm}      

\usepackage{makeidx}         
\makeindex
\usepackage{graphicx}        
\usepackage{multicol}        
\usepackage[bottom]{footmisc}
\usepackage{newtxtext}       %
 \usepackage{newtxmath}       

\title{\vskip -5cm{\normalsize\bf Luciano Pandolfi}\\
\vskip 2cm
\noindent
\bf \Large The Lebesgue Integral\\
\noindent  {\normalsize Via The Tonelli Method}\\
\vskip1cm \noindent
 }
\begin{document}
 
  \frontmatter   
  \pagestyle{empty}
  \pagenumbering{roman}
\vskip -5cm

 \maketitle
 
 \pagestyle{myheadings}
  \markboth{{PREFACE}}{{PREFACE}}

 \chapter*{Preface}
  \addcontentsline{toc}{chapter}{Preface}

The theory of Lebesgue integration was developed by H.-L. Lebesgue in his doctoral dissertation       
     published  in 1902  in the italian Journal \emph{Annali di Matematica Pura ed Applicata} 
(see~\cite{Lebesgue1902Annali}). The theory of Lebesgue integration was soon fully appreciated by the leading mathematicians in Italy,  also thanks to the activity of Vitali, and important results      soon proved: Fubini Theorem was proved in 1907 (in~\cite{Fubini07TeoRiduz}). Leonida Tonelli first used Lebesgue integral in 1908, in~\cite{Tonelli08RETTIFICAZIONE}, while studying  rectifiable curves\footnote{the study of the length of a curve and of the area of a surface using Lebesgue integration was initiated by Lebesgue himself in his thesis.}.

 As soon as the importance of Lebesgue integration was recognized,      different approaches   were proposed  both to put the   theory in a more general framework, as in~\cite{Daniell17INTEGRAL}, and also in order to speed up the presentation of the crucial ideas.   In fact, for didactic reasons, several authors tried to bypass or reduce to a minimum the preliminary study of the measure of sets and to lead students to appreciate the key ideas of Lebesgue integration in the shortest possible time. 
 Among these approaches, the one which is likely the  most well  known   is   due to   F. Riesz. This approach is the one used    in~\cite{RieszNAGYfunzANAL} and recently simplified in~\cite{KomornikSIMPLIFIEDrieszInteg}. Instead, the Daniell's approach, introduced in~\cite{Daniell17INTEGRAL},  intends to be more general and abstract.  

  Leonida   Tonelli  devised\footnote{S. Cinquini, one of Tonelli's students, asserts that    this approach was devised to quickly introduce the Lebesgue integral in a talk given at a conference (see~\cite{CinquiniCommemorTONELLI}).} an efficient approach to Lebesgue integration,
  based on the notion of ``quasicontinuous functions'', which requires only a basic knowledge of calculus together with a certain level of mathematical ingenuity. 
  
  Quasicontinuous functions are precisely the functions which are measurable according to Lebesgue, but defined in terms of elementary   continuity properties.
  
It seems to me that this approach proposed by Tonelli, and which is virtually unknown to young people, has its merits since  most of the methods used to construct the Lebesgue integral  hide the real reasons for the introduction of a new kind of 
integral. We cite from~\cite[p.~29]{RieszNAGYfunzANAL}: ``The reason for such a change [the shift from Riemann to Lebesgue integration] and, with them, the usefulness and beauty of the Lebesgue theory, will be seen in the course of the following chapters; there is no point in speaking of them in advance". Instead,   the  reason  is clearly displayed in the Tonelli approach, since it constitutes the foundation over which the integral is constructed.

   Lebesgue integral  is introduced by exploiting the method of Tonelli   in~\cite{CinquiniVARrel60} (for functions of one variable.   \cite[Chap.~VI]{CinquiniVARrel60} concisely defines the integral for   functions of two variables and  proves    Fubini theorem), 
and in part in the book~\cite{PiconeVIOLA52}. As stated by the authors, the first version of this book was prepeared for the lessons given by Picone when he gave the course of Tonelli, after Tonelli's death. 
Few books, as~\cite{HobsonTheoryFunctionsVol2,SansoneMerli}, present a sketch of Tonelli ideas.

For this reason we give here a quite complete account of Tonelli method.

We divide the presentation in three parts: in the first part we study integration  in the simpler case of the functions of one variable. The fundamental notion    that the student has to master is that of ``quasicontinuous function''.  Once this is done, Lebesgue integral appears as a straightforward extension of Riemann integral thanks to an ``exchange of limits and integrals'' which is precisely the goal for which the new integral has been introduced\footnote{in a sense, Tonelli construction can be see as an extension of Riemann integral ``via completion of a space'' but this abstract approach attach a number to equivalent classes of function while the concrete approach of Tonelli gives precisely the class of the functions which are Lebesgue integrable and their Lebesgue integral.}.

 This chapter ends with the statement of the theorems of the exchange of the limits and the integral, with a sketch of a proof. The details of the proofs     are   in Chap.~\ref{Cha1bis:INTElebeUNAvariaBIS}.
 
The first chapter can be used as an efficient introduction to  the Lebesgue integration even in courses for mathematically motivated engineering students, for example as a preliminary to a basic course on Hilbert spaces.

The second part extends the arguments of the first part to functions of several variables  (in the Chapters~\ref{Chap:2quasicontPIUvar} and~\ref{ch3:LIMITintegral}). Chap.~\ref{ch4CHAPFUBINI} is devoted to Fubini Theorem on the reduction of multiple integrals.

For the sake of completeness, the  third part, containing the sole Chap.~\ref{Chap:4Measure},  shows the well known fact that once the integral has been defined then Lebesgue and Borel measurable sets can be defined and the properties of the measure can be obtained from those of the integral.

Finally we repeat that the usual presentation of Lebesgue integration  proceeds from the study of the measure and measurable sets to the theory of measurable functions and then to the theory of the integral. This way the students have to study first abstract measure theory which has its independent importance, for example for the application to probability. Instead, in the Tonelli approach\footnote{as in the approaches proposed by
B. Levi, M. Picone and 
 F. Riesz.} measure theory is derived as a byproduct of the study of the integral and the abstract measure theory remains in the shadow. This fact, which can be seen   as an important limitation now, was instead the goal of Tonelli and others. For example, Tonelli states in the introduction of his 1924 paper~\cite{Tonelli24ANNALIsullaNOZinteg}: ``this paper is an attempt to put on a simple, let us say elementary, and so more acceptable,  basis the theory of Lebesgue integration by removing measure theory at all''.
Similar approaches to a direct introduction of the integrals done in the same period state explicitly the same goals.   Beppo Levi introduced his own approach  to Lebesgue integration in the same year 1924 and in the same journal as Tonelli did and for similar reasons (see~\cite{LeviBEPPO24ANNALIsullaNOZinteg}): ``It is a fact that in front of the importance of the new theory [i.e. the Lebesgue theory] there are obvious didactic and logical difficulties\dots due also to the need of a preliminary study of the theory of the measure   of sets''. The introduction of F. Riesz and 
B. Sz-Nagy book~\cite{RieszNAGYfunzANAL} states ``The two parts [of the book] form an organic unity centered about the concept of linear operator. This concept is reflected in the method by which we have constructed the Lebesgue integral; this method, which seems to us to be simpler and clearer then that based on the theory of measure\dots''

Finally we mention that Picone too proposed his approach to the Lebesgue integration for functions of Baire classes (see for example~\cite{PiconeAPPmatSUP}), by a repeated use of limiting processes.
 
 We conclude this introduction with a warning. The development of Lebesgue integration was stimulated by two main difficulties encountered with Riemann integration:
 \begin{enumerate}
 \item the fact that the pointwise limit of a sequence of continuous functions may not be Riemann integrable, a fact soon realized after that a rigorous definition of the integral (and also of the concept of  function) had been given. A consequence is that limits and Riemann integrals cannot be exchanged in the generality needed to study for example Fourier series (and  Riemann integral cannot be used to define Hilbert spaces).
 \item the existence of derivative functions which are not Riemann integrable, a fact
 discovered by V. Volterra in 1881, see~\cite{VolterraPrinCalINTbatt1881}.
 \end{enumerate}
In our exposition we concentrate on the first problem,   the problem of exchanging limits and integrals. The second problem, of recovering a function from its derivative, will not be considered here. Readers interested in the relation of the integral and the derivative can see for example~\cite{Benedetto76BOOK} or, for an approach based on Daniell ideas,~\cite{SHIlovUnifINTEgral}.
 
 People interested in the origin and the history of the Lebesgue integration can look at the book~\cite{HawkinsHISTORYlebe} and to the expositions~\cite{LebesgueEXPOsHISres,LebesgueSTORIArevMor} made by Lebesgue himself.
 \mainmatter
  \pagestyle{headings}

 {\footnotesize
  \tableofcontents
  \normalsize}
 \part{\ZLA{PART1OneVar}Functions of One Variable}
 
 \chapter{\ZLA{Ch1:INTElebeUNAvaria}
The Lebesgue Integral for Functions of One Variable} 
\chaptermark{Lebesgue Integral of Functions on $\zzr$}
This chapter is intended for general students af science and engineering with possibly just calculus courses but a sufficient motivation. 
In particular, we assume familiarity with Riemann integration whose key points are recalled in the Appendix~\ref{APPEch1SKETCH}.

This chapter can be used as the introductory chapter  of a course on Hilbert spaces. For this reason first we describe the limitations of the Riemann integral and then we use   Tonelli method as an efficient tool to introduce the key ideas   of the  Lebesgue integration. 

At the end of this chapter we state and shortly illustrate the theorems on the exchange of limits and integrals. The details  of the proofs are in Chap.~\ref{Cha1bis:INTElebeUNAvariaBIS}.

Our goal here is the presentation of key ideas and for this reason we confine ourselves to consider functions of one variable. The general case of functions of several variables is in Part~\ref{PART2MoreVar}.

\section{\ZLA{sect:Ch1LimitatRiemaINTEG}The Limitations of the Riemann Integral}

When presenting the definition of the Riemann integral, usually  instructors show the existence of functions which are not Riemann integrable. The standard example is the \defin{{Diriclet function}}\index{Dirichlet function}\index{function!Dirichlet} on $[0,1]$:
\begin{equation}
\ZLA{eq:Ch1:DiricFUNC}
d(x)=\left\{
\begin{array}{lll}
1&{\rm if}& x\in\zzq\\
0&{\rm if}& x\in\zzr\setminus\zzq\,.
\end{array}
\right.
\end{equation}
It is easily seen that this function is not Riemann integrable and it is often asserted that examples like this one   prompt for the definition of a more general kind of integral. 

We may
 ask ourselves whether there is any reason  to integrate  such kind of pathological functions. At first glance it seems that there is no reason at all. 
But, let us look at the problem from a different point of view.

Even the first elements of calculus make an essential use of the notion of limits and continuity. A function $f $ is {\sf continuous at $x_0\in\Dom \, f$}\index{continuity!at a point}\index{function!continuous} when
for every sequence $\{x_n\}$ with $x_n\in\Dom\, f$ and such that $x_n\to x_0$
we have
\[
\lim   f(x_n)=f(x_0)=f\left ( \lim  x_n\right ) \,.
\]
\emph{Continuity is the property that the limit and the function can be exchanged.}

Now we observe that Riemann integral is a function which associates a number to a set of functions. So, we introduce:
\begin{Definition}
{\rm
\begin{enumerate}
\item
Let $\mathcal D$ be a set of functions. We call {\sc functional}\index{functional} a transformation which associates a number to any element of $\mathcal D$.
\item Let it be possible to define the limit of sequences in $\mathcal D$ and let $\mathcal F$ be a functional defined on $\mathcal D$. We say that $\mathcal F$ is continuous at $f_0\in \mathcal{D}$ when the following holds for every sequence $\{f_n\}$ in $\mathcal D$, $f_n\to f_0$:
\[
\lim \mathcal{F}(f_n)=\mathcal{F}\left (\lim f_n\right )\,.
\]

\end{enumerate}
}
\end{Definition}

Riemann integral is a functional on $\mathcal{D}=C([h,k])$ (for every bounded interval $[h,k]$) and it is possible to define the limits of sequences in $C([h,k])$ as follows: {\sc the sequence $\{ f_n\}$ converges to $f $    uniformly on $[h,k]$}\index{uniform convergence}\index{convergence!uniform} when
the following holds:
\begin{equation}\ZLA{Ch1:DefUniConvPrim}
 \forall\ZEP>0\ \exists  N_\ZEP  \ \mbox{such that}\   n>N_\ZEP\ \implies\  
-\ZEP<f_n(x)-f(x)<\ZEP\quad \forall x\in [h,k]\,.
\end{equation}
I.e. we require that for $n>N_\ZEP$ we have
\[
-\ZEP<\inf _{[h,k]}(f_n-f)\leq \sup_{[h,k]}(f_n-f)<\ZEP\,.
\]
So we have the following result:
\begin{Theorem}\ZLA{Ch1TeoDefoCONVunif}
The sequence $\{f_n\}$ is uniformly convergent to $f$ if and only if we have
\[
\lim_{n\to +\ZIN} \inf _{[h,k]} (f_n-f)=0\quad \mbox{and}\quad \lim _{n\to +\ZIN} \sup _{[h,k]} (f_n-f)=0\,.
\]
\end{Theorem}

\begin{minipage}{2in}
 From a graphical point of view, uniform convergence is the property that for large $n$ the graphs of $f_n $ stay in a ``tube'' of width $\ZEP$ around that of $f $, as in the figure on the right.
\end{minipage}
\begin{minipage}{2in}
 
 \includegraphics[width=7cm,scale=.65]{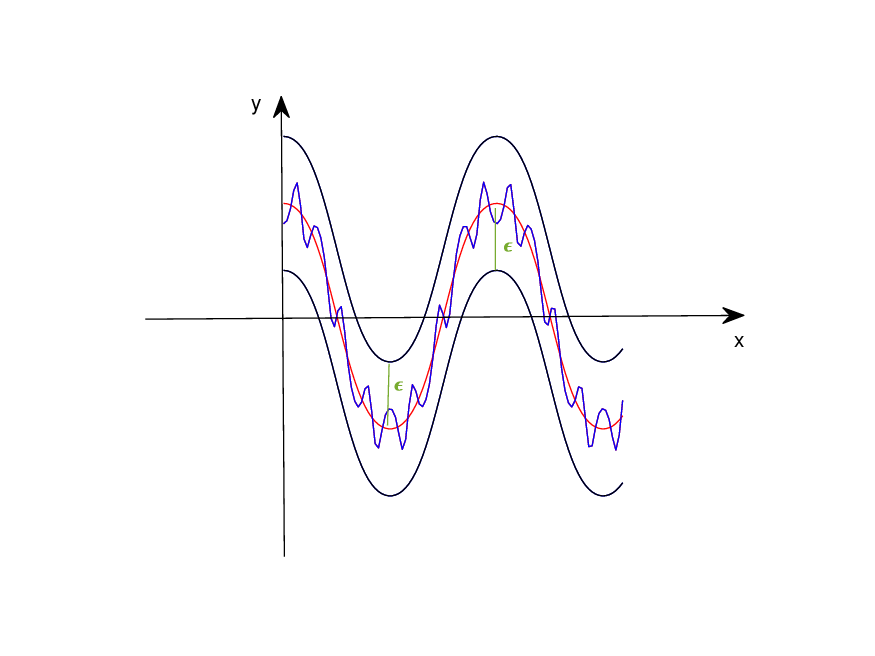}

\label{fig:ch1uniconve}       
 
  \end{minipage}

It turns out that Riemann integral as a functional on $C([h,k])$ is continuous:
\begin{Theorem}\ZLA{Teo:Ch1Sec1ScambLimRiem}
Let $[h,k]$ be a bounded interval and let $\{f_n\}$ be a uniformly convergence sequence in $C([h,k])$. Then:
\begin{equation}\ZLA{eq:ch1ContiRiemm}
\limn \int_h^k{ f_n(x)}\ZD x=\int_h^k\left (\limn f_n(x)\right )\ZD x\,.
\end{equation}
    \end{Theorem}

Now we observe:
\begin{itemize}
\item Riemann integral is defined on a set of functions which is larger then $C([h,k])$. In particular, any piecewise continuous function is Riemann integrable.
\item Uniform convergence is a very strong property.  
 Much too strong for most of the applications of mathematics. As an example, let us consider the following important equality:
 \begin{equation}\ZLA{eq:ch1SerFdiCHI}
\lim _{n\to+\ZIN} \left[\sum _{k=1}^n\frac{1}{k} \sin kx\right]=\chi(x) 
\end{equation}
where $\chi(x)$ is the extension of period $2\pi$ of the odd extension to $(-\pi,0)$  of  
 \smallskip
 
\hskip -.6cm\begin{minipage}{2.5in}
 
\[
[(\pi-x)/2]\quad 0<x<\pi 
\] 
(see the graph on the right).
The sequence of the functions
\[
f_n(x)= \sum _{k=1}^n\frac{1}{k} \sin kx 
\]
 is a sequence of continuous functions while $\chi$ has jumps. 
\end{minipage}
\begin{minipage}{2in}
 \includegraphics[width=7cm,scale=.65]{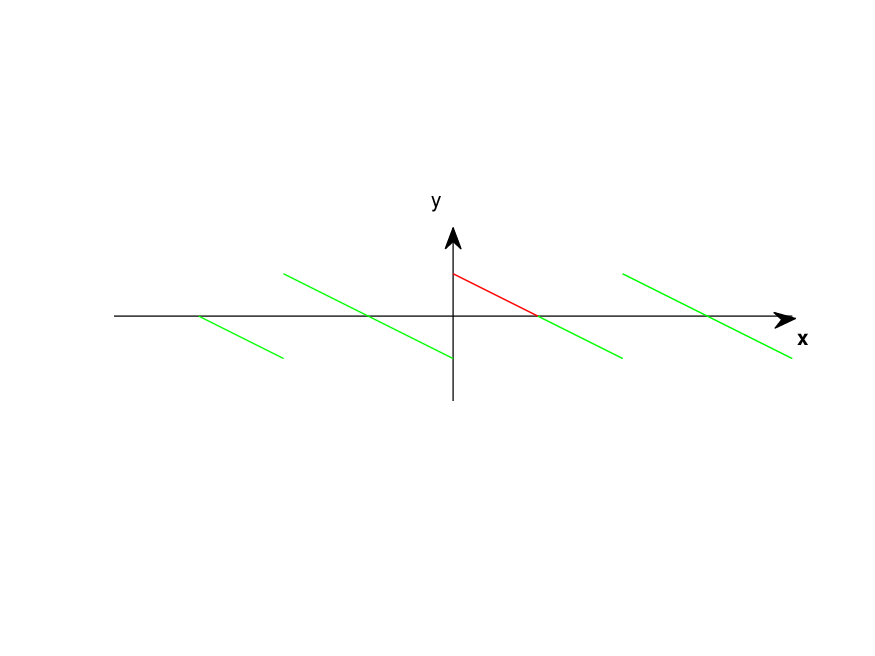}
  \end{minipage}
\end{itemize}

It is known that the uniform limit of a sequence of continuous functions is continuous and so the limit in~(\ref{eq:ch1SerFdiCHI}) cannot be uniform in an 
interval which contains jumps of $\chi$. \emph{On these intervals Theorem~\ref{Teo:Ch1Sec1ScambLimRiem} cannot be applied.}

The equality~(\ref{eq:ch1SerFdiCHI}) is a simple example of a Fourier series and it is seen in every elementary introduction to the Fourier series that a key point in the justification of the equality~(\ref{eq:ch1SerFdiCHI}) is the exchange of the series (i.e. of a limit) with suitable integrals.
So, in order to justify~(\ref{eq:ch1SerFdiCHI}) we must use an analogous of
Theorem~\ref{Teo:Ch1Sec1ScambLimRiem}   when
  the limit exists    but the convergence is not uniform.
Unfortunately,
the following example shows  that if the convergence is not uniform in general the limit and the Riemann integral cannot be exchanged.

\begin{Example}\ZLA{exCH1PrimoNoScLimINT}
{\rm 
The sole condition that each $f_n $ is Riemann integrable and that $\{f_n(x)\}$ converges to $f(x)$ for every $x\in [h,k]$ does not imply the equality~(\ref{eq:ch1ContiRiemm}). This can be seen as follows. We recall the following property: \emph{a bounded function which is equal $0$ a part finitely many points is Riemann integrable and its integral is $0$.}

It is known that the rationals in $[0,1]$ constitute  a  numerable set, i.e. they can be arranged to be the image of a sequence $\{q_n\}$.

For every $n$ we define 
\begin{equation}\ZLA{Ch1eq:sequeApproxDirichl}
f_n(x)=\left\{\begin{array}{lll}
1 &{\rm if}&   \mbox{$x=q_k$ with $k\leq n$}       \\
0 & &\mbox{otherwise}\,.
\end{array}\right.
\end{equation}
Each function $f_n $ differs from $ 0$ in $n$ points and so \emph{$f_n$ is integrable with   integral equal $0$.} The equality~(\ref{eq:ch1ContiRiemm}) does not hold since for every $x$ we have $f_n(x)\to d(x)$, the Dirichlet function,   and $d(x)$ is not Riemann integrable.\zdia

}
\end{Example}

This example however is not entirely satisfactory since we can wonder whether there is any concrete interest in the functions $f_n $ just constructed. So, let us see a second example.
\begin{Example}\ZLA{RseCH1ESEpeano}{\rm
We consider functions defined on $[0,1]$ and a limiting process in two steps.
Let us fix a natural number $m$ and let us consider the sequence of the continuous functions
\[
n\mapsto 
\left [\cos (m!)\pi x\right ]^{2n}\,.
\]
Then we define the function
\[
\phi_m(x)=\lim _{n\to+\ZIN}\left[\cos (m!)\pi x\right ]^{2n}\,.
\]

If $x\notin \zzq$  then $|\cos (m!)\pi x|<1$ and so
\[
x\notin \zzq\ \implies \ \phi_m(x)=0\,.
\]
Let instead $x\in \zzq$, $x=p/q$ with $p<q$. We have
\[
\begin{array}{l}
\cos (m!)\pi x= \cos  \dfrac{pm!}{q}\pi     =\pm1 \quad\mbox{if $q$ divides $pm!$}\\
|\cos (m!)\pi x|<1\quad \mbox{otherwise}\,.
\end{array}
\]
For a fixed value of $m$, the denominator $q$ divides $pm!$   only for finitely   rational number $p/q\in [0,1]$ and so

\[
\phi_m(x)=1 \ \mbox{for finitely many values of $x\in [0,1]$; otherwise $\phi_m(x)=0$}\,.
\]

Hence $\phi_m $ is Riemann integrable and its integral is equal zero.

Now we consider 
\[
\lim _{m\to+\ZIN} \phi_m(x)\,.
\]
 If $x\notin \zzq$ then the limit is zero since $\phi_m(x)=0$ for every $m$.

Let $x=p/q$. When $m$ is sufficient large, say $m>2q$, the number $(m!)p/q$ is an even number and $\phi_m(x)=1$.

It follows that $\lim_{m\to+\ZIN} \phi_m(p/q)=1$. Hence
\[
\lim _{m\to+\ZIN}\phi_m(x)=d(x)
\]
a non integrable function.

In conclusion, we are in the same situation as in the Example~\ref{exCH1PrimoNoScLimINT}: limit and integral cannot be exchanged.
 \emph{But now we have an example which  is significant  since it shows that the Dirichlet function can be obtained by computing limits of cosine functions; and sequences of sine and cosine functions are the basis of the Fourier analysis, a crucial tool in every application of  mathematics.}\zdia
}
\end{Example}
 
In fact, from the historical point of view it was Fourier analysis that stimulated the introduction of several variants of the Riemann integral in the course of the XIX century. This process then  culminated in 1902 with the introduction of the Lebesgue integral   in~\cite{Lebesgue1902Annali}.

\begin{Remark}\ZLA{Ch1RemaNelQualDEInteNulso}
{\rm
We repeat  again: the desired equality~(\ref{eq:ch1ContiRiemm}) is the reason for the construction of an integral which can integrate such patological functions like the Dirichlet function, and not the interest of such functions by themselves.

Note that Example~\ref{RseCH1ESEpeano} contains  also a second piece of information: if equality~(\ref{eq:ch1ContiRiemm}) has to hold then the Dirichlet function not only has to be integrable but its integral has to be zero.\zdia
}
\end{Remark}

Now we sum up our \emph{rather optimistic dream:}
we would like a new kind  of integral for which the following properties hold:
\begin{enumerate}
\item if $f(x)\equiv c$ on $[h,k]$ then it is integrable and  its integral is $c(k-h)$ i.e. the area of the rectangle identified by its graph;
\item additivity of  the integral: if $f(x)$ is defined on $A\cup B$  and it is integrable both on $A$ and on $B$ then  it is integrable on its domain $A\cup B$. Moreover we want:
\[
A\cap B=\emptyset\ \implies\ 
\int _{A\cup B} f(x)\ZD x=\int_A f(x)\ZD x+\int_B f(x)\ZD x\,.
\]
In concrete terms, if  $A=[h,\tilde h)$ and $B=[\tilde h,k]$ we want
\[
\int_h^k f(x)\ZD x=\int_h^{\tilde h} f(x)\ZD x+\int _{\tilde h}^k f(x)\ZD x\,.
\]
\item let  $\{f_n\}$ be  a sequence of integrable functions  on a set $A$ and let   $\lim f_n(x)=f(x)$ for every $x\in A$. Then we wish that
  $f $ be integrable too and 
  \begin{equation}
  \ZLA{eq:ch1ContiRiemmBIS}
  \lim _{n\to+\ZIN}\int_A f_n(x)\ZD x=\int_A\left ( 
   \lim _{n\to+\ZIN}  f_n(x)\right )\ZD x\,.
  \end{equation}
\end{enumerate}

\emph{The following examples show that these requirements are contradictory and cannot be achieved.}
\begin{Example}\ZLA{EseCH1EsempioNOlimOGNOcasI}
{\rm
Let $f_n $ be defined on $[0,1]$  as follows:
\[
f_n(x)=
\left\{\begin{array}{lll}
0&{\rm if}&  0\leq x\leq 1/n          \\
n&{\rm if}&   1/n<x    \leq 2/n     \\
0&{\rm if}&      2/n\leq x\leq 1\,.    
\end{array}\right.
\]
It is clear that
\[
\mbox{$ \lim f_n(x)=0=f(x) $  per ogni $x\in [0,1]$}\,.
\]
  The first and second requirements show  that the functions $f_n  $ and $f $ are integrable and give  the value of the integrals:
\[
  0=\int_0^1 f(x)\ZD x\neq \lim_{n\to+\ZIN} \underbrace{\int_0^1 f_n(x)\ZD x}_{=1\ \mbox{for every $n$}}=1\,.
\]

A similar example can be given for   integrals on unbounded domains. Let  the domain be $[0,+\ZIN)$. In this case we define
\[
f_n(x)=
\left\{\begin{array}{lll}
0&{\rm if}&  x\leq n          \\
1/n 
&{\rm if}&   n<x    \leq 2n   \\
0&{\rm if}&      x>2n\,.    
\end{array}\right.
\]
In this example $\limn f_n(x)=0$, even uniformly on the \emph{unbounded} set $[0,+\ZIN)$, but every $f_n $ has the integral equal to $1$.\zdia
}
\end{Example}

So, we reduce our goals and we require the properties listed in Table~\ref{TABLEch1RequestINTE}. It turns out that an integral with these properties exists, it is the Lebesgue integral, and it has a further bonus: the heavy boundedness assumption in the requirement~\ref{I3ZREFrequiABSTRinteg} can be weakened.

\begin{table}
\caption{\ZLA{TABLEch1RequestINTE}\bf Our request to the integral}
\smallskip
\ZREF{

\begin{enumerate}
\item\ZLA{I1ZREFrequiABSTRinteg} the function  $f(x)\equiv c$ on $[h,k]$ is integrable and  its integral is $c(k-h)$ i.e. the area of the rectangle identified by its graph;
\item\ZLA{I1bZREFrequiABSTRinteg}
 if $f$ and $g$ are defined on $[h,k]$ and integrable      then:
\begin{enumerate}
\item{\sc linearity:}\index{integral!linearity}\index{linearity of the integral} the function $\zaa f+\beta g$ is integrable for every 
real numbers 
$\zaa$ and $\beta$  and
\[
\int_h^k\left [\zaa f(x)+\beta g(x)\right ]\ZD x=\zaa \int_h^k f(x)\ZD x+\beta\int_h^k g(x)\ZD x\,.
 \]
 \item
{\sc monotonicity:}\index{integral!monotonicity}\index{monotonicity of the integral}  
  if $f(x)\geq g(x)$ for every $x\in [h,k]$ then 
 \[
\int_h^k f(x)\ZD x\geq \int_h^k g(x)\ZD x\,. 
 \]
\end{enumerate}
\item\ZLA{I2ZREFrequiABSTRinteg}  {\sc additivity:}\index{additivity of the integral}   if  $f $ is defined   and integrable both on $A$ and on $B$ then   it is integrable on   $A\cup B$ and  
\[
A\cap B=\emptyset\ \implies\ \int _{A\cup B} f(x)\ZD x=\int_A f(x)\ZD x+\int_B f(x)\ZD x\,.
\]
\item\ZLA{I3ZREFrequiABSTRinteg}  if   $\{f_n \}$ is is a sequence of functions and if:
 \begin{enumerate}
\item
the functions $f_n $ are defined   on the   set $A$ and  
 \[
\limn f_n(x)=f(x)  \ \mbox{ for every $x\in A$}\,;
\]
\item the set $A $ is  bounded; 
\item the sequence $\{f_n \}$ is bounded in the sense that there exists $M$ such that $|f_n(x)|<M$ for all $x\in A$ and every $n$;
 \item each $f_n $ is integrable on $A$;  
    \end{enumerate}
then    $f $ is integrable too and the equality~(\ref{eq:ch1ContiRiemm}) holds:
\[
\lim_{n\to+\ZIN} \int_A{ f_n(x)}\ZD x=\int_A{\left (\lim_{n\to+\ZIN} f_n(x)\right )}\ZD x=\int_A{  f (x)}\ZD x\,.
\]
 \end{enumerate}
 }
 \end{table}

  \section{\ZLA{sec:CH0REALline}Subsets of the Real Line and Continuity}
 
  The following definition has a crucial role:

\begin{Definition}\ZLA{CH1Parte1DefiMultiint}
A  {\sc multiinterval}\index{multiinterval} \ $\Delta$ is a \emph{finite or numerable sequence} of \emph{open  } intervals: $\Delta=\{(a_n,b_n)\}$. The intervals $(a_n,b_n)$ are the {\sc component intervals}\index{component!interval} of the multiinterval.  We associate to $\Delta$:
\[
\left\{\begin{array}{l}
\mbox{the set}\ \mathcal{I}_\Delta=\cup (a_n,b_n)\\
  \mbox{the number}\ 
L(\Delta)=\sum (b_n-a_n)\in (0,+\ZIN]\,.
\end{array}\right.
\]

The multiinterval is {\sc disjoint }\index{multiinterval!disjoint }\index{disjoint!multiinterval} when the component intervals are pairwise disjoint.
\end{Definition}
\begin{Remark}\ZLA{RemaCh1ImproperTERMS}{\rm 
We stress the following facts:
\begin{itemize}
 \item the term ``finite or numerable sequence'' is, strictly speaking, partly redundant and not strictly correct since a sequence (of intervals) is a map $n\mapsto (a_n,b_n)$ when the domain of the map is $\zzn$. We use this term since we consider also the case that the domain is finite, say $1\leq n\leq N$.
\item
The number $L(\Delta)$ cannot be zero since any open interval has a positive length and it  does not depend on the order in which the component intervals  $R_n$ are listed.
\item
We do not require that the \emph{component intervals} $(a_n,b_n)$ are disjoint.
Even more, we do not require that they are distinct intervals: the same interval can be listed several times. So, it can be $L(\Delta)=+\ZIN$ even if the set $\mathcal{I}_\Delta$ is bounded. Hence, in no way the number  $L(\Delta)$ can be considered as a ``measure'' of $\Delta$.

   \item
A bit pedantic observation is as follows. We defined $\Delta$ as a \emph{sequence} of intervals:
\[
  \mbox{$\Delta$ is the function $n\mapsto (a_n,b_n)$} 
\] 
but we noted that the number $L(\Delta)$ does not depend on the order in which the component intervals are listed. So, it might be tempting to define $\Delta$ as a \emph{family} of intervals. This is not quite correct since ``family'' is usually intended as a synonym of ``set'' and so the two families $\{(0,1)\,,\ (0,1)\}$ and $\{(0,1)\}$ are the same family, i.e. the same set, since they have the same element. But  they are different multiintervals, i.e. different sequences of intervals,  and 
\[
L\left (\{(0,1)\}\right )=1\,,\quad \mbox{while}\quad \mbox
L\left (\{(0,1)\,,\ (0,1)\}\right )=2. 
\]
\item In the definition of multiinterval we did not impose that the component intervals are bounded. In most of our applications we  use multiintervals $\Delta$ such that $\mathcal{I}_\Delta\subseteq [h,k]$ and so the component intervals are bounded. 
An exception is Theorem~\ref{teo:Ch2StrutturaAperti}.
\item the length of an interval does not depend on whether it contains the endpoints: the intervals $(a,b)$, $[a,b)$, $(a,b]$ and $[a,b]$ all have the same length $b-a$. So, we can associate the number $L$ also to multiintervals whose component intervals are not open. This will be done in Chap.~\ref{Cha1bis:INTElebeUNAvariaBIS}. \emph{ In the present chapter we assume that a 
multiinterval is composed of open intervals, as stated in Definition~\ref{CH1Parte1DefiMultiint}.}
 \end{itemize}
 }
 \end{Remark}
We  need the following observation:

\begin{Lemma}\ZLA{ch1:LemmaUnioMULTINT}
Let $\{\Delta_k\}$ be a sequence of multiintervals, $\Delta_k=\{I_{k,n}\}$.
There exists a multiinterval $\Delta$ whose component intervals are
the  intervals $I_{k,n}$ and such that
\[
L(\Delta)=\sum_{k=1}^{+\ZIN} L(\Delta_k)\,.
\]
\end{Lemma}
\zProof Let $p_k$ be the $k$--th prime number. The sequence $\Delta$ is the function $p_k^n\mapsto I_{k,n}$.\zdia

Now we define:
\begin{Definition}\ZLA{Ch1DefiNULLsetqO}
 The set $N\subseteq\zzr$ is a {\sc null set}\index{set!null!in $\zzr$}\index{null set!in $\zzr$} when for every $\ZEP>0$ there exists a multiinterval $\Delta $ such that
\[
L(\Delta)<\ZEP\,,\qquad N\subseteq \mathcal{I}_\Delta \,.
\]

 A function   whose points of discontinuity are a null set is   {\sc almost everywhere continuous (shortly, a.e. continuous).}\index{function!a.e. continuous}\index{continuity!a.e.}
\end{Definition}
\begin{Example}\ZLA{Exe:Ch0Qnullset}{\rm
The set $\zzq\cap (0,1)$ is a null set. In order to see this fact we recall that the set of the rational numbers  is  numerable: there exists a sequence $\{q_k\}_{1\leq k<+\ZIN}$ such that $q_k\neq q_j$ if $k\neq j$ and whose image is
 $\zzq\cap (0,1)$.
 
 We fix any $\ZEP>0$. To $q_n$ we associate the interval $I_n=(q_n-\ZEP/2^{n }, q_n+\ZEP/2^{n })$.
The sequence $\Delta=\{ I_n\}$ has the property that $\zzq\cap (0,1)\subseteq \mathcal{I}_\Delta$ and $L(\Delta)<\ZEP$.\zdia 
 
}\end{Example}
\begin{Remark}\ZLA{REMA:Ch0Qnullset}{\rm
The following observations have their interest:
\begin{enumerate}
\item\ZLA{I1Exe:Ch0Qnullset1}
the argument in Example~\ref{Exe:Ch0Qnullset} is very simple because we did not require that the component intervals   are  disjoint. Had we imposed this condition then the construction of the multiinterval  would   be    more delicate.

\item\ZLA{I2Exe:Ch0Qnullset1} Example~\ref{Exe:Ch0Qnullset} shows that in the definition of null set we cannot impose that the multiintervals have finitely many component intervals. 

\item\ZLA{I3Exe:Ch0Qnullset1} The method in Example~\ref{Exe:Ch0Qnullset} can be used to show that \emph{any numerable set is a null set.} There exists null sets which are not numerable. An example is the Cantor set, see Sect.~\ref{secCAP4:multirElebeMis}.\zdia

\end{enumerate}
}

\end{Remark}

 \begin{table}[h]\caption{Succinct notations and terminology}\ZLA{tableCH1shortNOTATIONS}
\ZREF{
In order to speed up the statements, it is convenient to introduce the following notations suggested by the previous considerations:
\begin{itemize}
\item let  $\{\Delta_n\}$ be a sequence of multiintervals and let the multiinterval $\Delta$ be constructed
from $\{\Delta_n\}$ with the procedure
 in Lemma~\ref{ch1:LemmaUnioMULTINT}. The multiinterval $\Delta$ is denoted $\cup\Delta_k$;
\item we say that a {\sc multiinterval covers a set $A$}\index{multiinterval!which covers $A$} when $A\subseteq \mathcal{I}_\Delta$. 
\item we say that a multiinterval {\sc is in a set $A$}\index{multiinterval!in a set} when $\mathcal{I}_\Delta\subseteq A$.
\item we say that $A$ is the set of $\Delta$ when $A=\mathcal{I}_\Delta$.
\item we say that a multiinterval $\tilde\Delta $ is {\sc extracted}\index{multiinterval!extracted} from $\Delta$ when any component interval of $\tilde \Delta$ is a component interval of $\Delta$.
\item a multiinterval which has finitely many component intervals is called a ``finite sequence'' (of intervals) or a {\sc finite multiinterval.}\index{multiinterval!finite} 
\end{itemize}

}
\end{table}

We use the terminology introduced in the table~\ref{tableCH1shortNOTATIONS} and we note the following simple observations.

\begin{Lemma}\ZLA{LemmaCAP1UnioNulliEnulla} Let $\{N_n\}$ be a sequence of null sets. Then $N=\cup N_n$ is a null set.
\end{Lemma}
\zProof
We fix $\ZEP>0$ and we construct a multiinterval $\Delta$ such that
 \[\mbox{$\Delta$ covers  $N$, i.e.}\
N\subseteq \mathcal{I}_\Delta\,,\ \mbox{and}\  L(\Delta)<\ZEP\,. 
 \]
 
The construction is as follows. For every $n$ there exists $\Delta_n$ such that
\[
N_n\subseteq \mathcal{I}_{\Delta_n}\,,\qquad L(\Delta_n)<\frac{\ZEP}{2^n}\,.
\]
The multiinterval $\Delta_n$ exists because $N_n$ is a null set.

The multiinterval $\Delta=\cup\Delta_n$ has the required properties.\zdia

We noted in Remark~\ref{REMA:Ch0Qnullset} that in general a multiinterval which covers a set $A$ has infinitely many component intervals  even if $A$ is bounded. The following observation has its interest: if $K$ is  compact, i.e. bounded and closed, then any family of open sets which covers $K$ has a \emph{finite} subfamily which still covers $K$. This fact can be recasted in terms of multiintervals as follows\footnote{the fact that a multiinterval is composed of \emph{open} intervals is crucial for Theorem~\ref{I2Exe:Ch0Qnullset} to hold.}:
\begin{Lemma}\ZLA{I2Exe:Ch0Qnullset} 
Let $\Delta$ be a multiinterval which covers a  compact set $K$. There exists a finite multiinterval $\tilde \Delta$ which covers $K$ and which is extracted from $\Delta$.

\end{Lemma}
 
Now we recast the property of being a.e. continuous given in Definition~\ref{Ch1DefiNULLsetqO}   in a more baroque style, which however suggests a   general definition: let $f$ be a function defined on $A$ and let $\mathcal{P}(x)$ be the following proposition which applies to the points  $x\in A$:   $\mathcal{P}(x) $ is true    when $f $ is continuous at $x$. We say that $f$ is a.e. continuous on $A$ if $\mathcal{P} $ is false on a null set. In symbols:
\[ \mbox{$f$ is a.e. continuous on $A$ }\ \iff\ 
\{ x\in A\ \ztc \ \neg \mathcal{P}(x)\}\quad \mbox{is a null set}\,.
 \]
 This baroque way of stating continuity a.e. shows that  the key ingredient is a property of the points of $A$. So we define: 

\begin{Definition} 
A property $\mathcal{P}$ of the points of $A$   holds {\sc almost everywhere}\index{almost everywhere} (shortly {\sc a.e.}\index{a.e.})   on $A$ when the subset where it is false is a null set.
\end{Definition}
 
 For example, a function defined on $A$ is a.e. positive when $\{ x\ztc f(x)\leq 0\}$ is a null set; a function is a.e. defined on $A$ when $A\setminus\Dom\, f$ is a null set.
 
 
The previous observations that we have seen in dimension $1$ have   counterparts in every space $\zzr^d$. The next result 
  holds only in dimension~1:

\begin{Theorem}\ZLA{teo:Ch2StrutturaAperti}
Any nonempty open   set $\mathcal{O}\subseteq \zzr$ is the union of the component intervals of a \emph{disjoint  multiinterval\footnote{we recall that according to our definition a multiinterval is composed of \emph{open} intervals.} :}   
there exists   $\Delta=\{  I_n\}$  such that 
\[
\mathcal{O}=\mathcal{I}_\Delta=\bigcup I_n\quad \mbox{and \ \ $I_n\cap I_k=\emptyset$ if $n\neq k$}\,.
\]

Two different multiintervals with this property differ solely for the order of the component intervals $I_n$.

\end{Theorem}

\begin{Corollary}\ZLA{CoroCH1PReMisOpe}
Let $\Delta$ be a multiinterval. There exists a \emph{disjoint } multiinterval $\hat  \Delta $ such that
\[
 \mathcal{I}_{\hat \Delta}=\mathcal{I}_{\Delta}\,.
\]

\end{Corollary}

 \begin{Remark}\ZLA{RemaCH1AlCoroCH1PReMisOpe}{\rm
The previous corollary shows that \emph{in dimension $1$}
it is equivalent to work with arbitrary multiintervals or with disjoint  multiintervals. Furthermore, it is the key observation for the following definition of the ``measure'' of an open set.\zdia
}\end{Remark}
\begin{Definition}\ZLA{DefiCH1DeFiMiSuOpen}
{\rm     
Le $\mathcal{O}\subseteq\zzr$ be a nonempty open set. We use Theorem~\ref{teo:Ch2StrutturaAperti} and we represent
 \begin{equation}\ZLA{Ch1DefiMISUaperti} 
 \mathcal{O}=\mathcal{I}_\Delta\,,
 \qquad \mbox{($\Delta=\{I_n\}$ is disjoint)}\,.
 \end{equation}  
 The {\sc measure}\index{measure!of an open set} of $ \mathcal{O}$ is the number
 \[
 \zl(\mathcal{O})=L(\Delta)\,.\zdiaform
 \]
}
\end{Definition}
The disjoint multiinterval $\Delta$ in~(\ref{Ch1DefiMISUaperti}) is not unique but two multiintervals differ only for the order of the component intervals $I_n$ and so $ \zl(\mathcal{O})$ is uniquely defined. 

The following property are clear:
\begin{Theorem}\ZLA{Teo:Ch1P1MeasMultirec}
We have:
\begin{enumerate}
\item Let $\mathcal{O}\neq \emptyset$ be an open set. We have:
\[
\zl(\mathcal{O})=\inf\{L(\Delta)\,:\ \mathcal{O}\subseteq \mathcal{I}_\Delta\}\,.
\]
\item if $\mathcal{O}_1\subseteq \mathcal{O}_2$ are nonempty open sets then $\zl(\mathcal{O}_1)\leq \zl(\mathcal{O}_2)\,.$
\end{enumerate}
\end{Theorem}
\begin{Remark}
{\rm
Note that $\zl(\mathcal{O})\leq +\ZIN$. When $\mathcal{O}$ is bounded its measure is finite but it is not difficult to construct examples of unbounded open sets with finite measure.\zdia
}
\end{Remark}
 \subparagraph{\em The proof of Theorem~\ref{teo:Ch2StrutturaAperti}}
 The proof follows from the following simple observations:
 \begin{enumerate}
 \item the union of open intervals which have a common point is an open interval.
 \item let $x_0\in\mathcal{O}$. The union of the open intervals $I$ such that $x_0\in I\subseteq \mathcal{O}$ is the \emph{maximal} open interval to which   $x_0$ belongs and which is contained in $\mathcal{O}$. We denote it $I_{x_0}$.
 \item let $x_0$ and $x_1$ be points of the open set $\mathcal{O}$. We have either $I_{x_0}=I_{x_1}$ or $I_{x_0}\cap I_{x_1}=\emptyset$.
 \end{enumerate}
 
 These properties imply that $\mathcal{O}$ is the union of disjoint open intervals. These intervals can be arranged to form a sequence by choosing one rational number from each one of them and recalling that the rational numbers are numerable.

  \subsection{Restriction and Extensions of Functions}

Let $A\subseteq \zzr  $ and let $f$ be a real valued function defined on $A$.
We define:
\begin{enumerate}
\item if $B\subseteq\zzr $ then $g=f_{|_{B}}$, the {\sc restriction}\index{function!restriction} of $f$ to $B$, is defined when $B\cap A\neq \emptyset$ and by definition $\Dom\, g=A\cap B$ and if  $b\in \Dom\, g$ then $g(b)=f(b)$. So, \emph{the graph of $g$ is obtained by that of $f$ by deleting the points $(a,f(a))$ for which $a\notin B$.}

\emph{The function $f$ has a unique restriction to $B$.}
\item Let instead $\Dom\, g\subseteq B$ and let $f$ be defined on    $A\supseteq \Dom\, g$. The function $f$ is an {\sc extension}\index{function!extension}\index{extension!function} of $g$ when $f_{|_{B}}=g$. So, a function defined on $B$ admits infinitely many extensions to $A$ (unless $A=\Dom\, g$!) and in practice the extensions which are used have additional special properties.

 \end{enumerate}
 It is obvious:
\begin{Lemma}\ZLA{Lemma:Ch1UNIlimAssocSeq} if $|f(x)|<M$ for every $x \in A $ then we have also $| f_{|_B}(x)|<M$ for every $x\in B$.
\end{Lemma} 
 
 It is important to be clear on the relation of continuity of a function and of its restrictions. In order to appreciate these relations we state explicitly the definitions of   continuity  at the points of $B\subseteq A$ and the definition of continuity 
of~$f_{|_{B}}$.

 Let $B\subseteq A$ and let $f$ be defined on $A$. Let
\[
g=f_{|_B}
\]
be the restriction of $f$ to the set $B$. 

\begin{description}
\item[\bf Continuity   of $f$ at  $x_0\in B\subseteq A$:] \index{continuity!at a point}\index{function!continuous} for every $\ZEP>0$ there exists $\ZDE>0$ such that  if \fbox{$x\in A$} and $|x-x_0|<\ZDE$ then $|f(x)-f(x_0)|<\ZEP$. \emph{Note that we do not impose $x\in B$.}

In terms of sequences, $f$ is continuous at $x_0$ when
\[
\fbox{$x_n\in A$}\ \mbox{and}\ x_n\to x_0\ \implies f(x_n)\to f(x_0)=g(x_0)\,.
\]
\item[\bf Continuity   of $g=f_{|_{B}}$ at  $x_0$:]
\index{continuity!of the restriction}\index{function!restriction!continuity}
 it must be $x_0\in B$ since 
$
\Dom\,f_{|_{ B}}=B 
$.
The definition of continuity is as follows: for every $\ZEP>0$ there exists $\ZDE>0$ such that  if \fbox{$x\in B$} and $|x-x_0|<\ZDE$ then $|f(x)-f(x_0)|<\ZEP$.
In terms of sequences, $f$ is continuous at $x_0\in B$ when
\[
\fbox{$x_n\in B$}\ \mbox{and} \ x_n\to x_0\ \implies \underbrace{f(x_n)}_{\begin{tabular}{c} $\rotatebox{90}{=}$\\
$g(x_n)$
\end{tabular}}\to \underbrace{f(x_0)}_{\begin{tabular}{c} $\rotatebox{90}{=}$\\
$g(x_0)$
\end{tabular}}\,.
\]
\end{description}

Furthermore we note:
if $f$ is continuous at $x_0\in B$ then $g=f_{|_{B}}$ is continuous too but
continuity of $g=f_{|_{B}}$  at $x_0\in B$ nothing says on the continuity of  $f$:  it is possible that $g$ is continuous while $f$ is not continuous. We give the following quite elaborate example (which is important for the definition of the integral) and we invite the reader to find a simpler one.
\begin{Example}\ZLA{Ch1Exe19DISCdirFUnaeCONT}{\rm
Let $A=[0,1]$ and  $f =d $ be the \defin{{Diriclet function}}\index{Dirichlet function}\index{function!Dirichlet} on $[0,1]$:
\begin{equation}
\ZLA{eq:Ch0:DiricFUNC}
f(x)=d(x)=\left\{
\begin{array}{lll}
1&{\rm if}& x\in\zzq\\
0&{\rm if}& x\in\zzr\setminus\zzq\,.
\end{array}
\right.
\end{equation}
\emph{The function $f=d $ is discontinuous at every point. }

Let $B=\zzq\cap(0,1)$, the set of the rational points of $(0,1)$. The function $g=f_{|_{B}}$   is continuous on $B$ because it is the constant function $0$.

Analogously we see that $f{_{|_C}}$ is continuous when $C=[0,1]\setminus\zzq$.

The message of this example is that when studying continuity of the restriction of a function on a set $B$, we only consider the part  of its graph which projects ortogonally to points of $B$. \emph{The remaining part of the graph is deleted.}\zdia
}\end{Example}

Finally, as an application of the definition of continuity,    we 
invite the reader to prove the following result\footnote{a proof is in Sect.~\ref{ch4SECDEFIprimPromeasurSETS}.}: 

\begin{Theorem}\ZLA{Ch1:TheoFuLimiTronca}
Let   $A\subseteq \zzr$ and let $f$
and $g$ be continuous on $A$. The functions
\[
\phi(x)=\max\{f(x),g(x)\}\,,\qquad \psi(x)=\min\{f(x),g(x)\}
\]
are continuous on $A$. In particular, let $N$ be any real number and let
\begin{equation}\ZLA{eq:DefiNOTafPM}
f_{+,N}(x)=\max \{f(x),N\}\,,\qquad f_{-,N}(x)=\min \{f(x),N\}\,.
\end{equation}
The functions $f_{+,N}$ and $f_{-,N}$ are continuous on $A$.
 \end{Theorem}
 
\subsubsection{Associated Multiintervals}

Let   $f$ be a function a.e. defined on an interval $R$. The interval can be unbounded, i.e. it can be a half line or it can be $\zzr$; the endpoints of the interval can belong to $R$ or not, i.e. $R$ can be open or closed or half open.

Let $\Delta$ be a multiinterval in $R$, i.e. such that $\mathcal{I}_\Delta\subseteq R$. We say that $\Delta $ is a {\sc multiinterval associated to $f$}\index{associated!multiinterval}
 when $f_{|_{R\setminus\mathcal{I}_\Delta}}$ is continuous. If $L(\Delta)<\ZEP$ we say that $\Delta$ is an {\sc associated multiinterval of order $\ZEP$}\index{multiinterval!associated!of order $\ZEP$}.

We note that   an   associated multiinterval of order $\ZEP$  needs not be unique.

We consider the following example:
 \begin{Example}\ZLA{eseCH1:Assocjumps}
 {\rm
Let $f$ be defined on $(h,k)$ and continuous on   $(h,\tilde h)\cup(\tilde h,k)$. Any interval $(\tilde h-1/n,\tilde h+1/n)$ is an associated interval of order $2\ZEP$ if $1/n<\ZEP$.

 The reader is invited to use the previous idea and to combine the examples~\ref{Exe:Ch0Qnullset}
 and~\ref{Ch1Exe19DISCdirFUnaeCONT} and to construct an associated multiinterval of order $\ZEP$ for the Dirichlet function (in case of difficulty see the Example~\ref{ESEch1DiriFUnzQUAScont}).\zdia
 }
 \end{Example}
 
 This concept of associated multiinterval is crucial in the Tonelli definition of the interval. So it is convenient to present few simple comments.

\begin{Lemma}\ZLA{LemmaCH1ProprieAssocMultintA}
let $R $  be an interval  and let $\ZEP>0$.  
The following properties hold:
\begin{enumerate}
 \item\ZLA{I1LemmaCH1ProprieAssocMultint}
let $R $    be not closed. Let 
 $f$ be defined on $R$  and let there exist  a multiinterval $\Delta$ of order $\ZEP$ associated to $f$.
 If $g$ is an extension of $f$ to ${\rm cl}\, R$ then there exists a multiinterval of order $\ZEP$ associated to $g$.
\item\ZLA{I2LemmaCH1ProprieAssocMultint} let $f$ and $g$ be a.e. defined on $R$ and let $N=\{ x\,:\ f(x)-g(x)\neq 0\}$ be a null set. If
there exists a multiinterval  $\Delta $ of order $\ZEP$ associated to $f$ then there exists also a
  multiinterval of order $\ZEP$ associated to $g$.
 \item  \ZLA{LemmaCH1ProprieAssocMultintB}
let $f$  and $g$ be a.e. defined on an interval $R$ and let 
  $\hat \Delta$ and $\tilde\Delta$ be multiinterval of order $\ZEP$ associated respectively to $f$ and to $g$. Then $\hat\Delta\cup\tilde\Delta$ is a multiinterval of order $2\ZEP$ associated to $f+g$.
\item  \ZLA{LemmaCH1ProprieAssocMultintC}
As in~(\ref{eq:DefiNOTafPM}) we define, for every $N\in\zzr$,
 \begin{equation}\ZLA{eqCh1DefiFpmN}
f_{+,N}(x)=\max\{f(x),N\}\,,\qquad  f_{-,N} (x) =\min\{f(x),N\}\,.
 \end{equation}
 If $\Delta$ is a multiinterval associated to $f$ it is also a multiinterval
 associated to $f_{+,N}$ and to $f_{-,N}$.
\end{enumerate}
\end{Lemma}
 \zProof 
 We prove statement~\ref{I1LemmaCH1ProprieAssocMultint}.
Note that $R\neq\zzr$ since it is not closed. So it has one or two (finite) endpoints. 
  The multiinterval associated to $g$ is obtained by adding to $\Delta$ one or two intervals which cover the end points of $R$ of length less then $\left (\ZEP-L(\Delta)\right )/2$.
 
 Statement~\ref{I2LemmaCH1ProprieAssocMultint}  is proved in a similar way: we construct a multiinterval $\Delta_1$ which covers $N$ and such that $L(\Delta_1)<\ZEP-L(\Delta)$. The required multiinterval is $\Delta\cup\Delta_1$.
 
 Statement~\ref{LemmaCH1ProprieAssocMultintB} is proved in a similar way and statement~\ref{LemmaCH1ProprieAssocMultintC} is a reformulation of
Theorem~\ref{Ch1:TheoFuLimiTronca}  when $A=R\setminus\mathcal{I}_\Delta$ and $\Delta$ is an associated multiinterval to $f$.\zdia
 \begin{Theorem}
 \ZLA{teo:CH1MultiASSOCsucceFUnz}
 Let $\{f_n\}$ be a sequence of functions. We assume that for every $\ZSI>0$ there exists a multiinterval $\Delta_{n,\ZSI}$ of order $\ZSI$ associated to $f_n$. Under this condition, for every $\ZEP>0$ there exists a multiinterval $\Delta_\ZEP$ which is associated to every $f_n$.
 \end{Theorem}
 \zProof The multiinterval $\Delta_\ZEP$ is 
 \[
 \Delta_\ZEP=\bigcup_{n=1}^{+\ZIN} \Delta_{n,\ZEP/2^n}\,.
 \zdiaform
 \]

\subsection{Tietze Extension Theorem: One Variable}
A key result which is used in the construction of the integral is
  Tietze extension theorem.
We introduce the following definition:

 \begin{Definition}
Let $f$ be continuous on the closed subset $K$ of $\zzr$.  We   call {\sc Tietze extension of $f$}\index{extension!Tietze}\index{Tietze extension} any continuous extension\footnote{the hypothesis that the function $f$ is defined on a closed set is crucial for the existence of continuous extensions. In general there is no continuous extension from an open set.} $f_e$  of $f$ to an interval $R\supseteq K$ such that 
  \begin{equation} \ZLA{eq:CarattDItietzeEXTEN}
  \begin{array}{l}
  \displaystyle 
\inf\left \{ f(x)\,,\  x\in K\right \}=\inf\left \{ f_e(x)\,,\  x\in R\right \}\,,\\[2mm]
\displaystyle
\sup\left \{ f(x)\,,\  x\in K\right \}=\sup\left \{ f_e(x)\,,\ x\in R\right \}\,.
\end{array}
  \end{equation}
  \end{Definition}

Tietze extensions exist in any normal topological space. Proofs in such generality can be found in books on general topology. Here we give a proof in the case of functions of one variable. The proof is intuitive thanks to Theorem~\ref{teo:Ch2StrutturaAperti} and furthermore it gives an extension which has an important additional property. 
  
  \smallskip
 
\hskip -.6cm\begin{minipage}{2.5in}
While reading the theorem, it may be helpful to look at the figure here on the right. 
The geometric definition of the extension and a look to the figure   easily shows the properties stated in the theorem. For completeness, the details of the proof are given at the end of this section. 

In the figure, the set $K$ and the graph of $f$ are red while the graph of the extension to $\zzr\setminus K $ is blue.
\end{minipage}
\begin{minipage}{2in}
 \includegraphics[width=7cm,scale=.65]{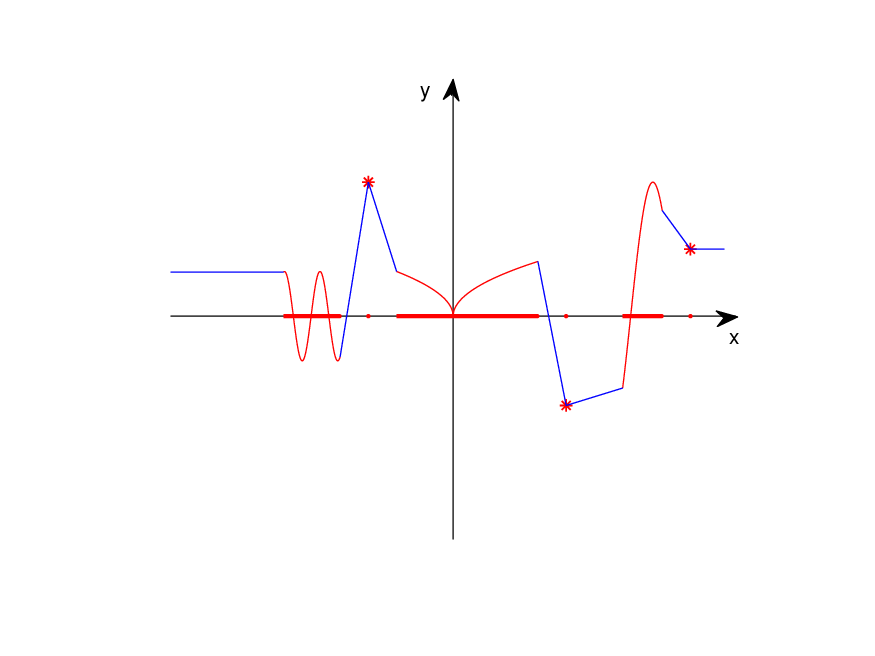}

  \end{minipage}


\begin{Theorem}[{\sc Tietze extension theorem}]\index{Theorem!Tietze!$1$ variable}\ZLA{teo:Ch1ESTEdaCHIUSO}
Let 
$K\subseteq\zzr $ be a  closed 
 set and let
  $f$ be defined on $K$ and continuous.
  

The extension $f_e$ of $f$   constructed with the procedure described below has the following properties:
\begin{enumerate}
\item\ZLA{I1teo:Ch1ESTEdaCHIUSO}  the function $f_e$ is continuous on $\zzr$;

  \item\ZLA{I2teo:Ch1ESTEdaCHIUSO}   we have
  \begin{eqnarray*}&&
\min\left \{ f(x)\,,\  x\in K\right \}=\min\left \{ f_e(x)\,,\  x\in \zzr\right \}\,,\quad \\
&&
\max\left \{ f(x)\,,\  x\in K\right \}=\max\left \{ f_e(x)\,,\ x\in \zzr\right \}\,.
  \end{eqnarray*}
\item\ZLA{I3teo:Ch1ESTEdaCHIUSO} let $g$ be a second continuous function defined on $K$ and let $g_e$ be its extension obtained with the procedure described below.
  If $f(x)\geq g(x)$ on $K$ then $f_e(x)\geq g_e(x)$.
\end{enumerate}

The function $f_e$ is defined as follows: $f_e(x)=f(x)$ if $x\in K$ while if $x\notin K$ we proceed as in the following steps:
\begin{description}
\item[\bf Step~0:]     we use Theorem~\ref{teo:Ch2StrutturaAperti}  and we represent        $\zzr\setminus K=\cup (a_n,b_n)$ (the open intervals $(a_n,b_n)$ are pairwise disjoint);
\item[\bf Step~1:]  we note that the endpoints $a_n$ and $b_n$ belong to $K$;
\item[\bf Step~2:] if $K$ is bounded above then (only) one of the interval $(a_n,b_n)$ has the form $(r,+\ZIN)$ with $r\in K$. On this interval we put $f_e(x)=f(r)$. Analogous observation and definition if $K$ is bounded below;
\item[\bf Step~3:] on the bounded interval $(a_n,b_n)$ the 
extension $f_e$ interpolates linearly among the points $(a_n,f(a_n))$ and $(b_n,f(b_n))$; i.e. the graph is the segment which joins these points. In analytic terms:
\begin{equation}
\ZLA{eqP1Ch1FormaAnalTietzEXTE}
\begin{array}{l}
\displaystyle\mbox{if $x\in (a_n,b_n)$ then $x=\zl a_n+(1-\zl)b_n$ (with $\zl\in  (0,1)$)}\,.\\
\displaystyle\mbox{By definition, $f_e(x)=\zl f(a_n)+(1-\zl) f(b_n)$}\,.
\end{array}
\end{equation}
\end{description}

\end{Theorem}

%
We note:
  \begin{itemize}
  \item the property in the statement~\ref{I3teo:Ch1ESTEdaCHIUSO} of Theorem~\ref{teo:Ch1ESTEdaCHIUSO} is a property of the special Tietze extension obtained with the procedure described in the theorem. It is not a property of any Tietze extension.
  \item
  it is clear that in general there are infinitely many extensions of a given function.  
  Fig.~\ref{fig:ch1:esedefiQUAScontFUN} in Example~\ref{eseCH1:ContFunzCONjumps} below shows two different Tietze extensions to $[h,k]$ of a  function which is continuous on $K=[h,\hat h-\ZEP)\cup(\hat h+\ZEP ,k]$  ($\ZEP>0$)
  while Example~\ref{ESEch1DiriFUnzQUAScont} below shows that the Dirichlet function admits a unique Tietze extension. More in general we have:

  \begin{Lemma}\ZLA{LemmaCH1TietzDIquasiCONSTfunz}
  A function which is a.e. constant on an interval, $f(x)=c$ a.e., admits the unique Tietze extension $f_e(x)\equiv c$.
  \end{Lemma}
  The proof (similar to that in Example~\ref{ESEch1DiriFUnzQUAScont})  is left as an exercise to the reader.
  \end{itemize}
  
We state:
\begin{Theorem}\ZLA{TheoCH1ConveTieSeq}
Let $\{f_n\}$ be a sequence of continuous function defined on the closed set $K$
and let 
let $(f_n)_e $ be    Tietze extensions of $f_n$.
 If $f_n\to 0$ \emph{uniformly} on $K$ then $(f_n)_e\to 0  $ \emph{uniformly} on $\zzr$.
\end{Theorem}
\zProof The statement  is an obvious consequence of Theorem~\ref{Ch1TeoDefoCONVunif} and the inequality~(\ref{eq:CarattDItietzeEXTEN}) which holds for every Tietze extension.\zdia
 
%
%
%

%
Theorem~\ref{teo:Ch1ESTEdaCHIUSO}   is extended to functions of several variables in Chap.~\ref{Chap:2quasicontPIUvar}. The proof is in the Appendix~\ref{AppeCH2} where we extend also Theorem~\ref{TheoCH1ConveTieSeq}.

Now we state further properties of the special Tietze extensions obtained with the procedure described in Theorem~\ref{teo:Ch1ESTEdaCHIUSO}. These properties, which do not hold for general Tietze extensions, will be used to give a simple proof of Egorov-Severini Theorem in Chap.~\ref{Cha1bis:INTElebeUNAvariaBIS}.

 \begin{Theorem} 
\ZLA{Theo:prprieSucceConvDIM1}
 Let $\{f_n\}$ be a sequence of continuous functions defined on a compact set $K$ and let $(f_n)_e$ be the Tietze extension of $f_n$ obtained with the algorithm described in Theorem~\ref{teo:Ch1ESTEdaCHIUSO}. 
The following properties hold: 
\begin{enumerate}
\item\ZLA{I1Theo:prprieSucceConvDIM1}
If $\{f_n\}$ is an increasing (decreasing) sequence on $K$ then $\{(f_n)_e\}$ is an increasing (decreasing) sequence on $\zzr$. 
\item \ZLA{I2Theo:prprieSucceConvDIM1}
  if $\{f_n(x)\}$ converges   for every  $x\in K$ then $\{(f_n)_e(x)\}  $ converges for every  $x\in \zzr$. 
  \end{enumerate}
\end{Theorem}
 \zProof
 Statement~\ref{I1Theo:prprieSucceConvDIM1} is an immediate consequence of the property~\ref {I3teo:Ch1ESTEdaCHIUSO} of 
 Theorem~\ref{teo:Ch1ESTEdaCHIUSO}. We prove Property~\ref{I2Theo:prprieSucceConvDIM1}. We recall~(\ref{eqP1Ch1FormaAnalTietzEXTE}):
Let $\zzr\setminus K=\cup _{k\geq 1}(a_k,b_k)$. 
Let $x\in (a_k,b_k)$, a bounded interval. There exists $\zl\in(0,1)$ such that $x=\zl a_k+(1-\zl)b_k$ and 
\[
(f_n)_e(x)=\zl f_n(a_k)+(1-\zl)f_n(b_k)\,.
\]
The end points $a_k$ and $b_k$ belong to $K$ so that $\{f_n(a_k)\} $ and  $\{f_n(b_k) \}$ both converge. So also
\[
(f_n)_e(x)=\zl f_n(a_k)+(1-\zl)f_n(b_k)
\]
converges.

The same hold on an unbounded interval, either $(a_{k_0},+\ZIN )$  or $(-\ZIN,b_{k_0})$, since on these intervals either $(f_n)_e(x)=f_n(a_{k_0})$
or $(f_n)_e(x)=f_n(b_{k_0})$.\zdia

\paragraph{The proof of Theorem~\ref{teo:Ch1ESTEdaCHIUSO}}

 We apply Theorem~\ref{teo:Ch2StrutturaAperti} to the open set
\[
\mathcal{O}=\zzr\setminus K:\qquad \zzr\setminus K =\cup I_n=\cup (a_n,b_n) 
\]
(the intervals $I_n$ are pairwise disjoint).

The end points $a_n$ and $b_n$ do not belong to $I_n$, which is open. Hence they belong to $K=\Dom\, f$ and $f(a_n)$ and $f(b_n)$ are defined.

It is convenient to denote $\hat f$ the restriction of $f_e$ to $\mathcal{O}$. So, $\hat f$ is defined as follows: 
\begin{itemize}
\item
 If $(r,+\ZIN) $ is one of the intervals which compose $\mathcal{O}$ then we define $\hat  f(x)=f(r)$ when $x>r$ and analogously when $(-\ZIN,r)  $ is one of the intervals which compose $\mathcal{O}$ we define $\hat  f(x)=f(r)$ when $x<r$. 
 \item we  define $\hat  f$ on the \emph{bounded} intervals $I_n=(a_n,b_n)$  by assigning its graph:
 the graph of $\hat f $ on $I_n$ is the segment which joins the two points $(a_n,f(a_n))$ and $(b_n,f(b_n))$.
 \end{itemize}
 
  The  function $\hat f$ is continuous on $\mathcal{O}$ and   
 
\begin{equation}\ZLA{eqteo:Ch1ESTEdaCHIUSOA}
\min _{K} f \leq \hat  f(x)\leq \max_K f\qquad \forall x\in \mathcal{O}\,.
\end{equation}

In order to complete the proof we must prove continuity of the function
\[
f_e(x)=\left\{\begin{array}{lll} 
\hat  f (x)&{\rm if}& x\in\mathcal{O}= \zzr\setminus K\\
f(x)&{\rm if}& x\in K
\,.
\end{array}\right.
\]

 Every $x_0\in \mathcal{O}$ has a neighborhood contained in $  \mathcal{O}$ and on this neighborhood $ f_e(x)=\hat  f(x)$, hence it is continuous.
 
 We must prove continuity at the points of $K$. 
 
 \emph{For most of clarity we use the following notation: a point of $K$ is denoted $k$ while a point of $\zzr\setminus K=\mathcal{O}$ is denoted $x$. In spite of this, we recall that the endpoints $a_n$ and $b_n$ of the intervals $I_n$ belong to $K$.}
 
We consider a point $k_0\in K$ and we prove continuity of $f_e$ at $k_0$. 
Continuity is clear if $k_0\in {\rm int}\, K$ and so we must consider solely the case   $k_0\in \partial K$   so that $k_0$ belongs to $K$ and it is an accumulation point of $\mathcal{O}$.

 We fix any $\ZEP>0$. Continuity of the function $f$ defined on $K$ shows the existence of $\ZDE>0$ such that
 
\begin{equation}
\ZLA{eq:ch1PropPERcompleTietzANTE}
\left\{\begin{array}{l}
k\in \mathcal{I}_\ZDE=\{|k-k_0|<\ZDE\}\,,\\ k\in K
\end{array}\right.\ \implies\ 
\underbrace{f(k_0)}_{=f_e(k_0)}-\ZEP < \underbrace{f(k )}_{=f_e(k )} < \underbrace{f(k_0)}_{=f_e(k_0)}+\ZEP\,.
\end{equation}
 
In order to prove continuity of $f_e$ at $k_0$ we must prove the existence of $\tilde \ZDE\leq\ZDE$ such that also the following property holds:
\begin{equation}
\ZLA{eq:ch1PropPERcompleTietz}
\left\{\begin{array}{l}
x\in \mathcal{I}_{\tilde \ZDE}=\{|x-k_0|<\tilde \ZDE\}\,,\\ x\in \mathcal{O}
\end{array}\right.\  \implies\ 
\underbrace{f_e(k_0)}_{=f (k_0)}-\ZEP<\underbrace{f_e(x )}_{=\hat f (x )} < \underbrace{f_e(k_0)}_{=f (k_0)}+\ZEP\,.
\end{equation}

The   intervals   $(a_n,b_n)$ belong to the following three classes:
\begin{enumerate}
\item\ZLA{cl1teo:Ch1ESTEdaCHIUSO} those intervals which do not intersect $\mathcal{I}_\ZDE$. Their points are not considered in~(\ref{eq:ch1PropPERcompleTietz}) since we choose $\tilde\ZDE\leq \ZDE$ and we do not need to consider them.

\item\ZLA{cl2teo:Ch1ESTEdaCHIUSO} The intervals $(a_{n_0},b_{n_0})\subseteq \mathcal{I}_\ZDE$. There can be infinitely many such intervals.
\item\ZLA{cl3teo:Ch1ESTEdaCHIUSO}
The intervals $(a_{n_0},b_{n_0})$ which are not contained in $ \mathcal{I}_\ZDE$ but which intersect $ \mathcal{I}_\ZDE$. There are at most two of such intervals, one on the right and one on the left: an interval $(a_{n_0},b_{n_0})$ such that
\[
k_0\leq a_{n_0}<k_0+\ZDE\leq b_{n_0}
\]
and the analogous intervals on the left. 

By reducing the value of $\ZDE$ we can assume 
\begin{equation}\ZLA{eq:Ch1NellaProTieTZreDdel}
k_0\leq a_{n_0}<k_0+\ZDE< b_{n_0}\,.
\end{equation}
\end{enumerate}

Now we show how $\tilde \ZDE$ can be chosen. We consider values of $x$ on the right of $k_0$. A similar procedure can be done  when $x<k_0$.

When $x\in \mathcal{O}$ there exists a unique $n_0$ such that 
 $x\in I_{n_0}=(a_{n_0},b_{n_0})$ and the construction of $\hat f$ is such that 
 \begin{equation}\ZLA{eq:ch1PropPERcompleTietz0ANTE}
 \begin{array}{l}
 \mbox{$\hat f(x)=f_e(x)$ belongs to the interval of end points}\\
 \mbox{  
$  f_e(a_{n_0})=  f (a_{n_0})$ and $  f_e(b_{n_0})=  f (b_{n_0})$ i.e.}\\
  f_e(a_{n_0})=  f (a_{n_0})\leq \hat f(x)=f_e(x)\leq   f (b_{n_0})=f_e(b_{n_0})  \,.
\end{array}
 \end{equation} 
 
 The interval $(a_{n_0}, b_{n_0})$ that has to be considered is either in the case~\ref{cl2teo:Ch1ESTEdaCHIUSO}
or in the case~\ref{cl3teo:Ch1ESTEdaCHIUSO}. 

First we consider the case that $(a_{n_0}, b_{n_0})$ is in the case~\ref{cl2teo:Ch1ESTEdaCHIUSO}  i.e. $x\in (a_{n_0},b_{n_0}) \subseteq \mathcal{I}_\ZDE$ then we have
\[
\underbrace{\underbrace{f(k_0)}_{=f_e(k_0)}-\ZEP < \overbrace{\underbrace{f(a_{n_0} )}_{=f_e(a_{n_0} )}
\leq \underbrace{\hat f(x)}_{=f_e(x)}
\leq  \underbrace{f(b_{n_0} )}_{=f_e(b_{n_0} )}}^{\tiny\mbox{from~(\ref{eq:ch1PropPERcompleTietz0ANTE})}}
 < \underbrace{f(k_0)}_{=f_e(k_0)}+\ZEP}_{\tiny\mbox{from~(\ref{eq:ch1PropPERcompleTietzANTE})}}
\]
 as wanted with $\tilde\ZDE=\ZDE$.
  
 Instead, let $(a_{n_0},b_{n_0})$
be in  the case~\ref{cl3teo:Ch1ESTEdaCHIUSO}. We consider the case of the interval on the right side. We have 
the following  cases:
\begin{enumerate}
\item the case that $k_0< a_{n_0}$. In this case we reduce the value of $\ZDE$ and we choose $\tilde\ZDE< a_{n_0}-k_0$. This way, 
\[
[k_0,k_0+\ZDE)=[k_0,a_{n_0})
\]
and
the points $x\in (a_{n_0}, b_{n_0})$ have not to be considered: the required inequality~(\ref{eq:ch1PropPERcompleTietz}) holds on $[k_0,k_0+\tilde\ZDE)$.

\item the case $ a_{n_0}= k_0<k_0+\ZDE<b_{n_0}$.    On  $ [ k_0,k_0+\ZDE) $ We have $f_e(x)=\hat f(x)$, a continuous function and by suitably reducing the value of $\ZDE $ the required inequality~(\ref{eq:ch1PropPERcompleTietz})
is achieved.
 \end{enumerate}
This last observation ends the proof of continuity.

The statement~\ref{I3teo:Ch1ESTEdaCHIUSO} follows from the following observation: if $f (a_n)\geq g(a_n)$ and $f(b_n)\geq g(b_n)$ then~(\ref{eqP1Ch1FormaAnalTietzEXTE})
gives $\hat f(x)\geq \hat g(x)$ for every $x\in (a_n,b_n)$.\zdia

\section{Tonelli Construction of the Lebesgue Integral}
 
 This section contains the fundamental ideas of the Tonelli method, presented in the simplest case of the functions of one variable.
We proceed as follows:  
 first we introduce and discuss  the main ingredients used by Tonelli. In particular we define the quasicontinuous functions defined on (bounded or unbounded) intervals. Then we define the Lebesgue integral of bounded quasicontinuous functions defined on bounded intervals. Finally we define quasicontinuous functions defined on a large family of domains and their integral. In this case we do not assume that the functions or its domain are bounded.

 \subsection{\ZLA{SubsCH1KeyINGRD}The Key Ingredients Used by Tonelli}
 
Tonelli construction of the Lebesgue integral is based on {\bf three main ingredients.} 

\begin{description}
\item[{\bf The first ingredient}] is the definition  of the null sets (see Definition~\ref{Ch1DefiNULLsetqO}). 

\item[{\bf The second ingredient}] is the class of the  quasicontinuous functions defined on (bounded or unbounded) intervals that we define now\footnote{the extension to functions defined in a larger class of domains is in Sect.~\ref{sec:CH1LebeGENEdom}.}.

Let $R $ be a  bounder or unbounded interval. A   function $f $ defined a.e. on $R$ is  a  {\sc  quasicontinuous function}\footnote{``quasicontinuous functions'' corresponds to ``funzioni quasi continue'' used by Tonelli. A better translation would be ``almost continuous functions'' but this term might be confused with ``almost everywhere   continuous functions'' and we prefer to avoid it.}\index{function!quasicontinuous!$1$ variable} when the following property holds: 
\emph{
for every $\ZEP>0$ there exists a multiinterval $\Delta_\ZEP$ of order $\ZEP$ which is associated to $f$ i.e. 
 such that
\begin{equation}\ZLA{EQCh1DefiQuasicontBpun}
\left\{\begin{array}{l}
\displaystyle
L\left (\Delta_\ZEP\right )<\ZEP\,,\\
\displaystyle
 f_{|_{R\setminus\mathcal{I}_{\Delta_\ZEP}}} \ \mbox{is continuous}\,.
\end{array}
\right.
\end{equation}
}

Note that
\begin{equation}\ZLA{Eq:P1C1IneqPreSecIngre}
\inf _{R} f\leq \inf _{R\setminus\mathcal{I}_{\Delta_\ZEP}} f \leq \sup _{R\setminus\mathcal{I}_{\Delta_\ZEP}} f \leq \sup _{R} f\,.
\end{equation}

A quasicontinuous function which is bounded is a {\sc bounded quasicontinuous function.}\index{function!bounded!quasicontinuous}

\begin{Remark}\ZLA{RemaP1Ch1RemaSecoIngre}
{\rm
We observe:
\begin{enumerate}
\item\ZLA{I1RemaP1Ch1RemaSecoIngre} if a function $f$ is quasicontinuous on the interval $R$ then   its restriction to the bounded intervals $R\cap [-k,k]$ is quasicontinuous  for every $k$. The converse implication holds too. Let $f$ be quasicontinuous on  $R\cap [-k,k]$    for every $k$. Let $\ZEP>0$. We associate to $f_{|_{R\cap [-k,k]}}$ a multiinterval $\Delta_k$ such that $L (\Delta_k)<\ZEP/2^k$. The multiinterval $\Delta_\ZEP=\cup _{k\geq 1}\Delta_k$
is associated to $f$ on $R$  and it is of order $\ZEP$.
\item
 in order to see whether a function is quasicontinuous it is sufficient to check that property~(\ref{EQCh1DefiQuasicontBpun}) holds solely for the sequence $\ZEP_n=1/n$.\zdia
 \end{enumerate}
}
\end{Remark}

Statement~\ref{LemmaCH1ProprieAssocMultintB} 
of Lemma~\ref{LemmaCH1ProprieAssocMultintA}
 has the following important consequence:
\begin{Theorem}\ZLA{TheoCH1LinearQC}
The classes of the quasicontinuous functions and  that of the bounded quasicontinuous functions (defined on a fixed interval $R$) are linear spaces.
\end{Theorem}

\item[{\bf The third ingredient.}] Let $f $ be a \emph{bounded  quasicontinuous} function on the interval $R $. 
 
We construct a sequence   of continuous functions on $R$ which we call an {\sc associated sequence of continuous functions.}\index{function!associated sequence of continuous functions}

We proceed in two steps:
\begin{description}
\item[\bf Step 1: when $f$ is defined on a closed interval]
\end{description}
Let  $\Dom\, f=R $. The interval $R$ can be unbounded but it is closed. Let $ \ZEP_n>0$, $\ZEP_n\to 0$.
Let\footnote{we use  the simplified notation $\Delta_n$ instead of the complete notation $\Delta_{\ZEP_n}$.} $\Delta_n$ be multiintervals associated to $f$, $\Delta_n$ of order $\ZEP_n$:
\[
\left\{\begin{array}{l}
\displaystyle
L\left (\Delta_n\right )<\ZEP_n\,,\\
\displaystyle
 f_{|_{R\setminus\mathcal{I}_{\Delta_n}}} \ \mbox{is continuous}\,.
\end{array}
\right.
\]
\emph{The set $R\setminus\mathcal{I}_{\Delta_n}$ is closed.} Theorem~\ref{teo:Ch1ESTEdaCHIUSO} asserts that $ f_{|_{R\setminus\mathcal{I}_{\Delta_n}}}$ admits a Tietze extension 
$\left ( f_{|_{R\setminus\mathcal{I}_{\Delta_n}}}\right )_e$ to $R$, i.e. a continuous extension with the same upper and lower bound as $f$.

The sequence    $n\mapsto \left ( f_{|_{R\setminus\mathcal{I}_{\Delta_n}}}\right )_e$ is an {\sc associated sequence {\rm  (to $f$)} of continuous functions   of order $\ZEP_n$.}\index{associated!continuous function of order $\ZEP_n$}\index{function!associated of order $\ZEP_n$}


\begin{description}
\item[\bf  Step 2:  $f$ is a.e. defined on a nonclosed interval $R$]
\end{description}
  The key observation is Statement~\ref{I2LemmaCH1ProprieAssocMultint} in   Lemma~\ref{LemmaCH1ProprieAssocMultintA} which can be reformulated as follows:

\begin{Theorem}\ZLA{Ch1TeoRIFOequivExteQCont}
Let $f$ be   a.e. defined on an interval $R$ which is not closed.
Let $g$ be one of its extension to ${\rm cl}\, R$. Then, $g $ is quasicontinuous if and only if $f$ is.
\end{Theorem}

An {\sc associated sequence of  continuous functions to $f$}\index{associated!continuous function of order $\ZEP_n$}\index{function!associated of order $\ZEP_n$}    is by definition a sequence associated to any of its extensions $g$.

Given $\ZEP_n\to 0$,
any sequence $\{\Delta_n,f_n\}$ where $f_n=\left (f_{|_{R\setminus\mathcal{I}_{\Delta_n}}}\right )_{e}$ is a sequence of {\sc associated multiintervals and continuous functions of order $\ZEP_n$.}\index{associated!multiintervals and continuous functions of order $\ZEP_n$}

\begin{Remark}{\rm
It is clear that  neither the multiinterval $\Delta_n$  nor the functions  $ \left ( f_{|_{R\setminus\mathcal{I}_{\Delta_n}}}\right )_e $
are uniquely identified by $\ZEP_n$ but when there is no risk of confusion it may be convenient to denote  $\left ( 
f_{|_{R\setminus\mathcal{I}_{\Delta_n}}}
\right )_e$ with the simpler notation $f_{n,e}$ or even $f_n$.\zdia  }
\end{Remark}

 \end{description}

\subsubsection{The Properties of the   Quasicontinuous Functions}

 First a simple example:
\begin{Example}\ZLA{eseCH1:ContFunzCONjumps}
{\rm
A function which is continuous on $[h,k]$ a part a jump at $\hat  h \in[h,k]$ is quasicontinuous. We show explicitly this fact in the case $\hat  h \in(h,k)$ and we leave to the reader the case that $\hat  h$ is one of the endpoints.

We fix any $\ZEP$ such that
\[
0<\ZEP<\min\{ \hat  h-h\,,\ k-\hat  h\}\,.
\]

The multiinterval $\Delta_\ZEP$ is composed by the sole interval $(\hat  h-\ZEP,\hat  h+\ZEP)$ (but of course we might choose a more complicated multiinterval).

The function $f $ is continuous both on $[h,\hat  h-\ZEP]$ and on $[\hat  h+\ZEP,k]$. The function $f_\ZEP$ is for example one of the functions whose graph is   in   Fig.~\ref{fig:ch1:esedefiQUAScontFUN}.

\begin{figure}[h]
\begin{center}
\includegraphics[width=8cm,scale=.65]{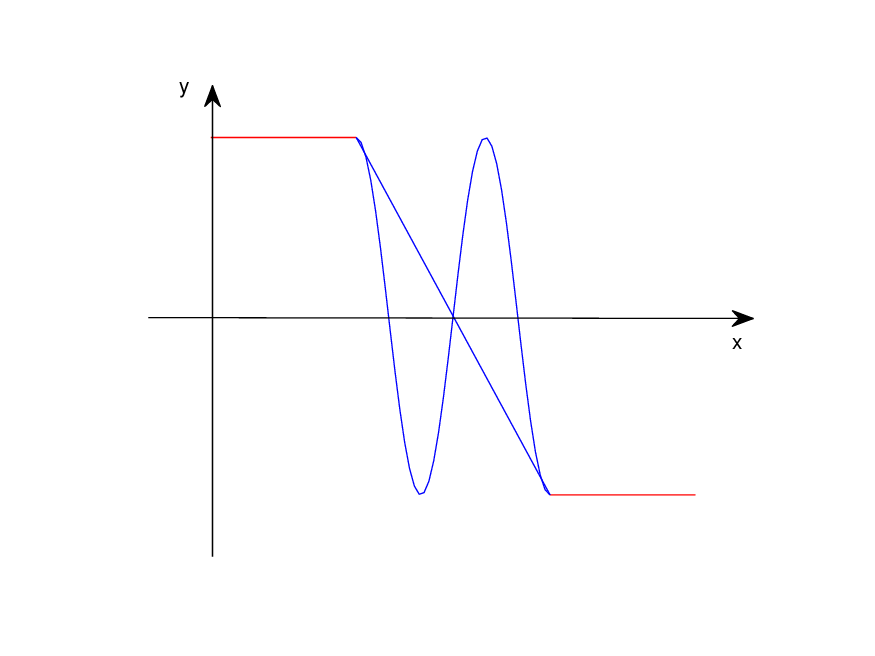}
 \end{center}
\caption{The function $f(x)$ and two different choices for  $\left (f_{|_{R\setminus\mathcal{I}_{\Delta_\ZEP}}}\right )_{e}$}

\label{fig:ch1:esedefiQUAScontFUN}       
\end{figure}

We leave the reader to write the analytic expression of the functions and we stress the fact that $\left (f_{|_{R\setminus\mathcal{I}_{\Delta_\ZEP}}}\right )_{e}$ is not uniquely identified by $\Delta_\ZEP$.

We observe that
\[
\lim_{\ZEP\to 0^+} \underbrace{\int_h^k \left (f_{|_{R\setminus\mathcal{I}_{\Delta_\ZEP}}}\right )_{e}(x)\ZD x}_{\tiny \begin{tabular}{l}
    Riemann      integral  \end{tabular}}
    =\underbrace{\int_h^k f(x)\ZD x}_{\tiny \begin{tabular}{l}
    Riemann      integral  \end{tabular}} 
\]
in spite of the fact that the convergence is not uniform.

This arguments is easily extended to functions with finitely many jumps, in particular to piecewise constant functions.
We recall that the approximation from below and from above of a function with piecewise constant functions is the key step in the definition of the Riemann integral.\zdia

}
\end{Example}

 Example~\ref{eseCH1:ContFunzCONjumps} can be much strengthened:

\begin{Lemma}\ZLA{lemmaCH1daqocontTOquasicont}
An a.e. continuous function which is a.e. defined on an interval $R $  is quasicontinuous.
\end{Lemma}
\zProof
We extend $f$    to ${\rm cl}\, R$ in an arbitrary fashion. The extension is still denoted $f$. By definition, for every $n$ there exists a multiinterval $\Delta$ such that $L(\Delta)<1/n$ and $f$ is continuous at the 
points of  $ R\setminus \mathcal{I}_\Delta$. So,   $f_{|_{R\setminus \mathcal{I}_\Delta}}$ is continuous too. It follows that $f$ is quasicontinuous.\zdia

The implication in Lemma~\ref{lemmaCH1daqocontTOquasicont} \emph{cannot} be inverted:
\begin{Example}\ZLA{ESEch1DiriFUnzQUAScont}
{\rm
We noted that the Dirichlet function on $[0,1]$ is discontinuous at every point. Hence \emph{it is not a.e.~continuous but it is quasicontinuous.} In fact, we proved in the Example~\ref{Exe:Ch0Qnullset} that the set of the rational points of $[0,1]$ is a null set. So, for every $n$ there exists a multiinterval $\Delta_n$
such that $L(\Delta_n)<1/n$ and $Q\subseteq\mathcal{I}_{\Delta_n}$.

The elements of the
set $K_n=[0,1]\setminus \mathcal{I}_{\Delta_n}$ are irrational points and so
\[
d_{|_{K_n}}(x)=0 \quad  \forall x\in K_n\,:
\]
the function $d_{|_{K_n}}$ is constant, hence continuous. So, \emph{the Dirichlet function is quasicontinuous in spite of the fact that it is not a.e. continuous.}

We note that the Tietze extension of $d_{|_{K_n}}$ is unique and it is identically zero:
\[
( d_{|_{K_n}} )_e=0 
\]
and this is a special instance of Lemma~\ref{LemmaCH1TietzDIquasiCONSTfunz}.\zdia
}
\end{Example}

The following simple observations have to be compared with Theorem~\ref{TheoCH1LinearQC}.
\begin{Theorem}\ZLA{TheoCH1PropELEMqcfunc}
Let $R$ be an interval. The following properties hold:
\begin{enumerate}
\item\ZLA{I1TheoCH1PropELEMqcfunc}
let   $f$ be quasicontinuous on $R$ and let   $R_1\subseteq R$ ($R_1$ is an  interval). Then $f_{|_{R_1}}$ is quasicontinuous on $R_1$.

\item\ZLA{I2TheoCH1PropELEMqcfunc} let $ \hat  h\in (h,k)$. If $f$ is  a.e. defined on $(h,k)$ and if it is   quasicontinuous both on $(h,\hat  h)$ and on $(\hat  h,k)$ then it is quasicontinuous on $(h,k)$.

\item\ZLA{I3TheoCH1PropELEMqcfunc}
the sum and the product of two (bounded)  quasicontinuous functions on an interval $R$ is a (bounded) quasicontinuous function on $R$. 

The quotient of quasicontinuous functions is quasicontinuous if the denominator is a.e. non zero.  

\item\ZLA{I3BisTheoCH1PropELEMqcfunc} let $f $ be quasicontinuouson the interval $R$ and let $g$ be a \emph{continuous} function on a domain which contains the image of $f$. The composition $x\mapsto g(f(x))$ is quasicontinuous.

 \item\ZLA{I4TheoCH1PropELEMqcfunc} let $f_n$ be quasicontinuous functions. For every $k$, the functions
 \begin{align*}
 &
\phi_k(x)=\max\{ f_1(x)\,,\ f_2(x)\,,\ \dots\,,\ f_k(x)\} \\
&\psi_k(x)=\min\{ f_1(x)\,,\ f_2(x)\,,\ \dots\,,\ f_k(x)\}
 \end{align*}
 are quasicontinuous.
\end{enumerate}
\end{Theorem}

These statements are either already noted or easily proved. In particular we examine the statement~\ref{I4TheoCH1PropELEMqcfunc}. It is sufficient to examine the maximum of two   functions $f(x)$ and $g(x)$. 

We use  Theorem~\ref{Ch1:TheoFuLimiTronca}: if both $f(x)$ and $g(x)$ are continuous at $x_0$ then $x\mapsto \max\{f(x)\,,\ g(x)\}$ is continuous at $x_0$ too.

Let the two functions $f$ and $g$ be quasicontinuous on $R$. For every $\ZEP>0$ there exists $\Delta_\ZEP$ such that both $f_{R\setminus\Delta_\ZEP}$ and $g_{R\setminus\Delta_\ZEP}$ are continuous. The restriction to $R\setminus\Delta_\ZEP$ of $\max\{f(x)\,,\ g(x)\}$ is equal to 
\[
\max\left \{f_{R\setminus\Delta_\ZEP}\,,\  g_{R\setminus\Delta_\ZEP}\right \}\,,
\]
the maximum of two continuous functions. Hence, it is a continuous function.

Let now $\carfunz_A$ be the {\sc characteristic function}\index{function!characteristic} of a set $A$:
\[
\carfunz_A(x)=\left\{\begin{array}
{lll}
1&{\rm if}& x\in A\\
0&{\rm if}& x\notin A\,.
\end{array}\right.
\]
The properties in the statements~\ref{I1TheoCH1PropELEMqcfunc} and~\ref{I2TheoCH1PropELEMqcfunc} of Theorem~\ref{TheoCH1PropELEMqcfunc} can be recasted as follows:
\begin{Corollary}
Let $f$ be a.e. defined on $[h,k]$  and let $\tilde h\in (h,k)$.  The function  $f$ is quasicontinuous if and only if both $ \carfunz_{[h,\tilde h]}f$ and $ \carfunz_{[\tilde h,k]}f$ are quasicontinuous.
 
\end{Corollary}

Now we recall the notations in~(\ref{eq:DefiNOTafPM})  and observe that

\[
f= \underbrace{f_{+,0}}_{=\max\{f,0\}}+\underbrace{f_{-,0}}_{=\min\{f,0\}}\,.
\]
The previous result can be applied separately to $f_{+,0} $ and to $f_{-,0}$. 
We apply Theorems~\ref{Ch1:TheoFuLimiTronca} and~\ref{TheoCH1PropELEMqcfunc}    to $f_{+,0} $ and to $f_{-,0}$  (when $A=R\setminus \mathcal{I}_{{\Delta_\ZEP}}$). We find:

\begin{Theorem}\ZLA{che1THEOFacFpmAC}
Let the function $f$ be a.e. defined on  $R $. Then we have:
\begin{enumerate}
\item\ZLA{I1che1THEOFacFpmAC} the function $f$  is quasicontinuous if and only if both $f_{+,0}$ and $f_{-,0}$ are quasicontinuous. 
\item\ZLA{I2che1THEOFacFpmAC}
if the function $f $ is quasicontinuous  on an interval $R$ then $|f|$ is quasicontinuous on $R$ too.
\end{enumerate}
\end{Theorem}

Statement~\ref{I2che1THEOFacFpmAC} follows  from
\[
|f(x)|=f_{+,0}(x)-f_{-,0}(x)\,.
\]

Theorem~\ref{che1THEOFacFpmAC} shows that it is not restrictive, when studying the properties of quasicontinuous functions, to assume that the functions do not change sign. This observation will be useful in the study of integration theory.

\begin{Remark}
{\rm
Observe that we  did not assert that the \emph{composition} of quasicontinuous function is quasicontinuous.
 In fact, in general it is not. See Remark~\ref{rema:AppeCap4NONcompos} in Appendix~\ref{APPEch4NullNONborel}.\zdia
}
\end{Remark}

\subsection{\ZLA{secCH1LebeBOUBound}The Lebesgue Integral under Boundedness Assumptions}
We define the Lebesgue integral of a bounded quasicontinuous  function a.e. defined on a bounded interval.  
 
We note a property of the Riemann integral:
 
\begin{Lemma}\ZLA{Lemma:Ch2LemmaPreliRiemaa}
Let $f $ be Riemann integrable (hence bounded) on $[h,k]$ and let $|f(x)|\leq M$ for all $x\in [h,k]$. Let there exists a multiinterval $\Delta$ such that $f(x)=0$ if $x\in [h,k]\setminus\mathcal{I}_\Delta$. Then the following inequality holds for the Riemann integral:
\[
\left | {\int_h^k f(x)\ZD x}\right |\leq ML(\Delta)\,.
\]
\end{Lemma}
\zProof 
The   Riemann integral is 
\begin{equation}\ZLA{Eq:Ch2LemmaPreliRiemaa}
\int_h^k f(x)\ZD x=\lim _{N\to+\ZIN} \sum_{i=0}^{N-1}f(x_{i,N}) (h_{i+1,N}-h_{i,N}) 
\end{equation}
where   $h_{i,N}=h+\dfrac{i}{N}(k-h)$
 and 
$ x_{N,i} \in [h_{i,N},h_{i+1,N})$  can be arbitrarily chosen
 (see Corollary~\ref{COROAppeCa1TheoDaPartApartEquiBIS} in Appendix~\ref{AppeTOch1IntegrRoemLebe}).
As the value of the integral does not depend on the choice of the  points $x_{i,N}$  we decide to choose $x_{i,N}\notin \mathcal{I}_\Delta$ provided that this is possible, i.e. when $[h_{i,N},h_{i+1,N})$ is not contained in $\mathcal{I}_\Delta$.

This way, the nonzero elements of the right sum in~(\ref{Eq:Ch2LemmaPreliRiemaa}) are those which corresponds to intervals  $[h_{i,N},h_{i+1,N})$  which are contained in $\mathcal{I}_\Delta$. The sum of the lengths of these intervals is less then $L(\Delta)$ so that
\[
\left |\sum_{i=0}^{N-1}f(x_{i,N}) (h_{i+1,N}-h_{i,N}) \right |\leq ML(\Delta)\,.
\]
This inequality is preserved in the limit.\zdia

We fix a bounded interval $R =[h,k]$ and a bounded quasicontinuous function $f $ a.e. defined on it:
\[
|f(x)|<M\,.
\]  Then we fix a sequence $(\Delta_n, f_n)$ of 
associated multiintervals and continuous functions of order $1/n$.
The functions $f_n$ are continuous, hence Riemann integrable, on $R$. 

Note that by the definition of the associated functions and from the inequality~(\ref{Eq:P1C1IneqPreSecIngre}) we have 
\begin{equation}\ZLA{eq:secCH1LebeBOUBoundSTIM}
|f_n(x)|<M\qquad \forall x\in[h,k]\quad\forall n\,.
\end{equation}

We prove:
\begin{Theorem}
The sequence of the Riemann integrals
\[
\int_h^k f_n(x)\ZD x
\]
is convergent and the limit does not depend
   on the particular associated sequence $(\Delta_n, f_n)$ of order $1/n$ in the sense that if $(\hat  \Delta_n, \hat  f_n)$ is a different associated sequence of order $1/n$ then
\[
\limn \int_h^k f_n(x)\ZD x=\limn \int_h^k \hat  f_n(x)\ZD x\,.
\]
\end{Theorem}
\zProof  
The continuous function  $f_n -f_m $ is different from zero solely on $\mathcal{I}_{\Delta_n}\cup \mathcal{I}_{\Delta_m}$. So, from Lemma~\ref{Lemma:Ch2LemmaPreliRiemaa}, we have
\[
\left |\int_h^k f_n(x)\ZD x-\int_h^k f_m(x)\ZD x\right |= 
\int_h^k |f_n(x)-f_m(x)|\ZD x\leq 2M\left (\frac{1}{n}+\frac{1}{m}\right )\,.
\]
It follows that the sequence
\[
\left \{\int_h^k f_n(x)\ZD x\right \}
\]
is a Cauchy sequence, hence it is convergent.

In a similar way we see that the limit does not depend of the special associated sequence.

First we note that  the inequalities~(\ref{eq:secCH1LebeBOUBoundSTIM}) holds both for $f_n$ and for $\hat f_n$: if $|f(x)|<M$ then we have $|f_n(x)|<M$  and     $|\hat  f_n(x)|<M$ for every $n$.
Then we observe 
\[
|f_n(x)-\hat  f_n(x)|=0\quad\forall x\in [h,k]\setminus\left \{\mathcal{I}_{\Delta_n}\cup \mathcal{I}_{\hat  \Delta_n}\right \}\,.
\]
Hence
\[
 \left |\int_h^k f_n(x)\ZD x-\int_h^k \hat  f_n(x)\ZD x\right |< \dfrac{4M}{n}
 \]
 so that 
 \[
 \limn  \left |\int_h^k f_n(x)\ZD x-\int_h^k \hat  f_n(x)\ZD x\right |=0\,.\zdiaform
\]

\emph{It is clear that the special sequence $1/n$ has no role and the same argument can be repeated for sequences  $(\Delta_n, f_n)$ of order $\ZEP_n$  and  $(\hat  \Delta_n, \hat  f_n)$  of order $\hat  \ZEP_n$
with $\ZEP_n\to 0$ and $\hat \ZEP_n\to 0$.} I.e. we have

\begin{Theorem}
Let $\{\ZEP_n\}$ and $\{\hat\ZEP_n\}$ be two sequences of positive numbers both convergent to zero. Let $ f_n $, $ \hat f_n $ be associated continuous functions (to $f$) of order respectively to $\ZEP_n$ and to $\hat\ZEP_n$. Then we have
\[
\limn \int_h^k f_n(x)\ZD x= \limn \int_h^k\hat f_n(x)\ZD x\,.
\]

\end{Theorem}

These observations justify the following definition:
let $f $ be a bounded quasicontinuous function a.e. defined on the bounded interval $R=[h,k]$. We fix any sequence $\{\Delta_n, f_n\}$ of 
associated multiintervals and continuous functions of order $1/n$.
The {\sc Lebesgue integral}\index{integral!Lebesgue!$1$ variable}\index{Lebesgue!integral!$1$ variable} of $f $ on $R$ is defined as follows
\begin{equation}\ZLA{eq:Ch1DEFIlebeFlimitAT}
\underbrace{\int_h^k f(x)\ZD x}_{\tiny \begin{tabular}{l}
    Lebesgue  \\    integral  \end{tabular}}=\limn\underbrace{\int_h^k f_n(x)\ZD x}_{\tiny \begin{tabular}{l}
    Riemann  \\    integral  \end{tabular}}\,.
\end{equation}
A consequence of the fact that $f$ is only a.e. defined on $R$ is that the integral does not change if the interval $R$ is open or half open.
\begin{Example}\ZLA{Ese:TWOeseInteCosTraDiri}{\rm
Let us see two examples:
\begin{enumerate}
\item
Let $f$ be a.e. defined on $[h,k]$, and let
\[
f(x)=f_1\quad h\leq x<\tilde h\,,\qquad f(x)=f_2\quad \tilde h<x<k\,.
\]
Tietze extensions of this functions have been considered in Example~\ref{eseCH1:ContFunzCONjumps}. From the arguments in Example~\ref{eseCH1:ContFunzCONjumps} we see that the \emph{Lebesgue} integral of $f$ is
\[
\int_h^k f(x)\ZD x= f_1(\tilde h-h)+ f_2(k-\tilde h)\,.
\]
It is known that the sum on the right is also the \emph{Riemann} integral of the function.

This observation can be extended to any piecewise continuous function. Let $h_i$, $0\leq i\leq N$ be points such that
\[
h=h_0\,,\quad h_i\leq h_{i+1}\,,\quad h_N=k
\]
and let $\chi(x)=\chi_i$ if $x\in [h_i,h_{i+1} )$. Then, its \emph{Lebesgue} integral is
\[
\int_h^k \chi(x)\ZD x=\sum _{i=0}^N \chi_i(h_{i+1}-h_i)
\]
\emph{and this number is also the \emph{Riemann} integral of $\chi$.}

\item
Example~\ref{ESEch1DiriFUnzQUAScont} shows that
\[
\underbrace{\int_h^k d(x)\ZD x}_{\tiny \begin{tabular}{l}
    Lebesgue      integral  \end{tabular}}=0\quad\mbox{($d$ is the Diriclet function)} 
\]
 as required in Remark~\ref{Ch1RemaNelQualDEInteNulso}.\zdia
 \end{enumerate}
 }\end{Example}

The definition of the   oriented  Riemann  integral    suggests to define also of the {\sc oriented  Lebesgue  integral:}\index{integral!Lebesgue!oriented}
\[
\underbrace{\int_h^k f (x)\ZD x
=-
\int_k^h f (x)\ZD x
}_{ \rm Lebesgue\  integrals }
\quad \mbox{since}\quad 
\underbrace{\int_h^k f_n (x)\ZD x
 =-
 \int_k^h f_n (x)\ZD x}_{ \rm Riemann\  integrals } \,.
\]
 
 When the interval of endpoints $h$ and $k$ (with $h\leq k$) is denoted $R$ then we use, both for the Riemann and the Lebesgue integrals,    
\[
\int_R f(x)\ZD x\quad \mbox{to denote}\quad \int_h^k f(x)\ZD x\qquad \mbox{(we repeat:   $h\leq k$)}\,.
\]

\ZREF{From now on, the integral sign will always  denote the Lebesgue integral, unless explicitly stated that it is a Riemann integral.
The fact that $\int$deno\-tes   a Lebesgue integral is explicitly indicated when convenient for  clarity.
}

Now we state the following obvious   result  (recall Example~\ref{Ese:TWOeseInteCosTraDiri} for the second statement and Lemma~\ref{LemmaCH1TietzDIquasiCONSTfunz} for the third):
\begin{Theorem}\ZLA{ch1TeoInteRIMlebContSTEP}
We have:
\begin{enumerate}
\item\ZLA{I1ch1TeoInteRIMlebContSTEP}
any $f\in C([h,k])$ is Lebesgue integrable and its Lebesgue integral coincide with its Riemann integral. In particular, if $f(x)\equiv c$ on $[h,k]$ then its Lebesgue integral is $c(k-h) $. 
\item\ZLA{I2ch1TeoInteRIMlebContSTEP} any piecewise continuous function is Lebesgue integrable and its Lebesgue integral coincides with its Riemann integral.
\item\ZLA{I3ch1TeoInteRIMlebContSTEP} two functions which are a.e. equal on $R$ have the same Lebesgue integral. In particular 
\[
f=0 \ {\rm a.e.} \ x\in [h,k]\ \ \implies \ \ \int_{h}^k f(x)\ZD x=0\,.
\]
\end{enumerate}
\end{Theorem}

The following result is a simple consequence of the corresponding result which holds for the Riemann integral of continuous functions:
\begin{Theorem}\ZLA{Ch1teoProprieINTEboundFun}
Let $f $ and $g $ be a.e. defined on $[h,k]$, bounded and quasicontinuous. Then:
\begin{enumerate}
\item
\ZLA{I1Ch1teoProprieINTEboundFun} {\sc linearity of the integral:}\index{integral!linearity}\index{linearity of the integral}  if $\zaa$ e $\beta$ are real numer then
\[
\int_h^k\left (\zaa f(x)+\beta g(x)\right )\ZD x=
\zaa\int_h^k f(x)\ZD x+\beta\int_h^k g(x)\ZD x\,.
\]
\item\ZLA{I5Ch1teoProprieINTEboundFun} {\sc additivity of the integral:}\index{integral!additivity}\index{additivity of the integral} if \ $\hat  h \in(h,k)$ then
\[
\int_h^k f(x)\ZD x=\int_h^{\hat  h} f(x)\ZD x+\int _{\hat  h}^k f(x)\ZD x\,.
\]

\item\ZLA{I3Ch1teoProprieINTEboundFun} {\sc monotonicity of the integral:}\index{integral!monotonicity}\index{monotonicity of the integral}  if $f(x)\leq g(x)$ then
\[
\int_h^k f(x)\ZD x\leq \int_h^k g(x)\ZD x\,.
\]
\item\ZLA{I4Ch1teoProprieINTEboundFun} inequality for the absolute value:  we have
\[ 
\left |\int_h^k f(x)\ZD x\right |\leq \int_h^k |f(x)|\ZD x\,.
\]

\item\ZLA{I2Ch1teoProprieINTEboundFun} integrability of the product and of the quotient: the product $f(x)g(x)$ is integrable; the quotient $f(x)/g(x)$ is integrable provided that $|g(x)|>\zaa>0$.
 \end{enumerate}
\end{Theorem}

\begin{Remark}\ZLA{RemaCH1MonotonInteTietze}{\rm 
The monotonicity property requires the following comment. By definition, the value of the Lebesgue integral is the limit of the sequence on the right side of~(\ref{eq:Ch1DEFIlebeFlimitAT}) and it does not depend on the sequence that it is used. So, monotonicity easily follows if we use the associated sequences to $f$ and $g$ used in the proof 
of Theorem~\ref{teo:Ch1ESTEdaCHIUSO} since in this case $f(x)\geq g(x)$  implies $f_n(x)\geq g_n(x)$ (see the statement~\ref{I1Theo:prprieSucceConvDIM1} of Theorem~\ref{Theo:prprieSucceConvDIM1}).}\zdia
\end{Remark}

We note the following equality which holds for the \emph{Riemann integral:}

\begin{equation}\ZLA{eqINVAtransaBOUfunz}
\int_{h +\zaa}^{k+\zaa}f(x-\zaa)\ZD x=\int_h^k f(x)\ZD x\,.
\end{equation}
Passing to the limit of sequences of associated functions we see that \emph{the property~(\ref{eqINVAtransaBOUfunz}) holds also for the Lebesgue integral.} 

Equality~(\ref{eqINVAtransaBOUfunz}) is   the {\sc  translation invariance}\index{transaltion invariance!of the integral}     of the (Riemann or Lebesgue) integral.

\ZREF{
As stated in Theorem~\ref{ch1TeoInteRIMlebContSTEP},   continuous functions and   piecewise constant functions (on a bounded interval)  are both Riemann and  Lebesgue integrable and the two integrals have the same value. In fact, the Lebesgue integral \emph{extends} the Riemann integral since:
 
\begin{Theorem}\ZLA{teoCH1:RiemESTESOdaLebe} Every Riemann integrable function is bounded quasicontinuous, hence    Lebesgue integrable, and   the values of its Riemann and Lebesgue    integrals coincide. 
\end{Theorem}

The proof is   in Appendix~\ref{leBiNteSUbsCh1RIEMann} where we prove also:
\begin{Theorem}
A bounded function defined on a bounded interval is Riemann integrable if and only if the set of its points of discontinuity is a null set.
\end{Theorem}

}

\subsection[Unbounded Functions and Domains]{\ZLA{sec:CH1LebeGENEdom}The  Integral of  Unbounded Functions on   Unbounded Domains}
 
First we investigate the definition of the integral of a function $f$ which is a.e. defined on $\zzr$ and which can be unbounded. Then we consider the integral of $f$ on a set $A$ provided that $A$ has a suitable property.

 So, we proceed in two steps:
\begin{description}

\item[\bf Step~1:]

We put 
 
\[
f_+(x)=f_{+,0}(x)=\max\{f(x),0\}\,,\qquad f_-(x)=f_{-,0}(x) =\min\{f(x),0\} 
\]
so that $f_+(x)\geq 0$  and $f_-(x)\leq 0$ for every $x$. Then,
  when $N\geq 0$, $K\geq 0$ and $R>0$, we put
\begin{eqnarray*}&&
f_{+;\, (R,N)}(x)=\left\{\begin{array}{lll}
\min\{f_+(x)\,,\ N\}&{\rm if}& |x|<R\\
0&{\rm if}& |x|\geq R\,.
\end{array}\right.\\  
&&
f_{-;\, (R,-K)}(x)=\left\{\begin{array}{lll}
\max\{f_-(x)\,,\ -K\}&{\rm if}& |x|<R\\
0&{\rm if}& |x|\geq R\,.
\end{array}\right.
\end{eqnarray*}
\begin{description}
\item[\bf Step~1A:]
Since $f $ is   quasicontinuous,  the function $f_{+;\,(R,N)}$ is bounded quasicontinuous   for every $R$ and every $N$. We define
\begin{equation}\ZLA{preliEq1DefinTE}
\underbrace{\int_\zzr f_{+ }(x)\ZD x}_{\tiny \begin{tabular}{l}
    Lebesgue  \\    integral  \end{tabular}} =\lim _{\stackrel{R\to +\ZIN}{N\to +\ZIN}}\underbrace{\int _{-R}^R 
f_{+;\,(R,N)}(x)\ZD x}_{\tiny \begin{tabular}{l}
   Lebesgue  \\    integral  \end{tabular}} \,.
\end{equation}
The limit has to be computed with $(R,N)\in \zzn\times\zzn$ (i.e. one independent   from the other). It exists since $f_{+;\,(R,N)}\geq 0$ and it can be~$+\ZIN$.
\item[\bf Step~1B:]
Since $f $ is quasicontinuous, the function $f_{-;\,(R,-K)}$ is bounded quasicontinuous for every $R$ and every $K$. We define
\begin{equation}\ZLA{preliEq2DefinTE}
\underbrace{\int_\zzr f_{- }(x)\ZD x}_{\tiny \begin{tabular}{l}
   Lebesgue  \\    integral  \end{tabular}}=\lim _{\stackrel{R\to +\ZIN}{K\to +\ZIN}}\underbrace{\int _{-R}^R 
f_{-;\,(R,-K)}(x)\ZD x}_{\tiny \begin{tabular}{l}
   Lebesgue  \\    integral  \end{tabular}}\,.
\end{equation}
The limit has to be computed with $(R,K)\in \zzn\times\zzn$ (i.e. one independent   from the other) and it can be $-\ZIN$.

\item[\bf Step~1C:] the quasicontinuous function $f $ is 
{\sc integrable}\index{function!integrable!$1$ variable}
on $\zzr$
 when at least one of the functions $f_+$ or $f_-$ has \emph{finite} integral.

In this case we define
\begin{equation}\ZLA{preliEq3DefinTE}
\underbrace{\int_\zzr f(x)\ZD x}_{\tiny \begin{tabular}{l}
   Lebesgue  \\    integral  \end{tabular}}=\underbrace{\int_\zzr f_+(x)\ZD x}_{\tiny \begin{tabular}{l}
   Lebesgue  \\    integral  \end{tabular}}+\underbrace{\int_\zzr f_-(x)\ZD x}_{\tiny \begin{tabular}{l}
   Lebesgue  \\    integral  \end{tabular}}\,.
\end{equation}
The integral can be a number (when both the integrals of $f_+$ and of $f_-$ are numbers) or it can be $+\ZIN$ or it can be $-\ZIN$. 

If  the integral is a number,  then we say that $f $ is {\sc summable\footnote{\ZLA{footnoteROYDENch1}we advise the reader that this distinction between summable and integrable functions (introduced by Lebesgue in his thesis) is not used in every text. For example,    in~\cite[p.~73]{Royden66LibroANALreale} the term ``integrable'' is used to intend that the integral is finite.}\index{summable (function)}\index{function!summable}}.
\end{description}
\item[\bf Step~2:]   let $A\subseteq \zzr$ satisfy the following assumption:
\begin{Assumption}\ZLA{ASSUch1QuascontINSa}{\rm 
The characteristic function of the set $A$ is quasicontinuous.\zdia}
\end{Assumption}
When $f$ is quasicontinuous on $\zzr$ and the set $A$ satisfies the Assumption~\ref{ASSUch1QuascontINSa} then the product $f\charfun _A$ is quasicontinuous. 
With an abuse of notation, even if $f$ is solely a.e. defined on $A$ we put
  \begin{equation}\ZLA{EqPar1Ch1AbuNOtaprodu}
f(x)\charfun_A(x)=
\left\{\begin{array}{lll}
f(x) &{\rm if}& x\in A \\
0&{\rm if}& x\notin A 
\end{array}\right.  \qquad\mbox{(abuse of notations)}
  \end{equation}
  and:
  \begin{enumerate}
  \item if $f\charfun_A$ (defined as in~(\ref{EqPar1Ch1AbuNOtaprodu})) is quasicontinuous  on $\zzr $ we say that $f$ is a {\sc quasicontinuous function on $A$;}\index{function!quasicontinuous!on a set $A$}
  \item
if $f\charfun_A$ is integrable on $\zzr $ then we say that $f$ is {\sc integrable on the set A;}\index{function!integrable!on a set $A$}\index{integrable!on a set $A$}   if it is summable then we say that $f$ is {\sc summable on the set A;}\index{function!summable!on a set $A$}\index{summable!on a set $A$}   
\item if $f$ is integrable or summable on $A$ then we introduce the notation
\[
\int_A f(x)\ZD x=\int_{\zzr} f(x)\charfun_A(x)\ZD x\,.
\]
\end{enumerate}
\end{description} 
\begin{Remark}\ZLA{Ch1RemeCONerrore}{\rm We note:
\begin{enumerate}

\item
sets which do not satisfy Assumption~\ref{ASSUch1QuascontINSa} exist. An example is in the Appendix~\ref{AppeCH4InseNONLmisur}.
\item any open set satisfies Assumption~\ref{ASSUch1QuascontINSa}. A \emph{faulty proof} is as follows: let ${\mathcal{O}}=\cup _{n\geq 1}(a_n,b_n)$ and let the intervals be pairwise disjoint. The function $\charfun_{\mathcal{O}}$ is discontinuous at $a_n$ and at $b_n$. The set of the points $a_n$ and $b_n$ is numerable so that  $\charfun_{\mathcal{O}}$ is discontinuous on a  numerable set. Hence,   $\charfun_{\mathcal{O}}$ is quasicontinuous. \emph{We invite the reader to discuss the error in this argument.} The error is discussed in Remark~\ref{RemaCH4InsOpeFrontPOSI}  while a correct proof is in  Theorem~\ref{Teo:Cha1bisFunCHarAPERTIquasCont}.\zdia

\end{enumerate}
}\end{Remark}
The following properties of the integral 
 clearly holds:
\begin{Theorem}\ZLA{ch1TeoInteRIMlebContSTEPSUMMA}
Let $A$ be a set which satisfies the Assumption~\ref{ASSUch1QuascontINSa}
and let $f $ and $g $ be summable on $A$. Then:
\begin{enumerate}
\item {\sc linearity of the integral:}\index{linearity of the integral}\index{integral!linearity} if $\zaa$ e $\beta$ are real numbers the function $\zaa f+\beta g$ is summable and
\[
\int_A\left (\zaa f(x)+\beta g(x)\right )\ZD x=
\zaa\int_A f(x)\ZD x+\beta\int_h^k g(x)\ZD x\,.
\]
\item if $g$ is bounded then the product $fg$ is summable; if $1/g$ is bounded then the quotient $f/g$ is summable.
\item  {\sc monotonicity of the integral:}\index{integral!monotonicity}\index{monotonicity of the integral}  if $f(x)\leq g(x)$ then
\[
\int_A f(x)\ZD x\leq \int_A g(x)\ZD x\,.
\]
 \item\ZLA{I3ch1TeoInteRIMlebContSTEPSUMMA} two functions which are a.e. equal on $R$ have the same Lebesgue integral. In particular the integral of a function which is a.e. zero is equal to zero.

\end{enumerate}
\end{Theorem}

The previous properties correspond to properties which hold  also for the improper integral. Instead,  the absolute value has a new and an important property. 
\begin{Theorem}\ZLA{Teo:Ch1teoValassoEsommab}
  We have:
\begin{enumerate}
\item\ZLA{I1Teo:Ch1teoValassoEsommab}
 If $f$ is integrable then $|f|$ is integrable and
\[
\int_A |f(x)|\ZD x=\int_A f_+(x)\ZD x-\int_A f_-(x)\ZD x\,.
\]
So, the usual inequality of the absolute value holds:
\[
\left |
\int_A f(x)\ZD x
\right |\leq \int_A |f(x)|\ZD x\,.
\]
\item\ZLA{I2Teo:Ch1teoValassoEsommab}
Let the  function $f $ be quasicontinuous. It is summable if and only if $|f |$ is summable.
 \end{enumerate}
\end{Theorem}

We  stress the assumption that $f$ is quasicontinuous in the statement~\ref{I2Teo:Ch1teoValassoEsommab}. This assumption cannot be removed since we shall see in the Remark~\ref{RemaFINAppeCH4InseNONLmisur} of the Appendix~\ref{AppeCH4InseNONLmisur} the existence of functions which are not quasicontinuous but whose absolute value is constant.

The last statement of Theorem~\ref{Teo:Ch1teoValassoEsommab} is crucial in functional analysis, since it permits to define integral norms and the spaces $L^p$.

Finally we state the translation invariance and the additivity of the integral.
\begin{Theorem}\ZLA{CH1P1TeoADDIinteGEBEiTRANSLA}
The following properties hold:
\begin{enumerate}
\item\ZLA{enuI1CH1P1TeoADDIinteGEBEiTRANSLA} {\sc  translation invariance:}\index{transaltion invariance!of the integral}
let $A$ be a set which satisfies the Assumption~\ref{ASSUch1QuascontINSa} 
and let $f$ be defined on $A$ and quasicontinuous. Let $\zaa\in \zzr$ and 
\[
A+\zaa=\{x+\zaa\,,\ x\in A\}\,.
\] 
We have:
\begin{enumerate}
\item\ZLA{I11CH1P1TeoADDIinteGEBEiTRANSLA} the set $A+\zaa$ satisfies the Assumption~\ref{ASSUch1QuascontINSa} and the function $x\mapsto f(x-\zaa)$ (defined on $A+\zaa$) is quasicontinuous.
\item\ZLA{I12CH1P1TeoADDIinteGEBEiTRANSLA}
if $f$ is integrable on $A$ then we have:
\[
\int_{ A+\zaa} f(x-\zaa)\ZD x=\int_A f(x )\ZD x\,.
\]
\end{enumerate}
\item\ZLA{I2CH1P1TeoADDIinteGEBEiTRANSLA} {\sc additivity of the integral:}\index{integral!additivity}\index{additivity of the integral} let $A_1$ and $A_2$ be two \emph{disjoint sets} which satisfy Assumption~\ref{ASSUch1QuascontINSa} and let $f$ be defined on $A_1\cup A_2$. We assume that $f_{|_{A_1}}$ and $f_{|_{A_2}}$ are summable (hence also quasicontinuous). We have:
\begin{enumerate}
\item\ZLA{I21CH1P1TeoADDIinteGEBEiTRANSLA} the set $A_1\cup A_2$ satisfy Assumption~\ref{ASSUch1QuascontINSa} and $f$ is summable  on $A_1\cup A_2$.
\item\ZLA{I22H1P1TeoADDIinteGEBEiTRANSLA} we have
\[
\int _{A_1\cup A_2} f(x)\ZD x=\int _{A_1} f(x)\ZD x+\int _{A_2} f(x)\ZD x \quad \mbox{(recall: $A_1\cap A_2=\emptyset$)}\,.
\]

\end{enumerate}
\item\ZLA{I3CH1P1TeoADDIinteGEBEiTRANSLA} if $A_1\cap A_2\neq \emptyset$ then:
\begin{enumerate}
\item
 the sets $A_1\cap A_2$,   $A_1\setminus A_2$ and $A_2\setminus A_1$ (if nonempty) satisfy Assumption~\ref{ASSUch1QuascontINSa};
 \item
if $ f$ is defined on $A_1\cup A_2$ and $f_{|_{A_1}}$ and $f_{|_{A_2}}$ are summable then $f $ is summable on $A_1\cup A_2$ and
\[
\int _{A_1\cup A_2} f(x)\ZD x\leq \int _{A_1} f(x)\ZD x+\int _{A_2} f(x)\ZD x\,.
\]  
\end{enumerate}
\end{enumerate}
\end{Theorem}

Note that the statement~\ref{I3CH1P1TeoADDIinteGEBEiTRANSLA} is consequence of the additivity of the integral since
\[
A_1\cup A_2=(A_1\setminus A_2)\cup (A_1\cap A_2)\cup (A_2\setminus A_1)\quad \mbox{(disjoint union)} 
\]
and
\[
\charfun_{A_1\cap A_2}(x)=\charfun _{A_1}(x)\charfun _{A_2}(x)\,,\qquad \charfun _{A_1\setminus A_2}(x)=\max\{\charfun _{A_1}(x)-\charfun _{A_2}(x),0\}\,.
\]
We conclude with two simple observations where we use the notation
\[
f_{-,N}(x)  =\min\{f(x),N\} \,.
\]
The theorem on the limits of monotone functions shows the  following criterion of summability:
  \begin{Theorem} 
  The quasicontinuous function $f$ is summable if and only if the function
  \[
(N,R)\mapsto \int _{-R	}^R  |f | _{-,N}(x)\ZD x  
  \]
(defined for $N>0$ and $R>0$)  is bounded.
  \end{Theorem}
  
  The second observation is a consequence of the linearity and additivity of the integral:
  \begin{Theorem}
  Let $f\geq 0$ be summable.  For every $\ZEP>0$ there exists $N_\ZEP$
  and $R_\ZEP$ such that:
  \begin{itemize}
  \item if $N>N_\ZEP$ we have
 $
\int_A \left [ f(x)-f_{-,N}(x)\right ]\ZD x<\ZEP$.
  \item
  if $R>R_\ZEP$ we have
  $
\int _{|R|>R_\ZEP} f(x)\ZD x<\ZEP$.
  \end{itemize}
  \end{Theorem}

%
%

\subsubsection{Lebesgue Integral and Improper Integral}

Theorem~\ref{teoCH1:RiemESTESOdaLebe} (yet to be proved) shows that Riemann integral   is
extended and  superseded by the Lebesgue integral.   Instead, the Lebesgue integral for unbounded functions or on unbounded intervals does not supersede  the improper integral since there exist  even continuous  functions which are important for the applications and which admit improper integral but which are not Lebesgue summable. Important examples are the improper integral of the $\sinc$ function  and the Fresnel integrals:
\begin{align*}&
\underbrace{\int_{-\ZIN}^{+\ZIN}\sinc\pi x\ZD x}_{\tiny \begin{tabular}{l}
   improper  \\    integral  \end{tabular}}=\underbrace{\int_{-\ZIN}^{+\ZIN}\frac{\sin\pi x}{\pi x}\ZD x}_{\tiny \begin{tabular}{l}
   improper  \\    integral  \end{tabular}}=1\,,
 \\&
 \underbrace{\int_{-\ZIN}^{+\ZIN} \sin^2 x}_{\tiny \begin{tabular}{l}
   improper  \\    integral  \end{tabular}}\ZD x=\underbrace{\int_{-\ZIN}^{+\ZIN}\cos x^2\ZD x}_{\tiny \begin{tabular}{l}
   improper  \\    integral  \end{tabular}}=\sqrt{\dfrac{\pi}{2}}\,.
\end{align*}
 \section{\ZLA{subsEXCHlimINTgenCASE}Limits and the  Integral}
 Both the definition of the Riemann and the Lebesgue integrals have a common idea: first we single out a class of ``simple'' functions whose integral can be defined in an obvious way. Then we single out classes of functions which can be ``approximated'' with simple functions. Let $f$ be a function in this class and $\{s_n\}$ be a sequence of ``approximating'' functions. We prove that the integrals of the functions $s_n$ converge to a number which is chosen as the (Riemann or Lebesgue) theorem of $f$.
 
 In the case of the Riemann integral, the simple functions are piecewise constant and the Riemann integrable functions  $f$ have to be uniformly approximated by piecewise constant function. So, it is natural to expect that limits and integrals can be exchanged when ``convergence'' is uniform convergence.
 
 In the case of the Lebesgue integrals, the ``simple'' functions are continuous functions but the ``approximation'' is more general then pointwise convergence on an interval: it is ``approximation'' in the weaker sense of the existence of sequences of associated functions. So, it is natural to expect that limits and integrals can be exchanged under more general conditions. In fact we have the results which we state here and we prove in Chap.~\ref{Cha1bis:INTElebeUNAvariaBIS}. 
 
The first result shows that the request~\ref{I3ZREFrequiABSTRinteg} in the Table~\ref{TABLEch1RequestINTE} is satisfied by the Lebesgue interval. The proof is in Chap.~\ref{Cha1bis:INTElebeUNAvariaBIS}  but, for the convenience of the readers who do not plan to study Chap.~\ref{Cha1bis:INTElebeUNAvariaBIS}, at the end of this section we sketch the proof under  restrictive assumptions.

\begin{Theorem}\ZLA{TeoCH1P1ConveUnderBound}
Let $\{f_n\}$ be a bounded sequence of  quasicontinuous function defined on a bounded interval $[h,k]$. 

Let $\{f_n(x)\}$ converges to $f(x)$ a.e. on $[h,k]$. Then the function $f$ is bounded quasicontinuous and we have

\[
\limn \left [\int_h^k f_n(x)\ZD x\right ]=\int_h^k \left [\limn f_n(x)\right ]\ZD x
=
\int_h^k f(x)\ZD x
\,.
\]
\end{Theorem}

 The heavy boundedness assumptions in this theorem can be much relaxed and in fact we have:

\begin{Theorem}[{\sc  Beppo Levi} or {\sc   monotone convergence}]\index{Beppo Levi Theorem}\index{Theorem!Beppo Levi} 
\index{monotone convergence theorem}\index{Theorem!monotone convergence}\ZLA{teo:ch1Beppo}
Let $\{f_n \}$ be a sequence of integrable functions defined a.e. on $\zzr$
and let us assume that the following properties hold  a.e. on $\zzr$:
\begin{eqnarray*}
{\bf 1)}\quad 0\leq f_n(x)\leq f_{n+1}(x)\,, && {\bf 2)}\quad \lim _{n\to +\ZIN} f_n(x)=f(x)\qquad a.e.\ x\in\zzr\,.
\end{eqnarray*}
Then, $f $ is integrable  and
\[
\lim_{n\to+\ZIN} \int_{-\ZIN}^{+\ZIN} f_n(x)\ZD x=\int_{-\ZIN}^{+\ZIN} f(x)\ZD x\,.
\]
\end{Theorem} 
 Note that in this theorem the functions $f_n$ need not be summable. The assumption is that they are integrable. And the conclusion is that $f$ is integrable, possibly not summable.  The function $f$ may not be summable even if each one of the functions $f_n$ is.
 
 Concerning summability we have:
 
\begin{Theorem}[{\sc Lebesgue} or {\sc     dominated convergence}]\index{Theorem!Lebesgue!}\index{Theorem!dominated convergence}\index{dominated convergence} \ZLA{teo:CH1TeoLebe}
Let $\{f_n \}$ be a sequence of summable functions a.e. defined on $\zzr$ and let $f_n\to f  $ a.e. on $\zzr$. If there exists a summable nonnegative function $g $ such that
\[
|f_n(x)|\leq g(x)\qquad a.e.\ x\in\zzr
\]
then $f(x)$ is summable and
\[
\lim \int_{-\ZIN}^{+\ZIN} f_n(x)\ZD x=\int_{-\ZIN}^{+\ZIN} f(x)\ZD x\,.
\]
\end{Theorem}

 \paragraph{Sketch of the proof of Theorem~\ref{TeoCH1P1ConveUnderBound}}
We assume that the assumptions of the theorem hold and furthermore that a seemingly   restrictive assumption holds too. So we assume:  
\begin{enumerate}
\item\ZLA{I1PARAproofteo:ch1:LebeCAROPartIntervLIMIT}
each function $f_n$ is a.e. defined on $[h,k]$ and it is quasicontinuous.
\item there exists $M$ such that $|f_n(x)|<M$ a.e. $x\in [h,k]$ and for every $n$.
\item\ZLA{I1PARAproofteo:ch1:LebeCAROPartIntervLIMITbis} the sequence converges a.e. $x\in [h,k]$, $f_n(x)\to f(x)$.
 \item\ZLA{I2PARAproofteo:ch1:LebeCAROPartIntervLIMIT}
\emph{for every $\ZEP>0$ there exists a multiinterval $\Delta_\ZEP$  such that $L(\Delta_\ZEP)<\ZEP$ and such that every $f_n$ is defined on $[h,k]\setminus\mathcal{I}_{\Delta_\ZEP}$
and furthermore:
\begin{enumerate}
\item   the restriction $(f_n)_{|_{[h,k]\setminus\mathcal{I}_{\Delta_\ZEP}}}$ is continuous.
\item the sequence
 $\{ f_n  \}$ is uniformly convergent on  $[h,k]\setminus\mathcal{I}_{ \Delta_\ZEP}$. 
 \end{enumerate}
 }
\end{enumerate}

 {Assumption~\ref{I2PARAproofteo:ch1:LebeCAROPartIntervLIMIT} looks unduly restrictive but Egorov-Severini Theorem (Theorem~\ref{Cha1bisTeoEgoSEVER} of Chap.~\ref{Cha1bis:INTElebeUNAvariaBIS}) shows that it is a consequence of pointwise convergence.}
 
The key points of the proof of Theorem~\ref{TeoCH1P1ConveUnderBound}   are the following ones.
\begin{description}
\item[{\bf The limit function $f$ is quasicontinuous.}]  We use the assumption~\ref{I2PARAproofteo:ch1:LebeCAROPartIntervLIMIT}: we fix $\ZEP>0$ and a corresponding multiinterval $\Delta_\ZEP$.  
In order to simplify the following notations, we put
\[
H_\ZEP=[h,k]\setminus \mathcal{I}_{\Delta_\ZEP }\quad\mbox{(note that $H_\ZEP$ is closed)}\,.
\]

We associate to the restriction of $f_n$ to $H_\ZEP$ the   particular Tietze extension constructed in Theorem~\ref{teo:Ch1ESTEdaCHIUSO}, via linear interpolation. This particular Tietze extension is denoted $ f_{n,\ZEP,e}$.

Now we use the following steps:
\begin{enumerate}
\item Assumption~\ref{I2PARAproofteo:ch1:LebeCAROPartIntervLIMIT} and  the statement~\ref{I2Theo:prprieSucceConvDIM1} of Theorem~\ref{Theo:prprieSucceConvDIM1} 
imply that the sequence $n\mapsto f_{n,\ZEP,e} (x)$ is convergent for every   $x\in [h,k]$. Let $\hat f$ be the limit function.

\item Assumption~\ref{I2PARAproofteo:ch1:LebeCAROPartIntervLIMIT} implies that the sequence $\{(f_n)_{|_{H_\ZEP}}\}$ is   uniformly convergent  and it is a sequence of continuous functions. So, its limit $\hat f_{|_{H_\ZEP}}$ is continuous.
\item from $(f_n)_{|_{H_\ZEP}}(x)=f_n(x)$ and $f_n(x)\to f(x)$ for every $x\in H_\ZEP$, it follows that $\hat f(x)=f(x)$ if $x\in H_\ZEP$.

So, $f_{|_{H_\ZEP}}$ is continuous 
and $\hat f$ is one of its Tietze extensions\footnote{by examining the limit process and the construction of $f_{n,\ZEP,e}$ we can see that $\hat f$ is precisely the Tietze extension of $f$, constructed with the method described in the statement of Theorem~\ref{teo:Ch1ESTEdaCHIUSO}.}.

This argument holds for every $\ZEP>0$ and so $f$ is quasicontinuous.

\end{enumerate}

\item[{\bf We prove convergence of the integrals.}] Once quasicontinuity of $f$ is known we can replace $f_n$ with $f_n-f$ and we can assume $f_n\to 0$. 

We prove that if $f_n\to 0$ a.e. then the sequence of the integrals converges to zero.

Assumption~\ref{I2PARAproofteo:ch1:LebeCAROPartIntervLIMIT} and Theorem~\ref{TheoCH1ConveTieSeq} imply that 
$\{ f_{n,\ZEP,e}\}$ converges \emph{uniformly} to~$0$.
\item We note 
 \[
\underbrace{\int_{h }^k f_n(x)\ZD x}_{\tiny \begin{tabular}{l}
     Lebesgue       integral  \end{tabular}}
     =\overbrace{\underbrace{\int_h^k f_{n,\ZEP,e}(x)\ZD x}_{A}}
^{\tiny \begin{tabular}{l}
    Riemann      integral  \end{tabular}}
+\overbrace{\underbrace{\int _h^k\left [ f_n(x)- f_{n,\ZEP,e}(x)\right ]\ZD x}_{B}}
^{\tiny \begin{tabular}{l}
    Lebesgue    integral  \end{tabular}}\,.
\]
The difference $\left [ f_n(x)- f_{n,\ZEP,e}(x)\right ]$ is zero when $x\notin \mathcal{I}_{\Delta_\ZEP}$. This fact implies\footnote{the proof is in Theorem~\ref{Teo:Cha1bisABSOluContiPERiSoliAperti} of Chap.~\ref{Cha1bis:INTElebeUNAvariaBIS}.}
\[
|B|\leq 2M\ZEP\,.
\]
The integrals $A$ are Riemann integrals and the integrands converge uniformly to zero. So, $|A|$ can be made as small as we wish by taking $n$ large. 

This implies that for every $\ZSI>0$ and every $\ZEP>0$ there exists $N=N_{\ZSI,\ZEP}$ such that if  $n>N_{\ZSI,\ZEP}$ we have
\[
\left |\int_{h }^k f_n(x)\ZD x\right |\leq |A|+|B|<\ZSI+\ZEP\,,
\]
as wanted.
\end{description}
 
 \section{\ZLA{AppeTOch1IntegrRoemLebe}\textbf{Appendix:}  Riemann  and Lebesgue Integrals}
 
   In this appendix we sketch the definition of the Riemann integral and we
characterize Riemann integrability of a function in terms of its points of discontinuity. Then we prove that every Riemann integrable function is Lebesgue integrable  and that the two integrals have the same value.

\subsection{\ZLA{APPEch1SKETCH}A Sketch of the Riemann Integral}
Functions which are Riemann integrable must be bounded and defined on a bounded interval which we denote $[h,k]$ but the fact that the end points belong to the interval has no effect on integrability of the function or on the value of its Riemann integral.

 We sketch the definition of the Riemann integral and the relation of the Riemann integral with the oscillation of the function.  

\subsubsection{\ZLA{secCHappe1OscillFUnz}The Oscillation of a Function}
We define
the  {\sc oscillation}\index{oscillation} 
 of a function on a set and at   a point. 

Let $f$ be a real valued function defined on a subset  $A\subseteq \zzr$.
Let $U$ be a set which intersects $A$:  $A\cap U\neq \emptyset$. The {\sc oscillation}\index{oscillation!of a function on a set} of $f$ on $U$ is
\[
\ZOMq(f,U)=\left [ \sup _{ x\in A\cap U} f(x)-\inf  _{x\in A\cap U} f(x)\right ]\,.
\]

We fix a point $x_0\in{\rm cl}\, A$ and we consider the intervals $B(x_0,r)=(x_0-r,x_0+r)$ i.e. the neighborhood of $x_0$ of radius $r$.
The function
\[
r\mapsto \ZOMq(f,B(x_0,r))
\]
is increasing, hence we can define the {\sc oscillation of $f$ at $x_0$}\index{oscillation!of a function at a point} as
\[
\ZOM(f;x_0)=\lim_{r\to 0^+}\ZOMq\left (f,B(x_0,r)\right )\,.
\]

The oscillation can be used to characterize continuity:
\begin{Theorem}
Let $x_0\in A=\Dom\, f$. The function $f$ is continuous at $x_0$ if and only if 
$\ZOM(f;x_0)=0$.
\end{Theorem}
\zProof Let $f$ be continuous at $x_0$ and let $\ZEP>0$. There exists $N$ such that for every $n>N$ the condition $|x-x_0|<1/n$  implies $|f(x)-f(x_0)|<\ZEP$. So,
\[
|x-x_0|<\dfrac{1}{n}\,,\quad |y-x_0|<\dfrac{1}{n}\ \implies\ |f(x)-f(y)|<2\ZEP
\]
and so  
\[
n>N_\ZEP \ \implies\ 
\ZOMq\left (f,B(x_0,1/n)\right )<2\ZEP\,.
\]
Then we have $\ZOM(f;x_0)=0$ since $\ZEP>0$ is arbitrary.

Conversely, let $ \ZOM(f;x_0)=0$. Then, for every $\ZEP>0$ there exists $\ZDE>0$ such that 
\begin{align*}&
0<r<\ZDE\ \implies\ 
\ZOMq\left (f,B(x_0,r)\right )<\ZEP\,. \\
&\mbox{In particular if $x\in B(x_0,r)$ then}\quad |f(x)-f(x_0)|\leq\ZOMq(f,B(x_0,r))<\ZEP\,.\zdiaform
\end{align*}

 We shall use the following result:
\begin{Theorem}\ZLA{APPE1TheoremCompaAzaa}Let
$f$ be any function defined on $[h,k]$ and let
$\zaa\geq 0$. The set   
\begin{equation}\ZLA{Eq:Ch1AppeDefiAfalpha}
A _\zaa =A _{f,\zaa}=\{x\in [h,k]\,:\ \ZOM(f;x)\geq \zaa\}  
\end{equation}
is closed.
\end{Theorem}
\zProof
 We prove that if $\hat x$ is an accumulation point of $A_\zaa $ then it  belongs to $A_ \zaa $ i.e. we prove that  
\[
\ZOM(f;\hat x)\geq \zaa\,.
\]
 The neighborhood $B(\hat x,1/n)$ contains a point $\tilde x\in A_\zaa$ and it contains also  a neighborhood  $U$ of $\tilde x$. So,
\begin{multline*}
\ZOMq\left (f,B(\hat x,1/n)\right )=
\left [ \sup _{x\in B(\hat x,1/n)  } f(x)-\inf  _{x\in B(\hat x,1/n) } f(x)\right ]\\
\geq 
\left [ \sup _{x\in U } f(x)-\inf  _{x\in U } f(x)\right ]\geq \ZOM(f; \tilde x)\geq \zaa\,.
\end{multline*}
The inequality is preserved by the limit so that
\[
\ZOM(f;\hat x)=\limn\left [ \sup _{x\in  B(\hat x,1/n)  } f(x)-\inf  _{x\in B(\hat x,1/n) } f(x)\right ]\geq\zaa\,.\zdiaform
\]

\subsubsection{Riemann Integral:  the Conditions of Integrability}
Riemann integral of the bounded function $f$ defined on the bounded interval $[h,k]$ can be defined with several slightly different but equivalent methods. We sketch one.

In the contest of Riemann integral, a {\sc partition}\index{partition} of $[h,k]$
is a finite set $\mathcal{P}=\{h_i\}_{0\leq i\leq N} $ such that
\[
 h_0=h\,,\quad   h_i<h_{i+1} \,,\quad h_N=k\,.
\]
We associate to the partition $\mathcal{P}$ the following numbers
\[
 \begin{array}{ll}
 \displaystyle
 
  l_i= h_{i+1}-h_i\,,&  \ZDE(\mathcal{P})=\max\{l_i\}
 \\[2mm]
 \displaystyle m_i=m_{i,\cal{P} }=\inf_{x\in [h_i,h_{i+1})} f(x)\,,& 
  M_i=M_{i,\mathcal{P}}=\sup_{x\in [h_i,h_{i+1})} f(x)\,.
  \end{array}
 \]
Then we introduce the two piecewise constant functions $\chi_ {+,\mathcal{P}}$ and  $\chi_{-,\mathcal{P}}$
\[
\mbox{if $x\in [h_i,h_{i+1})$ then}\  \left\{\begin{array}{l}
\chi_{+,\mathcal{P}}(x)=M_i\,, \\
\chi_{-,\mathcal{P}}(x)= m_i \,.
\end{array}\right.
\]
As seen in Example~\ref{Ese:TWOeseInteCosTraDiri}, their  \emph{Lebesgue} integrals are the numbers

\begin{equation}\ZLA{eqAppe1DefiInteLebeCHIpm}
\begin{array}{ll}\displaystyle
I_+(f,\mathcal{P})&=\underbrace{\int_h^k \chi_{+,\mathcal{P}}(x)\ZD x}_{\tiny \begin{tabular}{l}
    Lebesgue\\ integral  \end{tabular}}=\sum _{i=0}^{N-1} M_i (h_{i+1}-h_i)\,,\\[3mm]
    \displaystyle
I_-(f,\mathcal{P})&
=\underbrace{\int_h^k \chi_{-,\mathcal{P}}(x)\ZD x}_{\tiny \begin{tabular}{l}
    Lebesgue\\ integral  \end{tabular}}
=\sum _{i=0}^{N-1} m_i (h_{i+1}-h_i)\,.
\end{array}
\end{equation}
  It is clear that $I_-(f,\mathcal{P})\leq I_+(f,\mathcal{P})$
and in fact we have also
\[
I_-(f,\mathcal{P})\leq I_+(f,\mathcal{Q})
\]
even if $\mathcal{P}$ and $\mathcal{Q}$ are different partitions.

By definition, the function $f$ is {\sc Riemann integrable}\index{function!integrable!Riemann} when
\[
\sup\{I_-(f,\mathcal{P})\}=\inf\{I_+(f,\mathcal{P})\}
\]
 and this number is its {\sc Riemann integral:}\index{integral!Riemann}
 \[
\underbrace{\int_h^k f(x)\ZD x}_{\tiny \begin{tabular}{l}
    Riemann  \\    integral  \end{tabular}}=\sup\{I_-(f,\mathcal{P})\}=\inf\{I_+(f,\mathcal{P})\}\,.
 \]
 
\emph{It follows that every piecewise constant function is Riemann integrable and that its Lebesgue and Riemann integrals coincide\footnote{as stated in Example~\ref{Ese:TWOeseInteCosTraDiri} and in the statement~\ref{I2ch1TeoInteRIMlebContSTEP} of Theorem~\ref{ch1TeoInteRIMlebContSTEP}.}.} This observations hold in particular for the functions  $\chi_ {\pm,\mathcal{P}}$: \emph{the integrals in~(\ref{eqAppe1DefiInteLebeCHIpm}) are both Lebesgue and Riemann integrals.}
 
The following result is known:
\begin{Theorem}
\ZLA{TheoCH1AppePrimaCharactINTEGrIem}
Let $f$ be a bounded function defined on $[h,k]$. The function $f$
is Riemann integrable if and only if for every $\ZSI>0$ there exists a partition $\mathcal{P}$ of $[h,k]$ such that 
\begin{equation}\ZLA{eq:Ap1LEMMAinteRunaPart}
   I_+(f,\mathcal{P})-I_-(f,\mathcal{P})<\ZSI\,.
 \end{equation} 
 If inequality~(\ref{eq:Ap1LEMMAinteRunaPart}) holds for a partition $\mathcal{P}$ then it holds also for any partition $\mathcal{Q}\supseteq\mathcal{P}$. 
\end{Theorem} 
So:
\begin{Corollary}
Let $f$ be a bounded function defined on $[h,k]$. The function $f$
is Riemann integrable if and only if  there exists a sequence $\{\mathcal{P}_n\}$ of partitions of $[h,k]$ such that  
\[
\ZDE(\mathcal{P}_n)\to 0\,,\qquad   I_+(f,\mathcal{P}_n)-I_-(f,\mathcal{P}_n)\to 0\,.
\]
In this case
\begin{equation}
\ZLA{EQ:I2AppeCa1TheoDaPartApartEqui}
\underbrace{\int_h^k f(x)\ZD x}_{\tiny \begin{tabular}{l}
    Riemann  \\    integral  \end{tabular}}=\left\{\begin{array}{l}
\lim _{n\to +\ZIN} I_+(f,\mathcal{P}_n)=    
    \lim _{n\to +\ZIN} \underbrace{\int_h^k  \chi
     _{+, \mathcal{P}_n}
    (x)\ZD x}_{\tiny \begin{tabular}{l}
    Riemann and Lebesgue \\    integral  \end{tabular}}\\
 \lim _{n\to +\ZIN} I_-(f,\mathcal{P}_n)=   \lim _{n\to+\ZIN}\underbrace{\int_h^k  \chi
  _{-, \mathcal{P}_n}
 (x)\ZD x}_{\tiny \begin{tabular}{l}
    Riemann and Lebesgue  \\    integral  \end{tabular}}\,.
    \end{array}\right.
\end{equation}
\end{Corollary}
\begin{Remark}{\rm We stress the fact that the integrals on the right   of~(\ref{EQ:I2AppeCa1TheoDaPartApartEqui}) are both Riemann and Lebesgue integrals.\zdia}
\end{Remark}

A fact that has its interest is that the partitions $\mathcal{P}_n$ in~(\ref{EQ:I2AppeCa1TheoDaPartApartEqui}) can be taken composed by equispaced points: for every $n$ we can take a partition composed by the  $N=n+1$ equispaced points
\begin{equation}\ZLA{eq:ANTEnewAppeCa1TheoDaPartApartEqui}
h_i=h_{i,N}=h+\dfrac{i}{N}(k-h)\,,\quad 0\leq i\leq  N\,.
\end{equation}
 We have:
\begin{Corollary} \ZLA{COROAppeCa1TheoDaPartApartEquiBIS}
Let $f$ be Riemann integrable on $[h,k]$ and let $h_i=h_{i,N}$ be as defined in~(\ref{eq:ANTEnewAppeCa1TheoDaPartApartEqui}). Let $ x_{i,N}\in [h_{i,N},h_{i+1,N})$ be arbitrarily chosen. We have
\[
\underbrace{\int_h^k f(x)\ZD x}_{\tiny \begin{tabular}{l}
    Riemann  \\    integral  \end{tabular}}=
    \lim _{N\to+\ZIN} \sum _{i=0} ^{N-1} f(x_{i,N})(h_{i+1,N}-h_{i,N})\,.
\]
\end{Corollary}

\subsubsection{Riemann Integrability and Continuity} 
We noted in Sect.~\ref{secCHappe1OscillFUnz}   that continuity of a function can be formulated in terms of its oscillation and
Riemann integral has a relation with the oscillation of the function since
\[
M_i-m_i=\ZOMq(f,[h_i,h_{i+1}))
\,.
\]
We use this observation and we  reformulate Theorem~\ref{TheoCH1AppePrimaCharactINTEGrIem} as follows:

 \begin{Theorem}\ZLA{Ap1THEOinteRunaPart}
 The bounded function
 $f$, defined on $[h,k]$, is Riemann integrable if and only if for every $\ZSI>0$ there exists a partition $\mathcal{P}$ of $[h,k]$ such that
\begin{equation}\ZLA{eq:Appe1RelaRiemmEOSCILLA} 
 I_+(f,\mathcal{P})-I_-(f,\mathcal{P})=
 \sum _{i=0}^{N-1} \ZOMq\left (f,[h_i,h_{i+1})\right )(h_{i+1}-h_i)<\ZSI\,.
\end{equation}
 If $f$ is Riemann integrable then for every $\ZDE>0$ there exist partitions $\mathcal{P}$ for which~(\ref{eq:Appe1RelaRiemmEOSCILLA}) holds and such that $\ZDE(\mathcal{P})<\ZDE$.
 \end{Theorem}

This result suggests that we study in more details the relation of integrability and the oscillation of the function. 
 
 The integrability test in  Theorem~\ref{Ap1THEOinteRunaPart}
can be reformulated as in the   statement~\ref{I2ch1AppeLemmaIntegrRiemannOsci} of the following theorem.This reformulation of the test of integrability is often called the ``Dini test'' after its proof in~\cite[p.~242]{VDINIteoricaORIG} but it had already been stated by Riemann.  

We use the following notation. Let $\mathcal{P}$ be a partition of $[h,k]$ and let $\zaa>0$. We put
\[
\mathbf{I}_{\zaa,+}=\{ i\,:\ \ZOMq(f,[h_i,h_{i+1}))\geq \zaa\}\,,\qquad
\mathbf{I}_{\zaa,-}=\{ i\,:\ \ZOMq(f,[h_i,h_{i+1}))<\zaa\}\,.
\]

\begin{Theorem} \ZLA{ch1AppeLemmaIntegrRiemannOsci}
Let $f $ be a   function defined on $[h,k]$.
The following properties are equivalent:
\begin{enumerate}
\item\ZLA{I1ch1AppeLemmaIntegrRiemannOsci} the function $f$ is Riemann integrable.
\item\ZLA{I2ch1AppeLemmaIntegrRiemannOsci} the function $f$ is bounded and the following  {\sc  Dini  test}\index{Dini test}  holds: for every    $\ZEP>0$ and $\zaa>0$  there exists a partition $\mathcal{P}$ of $[h,k]$ such that
\begin{equation}\ZLA{eq:ch1AppeLemmaIntegrRiemannOsci}
\sum _{i\in\mathbf{I}_{\zaa,+} } (h_{i+1}-h_i)<\ZEP\,.
\end{equation}
\end{enumerate}  
\end{Theorem}
\zProof We prove that 
property~\ref{I1ch1AppeLemmaIntegrRiemannOsci} implies property~\ref{I2ch1AppeLemmaIntegrRiemannOsci}.  
We assign 
$\ZEP $ and $\zaa$ (both positive). 
We apply Theorem~\ref{Ap1THEOinteRunaPart}  with $\ZSI= \ZEP\zaa $ : there exists $\mathcal{P}$ such that  
\begin{multline*}
  \ZEP{\color{red}\cancel{\zaa}}\geq  \sum _{i=0}^{N-1} \ZOMq\left (f,[h_i,h_{i+1})\right )(h_{i+1}-h_i)\\
= 
  \sum _{i\in \mathbb{I}_{\zaa,+}}  \underbrace{\ZOMq\left (f,[h_i,h_{i+1})\right )}_{{\color{red}\geq \zaa}}(h_{i+1}-h_i)
+
 \sum _{i\in \mathbb{I}_{\zaa,-}}  
  \underbrace{\ZOMq\left (f,[h_i,h_{i+1})\right )}_{\geq 0}   \underbrace{(h_{i+1}-h_i)}_{>0}   
  \\
  \geq 
{\color{red}\cancel{\zaa}}\left [   \sum _{i\in\mathbf{I}_{\zaa,+}}  (h_{i+1}-h_i)\right ]
\,.
\end{multline*} 
So we have~(\ref{eq:ch1AppeLemmaIntegrRiemannOsci}).

We prove the opposite implication. We prove that the property in statement~\ref{I2ch1AppeLemmaIntegrRiemannOsci} implies the integrability condition in 
Theorem~\ref{Ap1THEOinteRunaPart}. 
We use boundedness of $f$, $|f(x)|<M$ for every $x\in [h,k]$. 
We fix any $\ZSI>0$ and we apply the property in the statement~\ref{I2ch1AppeLemmaIntegrRiemannOsci} with $\ZEP=\ZSI/4M$ and $\zaa=\ZSI/(2(k-h))$. Let $\mathcal{P}$ be a partition for which the inequality~(\ref{eq:ch1AppeLemmaIntegrRiemannOsci}) holds with these values of $\ZEP$ and $\zaa$. We have:
\begin{multline*}
 \sum _{i=0}^{N-1} \ZOMq\left (f,[h_i,h_{i+1})\right )(h_{i+1}-h_i)\\
 =
  \sum _{i\in \mathbf{I}_{\zaa,+}}  \underbrace{\ZOMq\left (f,[h_i,h_{i+1})\right )}_{<2M}(h_{i+1}-h_i)+
   \sum _{i\in \mathbf{I}_{\zaa,-}}  \overbrace{\ZOMq\left (f,[h_i,h_{i+1})\right )}^{< \zaa}(h_{i+1}-h_i)\\
   \leq 2M\ZEP+\zaa(k-h)<\ZSI.\zdiaform
\end{multline*}

 This result shows a relation between integrability and continuity: if 
 $x_0\in (h_i,h_{i+1})$   and if $\ZOM(f,x_0)>\zaa$ then   $i\in \mathbf{I}_{\zaa,+}$. More precisely, the next Theorem~\ref{TeoCH1AppeCONcrtDuBoisReimond} holds.

 \begin{Theorem}
 \ZLA{TeoCH1AppeCONcrtDuBoisReimond}
Let $f$ be defined on $[h,k]$. 
The following properties are equivalent:
\begin{enumerate}
\item \ZLA{I1TeoCH1AppeCONcrtDuBoisReimond}
the function is Riemann integrable 
\item \ZLA{I2TeoCH1AppeCONcrtDuBoisReimond} the function is bounded and 
the following version of Dini test holds: 
  for every $\zaa>0$ and $\ZEP>0$ there exists a finite disjoint multiinterval $\Delta=\{(a_i,b_i)\}_{1\leq i\leq K}$ such that
\[
  L(\Delta)<\ZEP\,,\qquad 
A_\zaa=\{x\,:\ \ZOM(f;x)\geq \zaa\} \subseteq  \mathcal{I}_\Delta  \,.
\] 
\end{enumerate}
 \end{Theorem}
\zProof The fact that Riemann integrability implies the property in the statement~\ref{I2TeoCH1AppeCONcrtDuBoisReimond} is easily seen from      Dini test. In fact, Dini test implies that  
\[
A_\zaa \subseteq \underbrace{\left [\bigcup_{i\in\mathbf{I}_{\zaa,+}}[h_i,h_{i+1})\right ]}_{\footnotesize\mbox{(sum of  lengths) $<\ZEP$}  }\bigcup\mathcal{P}\quad \mbox{($\mathcal{P}$ is a finite set, hence a null set)}\,.
\]
The multiinterval $\Delta$ is obtained by  slightly enlarging $[h_i,h_{i+1})$ to $(h_i-\ZDE,h_{i+1})$ with $\ZDE$ so small that the 
sum of the lengths is still less then $\ZEP$; and then by covering the endpoints $h_j$ where $\ZOM(f,h_j)\geq \zaa$, if not yet covered, with ``small'' open intervals.

Conversely, we prove that the property in the statement~\ref{I2TeoCH1AppeCONcrtDuBoisReimond} implies Dini test.

Let $\Delta$ be the finite disjoint multiinterval in the statement~\ref{I2TeoCH1AppeCONcrtDuBoisReimond}.
We have
\[
[h,k]\setminus\mathcal{I}_\Delta\subseteq \{x\,:\ \ZOM(f;x)<\zaa\}
\,.
\]

Let $x_0\in  \{x\,:\ \ZOM(f;x)<\zaa\}$. There exists an open interval $I_{x_0} $ such that
\[
\ZOMq(f,I_{x_0})<\zaa
\]
and $[h,k]\setminus\mathcal{I}_\Delta$ is a compact set such that
\[
[h,k]\setminus\mathcal{I}_\Delta\subseteq \bigcup_{x_0\in \tilde  A_{\zaa}} I_{x_0}\,.
\]
Compactness of $[h,k]\setminus\mathcal{I}_\Delta$ shows the existence of a finite number of intervals $I_{x_1}$,\dots, $I_{x_K}$ such that  
\[
[h,k]\setminus\mathcal{I}_\Delta\subseteq \bigcup _{1\leq j\leq K}I_{x_j}\,,\qquad 
\ZOMq(f,I_{x_j})<\zaa\,.
\]
By reordering the endpoints   $a_i$ and $b_i$ of the intervals which compose $\Delta$ and those of the intervals $I_{x_j}$ we get a partition of $[h,k]$ with the following property: if $[h_i,h_{i+1})$ is an interval determined by the partition then
\[
\ZOMq(f,[h_i,h_{i+1}))\geq \zaa\ \implies\ 
 [h_i,h_{i+1}) \subseteq \cup (a_i,b_i)
\]
and so
\[
\sum _{i\in\mathbf{I}_{\zaa,+} } (h_{i+1}-h_i)\leq \sum (b_i-a_i)<\ZEP\,.
\]
So, Dini test holds.\zdia

The property in the statement~\ref{I2TeoCH1AppeCONcrtDuBoisReimond} of Theorem~\ref{TeoCH1AppeCONcrtDuBoisReimond} implies:
 \begin{Corollary}
 If the function $f$ is Riemann integrable then 
 \[
 A_{ \zaa}  =\{x\,:\ \ZOM(f;x)\geq \zaa\}
\] 
  is a null set for every $\zaa$.
 \end{Corollary}

The main result of this section is the following one,   independently proved in~\cite {VITALIteoVITALIlebesgueGIOENIA,VITALIteoVITALIlebesgue} by Vitali  and in~\cite{Lebesgue04LibroPrimiRes} by Lebesgue:
 \begin{Theorem}[\sc Vitali-Lebesgue]\index{Theorem!Vitali-Lebesgue}\ZLA{Theo:Appe1VitaliLebe}
The bounded function $f$ defined on $[h,k]$ is Riemann integrable   if and only if it   is a.e. continuous. 
 \end{Theorem} 
\zProof The set of the points of discontinuity is the set
\[
\bigcup _{\zaa>0} A_{f,\zaa} =\bigcup _{n\geq 1}A_{f,1/n}\,.
\] 
We use statement~\ref{I2TeoCH1AppeCONcrtDuBoisReimond} of Theorem~\ref{TeoCH1AppeCONcrtDuBoisReimond}:
if $f$ is Riemann integrable then 
 each one of the sets $A_{f,1/n}$ is a null set, and their union is a null set too.
 
Conversely, let $f$ be a.e. continuous. We prove that
the property in the statement~\ref{I2TeoCH1AppeCONcrtDuBoisReimond} of Theorem~\ref{TeoCH1AppeCONcrtDuBoisReimond} holds.

The set
 $\{x\,:\ \ZOM(f;x)>0\}=\cup _{n\geq 1}A_{f,1/n}$ is a null set. It follows that $A_{f,\zaa}$ is a null set for every $\zaa>0$. 
So, $A_{f,\zaa}$     is a null set which is compact (see Theorem~\ref{APPE1TheoremCompaAzaa}):
 for every $\ZEP>0$ there exists a \emph{finite} and disjoint  multiinterval $\tilde \Delta$ such that\footnote{the fact that the multiinterval is composed by finitely many intervals is proved in Lemma~\ref{I2Exe:Ch0Qnullset}. The multiinterval can be chosen disjoint since if two open intervals intersect their union is still an open interval whose length is less then the sum of the two lengths.} 
\[
\begin{array}{l}\displaystyle
\tilde\Delta=\{ (\tilde a_1,\tilde b_1)\,, (\tilde a_2,\tilde b_2)\,, \cdots \,,(\tilde a_K,\tilde b_K)\}\,,   \quad  L(\tilde\Delta)<\ZEP\\
\displaystyle A_\zaa\subseteq  \mathcal{I}_{\tilde\Delta}\,.
\end{array}
\]
So, the integrability  test is verified.\zdia
%
%
%
%
%
%
%
%
%
%

In order to appreciate Theorem~\ref{Theo:Appe1VitaliLebe} the student should keep in mind that there exist  functions which are a.e. continuous but such that the set of the points of discontinuity is dense. An example is the following function, which is the    {\sc Riemann function}\index{function!Riemann} or {\sc Thomae function.}\index{function!Thomae} 
\begin{Example}{\rm
The function is defined on the interval $[0,1]$. We represent the rational points of $[0,1]$ as a   fraction in lowest terms, i.e. $p/q$ with $p$ and $q$ without common factors (different from $1$). This representation is unique if $p/q>0$. Then we define
\[
f(x)=\left\{
\begin{array}{lll}
0 &{\rm if}& x\in\zzr\setminus\zzq\\
1 &{\rm if}& x=0\\
\dfrac{1}{q}& {\rm if}& x=\dfrac{p}{q}\,,\qquad p>0\,.
\end{array}
\right.
\]
It is clear that this function is discontinuous at $x_0=p/q$. In fact, let $0<\ZEP<1/q$. Irrational points exist in every neighborhood of $x_0$ and at these
points the function is $0$. So, $\ZOMq\left (f;(x_0-\ZSI,x_0+\ZSI)\right )>\ZEP$ for every $\ZSI>0$ and the function is not continuous at $x_0$.

Let now $x_0$ be irrational. We prove continuity at $x_0$.

We fix any $\ZEP>0$ and we examine the inequality
\[
|f(x)-f(x_0)|=f(x_0)<\ZEP\,.
\]

We recall that $f(x)=0$ if $x\notin\zzq $ while 
\[
f(p/q)=1/q\,.
\]
The inequality
\[
\frac{1}{q} \geq   \ZEP 
\]
holds for finitely many values $q_1$,\,,\dots\, $q_K$ of the denominators.

The fractions $p/q\in (0,1]$ have $p<q$
and so the set $C$:
\[
C=\{ p/q\ \mbox{such that}\ 1/q>\ZEP \}
\]
is finite and
 \[
 {\rm dist}\, (x_0,C)>0\,.
 \]
 
 Let \[
\ZDE=\dfrac{1}{2}{\rm dist}\, (x_0,C)\,. 
 \]
When $|x-x_0|<\ZDE$ we have
 \[
f(x)-f(x_0)=\left\{\begin{array}{lll}
0&{\rm if}&x\in\zzr\setminus \zzq\\
\dfrac{1}{q}<\ZEP&{\rm if}&x\in\zzq\cap (x_0-\ZDE,x_0+\ZDE)\,.
\end{array}\right. 
 \]
 So, $f$ is continuous at every irrational point.\zdia
}
\end{Example}  
 
\subsection{\ZLA{leBiNteSUbsCh1RIEMann}Riemann and Lebesgue Integrability}

We invoke Lemma~\ref{lemmaCH1daqocontTOquasicont}: when the set of discontinuities of a function is a null set then $f$ is quasicontinuous. So, we can state the following corollary to Theorem~\ref{Theo:Appe1VitaliLebe}
\begin{Corollary}\ZLA{coroAPPE1LEBEdaRiem}
Any Riemann integrable function is quasicontinuous and (being bounded) it is   summable.
\end{Corollary}
  It remains to be proved that   the  Riemann and Lebesgue integrals have the same values. This we do now. First we prove   the following property of the functions $\chi_{+,N}$ and $\chi_{-,N}$:
\begin{Lemma}\ZLA{Ch1Appe1LemmaConveGenCHI}
If $f$ is continuous at $x_0$ then we have
\[
\lim _{N\to+\ZIN}\chi_{+,N}(x_0)=f(x_0)\,,\qquad \lim _{N\to+\ZIN}\chi_{-,N}(x_0)=f(x_0)\,.
\]
\end{Lemma}
\zProof we prove the statement for the functions $\chi_{+,N}$. We must prove:
\[
\forall \ZEP>0\ \exists N_\ZEP\,:\quad N>N_\ZEP\ \implies\ f(x_0)-\ZEP< 
\chi_{+,N}(x_0)< f(x_0)+\ZEP\,.
\] 
 The assumption is
 \[
\forall  \ZEP>0\ \exists \ZDE_{ \ZEP}\,:\quad |x-x_0|<\ZDE_{ \ZEP}\ \implies\ f(x_0)  -  \ZEP< 
f(x)< f(x_0)+  \ZEP\,.
\] 
  We choose $N_\ZEP$ such that
\[
\dfrac{k-h}{N_{\ZEP}}< \dfrac{1}{2}\ZDE_{\ZEP}\,.
\]
Let $N>N_{\ZEP}$ and let $x_0$ belong to the interval 
\[\left [  
h+
i_0\frac{k-h}{N}
, h+
(i_0+1)\frac{k-h}{N}
\right ) 
\,.
\] 
Then
\begin{multline*}
f(x_0)-\ZEP\leq
f(x_0)\leq\chi_{+,N}(x_0)\\
=\sup \left\{ f(x)\,, \ x\in
\left [  
h+
i_0\frac{k-h}{N}
, h+
(i_0+1)\frac{k-h}{N}
\right ) 
 \right\} 
 \leq f(x_0)+\ZEP\,.
\end{multline*}
This inequality verifies the thesis since $\ZEP>0$ is arbitrary.\zdia
 
Lemma~\ref{Ch1Appe1LemmaConveGenCHI} holds for any function $f$. If $f$ is Riemann integrable it is bounded and Vitali-Lebesgue Theorem implies:
\begin{Corollary}\ZLA{Appe1Ch1CoroConveqoCHI}
Lef $f$ be Riemann integrable on $[h,k]$. The sequences $\{\chi_{+,N}\}$
and $\{\chi_{-,N}\}$ are bounded and a.e. convergent to $f$.
\end{Corollary} 
 
 Now we can prove:
\begin{Theorem}\ZLA{teoAppe1Cap1DaRiemInteAlebeOnte}
If a function is Riemann integrable then   its Riemann and Lebesgue integrals have the same value.
\end{Theorem} 
\zProof
We use the definition~(\ref{EQ:I2AppeCa1TheoDaPartApartEqui})  of the Riemann integral, the fact that the integrals on the right sides are Lebesgue integrals and, thanks to Corollary~\ref{Appe1Ch1CoroConveqoCHI}, 
Theorem~\ref{TeoCH1P1ConveUnderBound} which shows that, under boundedness conditions, a.e. limits and Lebesgue integrals can be exchanged.

When $f$ is Riemann integrable we get:

\smallskip

\begin{center}
\begin{tabular}{ ccccc}
 
$\overbrace{\int_h^k \chi_{+,N}(x)\ZD x}^{\tiny \begin{tabular}{l}
    Riemann and Lebesgue  \\    integral  \end{tabular}}$ 
& 
& = 
& 
&$\overbrace{\int_h^k \chi_{+,N}(x)\ZD x}^{\tiny \begin{tabular}{l}
    Riemann and Lebesgue  \\    integral  \end{tabular}}$  \\
 $ 
 \rotatebox{-90}{$\to$} 
 $&&  &&$ \rotatebox{-90}{$\to$} $ \\
$\underbrace{\int_h^k f(x)\ZD x}_{\tiny\mbox{  } \begin{tabular}{l}
    Riemann     integral\\(by definition)  \end{tabular}}$&&
 && $\underbrace{\int_h^k f(x)\ZD x}_{\tiny \begin{tabular}{l}
    Lebesgue     integral\\
    (Theorem~\ref{TeoCH1P1ConveUnderBound})  \end{tabular}}$  
\end{tabular}
\end{center}
\medskip
So we have
\[
\underbrace{\int_h^k f(x)\ZD x}_{\tiny 
    \mbox{\tiny Riemann     integral}} = \underbrace{\int_h^k f(x)\ZD x}_{\tiny 
   \mbox{\tiny Lebesgue     integral}}
\]
 as wanted.\zdia

 \chapter[Limits and Integrals: One Variable]{\ZLA{Cha1bis:INTElebeUNAvariaBIS}
Functions of One Variable: the Limits and the Integral} 

In this chapter we prove the key theorems concerning the exchange of the limits and the integrals for functions of one variable. Similar theorems for functions of several variables are proved in Chap.~\ref{ch3:LIMITintegral} but we distinguish the two treatments since in the case of functions of one variable the treatment is  simpler because we can use Theorem~\ref{Theo:prprieSucceConvDIM1}.

In the course of our analysis in this chapter we introduce and use the important property of the absolute continuity of the integral.

\ZREF{
The key tool used in this chapter is Egorov-Severini Theorem which is proved here for function of one variable. The theorem holds for function of several variables too (see Chap~\ref{ch3:LIMITintegral}). Once the theorem is proved, its
  consequences   can be deduced with the same proofs regardless of the number of variables. 
In order to stress this fact, and in order to use the proofs give here also when the functions depend on several variables,
in this chapter an interval is denoted $R$ (initial of ``rectangle'' since the Tonelli construction of the integral of functions of several variables uses  rectangles  instead of intervals, see Ch.~\ref{Chap:2quasicontPIUvar}).
}

\ZREF{
In order to streamline certain statements, 
it is convenient to recall Theorem~\ref{teo:Ch2StrutturaAperti}:  if $\mathcal{O}\subseteq \zzr$ is an open set then there exists $\Delta$ (disjoint multiinterval) such that   $\mathcal{O}=\mathcal{I}_\Delta$. The corresponding number $L(\Delta)$ depends only on the set $\mathcal{O}$ and it is  denoted $\zl(\mathcal{O})$ (see the Definition~\ref{DefiCH1DeFiMiSuOpen}). 
}

\section{\ZLA{Cha1bisSECEgorovSev}Egorov-Severini Theorem and Quasicontinuity}

 It is convenient to note:
 \begin{Lemma}\ZLA{LemmaCha1bisSucceFunzQCinDiffeR}
Let $\{f_n\}$ be a sequence of functions each one a.e. defined on an interval $R$; i.e.,
\[
\Dom\, f_n=R\setminus N_n \qquad \mbox{($N_n$ is a null set)}\,.
\]
Then we have:
\begin{enumerate}
\item
\ZLA{I1LemmaCha1bisSucceFunzQCinDiffeR}
there exists a null set $N$ such that every $f_n$ is defined on $R\setminus N$.
\item\ZLA{I2LemmaCha1bisSucceFunzQCinDiffeR} if each $f_n$ is quasicontinuous then for every $\ZEP>0$ there exists a multiinterval $\Delta_\ZEP$ such that   $(f_n)_{|_{R\setminus \mathcal{I}_{\Delta_\ZEP}}}$ is continuous.
\end{enumerate}
 \end{Lemma}
  \zProof The set $N$ is $N=\cup N_n$. In fact, we proved in Lemma~\ref{LemmaCAP1UnioNulliEnulla} that $N$ is a null set.
  
The multiinterval $\Delta_\ZEP$ is constructed in a similar way: we associate a multinterval $\Delta_{n,\ZEP}$ of order $\ZEP/2^n$   to $f_n$ and we put $\Delta_\ZEP=\cup _{n\geq 1}\Delta_{n,\ZEP}$.\zdia

\ZREF{
 This observation shows that when working with sequences of  functions each one of them defined a.e. on $R$ we can    assume that they are all defined on $R\setminus N$ where $N$
 is a null set which does not depend on $n$. To describe this case we say (as   in Chap.~\ref{Ch1:INTElebeUNAvaria}) that \emph{the sequence $\{f_n\}$ is defined a.e. on $R$.} 
 }
 
Now we give a definition:
\begin{Definition}\ZLA{Cha1bisDEFINITalmUNIFconve}
{\rm Let $f_n$,   $f$ be  functions a.e. defined on a  set  $K\subseteq R$. We say that the sequence $\{ f_n\}$  {\sc converges almost uniformly}\index{convergence!almost uniform}\index{almost uniform convergence} to $f$  on $K$ when for every $\ZEP>0$     the following \emph{equivalent} statements hold:
\begin{enumerate}
\item there exists an open set $\mathcal{O}$ such that 
\[
\mbox{$\zl(\mathcal{O})<\ZEP $  and $\{ f_n\}$ converges \emph{uniformly} to $f$ on $K\setminus\mathcal{O}$}\,.  
\]
\item
there exists a multiinterval  $\Delta$   such that 
 \[ 
\mbox{ $
L(\Delta)<\ZEP$ and  
 $\{ f_n\}$ converges \emph{uniformly} to $f$ on $K\setminus{\mathcal{I}_\Delta}$}\,.  
\] 
\end{enumerate}

 }
\end{Definition}

\ZREF{
We must be clear on the content of this definition. We illustrate its content in terms of open sets and we invite the reader to recast it in terms of multiintervals.

  We fix any $\ZEP>0$ and we find an open set $\mathcal{O}_\ZEP$ such that
$\zl(\mathcal{O}_\ZEP)<\ZEP$ and such that
 the following property is valid: for every $\ZSI>0$ there exists a number $N$ \emph{which depends on $\ZSI$ and on the previously chosen set $\mathcal{O}_\ZEP$, $N=N_{\ZSI,\mathcal{O}_\ZEP }$}   such that
\[
\left\{\begin{array}{l}
n>N_{\ZSI,\mathcal{O}_\ZEP }\\
x\in R\setminus \mathcal{O}_\ZEP
\end{array}\right.\quad\implies \ |f_n(x)-f(x)|<\ZSI\,.
\]
\emph{The important point is that $\mathcal{O}_\ZEP$ does not depend on $\ZSI$.}

Finally we note: ``$\{ f_n\}$ converges uniformly  to $f$ on $K\setminus\mathcal{O} $'' 
is equivalent to ``$\{ (f_n)_{|_{K\setminus\mathcal{O}}}\}$ converges uniformly to $f_{|_{K\setminus\mathcal{O}}}$''. 

 }

Now we state the following preliminary result whose proof is   in Appendix~\ref{AppeCha1bisEgSEV}:

 \begin{Theorem}[Egorov-Severini: preliminary statement] 
\ZLA{Cha1bisTeoEgoSEVERANTE}Let $\{ f_n\}$ be a  sequence of  \emph{continuous} functions  \emph{everywhere defined} on the \emph{closed and bounded} interval $R=[h,k]$.
 If the sequence converges  on $R$
to a function $f$; i.e. \emph{if $f_n(x)\to f(x)$ for every $x\in R$,} then:
\begin{enumerate}
\item\ZLA{I1Cha1bisTeoEgoSEVERANTE}
the sequence    converges almost uniformly;
\item\ZLA{I2Cha1bisTeoEgoSEVERANTE} the limit function $f$ is quasicontinuous.
\end{enumerate}
\end{Theorem}

Theorem~\ref{Cha1bisTeoEgoSEVERANTE} is a special instance of the   Egorov-Severini Theorem\footnote{\ZLA{footP1Cap2EGSev}independently published by Egorov in~\cite{Egorov1911} and by Severini in~\cite{Severini1910}.} that we state now. Actually,   the two theorems are equivalent since the second is a consequence of the first.
  
\begin{Theorem}[{\sc Egorov-Severini}]\index{Theorem!Egorov-Severini}
\ZLA{Cha1bisTeoEgoSEVER}Let $\{ f_n\}$ be a  sequence of quasicontinuous functions a.e. defined on the \emph{bounded} interval $R$. 

We assume that for every $\ZEP>0$ there exists a multiinterval $\tilde \Delta_\ZEP$ such that  $L(\tilde \Delta_\ZEP)<\ZEP$ and such that the sequence $\{f_n\}$ is bounded on  $
R\setminus \mathcal{I}_{\tilde \Delta_\ZEP}$.

If the sequence converges a.e. on $R$
to a function $f$ (hence a.e. defined on $R$) then:
\begin{enumerate}
\item\ZLA{I1Cha1bisTeoEgoSEVER}
the sequence    converges almost uniformly;
\item\ZLA{I2Cha1bisTeoEgoSEVER} the limit function $f$ is quasicontinuous.
\end{enumerate}
\end{Theorem}
\zProof 
For clarity we split the   the proof in the following steps.  
\begin{description}
\item[\bf Step~1:] First we use the assumption that the sequence has to be   defined and   convergent \emph{a.e.} on $R$:  whether $R$ includes its end points has no effect. So, we can  assume that it is closed.
Then we use Lemma~\ref{LemmaCha1bisSucceFunzQCinDiffeR} and we find a null set $N\subseteq R$ such that each $f_n$ is defined on $N$ and converges to $f$ on $R\setminus N$.

\item[\bf Step~2:]
Each $f_n$ is quasicontinuous. We use again Lemma~\ref{LemmaCha1bisSucceFunzQCinDiffeR} and we see that for every $\ZEP>0$ there exists $\Delta_\ZEP$, the same multiinterval for every $n$, such that:
\begin{enumerate}
\item $L(\Delta_\ZEP)<\ZEP$;
\item $ \mathcal{I}_{\tilde\Delta _\ZEP}\subseteq \mathcal{I}_{\Delta_\ZEP}$;
\item for every $n$ the function  $(f_n)_{|_{R\setminus\mathcal{I}_{\Delta_\ZEP}}}$ is continuous;
\item the sequence $\{(f_n)_{|_{R\setminus\mathcal{I}_{\Delta_\ZEP}}}\}$ is bounded;
 
\item the sequence   $\{(f_n)_{|_{R\setminus\mathcal{I}_{\Delta_\ZEP}}}\}$  converges to $ f _{|_{R\setminus\mathcal{I}_{\Delta_\ZEP}}}$. 
\end{enumerate}

We denote   $f_{n,e}$    the 
Tietze extensions
of the functions $(f_n)_{|_{R\setminus\mathcal{I}_{\Delta_\ZEP}}}$
 \emph{defined by using the algorithm described in Theorem~\ref{teo:Ch1ESTEdaCHIUSO}.} So, for every $n$, $f_{n,e}$ is a continuous function defined on $\zzr$. We consider the restriction of $f_{n,e}$ to $R$.     \emph{Thanks to the property~\ref{I2Theo:prprieSucceConvDIM1} of Theorem~\ref{Theo:prprieSucceConvDIM1},} $\{f_{n,e}(x)\}$ converges \emph{for every $x\in R$.} So we can use Theorem~\ref{Cha1bisTeoEgoSEVERANTE}: we assign $\ZEP_1>0$ and we find a multinterval $\Delta_{\ZEP_1}$   such that $L(\Delta_{\ZEP_1})<\ZEP_1$ and
 $\{f_{n,e}\}$ converges uniformly on $R\setminus 
  \mathcal{I}_{\Delta_{\ZEP_1}}$.
 
 \item[\bf Step~3:]
 Uniform convergence of continuous functions implies continuity of the limit, so $\{f_{n,e}\}$ converges on  $R\setminus    \mathcal{I}_{\Delta_{\ZEP_1}}$ to a function which is continuous on this set. 
 
 \item[\bf Step~4:]
 Now we consider the set $R\setminus \mathcal{I}_{\Delta_\ZEP\cup    \Delta_{\ZEP_1}}$. The restrictions to this set of the functions $f_{n,e}$ are continuous and uniformly convergence holds. So, the limit is continuous on 
 $R\setminus \mathcal{I}_{\Delta_\ZEP\cup\Delta_{\ZEP_1}}$.  But, on 
 $R\setminus \mathcal{I}_{\Delta_\ZEP\cup\Delta_{\ZEP_1}}$ we have $f_{n,e}(x)=f_n(x)$ and we know that on this set $f_n(x)\to f(x)$.
 
 \item[\bf Step~5:]
 It follows that:
 \begin{itemize}
\item 
  the restriction of $f$ to  $R\setminus \mathcal{I}_{\Delta_\ZEP\cup\Delta_{\ZEP_1}}$ is continuous;
  \item the sequence $\{ f_n\}$ converges uniformly to $f$ on $ R\setminus\mathcal{I}_{\Delta_\ZEP\cup\Delta _{\ZEP_1}} $.
  \end{itemize}
 \end{description}
  
The result  follows since
 \[
L\left (  \Delta_\ZEP\cup\Delta_{\ZEP_1}\right )<\ZEP+\ZEP_1
 \]
 and both $\ZEP $ and $\ZEP_1$ can be arbitrarily assigned.\zdia

\begin{Remark}\ZLA{CHCha1bis:RemaESEunBEgorSEV}
{\rm  
Egorov-Severini Theorem does not assume that $\{f_n\}$ is   bounded while it is assumed that \emph{the interval $R$ is bounded.} This is a crucial assumption which cannot be removed  as the following example shows: 
 let $R=(0,+\ZIN)$ and let 
\[
f_n(x)=\left\{\begin{array}{lll}
0&{\rm if}& x\leq n-2 \\
x-(n-2)&{\rm if}& n-2\leq x\leq n-1\\
1 &{\rm if}&n-1\leq x\leq n+1\\
n+2-x &{\rm if}& n+1\leq x\leq n+2\\
0 &{\rm if}& x\geq n+2\,.
\end{array}\right.
\] 
We have $f_n(x)\to 0$ for every $x$ 
 but $\{ f_n(x)>1/2\}$ contains the interval $[n-1,n+1]$ of length $2$ for every $n$.\zdia
 
}
\end{Remark}
  In spite of this example we have:
  \begin{Corollary}\ZLA{CoroP1Ch2EsteEgSevReTTA}
  Let $\{f_n\}$ be a sequence of functions defined on $\zzr $ and let $f_n(x)\to f(x)$ a.e. $x\in\zzr$. The function $f$ is quasicontinuous.
  \end{Corollary}
\zProof We use the observation~\ref{I1RemaP1Ch1RemaSecoIngre} of Remark~\ref{RemaP1Ch1RemaSecoIngre}: the function $f$ is quasicontinuous if and only if its restrictions to every bounded interval $[-k,k]$ are quasicontinuous.  It is indeed so, by using Therem~\ref{Cha1bisTeoEgoSEVER} on the bounded intervals.\zdia

\subsection{\ZLA{SecP1CH2ConseEGOrSeve}Consequences of the Egorov-Severini Theorem}

 We prove several consequences of the Egorov-Severini Theorem.

 \begin{Corollary}\ZLA{Cha1bisCoroEGOSEveQuasContFunz}
 Let $\{f_n\}$ be a bounded sequence of quasicontinuous functions a.e. defined on a bounded interval $R$.
The following properties hold:
\begin{enumerate}
\item\ZLA{I1Cha1bisCoroEGOSEveQuasContFunz} let  $\{f_n\}$  be either a.e. increasing or   decreasing on $R$  and let
\[
f(x)=\limn f_n(x) \qquad a.e. \ x\in R\,.
\] 
The function $f$ is quasicontinuous.  
\item\ZLA{I2Cha1bisCoroEGOSEveQuasContFunz} Let $f$ be either
  \[
f(x)=\limsup_{n\to+\ZIN} f_n(x) \quad {\rm or}\quad 
f(x)=\liminf_{n\to+\ZIN} f_n(x)
     \qquad a.e. \ x\in R\,.  
  \]
  The function $f $ is quasicontinuous.
  \end{enumerate}
 \end{Corollary}
 \zProof 
Statement~\ref{I1Cha1bisCoroEGOSEveQuasContFunz} is just a restatement of Egorov-Severini Theorem in the special case of the monotone sequences of functions. Statement~\ref{I2Cha1bisCoroEGOSEveQuasContFunz} follows from statement~\ref{I1Cha1bisCoroEGOSEveQuasContFunz} since $\liminf$ and $\limsup$ are just limits of monotone sequences. In fact:

\medskip

\fbox{\parbox{2.4in}{
\[
\limsup_{n\to+\ZIN} f_n(x)=\limn \phi^{(s)}_n(x)
\]
where 
\begin{align*}
&\phi^{(s)}_n(x)=\lim_{m\to+\ZIN} \psi^{(s)}_{n,m}(x)\,,\\
& \psi^{(s)}_{n,m}(x)=\max_{0\leq i\leq m}  \{ f_{n+i}(x)  \}\,. 
\end{align*}

Statement~\ref{I4TheoCH1PropELEMqcfunc} of Theorem~\ref{TheoCH1PropELEMqcfunc} shows that the functions $\psi^{(s)}_{n,m}$ are quasicontinuous.

For every $n$, the sequence $m\mapsto  \psi^{(s)}_{n,m}$ is increasing so that $\phi^{(s)}_n$ is quasicontinuous; and $n\mapsto \phi^{(s)}_n$ is decreasing so that $f$ is quasicontinuous.\zdia
}
}
\fbox{
\parbox{2.4in}{
\[
\liminf_{n\to+\ZIN} f_n(x)=\limn \phi^{(i)}_n(x)
\]
where 
\begin{align*}
&\phi^{(i)}_n(x)=\lim_{m\to+\ZIN} \psi^{(i)}_{n,m}(x)\,,\\
& \psi^{(i)}_{n,m}(x)=\min_{0\leq i\leq m}  \{ f_{n+i}(x)  \}\,. 
\end{align*}

Statement~\ref{I4TheoCH1PropELEMqcfunc} of Theorem~\ref{TheoCH1PropELEMqcfunc} shows that the functions $\psi^{(i)}_{n,m}$ are quasicontinuous.

For every $n$, the sequence $m\mapsto  \psi^{(i)}_{n,m}$ is decreasing so that $\phi^{(i)}_n$ is quasicontinuous; and $n\mapsto \phi^{(i)}_n$ is  increasing so that $f$ is quasicontinuous.\zdia
}
}

\bigskip

It follows:
\begin{Corollary}
\ZLA{Coro:Cha1bisSUPinfPERch4}
Let $\{f_n\}$ be a bounded sequence of quasicontinuous functions defined on a interval $R$ and let   
 
\begin{align*}
\phi(x)&=\sup\{ f_n(x)\,,\quad n\geq 1\}\,\\
\psi(x)&=\inf\{ f_n(x)\,,\quad n\geq 1\}\,.
\end{align*}
The functions $\phi$ and $\psi$ are quasicontinuous.
\end{Corollary}
\zProof In fact
\begin{align*}
\phi(x)=\limn \phi_n (x)\,,&    &\phi_n(x)=\sup\{ f_k(x)\,,\quad 1\leq k\leq n\}\\
\psi(x)=\limn \psi_n (x)\,,& &\psi_n(x)=\inf\{ f_k(x)\,,\quad 1\leq k\leq n\}\,.
\end{align*}
and the functions $\phi_n$ and $\psi_n$ are quasicontinuous for every $n$.\zdia

Finally we prove:
\begin{Theorem}\ZLA{Teo:Cha1bisFunCHarAPERTIquasCont}
Let $\mathcal{O}$ be an open set. Its characteristic function $\charfun_\mathcal{O}$ is quasicontinuous.
 \end{Theorem}
 \zProof  
 We represent
  \[
\mathcal{O}=\bigcup _{n\geq 0} R_n\,, \quad    \mbox{$R_n=(a_n,b_n)$, and $R_n\cap R_j=\emptyset $ if $n\neq j$}\,.
 \]
 Then we have
 \[
\charfun_{\mathcal{O}}(x)=\sum _{n\geq 1} \charfun _{ (a_n,b_n)}(x)=\sum _{n\geq 1} \charfun _{ R_n}(x)
 \]
and the sum is either finite or a convergent series. The series converges for every $x$ since the intervals $R_n$ are disjoint  so that for every $x$ only one term of the series is different from zero.

We know that the the characteristic function of an interval is   quasicontinuous. So, if the sum is finite $\charfun_{\mathcal{O}}$ is quasicontinuous as the sum of a finite number of quasicontinuous functions. Otherwise it is quasicontinuous thanks to Egorov-Severini Theorem.\zdia

 
Theorem~\ref{Teo:Cha1bisFunCHarAPERTIquasCont} has the following consequence:
\begin{enumerate}
\item if $f$ is  summable   and if $\mathcal{O} $ is an open set then the integral of $f$ on $\mathcal{O}$ exists
\begin{equation}\ZLA{eq:Cha1bisAbsConTdefINTEopenSET}
\int _{\mathcal{O}} f(x)\ZD x=\int_R f(x)\charfun _{\mathcal{O}}(x)\ZD x\,.
\end{equation}

 \item
 Let $\mathcal{O}$ be a bounded open set,
 \[
\mathcal{O}=\bigcup _{n\geq 0} R_n\,, \quad    \mbox{$R_n=(a_n,b_n)$, and $R_n\cap R_j=\emptyset $ if $n\neq j$}\,.
 \]
we can associate two numbers to the open set $\mathcal{O}$:
 
\[
\begin{tabular}{ccccc}
$ \zl(\mathcal{O})=L\left (\{R_n\}\right )$
&
$=$
&
 $  \sum _{n\geq 1} \zl(R_n) =\sum _{n\geq 1} \int_R \charfun _{R_n}(x)\ZD x$\\
 &
\\ 
&&&&
\\
$\int_R\charfun_{\mathcal{O}}(x)\ZD x$
&
$=$
&
$ \int_R\left [\sum_{n\geq 1} \charfun_{R_n}(x)\right]\ZD x$\,.
&
&
 \end{tabular}
\]
If the sums are finite they can be exchanged with the integrals and   we have
\[
\zl(\mathcal{O})=\int _{R} \charfun_\mathcal{O}(x)\ZD x\,.
\]
We prove that this equality holds even if the series is not a finite sum:

 \end{enumerate}

 \begin{Theorem}\ZLA{Cha1bis:TeoPreAbsContPriLiMint}
We have:
\begin{equation}\ZLA{eq:Cha1bisPrimadelTEO127}
\zl(\mathcal{O})= 
  \int_R\charfun_{\mathcal{O}}(x)\ZD x
\,.
\end{equation}
 \end{Theorem}
\zProof
We note:
 \[
\charfun_{\mathcal{O}}(x)=\lim _{N\to+\ZIN} \charfun _{\mathcal{O}_N} (x)\,,
\quad
\mathcal{O}_N=\bigcup_{n=1}^N \mathcal{O}_n\,,
\qquad 
\charfun _{\mathcal{O}_N} (x)=\sum_{n=1}^{N} \charfun_{ R_n}(x) \,.
 \]
The function $\charfun_{\mathcal{O}_N}$ is quasicontinuous. 

We must study the two equalities~(\ref{eq:Cha1bisPSeCodelTEO127}) and~(\ref{eq:Cha1bisTeRzTEO127}) below:

\begin{subequations}
\begin{align}
\ZLA{eq:Cha1bisPSeCodelTEO127}
&\left\{\begin{array}{l}
\charfun_{\mathcal{O}}(x)=\lim _{N\to+\ZIN} \charfun _{\mathcal{O}_N} (x)=\sum_{n= 1}^{+\ZIN}  \charfun_{ R_n}(x) \\
{\rm where}\ \charfun_{{\mathcal{O}_N}} (x)=\sum_{n= 1}^N  \charfun_{ R_n}(x)
\end{array}\right.
\\
&
\ZLA{eq:Cha1bisTeRzTEO127}
 \int_R\charfun_{\mathcal{O}}(x)\ZD x
=  \int _R\left [\sum_{n\geq 1} \charfun_{{ R_n} }(x)\right]\ZD x\,.
\end{align} 
\end{subequations}
The inequality
 \[
 \charfun _{\mathcal{O}_N} (x)\leq \charfun_{\mathcal{O}}(x)\quad {\rm a.e.}\ x\in\mathcal{O}
 \]
implies
 \[
 \lim _{N\to+\ZIN} \sum _{n=1}^N
 \int _R  \charfun _{ R_n} (x)\ZD x
= 
  \lim _{N\to+\ZIN} 
 \int _R
 \sum _{n=1}^N 
   \charfun _{ R_n} (x)\ZD x 
 \leq \int_R \charfun_{\mathcal{O}}(x)\ZD x\,.
 \]
 So, for every $N$ we have
 \[
 \int_R \charfun_{\mathcal{O}}(x)\ZD x- \sum _{n=1}^N
 \int _R  \charfun _{ R_n} (x)\ZD x\,.
 \]
 
 \ZREF{
We prove that   
\begin{multline} \ZLA{eq:Cha1bisPSECONDAdelTEO127}
0\leq  \zaa\leq  \int_R \charfun_{\mathcal{O}}(x)\ZD x- \sum _{n=1}^N
 \int _R  \charfun _{ R_n} (x)\ZD x\\=
 \int_R \left [\charfun_{\mathcal{O}}(x)\ZD x-  \sum _{n=1}^N  \charfun _{ R_n} (x)\right] \ZD x 
 =\int _R\left [\sum _{n=N+1}^{+\ZIN}  \charfun _{ R_n} (x)\right] \ZD x  
 \end{multline}
   for every $N $  implies $\zaa=0$.} 
 
 Note that:
 \begin{enumerate}
 \item the integrals in~(\ref{eq:Cha1bisPSECONDAdelTEO127}) are Lebesgue integrals, i.e. limits of Riemann integrals of associated continuous functions.
 \item the integrand of the last integral in~(\ref{eq:Cha1bisPSECONDAdelTEO127}) is nonzero if $x\notin \cup _{n>N}R_n$.
 
 \end{enumerate}
 For every $  \nu  $ we construct an associated continuous function of order $1/  \nu  $
 of the integrand as follows:
 
\begin{enumerate}
\item we choose an associated multiinterval $\Delta_{1;  \nu  }$ of order $1/2  \nu  $ of the function $\charfun_{\mathcal{O}}$ and we denote \[
(\charfun_{\mathcal{O}})_  \nu  
\]
 an associated continuous function   of order $1/2  \nu  $.

 By definition the associated continuous function is a Tietze extension  and we know that it is possible to   choose an associated continuous function 
 which satisfies the monotonicity property    of Statement~\ref{I3teo:Ch1ESTEdaCHIUSO} of Theorem~\ref{teo:Ch1ESTEdaCHIUSO}. 
 
\ZREF{\begin{center} The difference $(\charfun_{\mathcal{O}})_  \nu  -\charfun_{\mathcal{O}}$ is nonzero on $\mathcal{I}_{\Delta_{1;  \nu  }} $ and $L(\Delta_{1;  \nu  })<1/2  \nu   $.
\end{center}
}
\item  
we choose an associated multiinterval $\Delta_{2; \nu }$ of order $1/2 \nu $
associated to $\sum _{n=1}^N  \charfun _{ R_n}$ and
we denote
\[
\left (\sum _{n=1}^N  \charfun _{ R_n}\right )_ \nu 
\]
 an associated continuous function   of order $1/2n$.

 Also in this case we choose   an associated continuous function which satisfies the monotonicity property of Statement~\ref{I3teo:Ch1ESTEdaCHIUSO} of Theorem~\ref{teo:Ch1ESTEdaCHIUSO}.
 
 \ZREF{
The difference 
\[
\left (\sum _{n=1}^N  \charfun _{ R_n}\right )_ \nu -\sum _{n=1}^N  \charfun _{ R_n}
\]
 is nonzero on $\mathcal{I}_{\Delta_{2; \nu }} $ and $L(\Delta_{2; \nu })<1/2 \nu $. 
 }
  
\end{enumerate} 
 
The function
\[
F_{N; \nu }(x)=(\charfun_{\mathcal{O}})_ \nu -\left (\sum _{n=1}^N  \charfun _{ R_n}\right )_ \nu 
\]
is an associated continuous function of order $1/ \nu $ to
\[
\charfun_{\mathcal{O}}-\sum _{n=1}^{N}  \charfun _{ R_n}
=
\sum _{n=N+1}^{+\ZIN}  \charfun _{ R_n}\,.
\]
This associated function takes values in $[0,1]$ thanks to the fact that both the associated functions we chosen satisfy the monotonicity assumption and because
\[
\charfun_{\mathcal{O}}(x)\geq  \sum _{n= 1}^{N}  \charfun _{ R_n}(x)\,.
\]

By definition, the Lebesgue integral in~(\ref{eq:Cha1bisPSECONDAdelTEO127})
is
\[
\lim _{ \nu \to+\ZIN} \underbrace{\int_R F_{N; \nu }(x)\ZD x}_{\tiny \begin{tabular}{l}
    Riemann integral  \end{tabular}}\,.
\]
We investigate where  $F_{N; \nu }$ can  possibly be   non zero: this is where the integrand in~(\ref{eq:Cha1bisPSECONDAdelTEO127}) is non zero and also  at the points of  $\mathcal{I}_{\Delta_ \nu}$ where $\Delta_ \nu =\Delta_{1; \nu }\cup\Delta_{2_; \nu }$:

\[
\{x\,:\ F_{N;\nu}(x)\neq 0\}\subseteq \left [\bigcup _{N+1}^{+\ZIN} R_n\right ]\cup \mathcal{I}_{\Delta_ \nu }
\]
and  
\[
\left [\bigcup _{N+1}^{+\ZIN} R_n\right ]\cup \mathcal{I}_{\Delta_ \nu }=\mathcal{I}_{\hat \Delta}
\]
where $\hat\Delta$ is a multiinterval such that
\[
L(\hat\Delta)\leq \dfrac{1}{\nu}+\sum _{n=N+1}^{+\ZIN} L( R_n)\,.
\]
We use Lemma~\ref{Lemma:Ch2LemmaPreliRiemaa} and we see that
\[
0\leq\zaa\leq  \underbrace{\int_R F_{N; \nu } (x)\ZD x}_{\tiny \begin{tabular}{l}
    Riemann integral  \end{tabular}}\leq \dfrac{1}{\nu}+\sum _{n=N+1}^{+\ZIN} L( R_n)\quad \mbox{(since $0\leq F_{N;n \nu } (x)\leq 1$)}\,.
\]
The   Lebesgue integral is obtained by taking the limit for $\nu\to+\ZIN$. So we have:  

\[
0\leq\zaa\leq \underbrace{\int _R\left [\sum _{n=N+1}^{+\ZIN}  \charfun _{ R_n} (x)\right] \ZD x
}_{\tiny \begin{tabular}{l}
    Lebesgue integral  \end{tabular}}
 =\lim _{ \nu \to+\ZIN} 
\underbrace{ \int_R F_{N; \nu }(x)\ZD x
}_{\tiny \begin{tabular}{l}
    Riemann integral  \end{tabular}} 
 \leq \sum _{n=N+1}^{+\ZIN} L( R_n)\,.
\]
This inequality holds for every $N$ and the limit for $N\to+\ZIN$ gives $\zaa=0$, as wanted.\zdia

\begin{Remark}{\rm
The statement of Theorem~\ref{Cha1bis:TeoPreAbsContPriLiMint} can be written
\[
\lim _{N\to+\ZIN}\left [\sum _{n=1 }^N \int_R \charfun _{ R_n}(x)\ZD x\right ]
= \int _R\left [\lim _{N\to +\ZIN}\sum_{n=1 }^N \charfun_{{ R_n} }(x)\right]\ZD x
\]
and it is a first instance of the exchange of limits and integrals, the main goal of this chapter.\zdia
}\end{Remark}

\subsection{\ZLA{Cha1bisSUBSabsoluConti}Absolute Continuity of the Integral}
We prove that   the integral is absolutely continuous, i.e. we prove the following theorem\footnote{this statement of absolute continuity is not the most general. The  general statement is 
 in Sect.~\ref{sectCH4MEASUREgeneral}, Theorem~\ref{TeoCH4ABsolContINgene}.}:
 \begin{Theorem}\ZLA{Teo:Cha1bisABSOluContiPERiSoliAperti}
 Let $f$ be summable on $\zzr $. The set valued function 
 \[
 \mathcal{O}\mapsto \int _{\mathcal{O}} f(x)\ZD x\qquad \mbox{($\mathcal{O}$ open)}
 \]
   is {\sc absolutely continuous}\index{set function!absolute continuity}\index{absolute continuity} in the sense that for every $\ZEP>0$ there exists $\ZDE>0$ such that 
 \begin{equation}\ZLA{EqTeo:Cha1bisABSOluContiPERiSoliAperti}
 \zl(\mathcal{O})<\ZDE\ \implies\ 
\left |\int _{\mathcal{O}} f(x)\ZD x \right |<\ZEP\,.
 \end{equation}
 \end{Theorem}
\zProof
It is sufficient to prove the theorem when $f\geq 0$.

First we consider the case that $f$ is a bounded quasicontinuous function, $0\leq f(x)\leq M$ on a bounded interval $R$. We combine~(\ref{eq:Cha1bisAbsConTdefINTEopenSET}) and~(\ref{eq:Cha1bisPrimadelTEO127}) and we find
\begin{equation}\ZLA{eq:Cha1bisSUBSabsoluConti0}
0\leq \int _{\mathcal{O}} f(x)\ZD x\leq M\int _{\mathcal{O}} 1\ZD x=M\int _R\charfun _{\mathcal{O}}
(x)\ZD x=M\zl(\mathcal{O})\,.
\end{equation}
So, absolute continuity holds when the integrand is bounded.

  We  consider the general case of summable functions on $\zzr $.  We fix to numbers $K$ and $N$ such that
 \[
\int _{\zzr } f(x)\ZD x-\int_{\zzr }f_{+;\,(K,N)}(x)\ZD x <\ZEP/2\,. 
 \]
 Then we use absolute continuity which holds when the integrand is the bounded function $f_{+;\,(K,N)}$: we fix $\ZDE>0$ such that
 \[
\zl(\mathcal{O})<\ \ZDE\ \implies 0\leq \int _{\mathcal{O}}f_{+;\,(K,N)}(x)\ZD x<\ZEP/2\,.
 \]
Then we have
 
\begin{multline*}
 \int_{\mathcal{O}}  f(x) \ZD x\leq \underbrace{\int_{\mathcal{O}} [f(x)-f_{+;\,(K,N)}(x)]\ZD x}_{
 \tiny \leq \int_{\tiny \zzr }[f(x)-f_{+;\,(K,N)}(x)]{\ZD x\, \leq\, \ZEP/2}
 }\\
+
 \int_{\mathcal{O}}  f_{+;\,(K,N)}(x) \ZD x<\ZEP\,.\zdiaform
\end{multline*}
 \section{\ZLA{Cha1bisSecDimoConve}The Limit of  Sequences and the Lebesgue Integral} In this section we prove the theorems concerning limits and integrals. We proceed as follows:
 \begin{enumerate}
 \item in Sect.~\ref{Cha1bisSecLIMINTboudSeqBoundINT}   we study the limits under boundedness assumptions.   
 \item
in Sect.~\ref{Cha1bisSecLIMINTposiSEQUE} we study the case of   sequences of nonegative functions.
\item the general case is in Sect.~\ref{Cha1bisSecLIMINTgeneSEQUELebTh}.
 \end{enumerate}
 
 \subsection{\ZLA{Cha1bisSecLIMINTboudSeqBoundINT}Bounded Sequences on a Bounded Interval}
 
 We restate and prove   Theorem~\ref{TeoCH1P1ConveUnderBound}:
 
 \begin{Theorem}\ZLA{teo:Cha1bis:LebeCAROPartIntervLIMIT}
 Let $\{f_n\}$ be a bounded sequence of quasicontinuous functions a.e. defined on a bounded   interval $R\subseteq\zzr $. If 
 \[
\limn f_n(x)=f(x)\qquad {\rm a.e. \ on}\  R 
 \]
 then we have also
 \[
\limn \int_R f_n(x)\ZD x=\int_R f(x)\ZD x\,. 
 \]
 \end{Theorem}
 \zProof
 Due to the fact that the assumed properties hold \emph{a.e.}, to fix our ideas we can assume that $R$ is closed.
 
 Note that $f$ is summable since $f$ is bounded and   quasicontinuous\footnote{quasicontinuity is   the   statement~\ref{I1Cha1bisTeoEgoSEVER} of Theorem~\ref{Cha1bisTeoEgoSEVER}.}. So, the integral on the right side is a number.
 
By replacing $f_n$ with $f_n-f$ we can prove:
 \begin{equation}\ZLA{Eq1:teo:Cha1bis:LebeCAROPartIntervLIMIT}
f_n\to 0\ {\rm a.e.}\ \implies\ \int_R f_n(x)\ZD x\to 0\,. 
 \end{equation}

 If $R=[h,h]$  then integrals on $R$ are zero  and the result is obvious. So we consider the case that $R=[h,k]$ with $k>h$.
 
 We rewrite the thesis~(\ref{Eq1:teo:Cha1bis:LebeCAROPartIntervLIMIT}) in explicit form:
 \begin{equation}\ZLA{Eq2:teo:Cha1bis:LebeCAROPartIntervLIMIT}
 \begin{array}{l}
 \displaystyle 
 \mbox{we fix any $\ZEP>0$. We must prove the existence of $N(\ZEP)$ such that}\\[2mm]
  \displaystyle 
 n>N(\ZEP)\ \implies\ \left | \int_R f_n(x)\ZD x\right |<\ZEP\,.
 \end{array}
 \end{equation}
 
   We invoke Egorov-Severini Theorem: there exists an open set $\mathcal{O}$ such that
   \[
\zl(\mathcal{O})\leq \dfrac{\ZEP}{4M}\,,\quad \mbox{$f_n\to 0$ uniformly on $R\setminus \mathcal{O}$}\,.   
   \]
   We use the fact that $R\setminus \mathcal{O}$ is closed and we construct the Tietze extension $f_{n,e}$ of~$f_n$.
   
 Theorem~\ref{TheoCH1ConveTieSeq} implies  
   \[
\limn f_{n,e}=0\qquad \mbox{uniformly on $R$}   \,.
   \]
   So there exists $N_\ZEP$ such that
   \[
n>N_\ZEP\ \implies\ \left |\int_R  f_{n,e}(x)  \ZD x\right|<\dfrac{\ZEP}{2}.
   \]
  \emph{ The integral here is a Riemann integral and we know from the Appendix~\ref{leBiNteSUbsCh1RIEMann} that it is also a Lebesgue integral.} 
  
   Now we observe:
   \begin{multline*} 
 \overbrace{\left |\int_R f_n(x)\ZD x\right |}^{\tiny \begin{tabular}{l}
    Lebesgue\\ integral  \end{tabular}}
  \leq \overbrace{\left |
   \int _R \left [ f_n(x)-f_{n,e}(x)\right ]\ZD x
  \right |}^{\tiny \begin{tabular}{l}
    Lebesgue  integral  \end{tabular}}
   +\overbrace{\left | \int_R f_{n,e}(x) \ZD x  \right |}^{\tiny \begin{tabular}{l}
    Lebesgue and Riemann\\ integral  \end{tabular}}\\
  =  \underbrace{
   \left |
  \int _\mathcal{O} \left [ f_n(x)-f_{n,e}(x)\right ]\ZD x
   \right |
   }_{
 \leq 2M\zl(\mathcal{O}) <\ZEP/2 \ 
  \mbox{(see~(\ref{eq:Cha1bisSUBSabsoluConti0}))} 
   }
 + 
 \underbrace{  \left |
 \int_R f_{n,e}(x)\ZD x   \right |
  }_{
  \leq \ZEP/2  \mbox{ (when $n>N_\ZEP$)}} 
  <\ZEP
   \end{multline*}
   as wanted.\zdia

 \subsection{\ZLA{Cha1bisSecLIMINTposiSEQUE}Sequences of Nonnegative Functions}
Let $A$ be a set which satisfies Assumption~\ref{ASSUch1QuascontINSa}, i.e. such that $\charfun_A$ is quasicontinuous.

Let $\{f_n\}$ be a sequence of \emph{integrable} functions on the set $A$. We assume $f_n(x)\to f(x)$  on $A$.   \emph{Furthermore we assume that the functions are  nonnegative:}
\[
f_n(x)\geq 0\quad \mbox{so that $f(x)\geq 0$ too}\,.
\]

  The set $A$ can be unbounded and the functions can be unbounded too. So, the functions   are integrable on $A$, possibly not summable. 
  
  We extend the functions with zero. We get nonegative integrable functions defined on $\zzr $: 
\begin{equation}\ZLA{eq:Cha1bisLIMIposiSEQUE0}
f_n(x)\geq 0\,,\qquad
\limn f_n(x)= f(x)\qquad \forall x\in\zzr \,.
\end{equation}
Corollary~\ref{CoroP1Ch2EsteEgSevReTTA} shows that  $f$ is quasicontinuous and it is integrable since it is nonnegative. We investigate the relations among the integrals. This can be done by lifting to this case the result we proved in the case of bounded sequences on  a bounded rectangle. 
 
We consider the interval $B_R=[-R,R]$ and
we define
\[
  f_{ (R,N)}(x)=  
\left [\min\{f_+(x)\,,\ N\}\right]\charfun_{B_N}(x)\,. 
\]

The support of $   f_{  (R,N)} $ is compact 
and  $|f_{  (R,N)}(x)|\leq N$. Then we define
\[
f_{n;\, (R,N)}(x)=\min\{ f_n(x),f_{  (R,N)}(x)\} \,.
\] 
 
So we have:
\[
  \limn f_{n;\, (R,N)}(x)= f_{  (R,N)}(x)\qquad \forall x\in\zzr \,.
  \]
 The support of $f_{n;\, (R,N)}$ is compact, contained in the bounded interval $B_R$. So, 
for every fixed $R$ and $N$ we can apply the result  in Sect.~\ref{Cha1bisSecLIMINTboudSeqBoundINT}. We have
\begin{multline}\ZLA{Cha1bisEquaPreFatou}
\overbrace{\int _{\zzr } f_{  (R,N)}(x)  \ZD x=\limn \int_{\zzr } f_{n;\, (R,N)}(x)\ZD x}^{\footnotesize\mbox{(Theorem~\ref{teo:Cha1bis:LebeCAROPartIntervLIMIT})}}
=  \liminf_{n\to+\ZIN} \int_{\zzr } f_{n;\, (R,N)} (x)\ZD x \\
\leq \liminf_{n\to+\ZIN} \int_{\zzr } f_n(x)\ZD x=\liminf_{n\to+\ZIN}\int_A f_n(x)\ZD x \,.
\end{multline}
The  inequality in the second line follows from the monotonicity of the integral, since $f_{n;\, (R,N)}(x)\leq f_n(x) $ for every $x$ while the last equality, $\int_{\zzr }=\int_A$,
holds because the functions are originally defined on $A$, then extended to $\zzr $ with zero.

Inequality~(\ref{Cha1bisEquaPreFatou}) implies the following result, which is known as {\sc Fatou Lemma,} proved in~\cite{FATOU06Lemma}:\index{Fatou Lemma}\index{Lemma of Fatou}
\begin{Lemma}[Fatou]
If $\{f_n\}$ is a sequence of nonnegatrive functions which are integrable on a set $A$ and if $f_n(x)\to f(x)$ a.e. on $A$ then we have
\begin{equation}\ZLA{eq:DiseqFATOUCha1bis}
\int _A f(x)\ZD x\leq  \liminf_{n\to+\ZIN} \int_{A} f_{n}(x)\ZD x \,.
\end{equation}
\end{Lemma}
\zProof Inequality~(\ref{eq:DiseqFATOUCha1bis}) follows by using again the fact that the integral on $A$ is the integral on $\zzr $ of the extension with zero, and  from the definition of the integral given in Sect~\ref{sec:CH1LebeGENEdom}:
\[
\int _{A} f(x)\ZD x=
  \int_{\zzr } f(x)\ZD x=\lim _{\stackrel{N\to+\ZIN}{R\to+\ZIN}}\int _{\zzr }f_{  (R,N)}(x)\ZD x\,.\zdiaform
\]

If it happens that the sequence of nonnegative functions $\{f_n\}$ is \emph{increasing} then we have Beppo Levi Theorem~\ref{teo:ch1Beppo} (first proved in~\cite{LeviBeppoTeoMonot}) that we restate:  
\begin{Theorem}[{\sc Beppo Levi} or {\sc monotone convergence}]\ZLA{TeoCha1bisBeppo}\index{Beppo Levi Theorem}\index{Theorem!Beppo Levi} 
\index{monotone convergence theorem}\index{Theorem!monotone convergence}
Let $\{f_n\}$ be a sequence of integrable nonnegative functions on $A$ and let
\[
0\leq f_n(x)\leq f_{n+1}(x)\quad {\rm a.e.}\ x\in A\\,\qquad \limn f_n(x)=f(x)\ {\rm a.e}\  x\in A\,.
\]
Then we have
\begin{equation}\ZLA{eq:TeoCha1bisBeppo}
\limn \int_A f_n(x)\ZD x=\int_A f(x)\ZD x\,.
\end{equation}
\end{Theorem}
\zProof  
Monotonicity of the integral gives

\[
\int_A f_n(x)\ZD x\leq \int _A f(x)\ZD x\ \ \mbox{for every $n$ so that }\ \limn \int _A f_n(x)\ZD x\leq \int_A f(x)\ZD x\,.
\] 
Fatou Lemma gives
\[
\int_{A} f(x)\ZD x\leq \liminf_{n\to +\ZIN} \int _A f_n(x)\ZD x=\limn \int _A f_n(x)\ZD x\,.
\]
The required equality follows.\zdia

\begin{Remark}\ZLA{Chap1BisRemaLemmaBeppoeFatou}
{\rm
Note that the assumption that the sequence is monotone increasing is crucial in the Beppo Levi Theorem. Let us consider the sequence of the functions $f_n =\carfunz _{[n,+\ZIN)}$. This sequence is decreasing and
\[
\limn f_n(x)=0\qquad \forall x
\]
but the integral of every $f_n$  is $+\ZIN$: Beppo Levi Theorem cannot be extended to decreasing sequences.

\emph{Note that this example shows also that in general the inequality in Fatou Lemma   is strict.}\zdia
}
\end{Remark}

 \subsection{\ZLA{Cha1bisSecLIMINTgeneSEQUELebTh}The General Case: Lebesgue Theorem}

Now we consider the case that the functions $f_n$ do not have fixed sign
and we prove Theorem~\ref{teo:CH1TeoLebe}, which we restate:

\begin{Theorem}[{\sc Lebesgue} or {\sc    dominated convergence}]\index{Theorem!Lebesgue!}\index{Theorem!dominated convergence}\index{dominated convergence}\ZLA{teo:Cha1bisTeoLebeRlimitato}
Let $\{f_n \}$ be a sequence of summable functions a.e. defined on    $A\subseteq \zzr $ and let $f_n\to f  $ a.e. on $A$. If there exists a summable nonnegative function $g $ such that
\[
|f_n ( x)|\leq g(x)\qquad a.e.\ x\in A
\]
then $f $ is summable and
\begin{equation}
\ZLA{eq:Cha1bisLebeDOMIconve0}
\lim \int_A f_n ( x)\ZD x=\int_A f(x)\ZD x\,.
\end{equation}
\end{Theorem}   
\zProof 
The function $f$ is quasicontinuous thanks to Corollary~\ref{CoroP1Ch2EsteEgSevReTTA}.

The assumption imply that every $f_n$ is summable and that
\[
|f(x)|\leq g(x)
\]
so that also $f$ is summable. Hence we can consider the  sequences $\{ g +f_n\}$ and  $\{g -f_n\}$ which are both sequences of  nonnegative functions.

We consider the sequence $\{ g +f_n\}$. Fatou Lemma gives
\begin{multline*}
\int_A g(x)\ZD x+\int_A f(x)\ZD x=
\int_A[ g(x)+f(x)]\ZD x\\\leq \liminf_{n\to+\ZIN}\int_A[ g(x)+f_n(x)]\ZD x
=\int_A  g(x)\ZD x+ \liminf_{n\to+\ZIN}\int_A f_n(x) \ZD x
\end{multline*}
so that
\begin{equation}
\ZLA{eq:Cha1bisLebeDOMIconve1}
\int_A f(x)\ZD x\leq  \liminf_{n\to+\ZIN}\int_A f_n(x) \ZD x\,.
\end{equation}

We consider the sequence $\{ g -f_n\}$. Fatou Lemma gives
\begin{multline*}
\int_A g(x)\ZD x-\int_A f(x)\ZD x=\int_A [g(x)-f(x)]\ZD x\\
\leq \liminf _{n\to +\ZIN}
\int_A [g(x)-f_n(x)]\ZD x=\int_A g(x)\ZD x-  \limsup _{n\to +\ZIN}
\int_A  f_n(x) \ZD x\,.
\end{multline*}
So we have
\begin{equation}
\ZLA{eq:Cha1bisLebeDOMIconve2}
\int_A f(x)\ZD x\geq \limsup _{n\to +\ZIN}
\int_A  f_n(x) \ZD x\,.
\end{equation}
The required equality~(\ref{eq:Cha1bisLebeDOMIconve0}) follows from~(\ref{eq:Cha1bisLebeDOMIconve1}) and~(\ref{eq:Cha1bisLebeDOMIconve2}).\zdia

\begin{Remark}{\rm
The result which holds for bounded sequences on bounded sets has been used to  prove Fatou Lemma, which is then used to prove   Beppo Levi and Lebesgue theorems. Actually, these three results, by Fatou, by  Beppo Levi  and by Lebesgue can be proved in different order. The reader can   see for example the books~\cite{DOOBbookMEASURE94,MalliavinBOKKintegr95,Royden66LibroANALreale}. From an historical point of view, Beppo Levi theorem was proved first, after the proof of Lebesgue concerning the bounded case.\zdia
}
\end{Remark}

\section{\ZLA{AppeCha1bisEgSEV}\textbf{Appendix:}  Egorov-Severini:  Preliminaries  in One Variable}
In this section we prove Theorem~\ref{Cha1bisTeoEgoSEVERANTE} for a sequence of functions of one variable.  
  
  We split the proof  in several parts.

 \subsection{Convergent and Uniformly Convergent Sequences}

Let $\{f_n \}$ be a sequence of functions defined on a set $R$.
We recall Cauchy theorem and we recast convergence of $\{f_n(x)\}$  without the explicit use of the limit function.   The sequence $\{ f_n(x)\}$ converges if and only if for every $\ZEP>0$ there exists $M=M(\ZEP,x)$ such that  
\begin{equation}\ZLA{Ch3APPEtestCauchy}
m\geq M(\ZEP,x)\,,\ r\,,\ s\in\zzn\ \implies\ 
|f_{m+r}(x)-f_{m+s}(x)|<\ZEP\,.
\end{equation}

The sequence $\{f_n\}$ \emph{converges uniformly on a set $S\subseteq R$} when 
 for every $\ZEP>0$ there exists $M=M(\ZEP) $  such that  
\begin{equation}\ZLA{Ch3APPEtestCauchyUNIF}
m\geq M(\ZEP )\,,\ r\,,\ s\in\zzn\ \implies\ 
|f_{m+r}(x)-f_{m+s}(x)|<\ZEP\qquad \forall x\in S\,.
\end{equation}
We stress the fact that $M(\ZEP)$ does not depend on $x$.
 
This observation suggests the introduction of the following functions: 
\begin{equation}\ZLA{Ch3AppeDifivNM}
\begin{array}{l}\displaystyle
 v_{n,m} (x)=\max\{
  |f_{m+r}(x)-f_{m+s}(x)| \quad 1\leq r<n\,,\qquad 1\leq s<n\}\,,\\[2mm]
\displaystyle
  w_m(x)=\sup\{v_{n,m}(x)\quad n\geq 1\}\leq +\ZIN\,.
  \end{array}
\end{equation}
These functions have the properties described in the following lemma:

\begin{Lemma}\ZLA{LemmaAppe3ProprVnmWm}
Let $\{f_n\}$ be a sequence of functions defined on a set $R$ and let $v_{n,m}(x)$, $w_m(x)$ be the functions defined in~(\ref{Ch3AppeDifivNM}). We have:

\begin{enumerate}
\item\ZLA{I1LemmaAppe3ProprVnmWm} monotonicity properties:
\begin{enumerate}\ZLA{I1VdipNLemmaAppe3ProprVnmWm}
\item
for every $x\in R$ and every $m$, the sequence $n\mapsto  v_{n,m}(x)$ is increasing.

\item \ZLA{I1WdipMLemmaAppe3ProprVnmWm}  the sequence $m\mapsto w_m(x)$ is decreasing.
\end{enumerate}
(monotonicity of the two sequences needs not be strict).
 
\item \ZLA{I2LemmaAppe3ProprVnmWm} the convergence of the sequence $n\mapsto  f_n(x) $ for a fixed value of $x$:

\begin{enumerate}
\item\ZLA{I2DIRELemmaAppe3ProprVnmWm} 
let   $n\mapsto f_n(x)$ converge. Then the sequence $n\mapsto v_{n,m}(x)$ is bounded so that  
\begin{equation}\ZLA{EqAI2DIRELemmaAppe3ProprVnmWm} 
w_m(x)=\sup _{n>0}v_{n,m}(x)=\lim _{n\to+\ZIN}v_{n,m}(x)   \in\zzr
\end{equation}
and we have also
\begin{equation}\ZLA{EqABI2DIRELemmaAppe3ProprVnmWm}
\lim_{m\to+\ZIN  } w_m(x)=0\,.
\end{equation}

 \item \ZLA{I2VICEVLemmaAppe3ProprVnmWm} 
let $\lim _{m\to +\ZIN} w_m(x)=0$. Then the sequence $n \mapsto f_n(x)$ converges.
 \item\ZLA{I3VICEVLemmaAppe3ProprVnmWm}let $S$ be a subset of $R$. The sequence $\{f_n\}$ converges uniformly on $S$ if and only if the sequence $\{w_m\}$ converges to $0$ uniformly on $S$.
\end{enumerate}

\item\ZLA{I3LemmaAppe3ProprVnmWm} 
 let  $S\subseteq R$. Let us assume\footnote{we are not assuming continuity on $ S $ of the functions $f_n$. We assume the weaker property that their restrictions are continuous.} that each     function $\left (f_n\right )_{|_{S}}$ be continuous on $S$ and let $\ZEP>0$.  
 The set
\[
A_{m,\ZEP}=\{x\in S\,, w_m(x)>\ZEP\}\,.
\]
  is relatively open in $S$.
\end{enumerate}
\end{Lemma}

\zProof The monotonicity property of $v_{n,m}(x)$ is obvious. We prove that  $m\mapsto w_m(x)$ is decreasing. We note:
\begin{multline*}
w_{m+1}(x)=
\sup _{n\geq 1}
\{ 
 \sup
 \{
 \left | f_{m+1+r}(x)-f_{m+1+s}(x)\right| \quad 1\leq r<n\,,\ 1\leq s<n\}
\,
\}
\\
 =
 \sup _{n\geq 1}
\{ 
 \sup
\{
 \left | f_{m +r}(x)-f_{m +s}(x)\right| \quad 2\leq r<n+1\,,\ 2\leq s<n+1\}
 \,
\}\\
\leq
 \sup _{n\geq 1}
\{ 
 \sup
\{
 \left | f_{m +r}(x)-f_{m +s}(x)\right| \quad 1\leq r<n+1\,,\ 1\leq s<n+1\}
 \,
\}
\\
=
 \sup _{n\geq 1}
\{ 
 \sup
\{
 \left | f_{m +r}(x)-f_{m +s}(x)\right| \quad 1\leq r<n \,,\ 1\leq s<n \}
 \,
\}=
  w_m(x)\,.
\end{multline*}

We prove statement~\ref{I2DIRELemmaAppe3ProprVnmWm}. 
Boundedness holds because
  
  \[
v_{n,m} (x)\leq 2\sup\{| f_k(x)|\,,\ k\geq 1\}<+\ZIN\quad \mbox{since \{$f_k(x)$\} is convergent} \,.
\]
So,~(\ref{EqAI2DIRELemmaAppe3ProprVnmWm}) holds. Property~(\ref{EqABI2DIRELemmaAppe3ProprVnmWm}) is seen by contradiction:
 if the limit is $l_0>0$ then  \emph{for every
$m $} we have
\[
w_m(x)=\limn v_{n,m}(x)>l_0>0
\]
and \emph{for every
$m $}   there exists 
 $N=N_m$ such that when $n>N_m$ we have
\[
v_{n,m}(x)>\dfrac{l_0}{2} \,.
\]
So, \emph{for every $m $} there exist $r$ and  $s$ such that
\[
|f_{m+r}(x)-f_{m+s}(x)|>\dfrac{l_0}{4}\,.
\]
This is not possible since $\{f_n(x)\}$ is a convergent sequence, see~(\ref{Ch3APPEtestCauchy}).

We prove statement~\ref{I2VICEVLemmaAppe3ProprVnmWm}. If $w_m(x)\to 0$ then for every $\ZEP>0$ there exists $M=M(\ZEP,x)$ such that  
\[
m>M(\ZEP,x)\ \implies w_m(x)<\ZEP
\]
so that, when $m>M(\ZEP,x)$ we have also $v_{n,m}(x)<\ZEP$ for every $n$ i.e.
\[
m>M(\ZEP,x)\ \implies |f_{m+r}(x)-f_{m+s}(x)|<\ZEP\quad \forall r\,, \ s\,.
\]
This is Cauchy condition of convergence.
 
 We prove the statement~\ref{I3VICEVLemmaAppe3ProprVnmWm}. We note that a sequence which converges uniformly to $0$ is bounded. For this reason boundedness of $\{w_m\}$ has not been explicily stated.
 
 Let $\{f_n\}$ be uniformly convergent.
The definition of $w_m$ and condition~(\ref{Ch3APPEtestCauchyUNIF}) shows that for $m>M(\ZEP)$ we have $0\leq w_m(x)<\ZEP$ for every $x\in S$. Hence, if $\{f_n\}$ converges uniformly on $S$ then $\{w_m\}$ converges to zero uniformly on $S$.
 
Conversely, let $\{w_m\}$ converge uniformly to zero. Then for every $\ZEP>0$ there exists $M(\ZEP)$ such that
\[
m>M_\ZEP\ \implies 0\leq w_m(x)<\ZEP\qquad \forall x\in S\,.
\]
And so when $m>M(\ZEP)$ we have also \[
 0\leq v_{n,m} (x)\leq \ZEP\qquad \forall x\in S\,,\qquad \forall n>0
\]
and this is the property that $\{f_n\}$ is uniformly convergent on $S$.

Finally we prove statement~\ref{I3LemmaAppe3ProprVnmWm}. 
   Let $x_0\in A_{m,\ZEP}$ so that
\[
w_m(x_0)=\ZEP_1>\ZEP\,.
\]
We prove the existence of an open set $B$ such that $x_0\in   A_{m,\ZEP}\cap B$.

We fix any $\ZSI_1\in (\ZEP,\ZEP_1)$. There exists $n_1$ such that $v_{n_1 ,m}(x_0)>\ZSI_1$ and so there exist $r(n_1)$ and $s(n_1)$ such that
\[
|f_{m+r(n_1)} (x_0)-f_{m+s(n_1)} (x_0)|> \ZSI_1\,.
\]
 
 The restriction to $S$ of the  function $x\mapsto |f_{m+r(n_1)} (x )-f_{m+s(n_1)} (x )|$ is continuous. So, for every $\ZSI_2\in (\ZEP,\ZSI_1)$ there exists an open ball $B$ of center $x_0$ such that  
 \[
|f_{m+r(n_1)} (x )-f_{m+s(n_1)} (x )|> \ZSI_2>\ZEP\qquad \forall x\in B\cap S \,.
\]
Hence,
\[
x\in B\cap S\ \implies\ w_m(x)>\ZSI_2>\ZEP
\]
and so $x\in A\cap B$: any point   $x_0\in A$ belongs to its relative interior, as wanted.\zdia
  
\begin{Remark}\ZLA{RemaP1C2AppeSULLAnonNECEinte}  
  {\rm
  Note that nowhere in this proof we used the fact that $R\subseteq \zzr$. The result holds and with the same proof (a part interpreting $|\cdot|$ as a norm) also if $R\subseteq \zzr^d$, any $d\geq 1$ (and in fact if $R$ is any normed space).\zdia
  }\end{Remark}
\ZREF{  
We stress the fact that the statements~\ref{I2DIRELemmaAppe3ProprVnmWm},~\ref{I2VICEVLemmaAppe3ProprVnmWm}  and~\ref{I3VICEVLemmaAppe3ProprVnmWm} of Lemma~\ref{LemmaAppe3ProprVnmWm} recast pointwise convergence of the sequence $\{f_n(x)\}$ in terms of the convergence to $0$ of the sequence $\{w_m(x)\}$ and uniform convergence of $\{f_n\}$ in terms of uniform convergence to zero of $\{w_m\}$.

}

\subsection{\ZLA{P1Ch2AppeSEctSulMultiint}The Two Main Lemmas}

We recall that any open set in $\zzr$ is the union of a sequence of disjoint open intervals,

\begin{equation}\ZLA{eq:CH2AppeRichiamAPErti}\mathcal{O}=\bigcup_{n\geq 1} R_n\,,\quad \mbox{$R_n$ pairwise disjoint open intervals}\,.
\end{equation}
Thanks to this observations, we can use Definition~\ref{DefiCH1DeFiMiSuOpen}: we represent   $\mathcal{O}$ as in~(\ref{eq:CH2AppeRichiamAPErti}) and we put
\begin{equation}\ZLA{eq:CH2AppeRichiamAPErtimeas}
\zl(\mathcal{O})=\sum_{n\geq 1}L(R_n)=L(\{R_n\})\,.
\end{equation}

In this section we prove two lemmas.
The first lemma concerns decreasing sequences   of open sets. In order to understand this lemma we must keep in mind that a decreasing sequence of open sets may have empty intersection, as in the example $\mathcal{O}_n=(0,1/n)$. But, in this example we have $\zl(\mathcal{O}_n)\to 0$. Instead we have the following lemma which we formulate both in terms of multiintervals and in terms of open sets:
\begin{Lemma}\ZLA{LemmaiAppeCh3SUOmultirettANNID}
The following (equivalent) statements hold:
 \medskip\begin{center}
\fbox{
\parbox{2.5in}{Let  
  $\{\Delta_n\}$ be a sequence of disjoint     multiintervals in $\zzr^d$. We assume:
\begin{enumerate}
\item\ZLA{I1LemmaiAppeCh3} the existence of a \emph{bounded} interval $R$ such that $\mathcal{I}_{\Delta_1}\subseteq R$; 
\item\ZLA{I0LemmaiAppeCh3}    every multiinterval  $\Delta_n$ is   disjoint; 
 
\item\ZLA{I2LemmaiAppeCh3} 
$\mathcal{I}_{\Delta_{n+1}}\subseteq \mathcal{I}_{\Delta_n}$ for every $n$, 
 \item\ZLA{I3LemmaiAppeCh3}   there exists     $l >0$  such that $L(\Delta_n)>l $ for every $n$.
  
\end{enumerate}
under these conditions $\cap _{n\geq 1} \mathcal{I}_{\Delta_n}=\emptyset$.
} 
}
 \fbox{
 \parbox{2in}{Let $\{\mathcal{O }_n\}$ be a   sequence of open sets. We assume:
 \begin{enumerate}
 \item  there is a \emph{bounded} interval $R$ such that $\mathcal{O}_1\subseteq R$;
 \item the sequence is decreasing, i.e. , $\mathcal{O}_{n+1}\subseteq \mathcal{O}_n$ for every $n$;
 \item there exists $l>0$ such that $\zl(\mathcal{O})>l$.
 \end{enumerate}
under these conditions $\cap_{n>1}\mathcal{O}_n\neq \emptyset$.
 } 
 }
 \end{center}
\end{Lemma} 
\zProof 
We divide the proof in two parts: first we present few preliminary observations and then we use them to prove the lemma.

\paragraph{Preliminary Observations on Sequences of Intervals}
 We assume that $\Delta$ is a \emph{disjoint} multiinterval (which is the case of interest in the proof).
\begin{enumerate}
\item\ZLA{P1LemmaiAppeCh3SUOmultirettANNID} The length of an interval does not change if we add (or remove) the endpoints. So we can define \[
L([h,k])=L((h,k])=L([h,k))=L((h,k))=k-h\,.
\]
\item\ZLA{P1bisLemmaiAppeCh3SUOmultirettANNID} If $R$ is an open interval and $\ZEP>0$ there exists a closed interval $S\subseteq R$ such that $L(S)>L(R)-\ZEP$.
\item\ZLA{P1LemmaiAppeCh3SUOmultirettANNIDa} if $S$ is an open interval and $\Delta=\{R_n\}$ is a multinterval, the multiinterval $\{S\cap R_n\}$ is denoted $S\cap \Delta$.

We use the notation $S\cap\Delta$   to denote the sequence $\{S\cap R_n\}$ also in the case that $S$ is closed or half closed. Then we can extend the definition of $L$ and we can define
\[
L(S\cap \Delta)=\sum _{n\geq 1} L(S\cap R_n)\,.
\]
 It is clear that
 \[
L(S\cap\Delta)\leq \min\{L(S)\,, \ L(\Delta)\}\,. 
 \]
\item\ZLA{P3LemmaiAppeCh3SUOmultirettANNID}
Let $\Delta=\{R_n\}$ be a multinterval such that $L(\Delta)<+\ZIN$. and
let $\ZEP>0$. There exists a sequence of closed intervals $S_n=[h_n,k_n]\subseteq R_n$ such that
\[
\sum _{n\geq 1} (k_n-h_n)>L(\Delta)-\ZEP\,.
\]
\item\ZLA{P3BISLemmaiAppeCh3SUOmultirettANNID} Let $\tilde\Delta=\{\tilde R_n\}$ and $ \Delta=\{  R_n\}$ be two disjoint multiintervals. We assume
\[
\mathcal{I}_{\tilde \Delta}=\bigcup _{n\geq 1}\tilde R_n
\subseteq \bigcup _{n\geq 1}  R_n=\mathcal{I}_{  \Delta}\,.
\]
This inclusion shows
\begin{equation}
\ZLA{P3eqBISLemmaiAppeCh3SUOmultirettANNID} 
L(\tilde\Delta)=\sum _{n\geq 1}\left (\sum _{k\geq 1} L(\tilde R_k\cap R_n)\right )=
\sum _{n\geq 0} L(R_n\cap\Delta)\,.
\end{equation}
\item\ZLA{P4LemmaiAppeCh3SUOmultirettANNID} If $S\subseteq R$ and $\Delta=\{R_n\}$ then we have
\[
\bigcup _{n\geq 1} (R\cap R_n) =
\left (\cup_{n\geq 1} \left (  S\cap R_n\right )\right ) \bigcup \left ( \cup_{n\geq 1} (R\setminus S) \cap R_n\right )
\]
and the intervals which appears in these expressions are pairwise disjoint (since $\Delta$ is disjoint). So we have
\begin{equation}\ZLA{eq:3APPeSec1Punto5}
L(R\cap\Delta_n)=L(S\cap \Delta_n)+L((R\setminus S)\cap\Delta_n)
\end{equation}
Note that the intervals in $S\cap \Delta_n$ and in $(R\setminus S)\cap\Delta_n$
need not be open and so we used the extended definitions introduced at the point~\ref{P1LemmaiAppeCh3SUOmultirettANNID}.
\end{enumerate}
 
 \ZREF{
 In order to facilitate the use the previous observations \emph{in the proof of Lemma~\ref{LemmaiAppeCh3SUOmultirettANNID}} it is conveniente to call ``multiinterval'' any sequence of intervals, possibly not open.}

 After these preliminaries we prove the Lemma.
 
 \paragraph{The proof of the Lemma~\ref{LemmaiAppeCh3SUOmultirettANNID}}

We note that Assumption~\ref{I2LemmaiAppeCh3} and the monotonicity of the measure  imply that the sequence $\{L(\Delta_n)\}$ is decreasing. We put
\[
\limn L(\Delta_n)=l_0 \quad\mbox{(assumption~\ref{I3LemmaiAppeCh3} implies $l_0\geq l>0$)}\,.
\]
 We introduce the notation
\[
\Delta_n=\{R_{n,k}\}_{k\geq 1}\,.
\]
and we proceed with the following steps.
\begin{description}
\item[{\bf Step~1:}] We prove the existence of $R_{1,k_1}$ and of~$\tilde l_1>0$ such that
\begin{equation}\ZLA{eq:CH3APPE1k1INTEmaggio}
L(R_{1,k_1}\cap \Delta_n)>  \tilde l_1>0 \qquad \forall n\geq 2\,.
\end{equation}
 Then we prove the existence of a \emph{closed} interval $S_1$ such that
\[
\begin{array}{l}S_1 \subseteq {\rm int}\, R_{1,k_1} \quad \mbox{such that}\quad
L(S_1\cap\Delta_n)\geq \hat l_1>0\quad \forall n\geq 2\,.   \\
 \mbox{So:}\quad \limn L(S_1\cap \Delta_n)=l_1\geq \hat l_1>0  \,.
\end{array}
\]

Note that these conditions imply
\[
S_1\neq \emptyset\,,\qquad L(S_1)>\hat l_1>0
\]
(it will be  $\hat l_1=\tilde l_1/2=l/4$ but the actual value has no role. The important point is that it is positive).

\emph{The proof  that these intervals exist is  postponed.}

\item[{\bf Step~2:}] We consider the new sequence of multiintervals $\{S_1\cap \Delta_n\} _{n\geq 2}$. The properties~\ref{I1LemmaiAppeCh3}-\ref{I3LemmaiAppeCh3} hold for this sequence and so we can apply the procedure in the {\bf Step $1$ } to the sequence $\{S_1\cap \Delta_n\}_{n\geq 2}$: there exists a interval $R_{2,k_2} $ such that
\[
L\left ((S_1\cap R_{2,k_2}  )\cap   (S_1\cap  \Delta_n)\right )>  \tilde  l_2>0 \qquad \forall n\geq 3\,.
\]
 Then we choose a \emph{closed} interval
$S_2$ such that 
  $S_2 \subseteq {\rm int}\, S_1\cap R_{2,k_2}$ and such that
\[
\begin{array}{l}
L\left (S_2\cap(S_1\cap \Delta_n )  \right )\geq \hat l_2>0\quad \forall n\geq 3\quad \mbox{and we put}\\
\limn L\left (S_2\cap (S_1\cap \Delta_n)\right )=  l_2>0\,.
\end{array}
\]
  Note that  
\[
\begin{array}{ll}
\mbox{$S_1$ closed nonempty}\,,\quad  L(S_1)>\hat l_1>0\,, 
S_1 \subseteq {\rm int}\, R_{1,K_1}\,,\\
\mbox{$S_2$ closed nonempty}\,, \quad  L(S_2)>\hat l_2>0\,, \\
 \left\{
\begin{array}{l}
S_2 \subseteq {\rm int}\, S_1 \cap {\rm int}\, R_{2,k_2}\subseteq {\rm int}\, R_{2,k_2} \,,\\
S_2 \subseteq {\rm int}\, S_1 \cap {\rm int}\, R_{2,k_2}\subseteq {\rm int}\, S_{1}
\subseteq R_{1,k_1} 
\end{array}\right.\quad {\rm i.e.}\ S_2\subseteq \left ({\rm int}\, R_{1,k_1}\right )\cap\left ({\rm int}\,  R_{2,k_2 }\right)\,.
\end{array}
\]

The assumptions~\ref{I1LemmaiAppeCh3}-\ref{I3LemmaiAppeCh3} of the lemma hold for the sequence   $\left \{S_2\cap (S_1\cap \Delta_n)\right \}$ and \emph{the procedure can be iterated.}

\item[In conclusion:]
  We 
single out a sequence of intervals $\{R_{n,k_n}\}$ and we
construct a sequence $\{S_n\}$ of \emph{closed} intervals such that 
\[
\begin{array}{l}
L(S_n)>0 \ \mbox{for every $n$, hence  $S_n\neq \emptyset$ }\ 
\,,\\
 S_n \subseteq \bigcap _{j=1}^n \left ({\rm int}\, R_{j,k_j}\right )\,,\quad S_{n+1}\subseteq{\rm int}\, S_{n}\,.
\end{array}
\]

Assumption~\ref{I1LemmaiAppeCh3} implies that the interval $S_1$ is bounded and, as we noted,  the closed sets $S_n$ are nonempty. The sequence $\{S_n\}$ is decreasing. Cantor theorem implies the existence of $x_0\in \cap _{n\geq 1} S_n$.
 
In particular we have
\[
x_0\in  S_j\subseteq  {\rm int}\, R_{j,k_j}\subseteq {\rm int}\, \mathcal{I}_{\Delta_j} \ \mbox{for every $j$}
\]
as wanted.
\end{description}
\emph{In order to complete the proof it is sufficient to prove the existence of the intervals $R_{1,k_1}$ and $S_1$ described in the {\bf Step~1.}}

First we prove the existence of $k_1$ such that
\[
L(R_{1,k_1}\cap \Delta_n)>l/2=\tilde l_1 \qquad \forall n\geq 2\,.
\]
 The assumptions~\ref{I1LemmaiAppeCh3} and~\ref{I0LemmaiAppeCh3}      imply 
\[
\sum _{k\geq 1}  L(R_{1,k})<+\ZIN  
\]
so that there exists $m_1$ such that
\begin{equation}\ZLA{eq:CH3APPE1k1INTEmaggioBIS}
\sum _{k\geq m_1+1}  L(R_{1,k})< \dfrac{l}{2}\,.
\end{equation}
 
The assumptions~\ref{I0LemmaiAppeCh3} and~\ref{I2LemmaiAppeCh3} 
and the equality~(\ref{P3eqBISLemmaiAppeCh3SUOmultirettANNID}) 
imply 
 that for every $n\geq 2$ we have

\[
L(\Delta_n)=\sum _{k\geq 1}L(R_{1,k}  \cap \Delta_n)
\]
so that from~(\ref{eq:CH3APPE1k1INTEmaggioBIS}) and Assumption~\ref{I3LemmaiAppeCh3}, we have 
 \begin{subequations}
 \begin{align}
 \ZLA{eq:subeqAppe3Sub1}
 &\sum _{k\geq m_1+1}L(R_{1,k}\cap \Delta_n) <l/2 \quad \forall n\geq 2\\
 \ZLA{eq:subeqAppe3Sub2}
 &\sum _{k=1}^m L(R_{1,k}\cap \Delta_n) >l/2\quad \forall n\geq 2\,.
 \end{align}
 \end{subequations}
We show:
\begin{equation}\ZLA{eq1APPELemmaiAppeCh3}
\exists k_1\leq m\,:\ \forall n\geq 2\ \mbox{we have}\ L\left (R_{1,k_1}\cap \Delta_n\right )>   \frac{l}{4m}=\tilde l_1\,.
\end{equation} 
The proof is by contradiction. If this property does not hold then for every $k\leq m $ there exists $n_k\geq 2$ such that
\[
 L\left (R_{1,k }\cap \Delta_{n_k}\right )\leq  \frac{l}{4m}\,.
\]
Let 
\[
N=\max\{ n_1\,,\dots\,, n_m\}\,.
\]
 
The inclusion $\mathcal{I}_{\Delta_N}\subseteq \mathcal{I}_{\Delta_{n_j}} $ and the monotonicity of the measure  imply
  
\[
\sum _{k=1} ^m L(R_{1,k}\cap \Delta_N)
\leq \sum _{k=1}^m \underbrace{L(R_{1,k}\cap \Delta_{n_k})}_{\leq l/4m}
 \leq \dfrac{l}{4}
\]
in contrast with~(\ref{eq:subeqAppe3Sub2}).
 
So, property~(\ref{eq1APPELemmaiAppeCh3}), i.e.~(\ref{eq:CH3APPE1k1INTEmaggio}), holds.
 
We use the {\bf Preliminary observation}~\ref{P1bisLemmaiAppeCh3SUOmultirettANNID}  and we  choose any \emph{closed} interval $S_1$ such that  
 
\[
S_1 \subseteq {\rm int}\, R_{1,k_1} 
\,,\quad
L(R_{1,k_1}\setminus S_1)< \frac{\tilde l_1}{2}\,.
\] 
We use~(\ref{eq:3APPeSec1Punto5}) and we get
 
\begin{multline*}
\tilde l_1 \leq L(R_{1,k_1}\cap \Delta_n) =L(S_1\cap \Delta_n)+L((R_{1,k_1}\setminus S_1)\cap \Delta_n)   
\leq L(S_1\cap \Delta_n)+\frac{\tilde l_1}{2}
\end{multline*} 
i.e.
\[
L(S_1\cap \Delta_n) \geq \dfrac{\tilde l_1}{2}=\hat l_1 >0\qquad \forall n\geq 2\,.
\]
The proof is finished.\zdia
 
\begin{Remark}\ZLA{Rema:AppeP1Ch2RemaRSullaNeceLimit}
{\rm
We observe that the lemma does not hold if the   boundedness assumption~\ref{I1LemmaiAppeCh3} is removed, as it is seen by considering the sequence of the intervals $[n,+\ZIN)$.\zdia
}
\end{Remark}

The second lemma is a weaker version of Theorem~\ref{Cha1bisTeoEgoSEVERANTE} and in fact it is the core of its proof.

\begin{Lemma}\ZLA{LemmaiAppeCh3PRELIsucceFUNcont}
  Let $R=[h,k]$ be a bounded and closed interval and let $\{f_n\}$ be a sequence of \emph{continuous} functions   defined on $R$. We assume that the sequence $\{ f_n(x)\}$ converges for every $x\in R$. 
 
We prove that for every pair of positive numbers $\zg >0$ and $\eta >0$ there exist 
\begin{description}
\item[\bf 1:] a  disjoint  multiinterval $\Delta_{\zg ,\eta } $ (which is composed of open intervals)      such that
\[
L(\Delta_{\zg ,\eta })=\zl(A_{\zg,\eta})<\zg\quad \mbox{(where we put  $A_{\zg,\eta}=\mathcal{I}_{\Delta _{\zg,\eta}}$)}\,.
\]
\item[\bf 2:] a number $M_{\zg ,\eta }$  such that 
  
\[
\left\{\begin{array}{l}
 x\in R\setminus\mathcal{I}_{ \Delta_{\zg ,\eta }}= R\setminus A_{\zg,\eta}\\
 m> M_{\zg ,\eta } \\ 
 r>0\,,\ s>0
 \end{array}\right.\quad \ \implies \ |f_{m+r}(x)-f_{m+s}(x)|<\eta  \,.
\]
\end{description}
\end{Lemma}
 \zProof   
We recast the thesis of the lemma   in terms of the functions  $w_m(x)$ defined in~(\ref{Ch3AppeDifivNM}) as follows:
  \emph{ for every $\zg >0$ and $\eta >0$ there exists a  disjoint  multiinterval $\Delta_{\zg ,\eta } $ (composed of \emph{open} intervals) such that
\begin{align*}
&L(\Delta_{\zg ,\eta })<\zg \quad\mbox{and}\\
&\left\{\begin{array}{l}
\displaystyle x\in R\setminus\mathcal{I}_{\Delta_{\zg ,\eta }}  \\
\displaystyle
m>M_{\zg,\eta}
\end{array}\right.\  \implies w_m(x)<\eta
\end{align*}
 }
This  we prove now.

Let $\eta >0$ and\footnote{note that in the definition below $A_{m,\eta}$ is defined in terms of $(h,k)$, the interior of the interval $[h,k]$.}
\begin{equation}\ZLA{ChP1Ch2eqRenaImpoNOunic}
A_{m,\eta }=\{x\in (h,k)\,, w_m(x)>\eta/2 \}\,.
\end{equation}
We proved in the statement~\ref{I3LemmaAppe3ProprVnmWm}  of Lemma~\ref{LemmaAppe3ProprVnmWm}  that the set $A_{m,\eta }$ is   open.    So, 
  there exists a  disjoint multiinterval $\Delta_{m,\eta }$ such that
$A_{m,\eta }=\mathcal{I}_{\Delta_{m,\eta }}$ (see Theorem~\ref{teo:Ch2StrutturaAperti}).

If $x\in (h,k)\setminus \mathcal{I}_{\Delta_{m,\eta }}= 
R\setminus A_{  m,\eta } 
$ we have
\[
  w_m(x)\leq \eta/2  \,.
\]
We   note
\begin{equation}\ZLA{eq2LemmaiAppeCh3PRELIsucceFUNcont}
\limm L( \Delta_{m,\eta } )=0\,.
\end{equation}
In fact, the sequence $m\mapsto  w_m(x)$ is decreasing so that
\[
\mathcal{I}_{\Delta _{{m+1},\eta }}=A_{m+1,\eta }\subseteq A_{m,\eta }\subseteq \mathcal{I}_{\Delta_{m,\eta }}\,.
\]
It follows that $\{ L(\Delta_{m,\eta })\}=\{\zl(A _{m,\eta })\}$ is decreasing too
and $\lim _{m\to+\ZIN}L(\Delta_{m,\eta })$ exists.

If the limit is positive then there exists $l >0$ such that
\[
L(\Delta _{m,\eta })>l>0 \qquad \forall m
\]
and, from Lemma~\ref{LemmaiAppeCh3SUOmultirettANNID} there exists $x_0\in  \mathcal{I}_{\Delta_m,\eta }$ for every $m$. So, for every $m$ we have
$w_m(x_0)>\eta/2 $.
Statement~\ref{I2DIRELemmaAppe3ProprVnmWm} of Lemma~\ref{LemmaAppe3ProprVnmWm} shows that the sequence $\{ f_n(x_0)\}$ is not convergent, in contrast with the assumption. So it must be
\begin{equation}\ZLA{ChP1Ch2eqBISRenaImpoNOunic}
\limm L(\Delta_{m,\eta })=\limm\zl(A_{m,\eta })=0\,.
\end{equation}
It follows that there exists $M_{\zg ,\eta } $ such that when $m>M_{\zg ,\eta }$
then we have
\[
\begin{array}{l}
\displaystyle
L(\Delta _{m,\eta })<\zg \,,\\
x\in R\setminus \mathcal{I}_{\Delta(m,\eta )}\  \implies\ w_m(x)\leq \eta/2<\eta \,.\zdiaform
\end{array}
\]
\begin{Remark}[Important observation]\ZLA{ChP1Ch2RenaImpoNOunic}{\rm
The sets $A_{m,\eta}$ can be chosen with different laws, for example by replacing $\eta/2$ in~(\ref{ChP1Ch2eqRenaImpoNOunic}) with $\eta/3$, or by replacing  $A_{m,\eta}$ in~(\ref{ChP1Ch2eqRenaImpoNOunic}) with larger open sets, provided that~(\ref{ChP1Ch2eqBISRenaImpoNOunic}) holds. So, the number $M_{\zg,\eta }$  does depend also the chosen set  $A_{m,\eta}$.\zdia

}
\end{Remark}

\subsection{\ZLA{AppeP2CH2FiNalStep}From Lemma~\ref{LemmaiAppeCh3PRELIsucceFUNcont}   to Theorem~\ref{Cha1bisTeoEgoSEVERANTE}: }

For clarity in the next box we report the statement of the weaker version  in Lemma~\ref{LemmaiAppeCh3PRELIsucceFUNcont} and
 we recast statement~\ref{I1Cha1bisTeoEgoSEVER} of  Theorem~\ref{Cha1bisTeoEgoSEVERANTE} with the notations in the Definition~\ref{Cha1bisDEFINITalmUNIFconve} but in terms of the functions $w_m$ defined 
  in~(\ref{Ch3AppeDifivNM}).
  
   {\footnotesize
  \begin{center}
  \ZREF{
	The functions $w_m$ are defined on a bounded closed interval $R$.
	\begin{center}
  \fbox{\parbox{2.3in}{ 
\begin{center}{\bf Under the assumptions of Lemma~\ref{LemmaiAppeCh3PRELIsucceFUNcont} we proved}\end{center}

For every $\zg>0$ and $\eta>0$ there exist an \emph{open} set $A_{\zg,\eta}$ and a
number $M_{\zg,\eta}$ such that 
\[
\left\{
\begin{array}{l}
\zl(A_{\zg,\eta})<\zg\\
\left\{\begin{array}{l}
\mbox{if $m>M_{\zg,\eta}$ and $x\notin A_{\zg,\eta}$}\\
\mbox{then $w_m(x)<\eta$}\,.
\end{array}\right.
\end{array}\right.
\]
 }} 
  \fbox{\parbox{2.3in}{
\begin{center}{\bf We must prove}\end{center}
 
 \vskip -.1cm
  \[\begin{array}{l}
\mbox{$\forall\ZEP>0\ \exists\  \mathcal{O}_\ZEP$   such that} \\
\mbox{$\mathcal{O}_\ZEP$ is open and  $\zl(\mathcal{O}_\ZEP)<\ZEP$ and} \\
\mbox{$\forall\ZSI>0\ \exists M>0$ such that}\\
\left\{\begin{array}{l}
\mbox{if $x\in R\setminus \mathcal{O}_\ZEP$, $m>M$}\\
\mbox{then  $  \ w_m(x)<\ZSI$ }\,.
\end{array}\right.
\end{array} 
  \]
  The number $M$ depends
on the previously chosen and fixed set $\mathcal{O}_\ZEP$ and  
   on~$\ZSI$. 
   }} 
   \end{center}
    \emph{The important fact to be proved is that $\mathcal{O}_\ZEP$ does not depend on $\ZSI$.}
  }
 {~}
  \end{center}
 }
 \smallskip

%
%
%
%
%
%
%

The proof   consists in this: we devise a procedure to replace the sets $A_{\zg,\eta}$, which depends on the two parameters $\zg$ and $\eta$, with a set $\mathcal{O}$ which depends solely on one parameter $\ZEP$.  The procedure uses the following steps.

\begin{description}
\item[\bf Step~1:] In this step we use the given number $\ZEP>0$. The number $\ZSI$ is not used. For every natural number $n$ we  do the following:
\begin{description}
\item[\bf Step~1A:] we apply Lemma~\ref{LemmaiAppeCh3PRELIsucceFUNcont} with
\[
\zg=\dfrac{\ZEP}{2^n}\,,\quad \eta=\dfrac{1}{n}\,.
\]
We find an \emph{open} set $A_{\zg,\eta}=A_{\ZEP,n} $ and a number $M_{\zg,\eta}=M_{\ZEP,n} $ such that
\[
x\notin A_{\ZEP,n}\,,\qquad m>M_{\ZEP,n}\ \implies w_m(x)<\dfrac{1}{n}\,.
\]
\item[\bf Important observation:] We recall from Remark~\ref{ChP1Ch2RenaImpoNOunic} that the number $M_{\ZEP,n }$ does depend also on the set $A_{\ZEP,n }$. So, $M_{\ZEP,n }=M_{n ,A_{\ZEP,n }}$.
 
\item[\bf Step~1B:] We define
\[
\mathcal{O}_\ZEP=\bigcup _{n=1}^{+\ZIN} A_{\ZEP,n}\quad \mbox{so that} \quad \left\{\begin{array}{l} \mbox{$\mathcal{O }_\ZEP$ is open}\,,\\
\zl(\mathcal{O}_\ZEP)<\ZEP\,.
\end{array}\right.
\]

For every $n_0$ the following holds: 
if $m>M_{\ZEP,n_0}$  and $x\notin\mathcal{O}_\ZEP$ then $w_m(x)<1/n_0$ since 
\[
x\notin\mathcal{O}_\ZEP\ \implies\ x\notin A_{\ZEP,n_0}\,.
\]
  And, we again note that $M_{\ZEP,n_0}$ does depoen on $\mathcal{O}_\ZEP$, 
  \[
M_{\ZEP,n_0}=M_{n_0,\mathcal{O}_\ZEP}\,.  
  \]
\end{description}
 
\item[\bf Sep~2:] 
Neither $\ZEP$ nor $\mathcal{O}_\ZEP$ are changed in this step.  In this step we
  take into account the number $\ZSI$ and we proceed as follows: we fix the least number $n$ such that\footnote{$\lfloor\cdot \rfloor$ denotes the integer part.}
\[
\dfrac{1}{n}<\ZSI\quad {\rm i.e.}\quad n=n_\ZSI=\left\lfloor \dfrac{1}{\ZSI}\right\rfloor+1\,.
\]
 
The number $M_{\ZEP ,n_\ZSI}$ does depend  on $\mathcal{O}_\ZEP$ and on  $\ZSI$: $M_{ \ZEP,n_\ZSI} =M_{ \ZSI,\mathcal{O}_\ZEP}$ but, as we stated, $\ZEP$ and $\mathcal{O}_\ZEP$ are kept fixed.    So, $M_{  \ZSI,\mathcal{O}_\ZEP}$ changes only if $\ZSI$ is changed.

If $m>M_{  \ZSI,\mathcal{O}_\ZEP} $ and if $x\notin \mathcal{O}_\ZEP$ we have
\[
w_m(x)<\dfrac{1}{n_\ZSI}<\ZSI
\]
since $x\notin \mathcal{O}_\ZEP$ implies $x\notin A_{\ZEP,n_\ZSI}$.
\item[\bf Step~3:] Once $\ZEP>0$ and $\mathcal{O}_\ZEP$  have been fixed, the set $\mathcal{O}_\ZEP$ and the number $M_{  \ZSI,\mathcal{O}_\ZEP}$
just constructed satisfy the required properties and so the Egorov-Severini Theorem is proved.\zdia
 
\end{description}
  \begin{Remark}\ZLA{P1REMACh2AppeNoEGsevSenzaLimit}
 {\rm 
The proof of Theorem~\ref{Cha1bisTeoEgoSEVERANTE} uses the assumption that $R$ is bounded, hidden in the    use 
of  Lemma~\ref{LemmaiAppeCh3SUOmultirettANNID}.\zdia
 
 }
 \end{Remark}
 

 
  \part{\ZLA{PART2MoreVar}Functions of Several Variables}
 \chapter{\ZLA{Chap:2quasicontPIUvar}Functions of Several Variables: the Integral}

In this chapter we define the Lebesgue integral for functions of several variables. We assume familiarity with the  elementary topological 
notions of sets in $\zzr^d$ and  with 
  Chap.~\ref{Ch1:INTElebeUNAvaria}     since, once the Tietze extension theorem in $\zzr^d$ is known, the procedure which leads to the construction of Lebesgue integral for functions of several variables is essentially the same as that for functions of one variable.
  
  In this chapter we do not discuss   the exchange of limits and integrals. These theorems are in Chap.~\ref{ch3:LIMITintegral}.

\section{\ZLA{secT:ch0:Rectangles}Rectangles, Multirectangles and Null Sets}

We call {\sc rectangle}\index{rectangle} a set $R\subseteq \zzr^d $
which is the cartesian product of $d$   intervals  of the real line
\[
R=\prod_{k=1}^d  I_k\,.
\]
When the intervals are open , $I_k=(a_k,b_k)$, the set $R$ is the set of the points \[
x=(x_1\,,\  x_2\,,\dots\,, x_d)\quad \mbox{such that}\quad a_k<x_k<b_k \,.
\]
When the intervals are bounded and closed,   $I_k=[a_k,b_k]$,
the set $R$ is the set of the points \[
x=(x_1\,,\  x_2\,,\dots\,, x_d)\quad \mbox{such that}\quad a_k\leq x_k\leq b_k \,.
\]
When $a_{k_0}=b_{k_0}$ for at least one value of the index then the  rectangle collapses to a ``face'', a rectangle in lower dimensions. If $a_k=b_k$ for every $k$ the rectangle collapses to a point.
The important observation is that   \emph{a rectangle is a nonempty set.}

\ZREF{From a geometrical point of view, the rectangles we are going to use are not arbitrary rectangles: they have faces parallel to the coordinate planes and sides parallel to the coordinate exes.

When using the term ``rectangle'', this fact is always intended and not explicitly repeated.
}

We recall that the   length of an interval  is the distance of its end points. The length of an interval whose endpoints are   $a_k$ and $b_k$ with $b_k\geq a_k$ is
\[
L(I_k)=b_k-a_k\,.
\]
Open, closed or half closed intervals with the same endpoints have the same length.
If the interval is open then $b_k>a_k$ and the length is positive.

The {\sc volume\footnote{i.e. the length in dimension $1$ and the surface in dimension $2$.}}\index{volume (of a rectangle)} of the rectangle is
\[
L\left (\prod_{k=1}^d  I_k\right )=\prod_{k=1}^d L(I_k)\,.
\]

Facts to be noted:
\begin{itemize}
\item if a rectangle collapses to a face, i.e. if one of the intervals collapses to a point, then  $L(R)=0$. Only in this case we may have $L(R)=0$. In particular, \emph{ the volume of an open rectangle is positive.  }
\item the volume of a rectangle does not depend on its topological nature\footnote{in the next equality there is a minor abuse of language: if the rectangle is \emph{degenerate,} i.e. if $L(R)=0$, then its interior is empty and we did not define the volume of an empty set. The equality holds also for degenerate rectangles provided that we complete the definition of the volume by imposing $L(\emptyset)=0$.}:    $L({\rm int}\,R)=L({\rm cl}\,R)$. 
\end{itemize}

The key notion we shall use is that of   multirectangle. 
\begin{Definition}\ZLA{DefiCh2DefiMultintGENER}{\rm A {\sc multirectangle}\index{multirectangle}\index{rectangle!multirectangle}  is any \emph{finite or numerable sequence} of  rectangles.   

  The rectangles of the sequence are   the {\sc component rectangles}\index{component!rectangle} of the multirectangle. \zdia
}\end{Definition}

\begin{Remark}{\rm
As in Remark~\ref{RemaCh1ImproperTERMS} we note: the term ``finite or numerable sequence'' is not strictly correct. We use it to intend that the domain of the index is either $\zzn$ or a finite set, say $1\leq n\leq N$.\zdia
}\end{Remark}

\ZREF{
\begin{center}
\fbox{Important observation}
\end{center}
A multiinterval is a multirectangle in dimension $1$ but we note an important difference. According to the definition in Chap.~\ref{Ch1:INTElebeUNAvaria}, a multiinterval is composed of \emph{open} intervals.  Definition~\ref{DefiCh2DefiMultintGENER} is more general since nothing is assumed on the topology of the component rectangles. 
}

As in the case $d=1$,
we associate a number and a set to any multirectangle $\Delta=\{R_n\}$:
\[
L(\Delta)=\sum L(R_n)\,, \qquad \mathcal{I}_{\Delta }=\bigcup R_n\,.
\]
 This number $L$ does not change if the rectangles are taken in different order. 
\emph{It is important to note that the number $L(\Delta)$ does not change if the component rectangles $R_n$ are changed by adding or removing parts of their boundaries.} 

As noted in dimension $1$, $L(\Delta)$ cannot be interpreted as an ``area'' or a ``volume'' since the component rectangles may not be pairwise disjoint.

The introduction of the notion of multirectangle and of the number $L$ allows us to define   null sets. 
\begin{Definition}\ZLA{DefiCH2DefiNULLset}{\rm 
A set $N$ is a {\sc null set}\index{set!null!in $\zzr^d$}\index{null set!in $\zzr^d$} when for every $\ZEP>0$ there exists a multirectangle $\Delta_\ZEP $ such that
\[
  L(\Delta_\ZEP)<\ZEP \quad {\rm and}\quad N\subseteq  \mathcal{I}_{\Delta_\ZEP}\,.
\]

A property of the points of a set $A$ which is false only when $x$ belongs to a null subset of $A$ is said to hold {\sc almost everywhere}\index{almost everywhere} (sortly {\sc a.e.}\index{a.e.}) on $A$.
} 
\end{Definition}
 
This definition looks different from the corresponding definition in dimension~$1$ since here the rectangles need not be open. In fact:
\begin{Theorem}\ZLA{TheoP2C1NullSETtwoDefI}
A set $N$ is a   null set  when for every $\ZEP>0$ there exists a multirectangle   $\Delta_\ZEP $ \emph{composed of open rectangles} such that
\[
  L(\Delta_\ZEP)<\ZEP \quad {\rm and}\quad N\subseteq  \mathcal{I}_{\Delta_\ZEP}\,.
\]
\end{Theorem}
\zProof Let $N$ be a null set and let $\ZEP>0$.
We construct a multirectangle $\hat\Delta_\ZEP$  such that
\[
\left\{\begin{array}{l}
\mbox{$\hat \Delta_\ZEP$ is composed by \emph{open} rectangles}\\
L(\hat\Delta_\ZEP)<\ZEP\\
N\subseteq  \mathcal{I}_{\Delta_\ZEP} \,.
\end{array}\right.
\]

By assumption, $N$ is a null set. Hence,   there exists a multirectangle $\Delta_\ZEP$ such that
 
\[
L(\Delta_\ZEP)<\frac{\ZEP}{2}\,,\qquad
N\subseteq  \mathcal{I}_{\Delta_\ZEP} \,.
 \] 
 Let $\Delta_\ZEP=\{R_n\}$. By slightly enlarging the sides of $R_n$ we construct an \emph{open} rectangle $\hat R_n$ such that
 \[
{\rm cl}\, R_n\subseteq  \hat R_n\,,\qquad L(\hat R_n)< L(R_n)+\frac{\ZEP}{2\cdot 2^n}\,.
 \]
 The required multirectangle is $\hat\Delta_\ZEP=\{\hat R_n\}$.\zdia
 
 Arguments like this will be further examined in Chap.~\ref{ch3:LIMITintegral}.
 
We state two lemmas:
\begin{Lemma}\ZLA{LemmaCAP2UniosequeMULTIINT}
Let $\{\Delta_n\}$ be a sequence of  multirectangles. There exists a multirectangle $\Delta$ wich is composed precisely by the rectangles which compose  the multirectangles $\Delta_n$. Hence we have
\[
L(\Delta)=\sum _{n=1}^{+\ZIN} L(\Delta_n)\,,\qquad 
\mathcal{I}_\Delta=\bigcup \mathcal{I}_{\Delta_n}\,.
\]

\end{Lemma}

The proof is similar to that of Lemma~\ref{ch1:LemmaUnioMULTINT}.

\begin{Lemma}\ZLA{LemmaCAP1UnioNulliEnullaDIMd} Let $\{N_n\}$ be a sequence of null sets. Then $N=\cup N_n$ is a null set.
\end{Lemma}
The proof  is similar to that of Lemma~\ref{LemmaCAP1UnioNulliEnulla}.

 \begin{Example}
{\rm 
An argument  analogous to the one in Example~\ref{Exe:Ch0Qnullset} can be used to prove that
any numerable set is a null set. For example, 
the points with rational coordinates of the square $(0,1)\times(0,1)\subseteq \zzr^2$ are a numerable set.    

A different argument is as follows. 
Let the set be $\{x_n\}$ where $x_n\in\zzr^d$. The set with the sole element $x_n$ is a null set and a numerable union of null sets is a null set.\zdia
}
\end{Example}
 
The previous definitions and observations  parallels the corresponding ones we have seen when $d=1$. Now we show an important difference. Theorem~\ref{teo:Ch2StrutturaAperti}  asserts that when $d=1$ any open set is the union of a disjoint sequence of open intervals.  \emph{This equivalence holds only in dimension $ 1$. It does not extend to higher dimension}   as the following example shows:
\begin{Example}\ZLA{Chap2ESEnoteo:Ch2StrutturaAperti} {\rm
Let $d=2$ and let $T$ be any open triangle.

Let   $\{R_n\}$ be a sequence of \emph{open  pairwise disjoint  rectangles} which are contained in $T$. Their union cannot fill the triangle $T$: let $x_0$ be a    point of $T$ which belongs to the boundary of $R_{n_0}$ (so that it does not belong to $R_{n_0}$ which  is open). Then we have also $x_0\notin R_n$ for every $n\neq n_0$ because $R_n$ is open: if $x_0\in R_n$ then  $R_n$ must intersect $R_{n_0}$ while the rectangles are disjoint.

Instead, it is easy to represent $T$ as the union of a sequence of (non pairwise disjoint) open rectangles.
 It is also the union of a sequence of closed rectangles with the property that two of them intersects only on the boundary. We shall see (in Theorem~\ref{TeoCH3STRUttAPERTI}) that this is a property of every open set.
\zdia

}\end{Example}


 The reader is already familiar with the rigorous   notations in Chap.~\ref{Ch1:INTElebeUNAvaria} 
and the concise terminology used to speed up the presentation (see the table~\ref{tableCH1shortNOTATIONS}). A similar terminology, collected in the table~\ref{tableCH2shortNOTATIONS}, we can use in any dimension.

 \begin{table}[h]\caption{Succinct notations and terminology}\ZLA{tableCH2shortNOTATIONS}
\ZREF{
 We recall that when 
 $\Delta=\{R_n\}$ ($R_n$ is a rectangle) we define
\[
L(\Delta)=\sum _{n\geq 1 }L(R_n)\,,\qquad 
\mathcal{I}_\Delta=\bigcup R_n\,.
\]
Then:
\begin{itemize}
\item the multirectangle $\Delta$ constructed in Lemma~\ref{LemmaCAP2UniosequeMULTIINT} is denoted $\cup\Delta_n$;

\item we say that a multirectangle $\tilde\Delta $ is {\sc extracted}\index{multirectangle!extracted} from $\Delta$ when any component interval of  $\tilde\Delta$ is a component interval  of $\Delta$.
\item a multirectangle which has finitely many component rectangles is called a ``finite sequence'' (of rectangles) or a {\sc finite multirectangle.}\index{multirectangle!finite} 
\item We say that $\Delta$ and $\Delta'$ are disjoint when
$\mathcal{I}_\Delta \cap \mathcal{I}_{\Delta'}=\emptyset $, i.e. when
 no component rectangle of one of them intersect a rectangle of the other.
 \item we say that a {\sc multirectangle covers a set $A$}\index{multirectangle!which covers $A$} when $A\subseteq \mathcal{I}_\Delta$. 

\item we say that a multirectangle {\sc is in a set $A$}\index{multirectangle!in a set} when $\mathcal{I}_\Delta\subseteq A$. 
 
\end{itemize}

Similar expressions we may use are self explanatory.
} 
 \end{table}
%
%

 \section{\ZLA{sectTIETZpiuVAR}The Tietze Extension Theorem:  Several Variables} 
 
Now we state   Tietze extension theorem in any dimension. 
In the proof of the monotonicity of the integral and in the study of the exchange of the limits and integral, it is convenient to know the existence of extensions which have the additional property stated in Theorem~\ref{Ch2TeoExteMonotSEQULebe} below. Most of the proposed proofs of Tietze theorem provide extensions for which Theorem~\ref{Ch2TeoExteMonotSEQULebe} holds.

Once the extension theorem (Theorem~\ref{teo:EsteDAcompattiinRn} below) is known, the definition of the Lebesgue integral for functions of $d$ variables  is essentially the same as that of functions of $1$ variable. So, in this section we confine ourselves to state the extension theorem which was first asserted by Lebesgue in~\cite{Lebesgue1907PALERMO} with a  hint to a possible proof. 
Later on and independently of Lebesgue, several simpler proofs have been proposed
(see  for example~\cite{DieudonneTRATTATOvol1}. See~\cite{BrudnyiBOOK2012} also for an interesting historical overview)
 and Tonelli proposed one  in~\cite{Tonelli42AnnSNSpiuVARIAB}.
For completeness, we reproduce this proof   in the Appendix~\ref{AppeCH2}  but   the reader can make reference to any of the   proofs that he may know, provided that the monotonicity property stated in Theorem~\ref{Ch2TeoExteMonotSEQULebe} holds for that extension.

We use the (standard)   notation $C(A)$  to denote the linear space of the functions which are continuous on the set $A$.
 
 \begin{Theorem}[{\sc  Tietze  extension theorem}]\index{Theorem!Tietze!$d$ variables}\ZLA{teo:EsteDAcompattiinRn}
Let  $K\subseteq \zzr^d$ be a closed set. There exists an algorithm which associates to every function $f\in C(K)$ a function $f_e\in C(\zzr^d)$ in such a way that the following properties hold:
 
\begin{enumerate}
\item\ZLA{I1teo:EsteDAcompattiinRn} the function $f_e$ is a continuous extension of $f$ to $\zzr^d$.
\item\ZLA{I2teo:EsteDAcompattiinRn} the following inequalities hold:
 \begin{multline}\ZLA{ch2:ineqExtensio}
 \inf\{f_e(  x)\,,\   x\in \zzr^d\}=\min\{f_e(  x)\,,\   x\in \zzr^d\}=  \min\{f (  x)\,,\   x\in K\} \\
 \leq \max\{f (  x)\,,\   x\in K\}=
  \max\{f_e(  x)\,,\   x\in \zzr^d\}=\sup\{f_e(  x)\,,\   x\in \zzr^d\}\,.
%
 \end{multline}
 
 \end{enumerate}
 \end{Theorem}

We call {\sc Tietze extension}\index{Tietze extension}\index{extension!Tietze} any continuous extension of $f$ from $K$ to $\zzr^d$ which enjoys the property~(\ref{ch2:ineqExtensio}).

Propery~(\ref{ch2:ineqExtensio}) has the following consequence:
 
\begin{Theorem}\ZLA{Ch2TeoConveUNIereditDAeste}
Let $\{f_n\}$ be a sequence of continuous functions defined on the closed set $K$ and let us assume that $\{f_n\}$ is uniformly convergent  to zero  on $K$:
\[
\mbox{$\forall \ZEP>0\ \exists N_\ZEP$ such that if $n>N_\ZEP$ then $|f_n(x)|<\ZEP$ for all $x\in K$}\,.
\]
 Let $\ f_{n,e} $ be a Tietze extension of $f_n$ to $\zzr^d$. 

The sequence $\{f_{n,e}\}$ is uniformly convergent   to zero  on\ $\zzr^d$.
\end{Theorem}
\zProof
Let
\[
m_n= \min\{f_n (  x)\,,\   x\in K\}\,,\qquad M_n= \max\{f _n(  x)\,,\   x\in K\}\,.
\]
The assumption can be reformulated as follows:
\[
\limn m_n=0\,,\qquad \limn M_n=0\,.
\]
This is the property that $\{f_{n,e}\}$ is uniformly convergent to $0$ since we have also 
\[
m_n= \min\{f_{n,e} (  x)\,,\   x\in \zzr^d\}\,,\qquad M_n= \max\{f _{n,e}(  x)\,,\   x\in \zzr^d\}\,.\zdiaform
\]

The construction    proposed by Tonelli, as most of the usual proofs of Theorem~\ref{teo:EsteDAcompattiinRn}, gives extensions which have a further property:
\begin{Theorem}\ZLA{Ch2TeoExteMonotSEQULebe}
We have:  
\begin{enumerate}
\item\ZLA{I1Ch2TeoExteMonotSEQULebe}
Let $f\in C(K)$ and $g \in C(K)$ be such that $g(x)\geq f(x)$ on $K$   and let   $ f_{ e}$, $g_e$ be their extensions  constructed as in Appendix~\ref{AppeCH2}. Then   $g_e(x)\geq f_e(x)$ on $\zzr^d$. 
\end{enumerate}
It follows:
\begin{enumerate}
\setcounter{enumi}{1}
\item\ZLA{I2Ch2TeoExteMonotSEQULebe}
Let $\{f_n\}$ be a sequence in $C(K)$ and let $ f_{n,e}$ be the extension of $f_n$ constructed as in Appendix~\ref{AppeCH2}. We have:
\begin{itemize}
\item  if $\{f_n\}$ is increasing   on $K$, i.e. if  $f_{n+1 }(x)\geq  f_{n }(x)$   for every $x\in K$,  then $\{f_{n,e}\}$ is increasing on $\zzr^d$.

\item
if $\{f_n\}$ is decreasing   on $K$, i.e. if  $f_{n+1 }(x)\leq  f_{n }(x)$   for every $x\in K$,  then $\{f_{n,e}\}$ is decreasing on $\zzr^d$.
\end{itemize}
 \end{enumerate}
\end{Theorem}

 \ZREF{
 Instead, it is important to note that the property in Theorem~\ref{Theo:prprieSucceConvDIM1} does not extend to functions of $d$ variables: it is possible that $\{f_n(x)\}$ converges for every $x\in K$ while the sequence of the extensions does not converges in $\zzr^d$. This is seen in Example~\ref{Exe:AppeCH2NONexteConve}.
 }

 \section{\ZLA{Ch2DefiInteUnbFUNzSet}The Lebesgue Integral in $\zzr^{\lowercase{d}}$}
 
 Now we define the Lebesgue integral for functions of $d$ variables. The procedure is the same as that seen in Chap.~\ref{Ch1:INTElebeUNAvaria} when $d=1$ and it is sufficient that we sketch the ideas. The definition is in three steps: first we define
the quasicontinuous functions    and then, with two steps, we    define their Lebesgue integral.  
 
  \begin{description}
\item[\bf Step A: quasicontinuous functions]
\end{description}
  
 We define:
 \begin{Definition}\ZLA{DefiVH2FunzQUasContLIMIT}
 {\rm Let $f$ be a   function a.e. defined on the ball\footnote{we stress the fact that it may be $R=+\ZIN$, i.e. that the ball can be the entire space $\zzr^d$.}
 \[
 D_R=\{x\,|\ \|x\|\leq R\}\,,\qquad R\leq +\ZIN\,.
 \] 
 The function is {\sc   quasicontinuous}\index{function!quasicontinuous!$d$ variables} when the following property holds:
let $\{\ZEP_n\}$ be a sequence of \emph{positive numbers} such that $\lim _{n\to+\ZIN}\ZEP_n=0 $. 
 For every $n$ there exists a  multirectangle $\Delta_n$ such that
 \[
L(\Delta_n)<\ZEP_n\,,\quad \mbox{$\mathcal{I}_{\Delta_n}$ is open}\,,\quad \qquad f_{|_{D_R\setminus\mathcal{I}_{\Delta_n}}}\quad  \mbox{is continuous}\,.\zdiaform
 \]

A quasicontinuous function which is bounded is a {\sc bounded quasicontinuous function.}\index{function!bounded!quasicontinuous}
 }
  \end{Definition}
 
 The multirectangle 
 $ \Delta_n $ is an   
 {\sc associated multirectangle of order $\ZEP_n$.}\index{multirectangle!associated of order $\ZEP$}\index{associated!multirectangle}
 
 A Tietze extension\footnote{i.e. a continuous extension with the 
 properties~(\ref{ch2:ineqExtensio}).} of $ f_{|_{D_R\setminus\mathcal{I}_{\Delta_n}}}$  is an  {\sc associated function of order $\ZEP_n$.}\index{associated!continuous function of order $\ZEP_n$}\index{function!associated of order $\ZEP_n$} 

The pair $(\Delta_n ,f_{|_{D_R\setminus\mathcal{I}_{\Delta_n}}})$ is a pair of
{\sc associated multiintervals and continuous functions of order $\ZEP_n$.}\index{associated!multiintervals and continuous functions of order $\ZEP_n$}\index{function!associated sequence of continuous functions}
 
 We note that, given $\ZEP_n$,   associated functions of order $\ZEP_n$  are not uniquely determined neither by $\ZEP_n$ nor by the choice of $\Delta_n$. 
 
 Results analogous to those seen in Chap.~\ref{Ch1:INTElebeUNAvaria} hold. In particular we state:
 
 \begin{Theorem}\ZLA{TheoCH2LinearQC}
 The following properties hold:
 \begin{enumerate}
 \item\ZLA{I1:TheoCH2LinearQC} a function which is a.e. continuous on $D_R$ is quasicontinuous;
 \item\ZLA{I2:TheoCH2LinearQC}
the classes of the quasicontinuous functions and  that of the bounded quasicontinuous functions on a ball $D_R$   are linear spaces.
\item the product of quasicontinuous functions is quasicontinuous and the quotient is quasicontinuous if the denominator is a.e. nonzero.
\item\ZLA{I3TheoCH2LinearQC} if $g$ is defined and \emph{continuous} on a domain which contains ${\rm im}\, f$ and if $f$ is quasicontinuous then the function $x\mapsto g(f(x))$ is quasicontinuous\footnote{we repeat  that in general the composition of quasicontinuous functions is not quasicontinuous, see See Remark~\ref{rema:AppeCap4NONcompos} in Appendix~\ref{APPEch4NullNONborel}.}.
\end{enumerate}
\end{Theorem}

A function defined on $D_R$ is {\sc piecewise constant}\index{function!piecewise constant on $\zzr^d$} when there exists a finite number  of   rectangles $R_n$, say  $ 1\leq n\leq N $,  such that:
\begin{enumerate} 
\item ${\rm int}\, R_i\cap{\rm int}\, R_j=\emptyset$ if $i\neq j$;
\item
$D_R\subseteq \cup _{n=1}^N {\rm cl}\, R_n$; 
\item  the function $f_{|_{{\rm int}\,R_n}}$ is constant  for every $n$.
\end{enumerate}

The boundaries of rectangles and balls are null sets.  
So, the Property~\ref{I1:TheoCH2LinearQC} of Theorem~\ref{TheoCH2LinearQC}  implies:
\begin{Corollary}
We have:
\begin{enumerate} 
\item
The characteristic functions of rectangles or balls of $\zzr^d$ are bounded quasicontinuous functions. 
\item
Piecewise constant functions are bounded quasicontinuous functions. 
\end{enumerate}
\end{Corollary}

Now we extend the definition of quasicontinuity.
Let $A$ be a set which satisfies the following condition:
\begin{Assumption}\ZLA{ASSUch2QuascontINSa}{\rm 
The characteristic function of the set $A$ is quasicontinuous.\zdia}
\end{Assumption} 

 Theorem~\ref{Teo:Ch3FunCHarAPERTIquasCont} proved in Chap.~\ref{ch3:LIMITintegral} shows that any   open set satisfies this assumption.
 
We note that if $f$ is quasicontinuous on $\zzr^d$ and $A $ satisfies Assumption~\ref{ASSUch2QuascontINSa} then $f\charfun _{A}$ is quasicontinuous on $\zzr^d$. So we define:

\begin{Definition}{\rm
 Let $f(x)$ be quasicontinuous on $\zzr^d$ and let $A$   satisfies the Assumption~\ref{ASSUch2QuascontINSa}. Under these conditions, we say that the  function $f$  is {\sc quasicontinuous on the set $A$.}\index{function!quasicontinuous!on a set $A$}   

If $f$ is defined on $A$ then we say that it is quasicontinuous on $A$ when its extension with $0$ to $\zzr^d$ is.\zdia
}\end{Definition}

Clearly:
 \begin{Theorem}\ZLA{TheoCH2LinearQCbIs}
The properties  stated Theorem~\ref{TheoCH2LinearQC}  for functions which are quasicontinuous on $D_R$ holds also for functions which are quasicontinuous on a set $A$ which satisfies Assumption~\ref{ASSUch2QuascontINSa}
 \end{Theorem}
 
 \begin{description}
\item[\bf Step B: Lebesgue integral   under boundedness assumptions]
\end{description}
 
 We state the following result which    extends Lemma~\ref{Lemma:Ch2LemmaPreliRiemaa}, and which has a similar proof:
 \begin{Lemma}\ZLA{Lemma:Ch2LemmaPreliRiemaaPIUvar}
 Let $f$ be \emph{Riemann integrable} on a bounded ball $D $ (so that $f$ is bounded) and let $\Delta$ be a multirectangle such that
 \[
f(x)=0 \quad {\rm if}\quad x\notin\mathcal{I}_\Delta\,. 
 \]
 We have
 \[
\underbrace{\int _{D } |f(x) |\ZD x}_{\tiny \begin{tabular}{l}
    Riemann       integral  \end{tabular}}\leq \left (\sup _{D }|f|\right )L(\Delta)\,.
 \]
 \end{Lemma}

The proof of the next statement is similar to that of Lemma~\ref{Lemma:Ch2LemmaPreliRiemaa}:
 \begin{Lemma}\ZLA{Ch2:TeoreINvarianLiMIT}
 Let $f$ be a \emph{bounded} quasicontinuous function   on a \emph{bounded} ball. 
 Let $\{\ZEP_n\}$ be a sequence of positive numbers such that $\lim _{n\to+\ZIN}\ZEP_n=0$. Let $\{f_n\}$ be a sequence of associated functions (to $f$) of order $\ZEP_n$.
 The sequence of the  \emph{Riemann} integrals
 \[
\left \{ \int_{D_R} f_n(x)\ZD x\right \} 
 \]
 is convergent and the limit does not depend  either on $\{\ZEP_n\}$  or on the particular chosen sequence of associated functions.  
 \end{Lemma}
 
 It is then legitimate to define:
 
 \begin{Definition}{\rm 
Let $f$ be a \emph{bounded} quasicontinuous function defined on the ball $D_R=\{x\,:\ \|x\|\leq R\}$. \emph{We assume that the ball is bounded, i.e. that $R<+\ZIN$.}
Let $\{\ZEP_n\}$ be any sequence of positive numbers such that $\lim_{n\to+\ZIN} \ZEP_n=0$ and let $\{f_n\}$ be any sequence of associated functions, $f_n$ of order $\ZEP_n$. We define:\index{Lebesgue!integral!$d$ variables} 
 
 \[
\underbrace{ \int_{D_R} f (x)\ZD x}_{\tiny \begin{tabular}{l}
    Lebesgue  \\    integral  \end{tabular}}=\lim _{n\to+\ZIN}
\underbrace{ \int_{D_R} f_n (x)\ZD x}_{\tiny \begin{tabular}{l}
    Riemann  \\    integral  \end{tabular}}\,.\zdiaform 
 \]
 }
 \end{Definition}

 As a consequence we have:
 \begin{Theorem}
  \ZLA{TEO:Ch2InteNULLOfunzQOnulla} 
 The following properties hold for functions a.e. defined on a bounded ball $D_R$:
 \begin{enumerate}
 \item
 \ZLA{I1TEO:Ch2InteNULLOfunzQOnulla} 
 if $f$ is defined on $D_R$ and it is either continuous or piecewise constant  then it is  both Riemann and Lebesgue integrable and the two integrals have the same values.
 \item
 \ZLA{I2TEO:Ch2InteNULLOfunzQOnulla}
 if $f$ is a.e. continuous  and bounded  then it is bounded quasicontinuous and so  it is  Lebesgue integrable.

 \item
 \ZLA{I3TEO:Ch2InteNULLOfunzQOnulla}
 If $f=0$ a.e. on $D_R$ then
 \[
\underbrace{\int _{D_R} f(x)\ZD x }_{\tiny \begin{tabular}{l}
    Lebesgue       integral  \end{tabular}}=0\,.
 \]
 \item \ZLA{I3BISTEO:Ch2InteNULLOfunzQOnulla} two bounded quasicontinuous functions which are a.e. equal have the same Lebesgue integral.
 \end{enumerate}
 \end{Theorem}  
 \zProof Statements~\ref{I1TEO:Ch2InteNULLOfunzQOnulla} is immediate from the definition of the integral  and the proof of statement~\ref{I2TEO:Ch2InteNULLOfunzQOnulla} is similar to that given in Lemma~\ref{lemmaCH1daqocontTOquasicont} in the case $d=1$.

 We prove statement~\ref{I3TEO:Ch2InteNULLOfunzQOnulla}. 
Let
\[
|f(x)|<M\,.
\]
 
 We fix   a multirectangle $\Delta_n$ with $L(\Delta_n)<1/n$ and such that $f=0$ on $D_R\setminus\mathcal{I}_\Delta$ and a corresponding associated function $f_n=\left (f_{|_{D\setminus \mathcal{I}_{\Delta_n}}}\right )_e$.
 Then, 
 \[
 \underbrace{\int_{D_R} f(x)\ZD x }_{\tiny \begin{tabular}{l}
    Lebesgue     integral  \end{tabular}} 
     =\limn\underbrace{
 \int_{D_R} f_n(x) \ZD x }_{\tiny \begin{tabular}{l}
     Riemann       integral  \end{tabular}} 
 \] 
and, from   Lemma~\ref{Lemma:Ch2LemmaPreliRiemaaPIUvar},
\[
\underbrace{
 \int_{D_R} |f_n(x)|\ZD x }_{\tiny \begin{tabular}{l}
     Riemann       integral  \end{tabular}}\leq \dfrac{M}{n}\,. 
\]
 The result follows by computing the limit for $n\to+\ZIN$.
 
Statement~\ref{I3BISTEO:Ch2InteNULLOfunzQOnulla} follows since the difference of the two functions is a.e. zero.\zdia 
  
\emph{We conclude by stating that
  any function which is Riemann integrable on $D$ is Lebesgue integrable too, and the values of the integrals coincide. The proof is analogous to that seen, when $d=1$, in the Appendix~\ref{AppeTOch1IntegrRoemLebe}.  }

 \begin{description}
\item[\bf Step C: the general case]
\end{description}

As in Chap.~\ref{Ch1:INTElebeUNAvaria}, first we define the Lebesgue integral of a quasicontinuous function on $\zzr^d$. Then we extend the definition when $f$ is quasicontinuous on a set $A$ which satisfies the Assumption~\ref{ASSUch2QuascontINSa}.

The function $f$ can be unbounded.

We define:
 
  \begin{subequations}
\begin{equation}\ZLA{EQ:stepBdefiINRfPM}
f_+(x)=\max\{f(x),0\}\,,\qquad f_-(x)=\min\{f(x),0\} 
\end{equation}
and, when $K>0$, $N>0$ and $R>0$,
\begin{equation}\ZLA{EQ:stepBdefiINRfPMrn}
\begin{array}{lll}
\displaystyle
  f_{+,\, (R,N)}(x)= 
\min\{f_+(x)\,,\ N\}& \  &\Dom\, f_{+,\, (R,N)}=\{x\,:\ \|x\|\leq R\}\,,\\
\displaystyle
  f_{-,\, (R,-K)}(x)= 
\max\{f_-(x)\,,\ -K\}&  \   &\Dom\, f_{-,\, (R,-K)}=\{x\,:\ \|x\|\leq R\}\, . 
\end{array}
\end{equation}
\end{subequations}
 
 The functions $f_+$ and $f_-$ are quasicontinuous and the functions  
 $f_{+,\,(R,N)}$ and $f_{-,\, (R,-K)}$ are bounded quasicontinuous on a bounded ball. So, their Lebesgue integrals exist.

We define
\begin{align*}
&
\underbrace{\int_{\zzr^d} f_{+ }(x)\ZD x}_{\tiny \begin{tabular}{l}
    Lebesgue  \\    integral  \end{tabular}}=\lim _{\stackrel{R\to +\ZIN}{N\to +\ZIN}}\underbrace{\int _{\|x\|\leq R} 
f_{+,\,(R,N)}(x)\ZD x}_{\tiny \begin{tabular}{l}
    Lebesgue  \\    integral  \end{tabular}} \,,
    \\
  &  \underbrace{\int_{\zzr^d} f_{- }(x)\ZD x}_{\tiny \begin{tabular}{l}
    Lebesgue  \\    integral  \end{tabular}}=\lim _{\stackrel{R\to +\ZIN}{K\to +\ZIN}}\underbrace{\int _{\|x\|\leq R}
 f_{-,\,(R,-K)}(x)\ZD x}_{\tiny \begin{tabular}{l}
    Lebesgue  \\    integral  \end{tabular}}\,.
\end{align*}
The limits have to be computed respectively with $(R,N)\in \zzn\times\zzn$
and   $(R,K)\in \zzn\times\zzn$  (i.e. one independent   from the other) and they can be respectively~$+\ZIN$ or~$-\ZIN$.

\begin{Definition}[{\sc Lebesgue integral   on $\zzr^d$}]\index{Lebesgue!integral!$d$ variables}\index{integral!Lebesgue!$d$ variables} 
{\rm 
{~}
\begin{enumerate}
\item
The function $f(x)$ is   {\sc integrable}\index{function!integrable!$d$ variables}\index{function!integrable!$d$ variables}
 when at least one of the function $f_+$ or $f_-$ has finite integral.
In this case we define
\[
\underbrace{\int_{\zzr^d} f(x)\ZD x}_{\tiny \begin{tabular}{l}
    Lebesgue      integral  \end{tabular}}=\underbrace{\int_{\zzr^d} f_+(x)\ZD x+\int_{\zzr^d} f_-(x)\ZD x}_{\tiny \begin{tabular}{l}
  both  Lebesgue    integrals  \end{tabular}}\,.
\]
The Lebesgue integral of $f$ can be a number (when both the integrals of $f_+$ and of $f_-$ are numbers) or it can be $+\ZIN$ or it can be $-\ZIN$.
 
The function $f$ is {\sc summable}\index{function!summable!}\index{summable}\footnote{as already noted in the footnote~\ref{footnoteROYDENch1} of Chap.~\ref{Ch1:INTElebeUNAvaria}, several books uses the term ``integrable''   to intend that the integral is finite.}\index{summable (function)}\index{function!summable} when its Lebesgue integral is finite.

 \item
  Let $A$ be a set which satisfies the Assumption~\ref{ASSUch2QuascontINSa} and let $f$ be defined on $A$. We use the notation $f \charfun_A$ to denote the product of $f$ and $\charfun _A$ when $f$ is defined on $ \zzr^d$. Otherwise, with a slight abuse of notations, we put
  \[
f(x)\charfun_A(x)=
\left\{\begin{array}{lll}
f(x) &{\rm if}& x\in A\,,\\
0&{\rm if}& x\notin A\,.
\end{array}\right.  
  \]
    
We say that $f$ is a {\sc quasicontinuous function on $A$}\index{function!quasicontinuous!on a set $A$} when $f\charfun_A$ is quasicontinuous on $\zzr^d$. In this case we define
\[
\underbrace{\int_A f(x)\ZD x}_{\tiny \begin{tabular}{l}
 Lebesgue    integral   \end{tabular}}=\underbrace{\int _{\zzr^d} f(x)\charfun_A(x)\ZD x}_{\tiny \begin{tabular}{l}
 Lebesgue    integral   \end{tabular}}
\]
when $f(x)\charfun_A$ is integrable (in particular, when it is summable) on $\zzr^d$ and correspondingly we say that $f$ is {\sc integrable (summable) on $A$.}
\index{function!integrable!on a set $A$}\index{summable!on a set $A$}\index{integrable!on a set $A$}\index{function!summable!on a set $A$}\zdia

\end{enumerate}
 }
\end{Definition}
 
 \ZREF{From now on, the integral sign will always  denote the Lebesgue integral, unless explicitly stated that it is a Riemann integral.}

 
%
%
%
%
%

\emph{We conclude  this section by stating that, as in the case $d=1$, The Lebesgue integral does not extend the improper integral.}
 
\subsection{The Properties of the Integral}

The properties of the set of the quasicontinuous functions and of the integral are the same as we listed in the case $d=1$, and with similar proofs. We repeat the statements for completeness.

\begin{Theorem}\ZLA{TheoCH2PropELEMqcfunc}
The sets which appear in the statements below satisfy    Assumption~\ref{ASSUch2QuascontINSa}. Under this condition
 the following properties hold:
\begin{enumerate}
\item\ZLA{I1TheoCH2PropELEMqcfunc}
let $f$ be quasicontinuous on $  A$  and let   $A_1\subseteq A$. Then $f$ is quasicontinuous on $A_1$.

\item\ZLA{I2TheoCH2PropELEMqcfunc} let $A_1$ and $A_2$ be disjoint and let   $f$ be  a.e. defined and quasicontinuous both on $A_1$ and on $A_2$. Then it is quasicontinuous on $A_1\cup A_2$.

\item\ZLA{I3TheoCH2PropELEMqcfunc}
the sum and the product of two    quasicontinuous functions   is a   quasicontinuous function. This statement holds also for the quotient  provided that the   denominator  is a.e. different from zero.

In particular, \emph{the set of the quasicontinuous functions a.e. defined on  $A$ is a linear space.}

\item Let $f$ be quasicontinuous on $A$.
Let $v \in \zzr^d$. We put $A+v=\{x+v\,,\ x\in A\}$. The
set $A+v$ satisfies the  Assumption~\ref{ASSUch2QuascontINSa}
and the  function $x\mapsto f(x-v)$ is quasicontinuous on $A+v$.

 \item\ZLA{I4TheoCH2PropELEMqcfunc} let $f_n$ be quasicontinuous functions. For every $k$, the functions
 \begin{align*}
 &
\phi_k(x)=\max\{ f_1(x)\,,\ f_2(x)\,,\ \dots\,,\ f_k(x)\} \\
&\psi_k(x)=\min\{ f_1(x)\,,\ f_2(x)\,,\ \dots\,,\ f_k(x)\}
 \end{align*}
 are quasicontinuous.
\end{enumerate}
\end{Theorem}

Similarly, we can  state the key properties of the integral:
\begin{Theorem}\ZLA{CH3ProprIntePiUvaRsuA}
Let $f(x)$ and $g(x)$ be integrable on $A$ (a set which satisfies Assumption~\ref{ASSUch2QuascontINSa}). Then:
\begin{enumerate}
 \item\ZLA{I1CH3ProprIntePiUvaRsuA} the integral of a function which is a.e. zero is zero. So, if $f=g$ a.e. then they have the same Lebesgue integrals.
 
\item {\sc monotonicity of the integral:}\index{integral!monotonicity}\index{monotonicity of the integral}  if $f(x)\leq g(x)$ then\footnote{as in dimension~1, the proof uses the
monotonicity property of the extensions given in Theorem~\ref{Ch2TeoExteMonotSEQULebe}: we associate 
  to $f$ and $g$ sequences $f_{n } $ and $g_{n } $ such that $f_{n }(x)\geq g_{n }(x)$. Compare with Remark~\ref{RemaCH1MonotonInteTietze}.}
\[
\int_A f(x)\ZD x\leq \int_A g(x)\ZD x\,.
\]
\item {\sc  translation invariance:}\index{transaltion invariance!of the integral}
let $v \in \zzr^d$. With $A+v=\{x+v\,,\ x\in A\}$ we have:
\[
\int_A f(x)\ZD x=\int _{A+v} f(x-v)\ZD x\,.
\]
 
\item the absolute value:
\begin{enumerate}
\item
 If $f$ is integrable  then $|f|$ is integrable and
\[
\int_A |f(x)|\ZD x=\int_A f_+(x)\ZD x-\int_A f_-(x)\ZD x\,.
\]
So, the usual inequality of the absolute value holds:
\[
\left |
\int_A f(x)\ZD x
\right |\leq \int_A |f(x)|\ZD x\,.
\]
\item
  the \emph{quasicontinuous} function $f $ is summable if and only if $|f |$ is summable.
 
 \end{enumerate}
\item let $f$ and $g$ be summable. We have:
 \begin{enumerate}  
 \item  {\sc linearity of the integral:}\index{linearity of the integral}\index{integral!linearity} the following equality holds for any real numbers $\zaa$ and $\beta$:
\[
\int_A\left (\zaa f(x)+\beta g(x)\right )\ZD x=
\zaa\int_A f(x)\ZD x+\beta\int_A g(x)\ZD x\,.
\]
 \item if   $g$ is bounded then the product $fg$ is summable. 
 \item   if  $1/g$ is bounded then the quotient $f/g$ is summable.
 
 \end{enumerate}
\end{enumerate}
\end{Theorem}

Now we consider   the additivity of the integral. 

\begin{Theorem}\ZLA{TeoCH2teoAbsANYdim}
Let $f$ be defined  on $A=A_1 \cup A_2$ and let $A_1$ and $A_2$ satisfy the Assumption~\ref{ASSUch2QuascontINSa}. Then:
\begin{enumerate}
   
\item \ZLA{I1TeoCH2teoAbsANYdim}
the set  $A_1\cup A_2$ and, if nonempty, the sets $A_1\cap A_2$,  $A_1\setminus A_2$, $A_2\setminus A_1$
satisfy the Assumption~\ref{ASSUch2QuascontINSa};
\item\ZLA{I2TeoCH2teoAbsANYdim} let  $f $ be summable on $A_1$ and on $A_2$. Then:
\begin{enumerate}

\item\ZLA{I2aTeoCH2teoAbsANYdim}   
 the function $f$ is summable on  $A_1\cup A_2$;
 \item\ZLA{I2bTeoCH2teoAbsANYdim}   we have
 \[
  \int _{A_1\cup A_2} f(x)\ZD x\leq \int _{A_1} f(x)\ZD x+\int _{A_2} f(x)\ZD x\,;
  \]
\item\ZLA{I2cTeoCH2teoAbsANYdim}{\sc additivity of the integral:}\index{integral!additivity}\index{additivity of the integral} 
\[
A_1\cap A_2=\emptyset\ \implies\ \int _{A_1\cup A_2} f(x)\ZD x=\int _{A_1} f(x)\ZD x+\int _{A_2} f(x)\ZD x\,.
\]

\end{enumerate}
\end{enumerate}
\end{Theorem}
\zProof
The result follows from the linearity of the integral and the equality
\[
f(x)=   f(x)\carfunz_{A_1}(x)+  f(x)\carfunz_{A_2}(x) \,.\zdiaform
\]
  
  Finally we extend Lemma~\ref{Lemma:Ch2LemmaPreliRiemaaPIUvar} from the Riemann to the Lebesgue integrals. We recall that the {\sc support}\index{support (of a function)} of a function $f$ is
\[
{\rm supp}\, f={\rm cl}\,\{x\,:\ f(x)\neq 0\}\,.
\]
So, the support is always closed and, if we consider functions defined in a bounded ball, it is compact.

Heine-Borel Theorem holds in every dimension. In particular, let $\Delta$ be a  multirectangle composed by open rectangles which \emph{covers} a compact set $K$; i.e. we assume   $K\subseteq \mathcal{I}_\Delta$. Then
$K$ is covered by finitely many of the rectangles which compose $\Delta$ (we already stated this observation when $d=1$ in Lemma~\ref{I2Exe:Ch0Qnullset}).

We use this observation and we extend Lemma~\ref{Lemma:Ch2LemmaPreliRiemaaPIUvar} as follows:
\begin{Lemma}\ZLA{Lemma:Ch3LemmaPreliLEBEaPIUvar}
Let $f$ be a bounded quasicontinuous function defined on a bounded ball $D$, $|f(x)|\leq M$ for every $x\in D$. Let $\Delta=\{R_n\}$ be a sequence of open  rectangles which covers ${\rm supp}\, f$:
\[
{\rm cl}\,\{x\,:\ f(x)\neq 0\}\subseteq \mathcal{I}_\Delta\,.
\]
Then we have
\begin{equation}\ZLA{eq:Lemma:Ch3LemmaPreliLEBEaPIUvar}
\left |\int_D f(x)\ZD x\right |\leq ML(\Delta)\,.
\end{equation}
\end{Lemma}
\zProof We use Heine-Borel Theorem and we reorder the rectangles of $\Delta $ so to have
\[
{\rm cl}\,\{x\,:\ f(x)\neq 0\}\subseteq  \bigcup _{i=1}^N R_i 
\] 
 (note that the rectangles $R_i$ need not be disjoint).
 
Statement~\ref{I2bTeoCH2teoAbsANYdim} of Theorem~\ref{TeoCH2teoAbsANYdim} gives
\[
\int_D f(x)\ZD x\leq \int _{\cup _{i=1}^N R_i} f(x)\ZD x \quad \mbox{so that}\quad
\int_D |f(x)|\ZD x\leq \sum_{i=1}^n \int_{R_i} |f(x)|\ZD x 
 \] 
Inequality~(\ref{eq:Lemma:Ch3LemmaPreliLEBEaPIUvar}) holds since for every $i$ we have
\[
  \int_{R_i} |f(x)|\ZD x\leq \int_{R_i} M\ZD x=ML(R_i)\zdiaform
\]
 
 Note that this result is quite weak. In particular it cannot be used to prove statement~\ref{I1CH3ProprIntePiUvaRsuA} of Theorem~\ref{CH3ProprIntePiUvaRsuA}.
%

 
\section{\ZLA{AppeCH2}\textbf{Appendix:}  Tonelli and the Tietze Extension}

 Here we reproduce the proof   of the Tietze extension theorem (Theorem~\ref{teo:EsteDAcompattiinRn})   presented  by Tonelli   in~\cite{Tonelli42AnnSNSpiuVARIAB}.

 \subsection{Few Notations}

For every $r\geq 0$ and every $  x\in\zzr^d$ let
 \[
D(  x,r)=\{  y\ztc \|   y-  x\|\leq r\}\qquad \mbox{($D(  x,r)$ is a \emph{closed} ball)}\,.
 \]
 Let $K$ be a compact set. We put\footnote{note in this definition: ``$\min$'' since $K$ is closed. For general sets the distance is an infimum.} 
\[
\rho(  x)={\rm dist}(  x,K)=\min\{\|x-k\|\,,\ k\in K\}\,.
\]

Note that  
\[
K\cap D(x,0)=\left\{\begin{array}{lll}
\emptyset &{\rm if}& x\notin K\\
x &{\rm if}& x\in K\,.
\end{array}\right.
\]

The set $K\cap D(  x,r) $ is empty when $r<\rho(x)$. If $r\geq \rho(x)$   it is a compact set over which $f$ is continuous. So, when $r\geq \rho(x)$,  the function $f$ reaches its maximum (and minimum) on $K\cap D(  x,r) $.

 We introduce the   function
 \[
M(  x,r)=\left\{\begin{array}{lll}
\max_{ D(  x,r)\cap K}f&{\rm if} &D(  x,r)\cap K\neq \emptyset\\
0&{\rm if}& D(  x,r)\cap K = \emptyset\,.
\end{array}\right. 
 \]
The function $r\mapsto M(x,r)$ is monotone non decreasing for every $x$ and both the functions $x\mapsto M(x,r)$ (with $r$ fixed) and $r\mapsto M(x,r)$ (with $x$ fixed) are discontinuous. In fact, let $d=1$, $K=[1,2]$ and $f(x)\equiv 1$ on $K$. Then\footnote{recall that $D(x,r)$ is a \emph{closed} ball.},
\[
M(0,r)=\left\{\begin{array}{lll}
0 &{\rm if}& r<1\\
1 &{\rm if}& r\geq 1\,,
\end{array}
\right.
\qquad  M(x,1/2)=
\left\{\begin{array}{lll}
0&{\rm if}& x<1/2\\
1&{\rm if}& x\geq 1/2\,.
\end{array}
\right.
\]

\subsection{The Proof of Theorem~\ref{teo:EsteDAcompattiinRn}}

For every $n$ and every $x\in\zzr^d$ we define
\begin{equation}\ZLA{eq:Appe2DefiTONELLIdiFn}
F_n(x)=\frac{1}{2^n} \sum _{k=0}^{2^n-1} M\left (  x, \left (1+\dfrac{k}{2^n}\right )\rho(  x)     \right ) \,.
\end{equation}

The number $F_n(x)$ is the arithmetic mean of the numbers
\[
\max\left \{\,
f(x)\,,\ x\in K\cap D\left (x,     (1+k/2^n  )\rho(  x) \right )   \,
\right \}\quad 0\leq k<2^{n}-1\,.
\]

\emph{We prove that
\begin{equation}\ZLA{eqESIlimTietze}
f_e(x)=\limn F_n(x)
\end{equation}
exists for every $x\in\zzr^d$, it is an extension of $f $ which satisfies~(\ref{ch2:ineqExtensio}) and it is continuous.}  
The proof is  in three steps. 
After that, it is a simple observation to note that the monotonicity 
 Theorem~\ref{Ch2TeoExteMonotSEQULebe} holds for the extension $f_e$ in~(\ref{eqESIlimTietze}).  See the statement~\ref{I3Rema:ch2RimuovBoundedDAtietze} of Remark~\ref{Rema:ch2RimuovBoundedDAtietze}.

Now we prove the theorem. 

\ZREF{
\begin{center}
{\bf A warning}
\end{center}

\smallskip
Note that the extension theorem has been used only in the special case that $K$ is compact and only this case  will be used  in the following. So, the reader can confine himself to consider the proof in the case that $K$ is compact. 

The proof is slightly more transparent when $K$ is compact and,  in order to help the reader, we prove the theorem   in this case.   
The points where \emph{boundedness of $K$ is used} are clearly stated by using a \fbox{box} and the way to remove  boundedness of $K$ is clearly indicated in \fbox{boxes.} But when first reading the proof it may be convenient to ignore these boxes.
 } 

\begin{description}
\item[\bf Step~1: the function $f_e$ is defined on $\zzr^d$.] We prove that 
for every fixed $x$ the sequence of real numbers $\{F_n(x)\}$ is bounded and nondecreasing. This implies that 
  the
limit~(\ref{eqESIlimTietze}) exists and that it is finite.

Boundedness is clear since for every $k$ and every $n$ we have
\[
\min_K f\leq  M\left (  x, \left (1+k/2^n\right )\rho(  x)     \right )\leq \max _K f
\]
and so we have also
\[
\min_K f\leq    F_n(x)\leq \max _K f\,.
\]
In order to prove monotonicity we prove $F_{n+1}(x)\geq F_n(x)$. Here we use monotonicity of $r\mapsto M(x,r)$ and the fact that the sum which defines $F_{n+1}(x)$ has \emph{twice as many addenda} as  that of  $F_n(x)$.
\emph{This is the reason for choosing the sum of \, $2^n$ terms.}

We write   the expression of $F_{n+1}(x)$ and we associate any term with even index   with its subsequent term with odd index:
 \begin{multline*}
 F_{n+1}(x)=
 \frac{1}{2^{n+1}} \sum _{k=0}^{2^{n+1}-1} M\left (  x,\left  (1+\frac{k}{2^{n+1}}\right )\rho(  x)     \right ) \\
  =\dfrac{1}{2^{n+1}}\left \{
   \left [
 M(x,\rho(x))+M\left (  x,\left  (1+\frac{1}{2^{n+1}}\right )\rho(  x)     \right )
  \right ]\right.\\
  +\left [
  M\left (  x,\left  (1+\frac{2}{2^{n+1}}\right )\rho(  x)     \right )+M\left (  x,\left  (1+\frac{3}{2^{n+1}}\right )\rho(  x)     \right )\right ]\\[3mm]
  +\cdots\\[3mm]
\left. +
 \left [M\left (  x,\left  (1+\frac{2\cdot2^n-2}{2^{n+1}}\right )\rho(  x)     \right )+
 M\left (  x,\left  (1+\frac{2\cdot 2^n-1}{2^{n+1}}\right )\rho(  x)     \right )
 \right ]
  \right \}\\
 =\frac{1}{2^n} \sum _{k=0}^{2^n-1}\frac{1}{2}\left [M\left (  x,\left (1+\dfrac{2k}{2\cdot2^n}\right )\rho(x)\right )+M\left (  x,\left (1+\dfrac{2k+1}{2\cdot2^n}\right )\rho(x)\right )\right ]\,.
 \end{multline*}
 
We compare the addenda of $F_{n+1 } (x)$ and those of $F_n(x)$. We see that
 \begin{multline*}
\frac{1}{2}\left [M\left (  x,\left (1+\dfrac{2k}{2\cdot2^n}\right )\rho(x)\right )+M\left (  x,\left (1+\dfrac{2k+1}{2\cdot2^n}\right )\rho(x)\right )\right ]\\
=
 \frac{1}{2}\left [M\left (  x,\left (1+\dfrac{ k}{ 2^n}\right )\rho(x)\right )
 +
 M\left (  x,\left (1+ \dfrac{ k }{  2^n}+\dfrac{ 1}{2\cdot2^n}\right )\rho(  x)\right )  \right]
\\[3mm]
 \geq  M\left (  x, \left (1+k/2^n\right )\rho(  x)     \right )\quad \mbox{\small (since $r\mapsto M(x,r)$ is nondecreasing)}\,.
 \end{multline*}
 It follows that
 \[
 \mbox{$
F_{n+1}(x)\geq F_n(x)$ 
 $ \forall n$  and so 
$f_e(x)=\lim F_n(x)\in\zzr$   for every $x$}\,.
 \]
 \item[\bf Step 2: the function $f_e$ extends $f$.] If $x\in K$ then $\rho(x)=0$ and $M(x,0)=f(x)$. Hence, when $x\in K$, $F_n(x)=f(x)$ for every $n$ so that we have also $f_e(x)=f(x)$.
 \item[\bf Step 3: the function $f_e$ is continuous on $\zzr^d$.] \fbox{In this step we use  $K$ compact}\\ (and we indicate how boundedness can be removed).
 
 If $x\in {\rm int}\, K$ then $f_e =f $ is continuous at $x$. We must prove continuity at the boundary points of $K$ and at the exterior points of $K$.

 \begin{description}
\item[Substep 3A: continuity at $x_0\in\partial K$.]  

We must prove
\[
\lim _{x\to x_0} f_e(x)=f_e(x_0)=f(x_0)
\]
(and we can confine ourselves to consider the limit from $\zzr^d\setminus K$ since $f$ is continuous on $K$ by assumption).

\fbox{A continuous function on a compact set is uniformly continuous.} Hence, for every $\ZEP>0$ there exists $\ZDE>0$ (which does not depend on $x_0$) such that
 \begin{equation}\ZLA{eq:CH2AppePriPerUniCont}
 y\in K\,,\ \|y-x_0\|<\ZDE\ \implies\ 
\left | f(y)-f(  x_0)\right |<\ZEP\,.
 \end{equation}
We note
\[
(1+k/2^n)\rho(  x)<2\rho(  x)\qquad 0\leq k <2^n\,.
\]
When $  x\in D(  x_0,\ZDE/4)$ then $\rho(x)\leq \ZDE/4$ and 
 \[
M\left (  x, (1+k/2^n)\rho(  x)\right )=f(y)\ \mbox{where}\ y\in D(x_0,\ZDE/4)\subseteq D(x_0,\ZDE)
 \]
  so that
 \[
\left | M\left (  x, (1+k/2^n)\rho(  x)\right )-f(  x_0)\right |<\ZEP\,.
 \]
Hence, for every  $n$,
 \begin{multline*}
\left | F_n(  x)-f(  x_0)\right |
=\left |
 \frac{1}{2^n} \sum _{k=0}^{2^n-1}\left [ M\left (  x, \left (1+\dfrac{k}{2^n}\right )\rho(  x)     \right )-f(x_0)\right ]
\right |\\
\leq \frac{1}{2^n} \sum _{k=0}^{2^n-1} \ZEP=\ZEP\,.
 \end{multline*}
The inequality is preserved by the limit, so it holds for $f_e$ as wanted.

 \ZREF{ 
\begin{center} 
{\bf How to remove   the assumption that $K$ is bounded} 
 \end{center} 
 
\smallskip

 We used boundedness of $K$ since we used uniform continuity but note that the point $x_0$ has been fixed. Once $x_0$ has been fixed, the value of $y\in K$ which are used in~(\ref{eq:CH2AppePriPerUniCont}) are confined to a ball of center $x_0$, for example   $y\in K\cap D(x_0,100\rho(x_0))$. The previous estimates holds also if $K$ is unbounded by using   uniform  continuity of $f$ on  $ K\cap D(x_0,100\rho(x_0))$ and the corresponding value of $\delta$.
 }

\item[Substep 3B: continuity at $x_0\in\zzr^d\setminus  K$.]
We fix a point $x_0\in \zzr^d\setminus  K$ so that $ \rho(x_0)>0$. 

We must prove that for every $\ZEP>0$ there exists $\ZDE>0$ such that
\begin{equation}\ZLA{eq:ch2TietzTEOdiseDAsottoEsopra}
\|x-x_0\|<\ZDE\ \implies  f_e(x_0)-\ZEP  <f_e(x)<  f_e(x_0)+\ZEP \,.
\end{equation}

First we note the following facts which holds for every   $\ZDE>0$:
\begin{itemize}
\item[\bf Fact 1:] we have 
\begin{equation}\ZLA{eqTeoTITZstep3b2}
\|x-x_0\|<\ZDE \ \implies\ \rho(x)\leq \rho(x_0)+\ZDE\,.
\end{equation} 
In fact let $k\in K $ be one of the points for which
\[
\rho(x_0)=\| x_0-k\|\,.
\]
Then we have
\[
\rho(x)\leq \|x-k\|\leq\|x-x_0\|+\|x_0-k\|=\|x-x_0\|+\rho(x_0)\leq \ZDE+\rho(x_0)\,.
\]
\item[\bf Fact 2:] for every $r>0$ we have
\begin{equation}\ZLA{eqTeoTITZstep3b1}
\begin{array}{l}
\displaystyle
\|x-x_0\|<\ZDE\ \implies
D(x,r)\subseteq D(x_0,r+\ZDE) \\[2mm]
\displaystyle \mbox{so that} \quad
  M(x,r)\leq M(x_0,r+\ZDE) \,.
\end{array}
\end{equation}
In fact:  
\[
y\in D(x,r)\ \implies\ 
\|y-x_0\|\leq \|y-x\|+\|x-x_0\|\leq r+\ZDE\,.
\]

\end{itemize}

 We use {\bf Fact~1} and {\bf Fact~2} to derive the following consequence:
  when  
\[
\|x-x_0\|<\ZDE 
\]
 we have:
\begin{multline}\ZLA{eq:LaCh2TietzeStmaPERm}
 M\left (x,(1+k/2^n)\rho(x)\right )\underbrace{\leq}_{\stackrel{\mbox{from}}{~(\ref{eqTeoTITZstep3b1})}} 
M\left (x_0,(1+k/2^n)\rho(x )+\ZDE\right ) \\
\underbrace{\leq}_{\stackrel{\mbox{from}}{(\ref{eqTeoTITZstep3b2})}} M\left (x_0,\left [(1+k/2^n)(\rho(x_0) +\ZDE)\right ]+\ZDE\right ) \\
\underbrace{\leq}_{\stackrel{\mbox{use}}{   (1+k/2^n)<2 }} M\left (x_0, (1+k/2^n) \rho(x_0)  +3\ZDE\right )\,.
\end{multline}

If $h$ is any number such that 
\begin{equation}\ZLA{eqCH2ProvaTIeSCLATh1}
3\ZDE< \dfrac{h\rho(x_0)}{2^n}\quad {\rm i.e.}\quad
h>2^n \dfrac{3\ZDE}{\rho(x_0)}
\end{equation}
then we have
\[
\left ( 1+\dfrac{k}{2^n}\right )\rho(x_0)+3\ZDE\leq 
\left ( 1+\dfrac{k+h }{2^n}\right )\rho(x_0) \,.
\]
We choose the \emph{smallest integer $h$ such that the   inequalities in~(\ref{eqCH2ProvaTIeSCLATh1}) hold:  }  
 
\begin{equation}
\ZLA{eqTeoTITZstep3b3hN}
h=h_n=\left\lfloor
2^n\dfrac{3 \ZDE}{\rho(x_0)}\right\rfloor+1 
\end{equation}
where $\lfloor
\zaa\rfloor$ denotes the integer part of the number $\zaa$.

Note that $h_n$   does not depend on $k<2^n$ so that for every $k<2^n$ we have
\begin{equation}
\ZLA{eqTeoTITZstep3b3}
\left ( 1+\dfrac{k}{2^n}\right )\rho(x_0)+3\ZDE\leq
\left ( 1+\dfrac{k+h_n }{2^n}\right )\rho(x_0)\,.
\end{equation}

The inequalities~(\ref{eq:LaCh2TietzeStmaPERm}) and~(\ref{eqTeoTITZstep3b3}) give, when $\|x-x_0\|\leq \ZDE $,
\begin{equation}\ZLA{eqTeoTITZstep3b3BBIS}
 M\left (x,\left (1+\dfrac{k}{2^n}\right )\rho(x)\right )\leq  M\left (x_0,\left ( 1+\dfrac{k+h_n }{2^n}\right )\rho(x_0)\right )\quad 0\leq k<2^n\,.
\end{equation}

\emph{After these preliminaries we prove separately the inequalities above and below in~(\ref{eq:ch2TietzTEOdiseDAsottoEsopra}). }

First we prove
the inequality above: we prove
 the existence of $\ZDE$ such that
\begin{equation}\ZLA{eq:ch2TietzTEOdiseDAsopra}
 \|x-x_0\|<\ZDE   \ \implies f_e(x)<f_e(x_0)+\ZDE  \,.
\end{equation}

We use~(\ref{eqTeoTITZstep3b3BBIS}) and \fbox{$M(x_0,r)\leq \max_K f$.} We have:
 
\begin{multline}\ZLA{eq:CH2TietFINEprimaDisu}
F_n(x)= \frac{1}{2^n} \sum _{k=0}^{2^n-1} M\left (x,\left( 1+\frac{k}{2^n}\right )\rho(x)\right )\\
\leq
\frac{1}{2^n} \sum _{k=0}^{2^n-1} M\left (x_0,\left( 1+\frac{k+h_n}{2^n}\right )\rho(x_0)\right ) 
\\
=
\frac{1}{2^n} \sum _{\nu =h_n}^{2^n-1+h_n} M\left (x_0,\left( 1+\frac{\nu}{2^n}\right )\rho(x_0)\right )\\
%
= \frac{1}{2^n} \sum _{\nu =0}^{2^n-1 } M\left (x_0,\left( 1+\frac{\nu}{2^n}\right )\rho(x_0)\right )\\
+\frac{1}{2^n}\left [
  \sum _{\nu =2^n}^{2^n-1+h_n } \underbrace{M\left (x_0,\left( 1+\frac{\nu}{2^n}
  \right )\rho(x_0)\right )}_{\tiny \fbox{$\leq \hat M= \max_K |f|$}}
-
  \sum _{\nu =0}^{ h_n -1} \underbrace{M\left (x_0,\left( 1+\frac{\nu}{2^n}\right )
  \rho(x_0)\right )}_{\tiny\fbox{$\leq \hat M=  \max_K | f|$}}\right ]
\\
\leq 
F_n(x_0)
+2\dfrac{1}{2^n} h_n\hat M \quad \fbox{where $\hat M=  \max_K |f| $}\,.
\end{multline}
 \end{description}
Equality~(\ref{eqTeoTITZstep3b3hN}) shows that
\[
\dfrac{h_n}{2^n}\leq \dfrac{1}{2^n}+\dfrac{3\ZDE}{\rho(x_0)}
\]
so that, when $\|x-x_0\|<\ZDE $,  we have
\begin{equation}\ZLA{semiFineTietzS}
F_n(x)\leq F_n(x_0)+\dfrac{2}{2^n}\hat M+\dfrac{6\ZDE}{\rho(x_0)}\hat M\quad\mbox{and so}\quad 
f_e (x)\leq f_e (x_0)+ \dfrac{6\ZDE}{\rho(x_0)}\hat M \,.
\end{equation}
The required inequality~(\ref{eq:ch2TietzTEOdiseDAsopra}) holds provided that we choose
\[
  \ZDE<\dfrac{\rho(x_0)}{6\hat M}\ZEP\quad \fbox{where, we recall, $\hat M=\max_K |f|$ }\,.
\]
 \ZREF{ 
\begin{center} 
{\bf How to remove   the assumption that $K$ is bounded} 
 \end{center}
 
\smallskip

 In this step of the proof boundedness of $K$ was used when we defined $\hat M=\max_K |f|$. But, it is still true that $x_0$ is fixed and that the values of $f$ which are used in this computation are the values $f(y)$ when $y\in   K\cap D(x_0,100\rho(x_0))$. So, boundedness of $f$ is easily removed by redefining 
 $\hat M=\max_{ K\cap D(x_0,100\rho(x_0)) } |f|$.
 }

In a similar way we prove the inequality from below in~(\ref{eq:ch2TietzTEOdiseDAsottoEsopra}), i.e. we prove
\begin{equation}\ZLA{eq:ch2TietzTEOdiseDAsottoFiN}
\forall\ZEP>0\ \exists\ZDE>0\,:\ \|x-x_0\|<\ZDE\ \implies   f_e(x_0)-\ZEP  < f_e(x).
\end{equation}
 
We sketch the steps in order to see a (very) minor difference.

  \emph{Note that it is not restrictive to assume from the outset $\ZEP<1$ and   $\ZDE\in (0,\rho(x_0)/4)$.}  The value of $\ZDE$ will be further reduced later on.

 First we use the inequalities~(\ref{eqTeoTITZstep3b2}) and~(\ref{eqTeoTITZstep3b1}) with the roles of $x$ and $x_0$ exchanged, i.e.
\begin{equation}\ZLA{che2EqDellaIIparteTiet}
\|x-x_0\|<\ZDE\ \implies\left\{
\begin{array}{l}
\rho(x_0)<\rho(x)+\ZDE\\
r>0\ \implies\ \left\{\begin{array}{l}
D(x_0,r)\subseteq D(x,r+\ZDE)\\
M(x_0,r)\leq M(x,r+\ZDE)\,.
\end{array}\right.
\end{array}
\right.
\end{equation}
So we have
\begin{multline}\ZLA{eqCH2TietMULTI1IIparte}
M\left (x_0,\left (1+\dfrac{k}{2^n}\right )\rho(x_0)\right )\leq 
M\left (x ,\left (1+\dfrac{k}{2^n}\right )\rho(x )+3\ZDE\right )\\
\leq M\left (x ,\left (1+\dfrac{k+\tilde h}{2^n}\right )\rho(x )\right )
\end{multline}
provided that
\begin{equation}\ZLA{eqCH2ProvaTIeSCLATh2}
3\ZDE<\dfrac{\tilde h\rho(x)}{2^n}\,.
\end{equation}
\emph{Here is the point of difference: the inequality~(\ref{eqCH2ProvaTIeSCLATh1}) does not depend on $x$ while the right side of~(\ref{eqCH2ProvaTIeSCLATh2}) does depend on $x$.} But, this difficulty is easily overcame since  $\rho(x)>\rho(x_0)-\ZDE$ from~(\ref{che2EqDellaIIparteTiet}) and we did impose $ \ZDE\in (0,\rho(x_0)/4)$. Hence, inequality~(\ref{eqCH2ProvaTIeSCLATh2}) holds if we impose
\[
3\ZDE< \frac{\tilde h(\rho(x_0)-\ZDE)}{2^n}\,.
\]
This condition is satisfied if we choose $\tilde h$ such that
\[
  \ZDE<\dfrac{\tilde h\rho(x_0)}{4\cdot 2^n}\quad{\rm i.e.}\quad \tilde h>2^n\frac{4\ZDE}{\rho(x_0)}\,.
\]

We write the inequality~(\ref{eqCH2TietMULTI1IIparte}) with $h=h_n$:
\[
h_n=\lfloor  2^n\frac{4\ZDE}{\rho(x_0)}\rfloor +1\,.
\]
The same computations as in~(\ref{eq:CH2TietFINEprimaDisu}) with the roles of $x$ and $x_0$ exchanged  give
\[
F_n(x_0)\leq F_n(x)+\dfrac{2\hat M}{2^n}+\dfrac{8\hat M\ZDE}{\rho(x_0) }\quad \mbox{hence}\quad
f_e(x_0)-\dfrac{8\hat M\ZDE}{\rho(x_0) }\leq f_e(x) 
\]
\fbox{where $\hat M=\max_K |f|$. }
It follows that the required inequality~(\ref{eq:ch2TietzTEOdiseDAsottoFiN}) holds if we further reduce the value of $\ZDE$ and we impose  
\[
 \ZDE<\min\left \{
\dfrac{\rho(x_0)}{4}\,,\   \dfrac{\rho(x_0)}{8\hat M}\ZEP\right  \}\,.
\]

 \ZREF{ 
\begin{center} 
{\bf How to remove   the assumption that $K$ is bounded} 
 \end{center}
 
\smallskip

Also in this step of the proof boundedness of $K$ is used  since $\hat M=\max_K |f|$. But, it is still true that $x_0$ is fixed and that the values of $f$ which are used in this computation are the values $f(y)$ when $y\in   K\cap D(x_0,100\rho(x_0))$. So, boundedness of $f$ is easily removed by redefining 
 $\hat M=\max_{ K\cap D(x_0,100\rho(x_0)) } |f|$.
 }

\end{description}

\emph{Te proof is now complete. }
\begin{Remark}\ZLA{Rema:ch2RimuovBoundedDAtietze}{\rm 
We note:
\begin{enumerate}
\item\ZLA{i1Rema:ch2RimuovBoundedDAtietze}
We described how boundedness of $K$
 can be   removed. Instead, the assumption that $K$ is closed cannot be removed.
 
\item\ZLA{I2Rema:ch2RimuovBoundedDAtietze}
Inequality~(\ref{ch2:ineqExtensio}) holds since it holds for every addendum $M(x,(1+k/2^n)\rho(x))$ and since $F_n(x)$ is the arithmetic mean of these numbers.
\item\ZLA{I3Rema:ch2RimuovBoundedDAtietze}  Statement~\ref{I1Ch2TeoExteMonotSEQULebe} of Theorem~\ref{Ch2TeoExteMonotSEQULebe} holds for the extension $f_e(x)$ proposed by Tonelli  since if $f\leq g$ on a compact set $A$ then $\max_A f\leq \max_A g$ so that
\[
\underbrace{ M(x,(1+k/2^n)\rho(x))}_{\tiny \mbox{computed from $f$}}
\leq \underbrace{ M(x,(1+k/2^n)\rho(x))}_{\tiny \mbox{computed from $g$}}\,.
\] 
Statement~\ref{I2Ch2TeoExteMonotSEQULebe} of Theorem~\ref{Ch2TeoExteMonotSEQULebe} is an obvious consequence of Statement~\ref{I1Ch2TeoExteMonotSEQULebe}.\zdia
\end{enumerate}
 }
\end{Remark}

Finally we prove that statement~\ref{I2Theo:prprieSucceConvDIM1} of Theorem~\ref{Theo:prprieSucceConvDIM1}   does not hold for the extension of functions of $d$ variables obtained with the Tonelli method: \emph{the fact that a sequence $\{f_n\}$   converge pointwise on $K$ \emph{does not imply} that the sequence if the extensions converges.}

\begin{Example}
\ZLA{Exe:AppeCH2NONexteConve}
{\rm
We construct a sequence $\{f_j\}$ such that  $f_j(x)\to 0$ for every $x\in K$  while the extensions do not converge.

We consider $d=2$ and we denote $(x,y)$ the points of $\zzr^2$. The set $K$ is the boundary of the disk whose center is the origin and  of radius $1$:
\[
K=\{(x,y)\,,\quad x=\cos\zthe\,,\quad y=\sin\zthe\qquad \zthe\in[0,2\pi)\}\,.
\]
So,
\[
\rho(0,0)=1\,,\quad K\subseteq D\left ((0,0),(1+k/2^n)\rho(0,0)\right )\qquad \forall k\leq 2^{n-1} 
\]
and
\[
M\left (  (0,0), \left (1+\dfrac{k}{2^n}\right )\rho(  0,0)     \right ) =\max_K f\qquad \forall k\leq 2^{n-1} 
\]
so that
\[ F_{n}(0,0)=\max_K f\,,\qquad f_e(0,0)=\max_K f\,.
\]

Now we consider the sequence of the functions
  $f_j$   defined as follows:  
\[
\left\{\begin{array}{lll}
f_{2\nu+1}(\cos\zthe,\sin\zthe)&=&0\\[3mm]
f_{2\nu}(\cos\zthe,\sin\zthe)&=&\left\{\begin{array}{lll}
0 &{\rm if}& \zthe\notin(0,2/\nu)\\
 \nu\zthe &{\rm if}&0<\zthe\leq 1/\nu\\
2-\nu\zthe &{\rm if}& 1/\nu\leq \zthe\leq 2/\nu\,.
\end{array}\right.
\end{array}\right.
\]
  \begin{minipage}{2in}
The graphs of few of the functions $f_{2\nu} $ ora represented in the figure on the right.

The sequence $\{f_{j}(x)\}$ converges to zero for every $x\in K$. In spite of this, $(f_j)(0,0)$ oscillates, $f_{2\nu}(0,0)=1$  while $f_{2\nu+1}(0,0)=0$.\zdia 
    \end{minipage}
  \begin{minipage}{2in}
 
 \includegraphics[width=7cm,scale=.65]{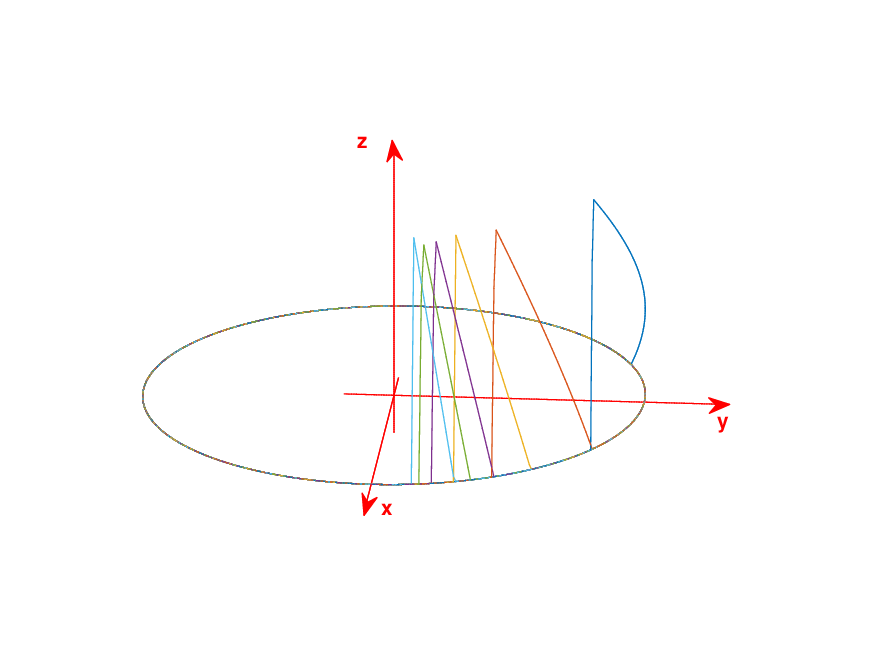}

%
  \end{minipage}
}\end{Example}

 \chapter{\ZLA{ch3:LIMITintegral}The Limits and the Integral}

The theorems concerning the exchange of  limits and   integrals are the main object of this chapter. We need few observations on multirectangles and we represent open sets as union of rectangles. This representation opens the way to   the proof of a version of    the absolute continuity of the integral.

{
\begin{table}[h]\caption{Succinct notations and terminology}
 \ZLA{tableCH3ADDITIONALshortNOTATIONS}
\ZREF{ 
When convenient, we use the informal terminology introduced in the table~\ref{tableCH2shortNOTATIONS}.   
 Moreover we introduce the following definitions.
\begin{itemize}
\item
 A {\sc  c-rectangle}\index{rectangle!c-rectangle}\index{c-!rectangle} is a closed rectangle and  a {\sc c-multirectangle}\index{c-!multirectangle}\index{multirectangle!c-multirectangle} is a multirectangle composed of  c-rectangles.
 \item A multirectangle composed of open rectangles is called an {\sc open multirectangle}\index{multirectangle!open}
 \item A multirectangle is {\sc disjoint }\index{disjoint!multirectangle}\index{multirectangle!disjoint }  when the component rectangles are pairwise disjoint.
 \item A multirectangle   is {\sc almost disjoint }\index{multirectangle!almost disjoint}\index{almost disjoint!multirectangle} when any two of its component rectangles do not have common points which are interior to at least one of them. So,  if the rectangles have interior points and intersect,   the common points belong to   faces of both the rectangles.
 Hence, \emph{an open almost disjoint  multirectangle is disjoint.}
 \end{itemize}
}
\end{table}
 \section{\ZLA{CH3MoreonMULTIREC}Special Classes of Sets and Multirectangles}
The following observations comment the definition in the table~\ref{tableCH3ADDITIONALshortNOTATIONS}.
 \begin{enumerate}
\item
 a c-rectangle is a closed set. If instead $\Delta$ is    a c-multirectangle then $\mathcal{I}_\Delta$ may not be closed. 
\item
  if a multirectangle $\Delta$ is composed of open rectangles then $\mathcal{I}_\Delta$ is open. For this reason a multirectangle composed of open rectangles is simply called an ``open multirectangle''. No confusion can arise provided that we keep in mind the following facts:
  \begin{enumerate}
  \item there exists open sets which cannot be represented as $\mathcal {I} _{\Delta}$ if we impose that $\Delta$ is disjoint (see Example~\ref{Chap2ESEnoteo:Ch2StrutturaAperti}
of Chap.~\ref{Chap:2quasicontPIUvar}).
  \item as we shall prove below, any open set is $\mathcal{I}_\Delta$
where $\Delta$ is an almost disjoint \emph{ c-multirectangle.  }
 \item in particular,   \emph{the representation of $\mathcal{I}_\Delta$ as the union of rectangles is never unique.}
   \end{enumerate}
 \item to understand the definition of almost disjointness   consider the following example: $R_1=[-1,1]$, $R_2=[0,0]$. We do not want to say that they are almost disjoint and they are not since the common point $0$ is interior to one of them.
  \item finite or numerable sets of rectangles can be ordered to form a (finite or numerable) multirectangle. Thanks to this observation we extend the term ``almost disjoint'' to finite or numerable sets $\{ R_i\}$ of multirectangles: we say that   the   {\sc rectangles are 
  almost disjoint}\index{rectangles!almost disjoint}\index{almost disjoint!rectangles} when the intersection of any pair of the rectangles either is empty or the common points belong to the boundary of both of them.
\end{enumerate}
 
 The following observations are obvious and constitute a more precise elaboration of the argument used in the proof of Theorem~\ref{TheoP2C1NullSETtwoDefI}: 
 \begin{itemize}
 
 \item
Let $R$ and $S$ be rectangles. Their union, their intersection   and their difference   are either empty or
a finite union of almost disjoint rectangles. Hence they are   sets $\mathcal{I}_\Delta$ where $\Delta$ is a finite almost disjoint multirectangle. 

Let $R$ and $S$ be rectangles and let 
 $\Delta$ and $\hat\Delta$ be \emph{finite almost disjoint} multirectangles such that  $\mathcal{I}_\Delta=R\setminus S$ and   $\mathcal{I}_{\hat\Delta}= R\cap S$. Then we have
\begin{equation}
\ZLA{CH3:EqMisuDIFFEinsie}
  L(R)=L(\Delta)+L(\hat\Delta)\,.
\end{equation}
\item
  let $R$ be a   rectangle such that $L(R)>0$. For every $\ZEP>0$ there exist  open  multirectangles $R_- $
and $R_+$
such that
\begin{equation}\ZLA{eq:CH3incluRettang}
R_-\subseteq{\rm cl}\, R_-\subseteq {\rm int}\, R\subseteq {\rm cl}\, R\subseteq {\rm int}\,R_+\quad {\rm and}\quad\left\{\begin{array}{l}
  0< L( R)-L( R_-)  <\ZEP\,,\\
  0<L( R_+)-L( R)<\ZEP\,,\\
  0<L( R_+)-L( R_-)<\ZEP\,.
\end{array}\right.
\end{equation}
 
If $R=\prod_{i=1}^d (a_i,b_i)$    then the   rectangles $R_{\pm} $ are $\prod_{i=1}^d(a_i\mp\ZSI,b_i\pm\ZSI)$
with $\ZSI$ sufficiently small.
\end{itemize}
So we can state\footnote{observe that we already used this fact in the proof of Theorem~\ref{TheoP2C1NullSETtwoDefI}.}:
\begin{Lemma}\ZLA{Lemma:Ch3RacchiudRettINmultirett}
Let
$R=\prod_{i=1}^d (a_i,b_i)$ be a rectangle such that $L(R)>0$
and let
\[
R_{\pm}=\prod_{i=1}^d (a_i\mp \ZSI,b_i\pm\ZSI)\,.
\]
  For every $\ZEP>0$ there exists $\ZDE_\ZEP>0$ such that if $0<\ZSI<\ZDE_\ZEP$  there exists 
a finite almost disjoint multirectangle $\Delta$ such that  
  
\[
\mathcal{I}_{ \Delta}=R_+\setminus R_-\,,\qquad L( \Delta)<\ZEP\,.
\]
\end{Lemma}
This observation and~(\ref{eq:CH3incluRettang}) have the following consequences\footnote{the statement~\ref{coroCH2copriOPENNULLset}   is a reformulation of Theorem~\ref{TheoP2C1NullSETtwoDefI}.}:
\begin{Corollary}\ZLA{Ch3CoriFrontiRectNULLO}
we have:
\begin{enumerate}
\item\ZLA{I1Ch3CoriFrontiRectNULLO}
The boundary of a rectangle is a null set.
\item\ZLA{coroCH2copriOPENNULLset} A set $N$ is a null set when there exists a sequence $\{\Delta_k\}$ of \emph{open multirectangles} such that 
\[
N\subseteq \mathcal{I}_{\Delta_k}\,,\qquad \lim_{k\to+\ZIN} L(\Delta_k)=0\,.
\]
\end{enumerate}
\end{Corollary}
\subsection{\ZLA{Ch3SecAlmostDisgMultir}Almost Disjoint Multirectangles and their Measure}

Let   $\Delta=\{R_k\}$ with $L(R_k)>0$  and let $\ZEP>0$. Let $ (R_k)_{\pm}$ be rectangles with the properties~(\ref{eq:CH3incluRettang}) with $\ZEP/2^k$ in the place of $\ZEP$.
Let $(\Delta_k)_\pm=\{ (R_k)_{\pm}\}$. The properties of the rectangles $(R_k)_{\pm}$ and Lemma~\ref{Lemma:Ch3RacchiudRettINmultirett} give:
\begin{Lemma}\ZLA{ch3Lemma:approxMultrect}
Let $\Delta=\{R_k\}$ be an almost disjoint  multirectangle such that
$L(R_k)>0$ for every $k$ and $L( \Delta)<+\ZIN$. Under these conditions,
 for every $\ZEP>0$ there exists a disjoint  open multirectangle $  \Delta_{-}$ and an open multirectangle $\Delta_+$ such that
\[
\mathcal{I}_{ \Delta_-}\subseteq \mathcal{I}_\Delta\subseteq \mathcal{I}_{\Delta_+}
\qquad 0\leq |L( \Delta_\pm)-L(   \Delta)|<\ZEP\,,\quad 0<L( \Delta_+)-L( \Delta_-)<\ZEP\,.
\]
\end{Lemma}

\emph{Note that the \emph{open disjoint} multirectangle $\Delta_-$ does not exist if the condition $L(R_k)>0$ is removed. If the condition $L(R_k)>0$ is remove then the statement still hold, provided that we remove the condition that $\Delta_-$ is open.}
 
We noted that the number $L(\Delta)$ is associated to the \emph{sequence} $\Delta$ and not to the set $\mathcal{I}_\Delta$. If $\Delta$ is almost disjoint then we can be more precise. First we recall:
\begin{enumerate}
\item  $L(\Delta)=\ZSUno L(R_n)$ does not depend on the order of the rectangles;
\item $L(\Delta)=\ZSUno L(R_n)$   does not change when   its component rectangles   are changed by adding or removing parts of their boundary; in particular if $R_n$ is replaced by its closure or by its interior. This 
observation is a reformulation of the statement~\ref{I1Ch3CoriFrontiRectNULLO} of Corollary~\ref{Ch3CoriFrontiRectNULLO}. 
\item
$L(\Delta)=\ZSUno L(R_n)$ does not change if the rectangles $R_n$ are represented as   countable unions of sequences of almost disjoint rectangles.
\end{enumerate}
These observations imply:
\begin{Lemma}\ZLA{ch3LemmaPrimaDefiLunivDisjquasiMult}
Let $\Delta$ be an almost disjoint multirectangle, $\Delta=\{R_k\}$ and let $ L(R_k)>0$ for every $k$. We have\footnote{''ins'' for ''inside'' and ``out'' for ''outside''.}:
\begin{subequations}
\begin{equation}\ZLA{eq:CH3RappreMeasMISESTER-A}
L(\Delta) =\inf \{ L(  \Delta_{\rm out})\,,\quad \mbox{$ \Delta_{\rm out}$ \emph{open} and $\mathcal{I}_\Delta\subseteq \mathcal{I}_{  \Delta_{\rm out }}$}\}
\end{equation}
and also
\begin{equation}\ZLA{eq:CH3RappreMeasMISESTER-B}
L(\Delta)= \sup \{ L(  \Delta_{\rm ins})\,,\quad \mbox{$  \Delta_{\rm ins}$ \emph{(open)   disjoint} and $\mathcal{I}_{  \Delta_{\rm ins}}\subseteq \mathcal{I}_{  \Delta}$}\}
\end{equation}
\end{subequations}
\end{Lemma}
(the reason for putting ``open'' in parenthesis in the formula~(\ref{eq:CH3RappreMeasMISESTER-B})  is explained in Remark~\ref{P2Ch2RemaModiDefiMultRwSET} below).

It follows: 
\begin{Theorem}
Let $ \mathcal{I}_{\Delta_1}=\mathcal{I}_{\Delta_2}$ 
and let the multirectangles be almost disjoint. We have $L(\Delta_1)=L(  \Delta_2)$.
\end{Theorem}
 Thanks to this theorem, we associate a number to every set $A$ which can be represented as $\mathcal{I}_\Delta$ when $\Delta$ is almost disjoint\footnote{strictly speaking at this point when $\Delta=\{R_k\}$ is almost disjoint and such that   $L(R_k)>0$ for every $k$. Se  Remark~\ref{P2Ch2RemaModiDefiMultRwSET} below to see that the condition
 $L(R_k)>0$   can be removed.}:
 \begin{equation}\ZLA{eq:CH3DefiMisuRectangleseTS}
\zl(A)=L(\Delta)\ \mbox{where $\Delta $ is almost disjoint and $A=\mathcal{I}_\Delta$}\,.  
 \end{equation}
 We shall see in Chap.~\ref{Chap:4Measure} that this number $\zl(A)$ is the {\sc Lebesgue 
 measure}\index{measure!Lebesgue}\index{Lebesgue!measure} of the   set $A$.\index{measure!of an open set} 
 

\begin{Remark}\ZLA{P2Ch2RemaModiDefiMultRwSET}
{\rm It is easily seen that the formula~(\ref{eq:CH3RappreMeasMISESTER-A}) holds also if the condition  $L(R_k)>0$ is removed from the statement of Lemma~\ref{ch3LemmaPrimaDefiLunivDisjquasiMult}. If this condition is removed then formula~(\ref{eq:CH3RappreMeasMISESTER-B}) holds too 
provided that we remove the request that $\Delta_{\rm ins}$ is open.
In fact, the degenerate component rectangles $R_k$ do not contain open rectangles, but we have also $L(R_k)=0$.  

This is the reason why ``open'' has been written in parenthesis.

 Note also that    if $A$ is a sequence of points, $A=\{q_k\}$, then both the numbers in~(\ref{eq:CH3RappreMeasMISESTER-A}) and~(\ref{eq:CH3RappreMeasMISESTER-B}) are zero. Compare this observation with the item~\ref{I3Exe:Ch0Qnullset1} of Remark~\ref{REMA:Ch0Qnullset}.

}
\end{Remark}

In order to stress the previous observation we state explicitly:
\begin{Theorem}\ZLA{TEOeq:CH3RappreMeasMISESTER-C}
Let $A=\mathcal{J}_{\Delta}$ and let $\Delta=\{R_k\} $ be almost disjoint. The number $\zl(A)=L(\Delta)$ is given by both the formulas~(\ref{eq:CH3RappreMeasMISESTER-A}) and~(\ref{eq:CH3RappreMeasMISESTER-B}) with the following warning: if   $L(R_k)$ is not strictly positive   then formula~(\ref{eq:CH3RappreMeasMISESTER-B}) takes the following form:
\begin{equation}\ZLA{eq:CH3RappreMeasMISESTER-C}
L(\Delta)= \sup \{ 
L(  \Delta_{\rm ins})\,,\quad
 \mbox{$  \Delta_{\rm ins}\ $   disjoint  and $\mathcal{I}_{  \Delta_{\rm ins}}\subseteq \mathcal{I}_{  \Delta}$}\}\,.
\end{equation}
\end{Theorem}
 
\smallskip
\begin{table}[h]\caption{Multirectangle sets}\ZLA{tableCH3ADDITIONALshortMulti} 
\ZREF{
To streamline the  presentation, a set $A=\mathcal{I}_\Delta$ where $\Delta=\{R_k\}$ is almost disjoint  
will be called a {\sc    multirectangle set}\index{set!multirectangle set}\index{multirectangle!set} and, if convenient, we specify {\sc the multirectangle set of $\Delta$.} We say also that {\sc $A$ is the set of the multirectangle $\Delta$}\index{set!of a multirectangle} and, reciprocally, that {\sc $\Delta$ is a multirectangle of the set $A$.}\index{multirectangle!of a set}

\emph{We stress that when using these terms we always intend that $\Delta$ is almost disjoint.}
} 

\end{table}
\smallskip
%
%
%

The following result is clear  from~(\ref{eq:CH3RappreMeasMISESTER-A}):

 \begin{Lemma}
 \ZLA{CH3:lemmaMONOTONmisura} Let $A=\mathcal{I}_{\Delta_1}$ and $B=
 \mathcal{I}_{\Delta_2}$ be multirectangle sets. We have:
\begin{enumerate}
\item  \ZLA{I1CH3:lemmaMONOTONmisura}
 If $A\subseteq B$ then $\zl(A)\leq \zl(B)$.
 \item \ZLA{I2CH3:lemmaMONOTONmisura} Let $\Delta_1=\{R_{1,n}\}$ and $\Delta_2=\{R_{2,n}\}$ be almost disjoint multirectangles\footnote{we recall the notation $ \Delta_1\cup\Delta_2$ to denote the multirectangle composed by  the rectangles $R_{1,n}$ and $R_{2,k}$ (taken in any order).}.  Then:
 \begin{align*}
 R_{1,n}\cap R_{2,j}=\emptyset \ \ \forall n\,,\ j\quad \implies\quad 
\zl(\mathcal{I}_{{\Delta_1\cup\Delta_2}})=L(\Delta_1\cup\Delta_2)&=L(\Delta_1)+L(\Delta_2)\\
&=\zl(\mathcal{I}_{\Delta_1})+\zl(\mathcal{I}_{\Delta_2})
\,.
 \end{align*} 
 This property can be written as
 \[
A\cap B=\emptyset\ \implies \zl(A\cup B)=\zl(A)+\zl(B)\quad \mbox{where}\quad A=\mathcal{I}_{\Delta_1}\,,\ \  B=\mathcal{I}_{\Delta_2}\,. 
 \]
 \end{enumerate}
 \end{Lemma}
 
 The property in the statement~\ref{I1CH3:lemmaMONOTONmisura} of Lemma~\ref{CH3:lemmaMONOTONmisura} is the  
  {\sc monotonicity property of the measure}\index{measure!monotonicity}
  while that in the statement~\ref{I2CH3:lemmaMONOTONmisura} is the
  {\sc additivity property of the measure.}\index{measure!additivity} 
\begin{Remark}{\rm
We observe:
\begin{enumerate}
\item If $A$ is a multirectangle set then it does not identify $\Delta$ uniquely;
\item if $A$ is the multirectangle set of $\Delta$, there exist multirectangles   $\tilde\Delta$ which are not almost disjoint and such that $A=\mathcal{I}_{\tilde\Delta}$.
\item
The computation used in Example~\ref{Exe:Ch0Qnullset}   to see that the set $\{q_n\}$, the set of the rational points in $  (0,1)$, is a null set is an application of ~(\ref{eq:CH3RappreMeasMISESTER-A}).\zdia
\end{enumerate}
}\end{Remark}

 \subsection{\ZLA{Ch3SECTapertiABScont}Representation of Nonempty Open Sets}
 
Important sets to which Lemma~\ref{ch3LemmaPrimaDefiLunivDisjquasiMult} 
and formula~(\ref{eq:CH3DefiMisuRectangleseTS}) 
can be applied are the open sets. In fact, 
in this section we represent any nonempty    open set  in terms of almost disjoint multirectangles.

  Example~\ref{Chap2ESEnoteo:Ch2StrutturaAperti}
of Chap.~\ref{Chap:2quasicontPIUvar} 
shows  that, when $d>1$, an open set is not a disjoint  open multirectangle in the sense that it is not equal to any such $\mathcal{I}_\Delta$. Instead we have the following result, which we state in the case of bounded open sets since this is the case we shall need, but which can be extended also to unbounded open sets (in this case without the   condition $L( \Delta)<+\ZIN$):
\begin{Theorem}\ZLA{TeoCH3STRUttAPERTI}
Let $\mathcal{O}$ be a  nonempty  bounded open set. There exists an almost disjoint   c-multirectangle $\Delta$ such that $L( \Delta)<+\ZIN$ and   $\mathcal{O}=\mathcal{I}_\Delta$. It is possible to chose $\Delta=\{R_n\}$ such that $L(R_n)>0$               for every $n$.
\end{Theorem}
\zProof  
The construction of the multirectangle $\Delta$ is iterative, as follows:

\begin{description}
\item[\bf Step 0:] We divide $\zzr^d$ with \emph{closed} ``cubes''   whose vertices are the points with integer coordinates. This way we get a set $\Q_0$ of closed cubes   with sides parallel to the coordinate axes and   of length $1=1/2 ^0$.

Two different cubes either have empty intersection or they intersect along a face.

 Then we perform the following operations:
\begin{description}
\item[\bf Substep 0-a:] we single out the cubes of $\Q_0$ which are contained in $\mathcal{O}$. We call $\QP_0$ the set of these cubes\footnote{the index $_{\pmb P}$ of
$\QP$   is for ``present'', since these cubes are
retained in this step while    $_{\pmb N}$
in $\QN$  is for ``next'' since these cubes are elaborated in the next step.}:
\[
\QP_0=\{Q\in\Q_0\,:\quad Q\subseteq \mathcal{O}\}\,. 
\]
Then we define 
\[
\QN_0=\{Q\notin \QP_0 \,,\quad Q\cap \mathcal{O}\neq \emptyset\}\,.
\]
The set of cubes $\QP_0$ is either empty or finite while $\QN_0$ is always nonempty and finite 
(these sets are finite since  $ \mathcal{O}$ is bounded).

 \emph{Note that both $\QP_0$ and $\QN_0$  are sets of almost disjoint closed rectangles\footnote{of course this statement holds for $\QP_0$ when it is not empty. A similar warning holds also in the following steps and it is not repeated.}.}

\item[\bf Substep 0-b:]  
The elements of $\QP_0$  (if any) are retained as elements of $\Delta$.

\end{description}

\item[\bf Step 1:] We divide every cube $Q\in \QN_0$ in $2^d$  equal \emph{closed} cubes (with planes which are parallel to the coordinate planes and which cut the edges in the middle). This way we get a set $\Q_1$ of cubes  with sides parallel to the coordinate axes and   of length $ 1/ 2 ^1$. Then we perform the following operations:

\begin{description}
\item[\bf Substep 1-a:] we single out the cubes of $\Q_1$ which are contained in $\mathcal{O}$. We call $\QP_1$ the set of these cubes:
\[
{\QP}_1=\{Q\in\Q_1\,:\quad Q\subseteq \mathcal{O}\}\,.
\]
Then we define
\[
\QN_1=\{ Q\in \Q_1\setminus\QP_1 \,,\quad  Q\cap \mathcal{O}\neq \emptyset\}\,.
\]

\item[\bf Substep 1-b:]  
The set $\QP_1$ is empty or finite. If   nonvoid, its elements are retained as elements of $\Delta$.

 
  \end{description}
    \end{description}

  \begin{description}
 
\item[\bf Step j:] As in the previous steps, we divide every cube $Q\in \QN_{j-1}$ in $2^d$ equal \emph{closed} cubes with planes which cut the edges in the middle. This way we get a set $\Q_j$ of cubes, with sides parallel to the coordinate axes ad   of length $ 1/2^j$. Then we perform the following operations:

\begin{description}
\item[\bf Substep j-a:] we single out the cubes of $\Q_j$ which are contained in $\mathcal{O}$. We call $\QP_j$ the set of these cubes:
\[
\QP_j=\{Q\in\Q_j\,:\quad Q\subseteq \mathcal{O}\}\,.
\]
Then we define 
\[
\QN_j=\{Q\in \Q_j\setminus\QP_j\,,\quad Q\cap \mathcal{O}\neq \emptyset\}\,.
\]

\item[\bf Substep j-b:]  
 The elements of $\QP_j$ are retained as elements of~$\Delta$.

\end{description}
 
\item[\bf We iterate this procedure.]
\end{description}

This way we obtain a \emph{almost disjoint  c-multirectangle} $\Delta$
  such that $\mathcal{I}_\Delta\subseteq \mathcal{O}$:  the component rectangle of $\Delta$ are the ``squares'' $Q\in  \cup _{k=1}^{+\ZIN}\QP_k$.  The inclusion is clear and we have also $\mathcal{I}_\Delta=\mathcal{O}$  since any point   $x\in \mathcal{O}$ is an interior point and    it exists  a closed  cube whose sides have length less then $1/2^n$ ($n$ suitably large), which is contained in $\mathcal{O}$ and which contains $x$.\zdia

 \emph{An important consequence of Theorem~\ref{TeoCH3STRUttAPERTI} is that Lemma~\ref{ch3LemmaPrimaDefiLunivDisjquasiMult} can be applied and the number $\zl(\mathcal{O})$ is well defined for every (nonempty) open set $\mathcal{O}$.}
We combine Lemma~\ref{ch3Lemma:approxMultrect} and Theorem~\ref{TeoCH3STRUttAPERTI} and we get:
\begin{Theorem}
\ZLA{TeoCH3quasiSTRUttAPERTI}Let $\mathcal{O}$ be a nonempty bounded \emph{open} set. For every $\ZEP>0$ there exist a  disjoint   open
 multirectangle $  \Delta_{\rm ins}$ and an open multirectangle $ \Delta_{\rm out}$ such that
\[
    \Delta_{\rm ins} \subseteq \mathcal{O} \subseteq   \Delta_{\rm out}\,,\qquad L(  \Delta_{\rm out})-L(  \Delta_{\rm ins})<\ZEP\,.
\]
The multirectangle  $ \Delta_{\rm ins}$ can be chosen finite:
\[
 \Delta_{\rm ins}=\{R_1\,,\ R_n\,,\ \dots\,,\ R_K\}\quad \mbox{with $\left\{\begin{array}{l}
({\rm cl} \, R_i)\cap ({\rm cl} \, R_j)  =\emptyset\\
{\rm cl} \, R_i\subseteq \mathcal{O}\qquad 1\leq i\leq K\,,\\
L(R_i)>0\ \mbox{for every $i$}\,.
\end{array}\right.$} 
\]
\end{Theorem}

We put  \[
A=\bigcup _{i=1}^K R_i=\mathcal{I}_{ \Delta_{\rm ins}}\,,\qquad B=\mathcal{I}_{ \Delta_{\rm out}}\,.
\]
We use additivity of the measure. The previous theorem can be reformulates as follows:
\begin{Corollary}\ZLA{CoroCH3quasiSTRUttAPERTI}
Let $\mathcal{O}$ be a bounded open set. For every $\ZEP>0$ there exist  open sets $A$ and $B$ such that\footnote{observe that we are not asserting ${\rm cl}\, O\subseteq B$. This inclusion is false in general
   and so the inclusion $\mathcal{O}\subseteq B$ is trivial when studying open sets  since we can choose $B=\mathcal{O}$.} $A\subseteq{\rm cl}\,A\subseteq \mathcal{O} \subseteq B$ such that
\[
 \zl(\mathcal{O}\setminus {\rm cl}\,A)<\ZEP\,,\qquad \zl(B)-\zl( \mathcal{O})<\ZEP\,.
\]
\end{Corollary}

So, when $d>1$ an open set is not a disjoint  open multirectangle, but ``its difference    from a disjoint  open multirectangle is as small as we whish''.

We combine Lemma~\ref{ch3Lemma:approxMultrect} and Theorem~\ref{TeoCH3STRUttAPERTI} and we see that the formulas~(\ref{eq:CH3RappreMeasMISESTER-A}) 
and~(\ref{eq:CH3RappreMeasMISESTER-B}) hold for every open set:
 
\begin{equation}\ZLA{ch3DefiLebeMISaperti}
\begin{array}{l}
 \mbox{if $\Delta$ is almost disjoint and $\mathcal{O}=\mathcal{I}_\Delta$ is open and bounded then}\\
\zl(\mathcal{O})=L(\Delta)
 =\left\{\begin{array}{l}
\sup\{L( \Delta_{\rm ins})\,,\ \mbox{$ \Delta_{\rm ins}$ finite \emph{open} disjoint and $\mathcal{I}_{ \Delta_{\rm ins}}\subseteq  \mathcal{O}$}\} \\
\inf\{L(\Delta_{\rm out})
\,,\ \mbox{$ \Delta_{\rm out}$    \emph{open} multirectangle, $\mathcal{I}_{  \Delta_{\rm out}}\supseteq \mathcal{O}$\}
}\,.
\end{array}\right.
\end{array}
\end{equation}

\section{\ZLA{Cha3SECEgorovSev}Egorov-Severini Theorem and Quasicontinuity}

In this section we discuss the version of Egorov-Severini Theorem for functions of $d$ variables. Note that the statement is the same as that in Sect.~\ref{Cha1bisSECEgorovSev} but the proof is different for a reason we illustrate below\footnote{but note that the proof we are going to give here holds in any dimension, also in dimension~$1$.}.
 
 First of all,
as in the case of dimension $1$, we present a lemma and a definition:
 \begin{Lemma}\ZLA{LemmaCH3SucceFunzQCinDiffeR}
Let $\{f_n\}$ be a sequence of functions each one a.e. defined on a rectangle $R$; i.e.,
\[
\Dom\, f_n=R\setminus N_n \qquad \mbox{($N_n$ is a null set)}\,.
\]
Then we have:
\begin{enumerate}
\item
\ZLA{I1LemmaCH3SucceFunzQCinDiffeR}
there exists a null set $N$ such that every $f_n$ is defined on $R\setminus N$\,.
\item\ZLA{I2LemmaCH3SucceFunzQCinDiffeR} if each $f_n$ is quasicontinuous then for every $\ZEP>0$ there exists a multirectangle $\Delta_\ZEP$ such that $\mathcal{I}_{\Delta_\ZEP}$ is open and such that $(f_n)_{|_{R\setminus \mathcal{I}_{\Delta_\ZEP}}}$ is continuous.
\end{enumerate}
 \end{Lemma}
 
 The proof is similar to that of Lemma~\ref{LemmaCha1bisSucceFunzQCinDiffeR}.
 
\ZREF{
 This observation shows that when we have a sequences of  functions each one of them defined a.e. on $R$ we can    assume that they are all defined on $R\setminus N$ where $N$
 is a null set which does not depend on $n$. To describe this case we say (as   in Chap.~\ref{Ch1:INTElebeUNAvaria}) that \emph{the sequence $\{f_n\}$ is defined a.e. on $R$.} 
 }
 
Similar to   Definition~\ref{Cha1bisDEFINITalmUNIFconve}:
\begin{Definition}\ZLA{Ch3DEFINITalmUNIFconve}
{\rm Let $f_n$,   $f$ be  functions a.e. defined on a  set  $K$. We say that the sequence $\{ f_n\}$  {\sc converges almost uniformly}\index{convergence!almost uniform}\index{almost uniform convergence} to $f$  on $K$ when for every $\ZEP>0$     the following \emph{equivalent} statements hold:
\begin{enumerate}
\item there exists an open set $\mathcal{O}$ such that 
\[
\mbox{$\zl(\mathcal{O})<\ZEP $  and $\{ f_n\}$ converges \emph{uniformly} to $f$ on $K\setminus\mathcal{O}$}\,.  
\]
\item
there exists a multirectangle\footnote{observe a discrepancy between the present definition and Definition~\ref{Cha1bisDEFINITalmUNIFconve}: here we must explicitly state that~$\mathcal{I}_\Delta$ is open. This specification was not needed in Definition~\ref{Cha1bisDEFINITalmUNIFconve} since when $d=1$ multiinterval are composed by open intervals and because of
  Theorem~\ref{teo:Ch2StrutturaAperti}.}   $\Delta$   such that 
\[
\mbox{$
L(\Delta)<\ZEP$,    $\mathcal{I}_\Delta$ is open and $\{ f_n\}$ converges \emph{uniformly} to $f$ on $K\setminus{\mathcal{I}_\Delta}$}\,.  
\] 
The multirectangle can be chosen c-closed and almost disjoint.
\end{enumerate}

 }
\end{Definition}

 We note: ``$\{ f_n\}$ converges uniformly  to $f$ on $K\setminus\mathcal{O} $'' 
is equivalent to ``$\{ (f_n)_{|_{K\setminus\mathcal{O}}}\}$ converges uniformly to $f_{|_{K\setminus\mathcal{O}}}$''. 

\ZREF{
We must be clear on the content of this definition. 
The sequence $\{f_n\}$ converges almost uniformly to $f$ when the following holds:
   we fix any $\ZEP>0$ and we find an open set $\mathcal{O}_\ZEP$ such that
$\zl(\mathcal{O}_\ZEP)<\ZEP$ and such that
 the following property is valid: for every $\ZSI>0$ there exists a number $N$ \emph{which depends on $\ZSI$ and on the previously chosen set $\mathcal{O}_\ZEP$, $N=N_{\ZSI,\mathcal{O}_\ZEP }$,} and such that
\[
\left\{\begin{array}{l}
n>N_{\ZSI,\mathcal{O}_\ZEP }\\
x\in R\setminus \mathcal{O}_\ZEP
\end{array}\right.\quad\implies \ |f_n(x)-f(x)|<\ZSI\,.
\]
\emph{The important point is that $\mathcal{O}$ does not depend on $\ZSI$.}
 }
 
Now we can state a result similar to that in Theorem~\ref{Cha1bisTeoEgoSEVERANTE}:
\begin{Theorem}[Egorov-Severini: preliminary statement] \ZLA{P2CH1TeoEgoSEVERante} Let $\{ f_n\}$ be a  sequence of \emph{continuous} functions defined  on the \emph{closed and bounded} rectangle $R$. If the sequence converges   on $R$
to a function $f$:
\[
\limn f_n(x)=f(x) \qquad \fbox{$\forall x\in R$}
\]  
 then:
\begin{enumerate}
\item\ZLA{P2CH1I1TeoEgoSEVERante}
the sequence    converges almost uniformly;
\item\ZLA{P2I2TeoEgoSEVERante} the limit function $f$ is quasicontinuous.
\end{enumerate}
\end{Theorem}
 
The proof is  in Appendix~\ref{P2appeCH3EGOSEVE}.

A consequence (actually an equivalent formulation) is Egorov-Severini Theorem\footnote{the   proofs by Egorov in~\cite{Egorov1911} and by Severini in~\cite{Severini1910} concern functions of one variable.}:

\begin{Theorem}[{\sc Egorov-Severini}]\index{Theorem!Egorov-Severini}
\ZLA{P2CH3TeoEgoSEVER} Let $\{f_n\}$   be a sequence of \emph{quasicontinuous} function \emph{a.e. defined} on a  rectangle $R$.   We assume:
\begin{enumerate}
\item the rectangle is bounded;
\item for every $\ZEP>0$ there exists an open set $\mathcal{O}$ such that
$\zl(\mathcal{O})<\ZEP$ and
 $\{f_n\}$   is a bounded sequence on $R\setminus\mathcal{O}$;
\item we have
\[
\limn f_n(x)=f(x)\qquad \fbox{a.e. $x\in R$}\,.
\]
\end{enumerate}
Under these conditions:
\begin{enumerate}
\item\ZLA{I1P2TeoEgoSEVER}the convergence is almost uniform;
\item\ZLA{I2P2TeoEgoSEVER}
 $f$ is quasicontinuous.
 \end{enumerate}
\end{Theorem}

As seen in Remark~\ref{CHCha1bis:RemaESEunBEgorSEV} \emph{the assumption that $R$ is bounded cannot be removed.}

 The two theorems~\ref{P2CH1TeoEgoSEVERante}    and~\ref{P2CH3TeoEgoSEVER} are equivalent since the first is a particular case of the second and, in its turn, implies the second. This fact precisely correspond to the fact which holds when $d=1$ but the proof now is less direct since we 
 cannot relay on a convergence property analogous to that in the statement
 \ref{I2Theo:prprieSucceConvDIM1} of Theorem~\ref{Theo:prprieSucceConvDIM1}. Instead, the proof relays on the the monotonicity   in Theorem~\ref{Ch2TeoExteMonotSEQULebe}. 
 
 Before proving that Theorem~\ref{P2CH1TeoEgoSEVERante} implies Theorem~\ref{P2CH3TeoEgoSEVER}, we state the following corollaries, which are analogous  to Corollary~\ref{Cha1bisCoroEGOSEveQuasContFunz} 
 and~\ref{Coro:Cha1bisSUPinfPERch4} seen when $d=1$. 
 
 Corollary~\ref{Coro:Ch3SUPinfPERch4} is a consequence of Corollary~\ref{Ch3CoroEGOSEveQuasContFunz} and the proof is analogous to that in Sect.~\ref{SecP1CH2ConseEGOrSeve} while   Corollary~\ref{Ch3CoroEGOSEveQuasContFunz} does not need an independent proof since it is a step in the proof of Egorov-Severini Theorem.

 \begin{Corollary}\ZLA{Ch3CoroEGOSEveQuasContFunz}
 Let $\{f_n\}$ be a  bounded sequence of quasicontinuous functions a.e. defined on a bounded rectangle $R$.
The following properties hold:
\begin{enumerate}
\item\ZLA{I1Ch3CoroEGOSEveQuasContFunz} let  $\{f_n\}$  be either a.e. increasing or   decreasing on $R$  and let
\[
f(x)=\limn f_n(x) \qquad a.e. \ x\in R\,.
\] 
The function $f$ is quasicontinuous.  
\item\ZLA{I2Ch3CoroEGOSEveQuasContFunz} Let $f$ be either
  \[
f(x)=\limsup_{n\to+\ZIN} f_n(x) \quad {\rm or}\quad 
f(x)=\liminf_{n\to+\ZIN} f_n(x)
     \qquad a.e. \ x\in R\,.  
  \]
  The function $f(x)$ is quasicontinuous.
  \end{enumerate}
 \end{Corollary}

Similar to Corollary~\ref{Coro:Cha1bisSUPinfPERch4}:
\begin{Corollary}
\ZLA{Coro:Ch3SUPinfPERch4}
Let $\{f_n\}$ be a bounded sequence of quasicontinuous functions defined on a rectangle $R$ and let   
 
\begin{align*}
\phi(x)&=\sup\{ f_n(x)\,,\quad n\geq 1\}\\
\psi(x)&=\inf\{ f_n(x)\,,\quad n\geq 1\}\,.
\end{align*}
The functions $\phi$ and $\psi$ are quasicontinuous.
\end{Corollary}

 \paragraph{The Proof that Theorem~\ref{P2CH1TeoEgoSEVERante} Implies Theorem~\ref{P2CH3TeoEgoSEVER}}

Due to the fact that $\partial R$ is a null set, we can assume that $R$ is closed. We proceed in 3 steps.

\begin{description}
\item[{\bf Step~1: the case that $\{f_n\}$ is a monotone sequence.}]
\end{description}
Let $\{f_n\}$ be increasing. We fix any $\ZEP>0$ and any set $\tilde{\mathcal{O}}_\ZEP$ such that $\zl(\tilde{\mathcal{O}}_\ZEP)<\ZEP$ and such that:
\begin{enumerate}
\item the sequence $\{f_n\}$ is bounded on $R\setminus\tilde{\mathcal{O}}_\ZEP$;
\item for every $n$, the restriction of $f_n$ to $R\setminus \tilde{\mathcal{O}}_\ZEP$
is continuous.
\end{enumerate}
We denote $f_{n,e}$ the Tietze extension of $f_n$ obtained from Tonelli algorithm\footnote{or, from any other algoritm provided that it preserves monotonicity of sequences.} in the Appendix~\ref{AppeCH2}. We have:
\begin{enumerate}
\item every function $f_{n,e}$ is continuous on $R$;
\item the sequence $\{f_{n,e}\}$ is bounded and increasing on $R$ and so $\{f_{n,e}(x)\}$ converges pointwise to a function $\tilde f$ for every $x\in R$.
\item on $R\setminus\tilde{\mathcal{O}}_\ZEP$ we have $f_{n,e}(x)=f(x)$ and so also $\tilde f(x)=f(x)$.
\end{enumerate}
The sequence $\{f_{n,e}\}$ satisfies the assumptions of
  Theorem~\ref{P2CH1TeoEgoSEVERante}. So, there exists an open set $\hat{\mathcal{O}}_\ZEP$ such that
\[
\zl\left (\hat{\mathcal{O}}_\ZEP\right )<\ZEP\,,\quad \mbox{$f_{n,e}\to \tilde f $ \emph{uniformly} on $R\setminus\hat{\mathcal{O}}_\ZEP $}\,.
\]
Uniform convergence implies continuity of $\tilde f _{|_{R\setminus\hat{\mathcal{O}}_\ZEP }}$ and so
\begin{align*}
&
 f _{|_{R\setminus \left (\tilde{\mathcal{O}}_\ZEP\cup \hat{\mathcal{O}}_\ZEP\right )}}=\tilde f _{|_{R\setminus \left (\tilde{\mathcal{O}}_\ZEP\cup \hat{\mathcal{O}}_\ZEP\right )}}
\quad \mbox{is continuous}\\
&\mbox{$f_n\to f$ uniformly on $R\setminus \left (\tilde{\mathcal{O}}_\ZEP\cup \hat{\mathcal{O}}_\ZEP\right )$}\\
 &\zl \left (\tilde{\mathcal{O}}_\ZEP\cup \hat{\mathcal{O}}_\ZEP\right )<2\ZEP\,.
\end{align*}
This argument shows that Egorov-Severini Theorem holds for \emph{increasing} sequences since $\ZEP>0$ is arbitrary.

Analogously it is proved that it holds for \emph{decreasing} sequences too.
\begin{description}
\item[{\bf Step~2: quasicontinuity of $f$   even if    $\{f_n\}$ is not monotone.}]
\end{description}
We give two proofs of this statement since the formulas we find in the two cases are both used in the third step.
\begin{description}
\item[\emph{The first proof.}]
\end{description}
We consider a sequence $\{f_n\}$ which satisfies the assumption of Egorov-Severini Theorem. 
We do not assume that $\{f_n\}$ is monotone as we did in the {\bf Step~1.}  In spite of this, we prove that $f$ is quasicontinuous.

We fix $\ZEP>0$ and the open set $\tilde{\mathcal{O}}_\ZEP$ as in the {\bf Step~1,} such that
\[
\zl\left (\tilde{\mathcal{O}}_\ZEP\right )<\ZEP\,,\qquad 
f_n(x)\to f(x)\qquad \forall x\in R\setminus \tilde{\mathcal{O}}_\ZEP\,.
\]

For every pair   $N $   and $m$ of natural numbers,  we define
\[
\psi^{(s)}_{N,m}=\max\{ f_N(x)\,,\ f_{N+1} (x)\,,\ \dots\,,\ f_{N+m}(x)\}\,.
\]
The functions $\psi^{(s)}_{N,m}$ satisfy the assumptions of Egorov-Severini Theorem and furthermore    sequence $m\mapsto \psi^{(s)}_{N,m}$ is \emph{increasing and bounded on $R\setminus\tilde{\mathcal{O}}_\ZEP$.} Hence, it is convergent and from the {\bf Step~1,}
\[
\phi^{(s)}_N =\lim _{m\to +\ZIN}\psi^{(s)}_{N,m}
\]
is quasicontinuous.

The sequence $N\mapsto\phi^{(s)}_N$ is \emph{decreasing and bounded $R\setminus\tilde{\mathcal{O}}_\ZEP$.}  It satisfies the Assumption of Egorov-Severini Theorem and
 
\[
\lim_{N\to+\ZIN} \phi^{(s)}_N=\underbrace{\limsup_{n\to+\ZIN} f_n(x)=f(x)}_{\tiny\rm\mbox{since $\{f_n\}$ is a.e. convergent by assumption}} \quad {\rm a.e.}\ x\in R\,.
\]
 We invoke again    the {\bf Step~1} 
and we see that 
 \[
f=\lim_{N\to+\ZIN} \phi^{(s)}_N\quad \mbox{is quasicontinuous} 
\]
as wanted.

\ZREF{We recapitulate what we found: we fix arbitrary numbers $\ZEP>0$ and $\ZSI>0$. Then:
\begin{multline}
\ZLA{EqP2Ch2QcontPriMet}
 \exists \mathcal{O}^{(s)}_\ZEP\,, \exists N^{(s)}_{\ZSI, \mathcal{O}^{(s)}_\ZEP}\,:
\ \zl( \mathcal{O}^{(s)}_\ZEP)<\ZEP\ \mbox{and}\\ 
\left\{\begin{array}{l}
x\in R\setminus \mathcal{O}^{(s)}_\ZEP\\
N>N^{(s)}_{\ZSI, \mathcal{O}^{(s)}_\ZEP}
\end{array}\right.
  \ \implies\ |\phi^{(s)}_N(x)-f(x)|<\ZSI \,.
\end{multline}
}
{
 \begin{description}
\item[\emph{The second proof.}]
\end{description}
In the Step~2 we proved quasicontinuity of $f$ by using $f=\limsup_{n\to+\ZIN} f_n$. We can equivalently use $f=\liminf_{n\to+\ZIN} f_n$. We give the details  which will be used in the {\bf Step~3.}
We fix $N>0$ and we define
\[
\psi^{(i)}_{N,m}=\min\{ f_N(x)\,,\ f_{N+1} (x)\,,\ \dots\,,\ f_{N+m}(x)\}\,.
\]
The functions $\psi^{(i)}_{N,m}$ satisfy the assumptions of the Egorov-Severini Theorem and the sequence $m\mapsto \psi^{(i)}_{N,m}$ is \emph{decreasing.} Hence, from the {\bf Step 1,}
\[
\phi^{(i)}_N =\lim _{m\to +\ZIN}\psi^{(i)}_{N,m}
\]
is quasicontinuous.

The sequence $N\mapsto\phi^{(i)}_N$ (which satisfies the assumptions of the Egorov-Severini Theorem) is \emph{increasing} and, from the {\bf Step 1,} \[
\lim_{N\to+\ZIN} \phi^{(i)}_N\quad \mbox{is quasicontinuous}\,.
\]
Now we use
\[
\lim_{N\to+\ZIN} \phi^{(i)}_N=\underbrace{\liminf_{n\to+\ZIN} f_n(x)=f(x)}_{\tiny\rm\mbox{since $\{f_n\}$ is a.e. convergent by assumption}} \quad {\rm a.e.}\ x\in R
\]
and we conclude that $f$ is quasicontinuous.

\ZREF{We recapitulate what we found: we fix arbitrary numbers $\ZEP>0$ and $\ZSI>0$. Then:
\begin{multline}
\ZLA{EqP2Ch2QcontSEcMet}
 \exists \mathcal{O}^{(i)}_\ZEP\,, \exists N^{(i)}_{\ZSI, \mathcal{O}^{(i)}_\ZEP}\,:
\zl(\mathcal{O}^{(i)}_\ZEP)<\ZEP\ \mbox{and} 
 \\
\left\{\begin{array}{l}
x\in R\setminus \mathcal{O}^{(i)}_\ZEP\\
N>N^{(i)}_{\ZSI, \mathcal{O}^{(i)}_\ZEP}
\end{array}\right.
  \ \implies\ |\phi^{(i)}_N(x)-f(x)|<\ZSI \,.
\end{multline}
}
}
\begin{description}
\item[{\bf Step~3: end of the proof:    $\{f_n\}$ converges almost uniformly
to $f$ on $R$.}]
\end{description}
We assign arbitrary positive numbers $\ZEP$ and $\ZSI$ and we combine both the result in the {\bf Step~2.} We put
\[
\mathcal{O}_\ZEP= \mathcal{O}^{(s)}_\ZEP\cup  \mathcal{O}^{(i)}_\ZEP\,, \
N_{\ZSI,\mathcal{O}_\ZEP}=\max\left \{N^{(s)}_{\ZSI, \mathcal{O}^{(s)}_\ZEP},N^{(i)}_{\ZSI, \mathcal{O}^{(i)}_\ZEP}\right \}\,.
\]
The formulas~(\ref{EqP2Ch2QcontPriMet}) and~(\ref{EqP2Ch2QcontSEcMet}) give: 
\begin{multline*}
\zl(\mathcal{O}_\ZEP)<2\ZEP \,, \ \mbox{and}\\
\left\{
\begin{array}{l}
x\in R\setminus \mathcal{O}_\ZEP\,,\\
N>N_{\ZSI,\mathcal{O}_\ZEP}
\end{array}\right.\ \implies\left\{\begin{array}{l}
f(x)-\ZSI< \phi^{(s)}_N(x)<f(x) +\ZSI\,,\\
 f(x)-\ZSI<\phi^{(i)}_N(x)<f(x) +\ZSI\,.
\end{array}\right.
\end{multline*}

Now we observe:
\begin{subequations}
\begin{equation}\ZLA{Ch3P2Ste3ProvEgSevA}
 f_N(x)\leq \phi_N^{(s)}(x) \ \mbox{so that}\ \left(\left\{\begin{array}{l}
x\in R\setminus \mathcal{O}_\ZEP\,,\\ 
 N>N_{\ZSI,\mathcal{O}_\ZEP}
 \end{array}\right.\ \implies 
 f_N(x)<f(x) +\ZSI\right )\,.
\end{equation}

Analogously:
\begin{equation}\ZLA{Ch3P2Ste3ProvEgSevB}
 f_N(x)\geq \phi_N^{(i)}(x) \ \mbox{so that}\ \left(\left\{\begin{array}{l}
x\in R\setminus \mathcal{O}_\ZEP\,,\\ 
 N>N_{\ZSI,\mathcal{O}_\ZEP}
 \end{array}\right.\ \implies f(x) -\ZSI<
 f_N(x)  \right )\,.
\end{equation}
\end{subequations}

\ZREF{
We combine~(\ref{Ch3P2Ste3ProvEgSevA}) and~(\ref{Ch3P2Ste3ProvEgSevB}): for every $\ZEP>0$ we found an open set $\mathcal{O}_\ZEP$
such that $\zl(O_\ZEP)<2\ZEP$ and with the following property: for every $\ZSI>0$ there exists $N$ \emph{which  depends on $\ZSI$ and on the previously chosen and fixed set $O_\ZEP$, $N=N_{\ZSI,\mathcal{O}_\ZEP}$,} such that
\[
N>N_{\ZSI,\mathcal{O}_\ZEP} \ \implies \Big (
|f_N(x)-f(x)|<\ZEP\ \ \forall x\in R\setminus\mathcal{O}_\ZEP
\Big )\,.
\]
This is the property that  the sequence converges uniformly on $R\setminus\mathcal{O}_\ZEP$. The proof is finished since $\ZEP>0$ is arbitrary.\zdia
} 
\section{\ZLA{Ch3SUBSabsoluConti}Absolute Continuity of the Integral}
As in the case of functions of one variable, the integral is an absolutely continuous set function, i.e. the following result holds\footnote{this statement of absolute continuity is not the most general. The  general statement is 
 in Sect.~\ref{sectCH4MEASUREgeneral}, Theorem~\ref{TeoCH4ABsolContINgene}.}:
 \begin{Theorem}\ZLA{Teo:CH3ABSOluContiPERiSoliAperti}
 Let $f$ be summable on $\zzr^d$. The set valued function 
 \[
 \mathcal{O}\mapsto \int _{\mathcal{O}} f(x)\ZD x\qquad \mbox{($\mathcal{O}$ open)}
 \]
is well defined and it  is {\sc absolutely continuous}\index{set function!absolute continuity}\index{absolute continuity} in the sense that for every $\ZEP>0$ there exists $\ZDE>0$ such that 
 \begin{equation}\ZLA{EqTeo:CH3ABSOluContiPERiSoliAperti}
 \zl(\mathcal{O})<\ZDE\ \implies\ 
\left |\int _{\mathcal{O}} f(x)\ZD x \right |<\ZEP\,.
 \end{equation}
 \end{Theorem}
 
 Of course, in order that this statement makes sense, we must know that the open set $\mathcal{O}$ satisfies Assumption~\ref{ASSUch2QuascontINSa}, i.e. that its characteristic function is quasicontinuous. This we prove first.

\subsection{\ZLA{Ch3:SEct:CarFunzOPENseta}The Characteristic Function of an Open Set}

We extend Theorem~\ref{Teo:Cha1bisFunCHarAPERTIquasCont} to functions of several variables. The proof uses the representation of open sets in terms of multirectangles, hence it is a bit more elaborated since we must relay on Theorem~\ref{TeoCH3STRUttAPERTI} while when $d=1$ we can use Theorem~\ref{teo:Ch2StrutturaAperti}.
 
\begin{Theorem}\ZLA{Teo:Ch3FunCHarAPERTIquasCont}
Let $\mathcal{O}$ be a nonempty bounded open set. Its characteristic function $\charfun_\mathcal{O}$ is quasicontinuous.
 \end{Theorem}
 \zProof 
Let $R$ be a bounded rectangle such that $\mathcal{O}\subseteq R$. 
 
 We represent $\mathcal{O}=\cup _{n\geq 1} R_n$ where $R_n$ are closed almost disjoint rectangles such that

 \[
 \Rpalla_n={\rm int}\, R_n \neq \emptyset\,.
 \]

We note that    $\partial R_n=\partial\Rpalla_n$ is a null set so that also $\cup _{n\geq 1}\partial R_n= \cup _{n\geq 1}\partial \Rpalla_n$ is a null set\footnote{we are not asserting that $\partial  [ \cup _{n\geq 1}\Rpalla_n ]$ is a null set. In fact, in general it is not, see Remark~\ref{RemaCH4InsOpeFrontPOSI}.}.

We consider the increasing sequence of  open sets 
 \[
 \mathcal{O}_N =\bigcup_{n=1}^N \Rpalla_n\,.
 \]
 The characteristic function   $\charfun_{\mathcal{O}_N}$ is bounded quasicontinuous since
 \[
 \charfun _{ \mathcal{O}_N}(x)=\sum _{n=1}^N \charfun _{\Rpalla_n}(x) 
 \]
 and the sequence $\{\charfun _{ \mathcal{O}_N}\}$ is increasing.
 We observe:
 \begin{itemize}
 \item when $x\in \mathcal{O}\setminus \left [\cup _{n\geq 1}\partial\Rpalla_n\right]$: in this case there exists $ N_x$ such that $x\in  \mathcal{O}_N$ for every $N>N_x$ and so  for $N>N_x$ we have $\charfun _{\mathcal{O}_N}(x)=1$. Hence
 \[
\lim_{N\to+\ZIN}\charfun _{\mathcal{O}_N}(x)=1=\charfun _{\mathcal{O} }(x)\,.
 \]
 \item when $x\notin\mathcal{O}$: in this case $x\notin{\mathcal{O}_N}$ for every $N$ and 
 $\charfun_{\mathcal{O}_N}(x)=0$. So, also in this case we have
 \[
 \lim_{N\to+\ZIN}\charfun _{\mathcal{O}_N}(x)=0=\charfun _{\mathcal{O} }(x)\,.
 \]
 \item nothing we assert when $x$ belongs to the null set $ \cup _{n\geq 1}\partial \Rpalla_n$.
 \end{itemize}
 It follows that $\charfun _{\mathcal{O} }$ is a.e. limit of a sequence of quasicontinuous functions and it is quasicontinuous from Egorov-Severini Theorem.\zdia
 
Theorem~\ref{Teo:Ch3FunCHarAPERTIquasCont} has the following consequence:
\begin{enumerate}
\item if $f$ is summable   and if $\mathcal{O}\subseteq R$ is an open set then the integral of $f$ on $\mathcal{O}$ exists
\begin{equation}\ZLA{eq:Ch3AbsConTdefINTEopenSET}
\int _{\mathcal{O}} f(x)\ZD x=\int_R f(x)\charfun _{\mathcal{O}}(x)\ZD x\,.
\end{equation}

 \item
we can associate two numbers to the open set $\mathcal{O}$:
 
\[
\begin{tabular}{ccccc}
$ \zl(\mathcal{O})$
&
$=$
&
 $\overbrace{\sum _{n\geq 1} \zl(R_n)}^{\mbox{\small \begin{tabular}{c}  $\sum _{n\geq 1} \int_R \charfun _{R_n}(x)\ZD x$\\ $\rotatebox{90}{=}$\end{tabular}}}$
 &
 $=$
 &
$\overbrace{ \sum _{n\geq 1} \zl(\Rpalla_n)}^{
\mbox{\small\begin{tabular}{c}  $\sum _{n\geq 1} \int_R \charfun _{\Rpalla_n}(x)\ZD x$\\ $\rotatebox{90}{=}$\end{tabular}
}
}$    
\\ 
&&&&
\\
$\int_R\charfun_{\mathcal{O}}(x)\ZD x$
&
$=$
&
$ \int_R\left [\sum_{n\geq 1} \charfun_{R_n}(x)\right]\ZD x$
&
$=$
&
$ \int _R\left [\sum_{n\geq 1} \charfun_{{\Rpalla_n} }(x)\right]\ZD x$\,.
 \end{tabular}
\]
 since,  we repeat, a.e. $x\in R$ we have\footnote{note that $\{R_n\}$ is not a finite sequence since a finite union of closed sets is not open, unless it is the entire space.}
 \[
\charfun_{\mathcal{O}}(x)=\lim _{N\to+\ZIN} \charfun _{\mathcal{O}_N} (x)=\sum_{n=1}^{+\ZIN} \charfun_{\Rpalla_n}(x)= \sum_{n=1}^{+\ZIN} \charfun_{R_n}(x)\,. 
 \]
 \end{enumerate}
 We  prove:
 
 \begin{Theorem}\ZLA{Ch3:TeoPreAbsContPriLiMint}
We have:
\begin{equation}\ZLA{eq:ch3PrimadelTEO127}
\zl(\mathcal{O})= 
  \int_R\charfun_{\mathcal{O}}(x)\ZD x
\,.
\end{equation}
 \end{Theorem}
\zProof
 We use the fact that $\charfun_{\mathcal{O}}$ is quasicontinuous and the  equalities
\begin{align}
\ZLA{eq:ch3PSeCodelTEO127}
&
\charfun_{\mathcal{O}}(x)=\lim _{N\to+\ZIN} \charfun _{\mathcal{O}_N} (x)=\sum_{n=1}^{+\ZIN} \charfun_{\Rpalla_n}(x)\qquad {\rm a.e.}\ x\in\mathcal{O}\,,
\\
&
\ZLA{eq:ch3TeRzTEO127}
 \int_R\charfun_{\mathcal{O}}(x)\ZD x
=  \int _R\left [\sum_{n= 1}^{+\ZIN} \charfun_{{\Rpalla_n} }(x)\right]\ZD x\,.
\end{align} 
Equality~(\ref{eq:ch3PSeCodelTEO127}) does not hold on the null set $\cup _{n=1}^{+\ZIN}\partial \Rpalla$ and 
 \[
 \charfun _{\mathcal{O}_N} (x)\leq \charfun_{\mathcal{O}}(x)\quad {\rm a.e.}\ x\in\mathcal{O}
 \]
 so that
 \[
 \lim _{N\to+\ZIN} \sum _{n=1}^N
 \int _R  \charfun _{\Rpalla_n} (x)\ZD x\leq \int_R \charfun_{\mathcal{O}}(x)\ZD x\,.
 \]

 \ZREF{
  We prove that equality holds: we prove that  if   for every $N $  we have
\begin{multline} \ZLA{eq:ch3PSECONDAdelTEO127}
0\leq  \zaa\leq  \int_R \charfun_{\mathcal{O}}(x)\ZD x- \sum _{n=1}^N
 \int _R  \charfun _{\Rpalla_n} (x)\ZD x\\=
 \int_R \left [\charfun_{\mathcal{O}}(x)\ZD x-  \sum _{n=1}^N  \charfun _{\Rpalla_n} (x)\right] \ZD x 
 =\int _R\left [\sum _{n=N+1}^{+\ZIN}  \charfun _{\Rpalla_n} (x)\right] \ZD x 
 \end{multline}
 then it must be $\zaa=0$.} 
 
 Note that the integrals in~(\ref{eq:ch3PSECONDAdelTEO127}) are Lebesgue integrals, i.e. limits of Riemann integrals of associated continuous functions.
 
 For every $n$ we construct an associated continuous function of order $1/n$
 of the integrand as follows:
 
\begin{enumerate}
\item we choose an associated multirectangle $\Delta_{1;\nu}$ of order $1/2\nu$ of the function $\charfun_{\mathcal{O}}$ and we denote \[
(\charfun_{\mathcal{O}})_\nu
\]
 an associated continuous function    of order $1/2\nu$.

 By definition the associated continuous function is a Tietze extension  and we know that it is possible to   choose an associated continuous function which satisfies the monotonicity assumptions of  Theorem~\ref{Ch2TeoExteMonotSEQULebe}.
 
 The difference $(\charfun_{\mathcal{O}})_\nu-\charfun_{\mathcal{O}}$ is nonzero on $\mathcal{I}_{\Delta_1} $ where $L(\Delta_1)<1/2n$.
\item  
we choose an associated multirectangle $\Delta_{2;\nu}$ of order $1/2\nu$
associated to $\sum _{n=1}^N  \charfun _{\Rpalla_n}$ and
we denote
\[
\left (\sum _{n=1}^N  \charfun _{\Rpalla_n}\right )_\nu
\]
 an associated continuous function   of order $1/2\nu$.

 Also in this case we choose   an associated continuous function which satisfies the monotonicity assumptions of  Theorem~\ref{Ch2TeoExteMonotSEQULebe}.
 
The difference 
\[
\left (\sum _{n=1}^N  \charfun _{\Rpalla_n}\right )_\nu-\sum _{n=1}^N  \charfun _{\Rpalla_n}
\]
 is nonzero on $\mathcal{I}_{\Delta_{2,\nu}} $ where $L(\Delta_{2,\nu})<1/2\nu$.

\end{enumerate} 
 
The function
\[
F_{N;\nu}(x)=(\charfun_{\mathcal{O}})_\nu-\left (\sum _{n=1}^N  \charfun _{\Rpalla_n}\right )_\nu
\]
is an associated continuous function of order $1/\nu$ to
\[
\charfun_{\mathcal{O}}-\sum _{n=1}^{N}  \charfun _{\Rpalla_n}
=
\sum _{n=N+1}^{+\ZIN}  \charfun _{\Rpalla_n}\,.
\]
This associated function takes values in $[0,1]$ thanks to the fact that both the associated functions we chose satisfy the monotonicity assumption and because
\[
\charfun_{\mathcal{O}}(x)\geq  \sum _{n=N+1}^{+\ZIN}  \charfun _{\Rpalla_n}(x)\,.
\]

By definition, the Lebesgue integral in~(\ref{eq:ch3PSECONDAdelTEO127})
is
\[
\lim _{\nu\to+\ZIN} \underbrace{\int_R F_{N;\nu}(x)\ZD x}_{\tiny \begin{tabular}{l}
    Riemann integral  \end{tabular}}\,.
\]
We investigate where  $F_{N;\nu}$ can  possibly be different non zero: this is where the integrand in~(\ref{eq:ch3PSECONDAdelTEO127}) is non zero and on $\mathcal{I}_\Delta$ where $\Delta_\nu=\Delta_{1,\nu}\cup\Delta_{2,\nu}$.

\[
\{x\,:\ F_{N;\nu}(x)\neq 0\}\subseteq \left [\bigcup _{N+1}^{+\ZIN}\Rpalla_n\right ]\cup \mathcal{I}_{\Delta_\nu}
\]
and  
\[
\left [\bigcup _{N+1}^{+\ZIN}\Rpalla_n\right ]\cup \mathcal{I}_{\Delta_\nu}=\mathcal{I}_{\hat \Delta}
\]
where $\hat\Delta$ is a multirectangle such that
\[
L(\hat\Delta)\leq \dfrac{1}{\nu}+\sum _{n=N+1}^{+\ZIN} L(\Rpalla_n)\,.
\]
We use Lemma~\ref{Lemma:Ch2LemmaPreliRiemaaPIUvar} and we see that
\[
0\leq \underbrace{\int_R F_{N;n} (x)\ZD x}_{\tiny \begin{tabular}{l}
    Riemann integral  \end{tabular}}\leq \dfrac{1}{\nu}+\sum _{n=N+1}^{+\ZIN} L(\Rpalla_n)\quad \mbox{since $0\leq F_{N;n} (x)\leq 1$}\,.
\]
The limit for $\nu\to +\ZIN$ gives

\[
0\leq\zaa\leq \underbrace{\int _R\left [\sum _{n=N+1}^{+\ZIN}  \charfun _{\Rpalla_n} (x)\right] \ZD x
}_{\tiny \begin{tabular}{l}
    Lebesgue integral  \end{tabular}}
 =\lim _{n\to+\ZIN} 
\underbrace{ \int_R F_{N;n}(x)\ZD x
}_{\tiny \begin{tabular}{l}
    Riemann integral  \end{tabular}} 
 \leq \sum _{n=N+1}^{+\ZIN} L(\Rpalla_n)\,.
\]
This inequality holds for every $N$ and the limit for $N\to+\ZIN$ gives $\zaa=0$, as wanted.\zdia

\begin{Remark}{\rm
The statement of Theorem~\ref{Ch3:TeoPreAbsContPriLiMint} can be written
as
\[
\lim _{N\to+\ZIN}\left [\sum _{n =1}^N \int_R \charfun _{\Rpalla_n}(x)\ZD x\right ]
= \int _R\left [\lim _{N\to +\ZIN}\sum_{n=1 }^N \charfun_{{\Rpalla_n} }(x)\right]\ZD x
\]
and it is a first instance of the exchange of limits and integrals, the main goal of this chapter.\zdia
}\end{Remark}

\subsection{Absolute Continuity: the Proof of Theorem~\ref{Teo:CH3ABSOluContiPERiSoliAperti}}
It is sufficient to prove the theorem when $f\geq 0$.

First we consider the case that $f$ is a bounded quasicontinuous function, $0\leq f(x)\leq M$ on $R$. We combine~(\ref{eq:Ch3AbsConTdefINTEopenSET}) and~(\ref{eq:ch3PrimadelTEO127}) and we find
\begin{equation}\ZLA{eq:Ch3SUBSabsoluConti0}
0\leq \int _{\mathcal{O}} f(x)\ZD x\leq M\int _{\mathcal{O}} 1\ZD x=\int _R\charfun _{\mathcal{O}}
(x)\ZD x=M\zl(\mathcal{O})\,.
\end{equation}
So, absolute continuity holds when the integrand is bounded.

  We  consider the general case of summable functions on $\zzr^d$. The procedure is similar to that used in the proof of Theorem~\ref{Teo:Cha1bisABSOluContiPERiSoliAperti} in Chap.~\ref{Cha1bis:INTElebeUNAvariaBIS}: first we fix $R$ and $N$ such that
 \[
\int _{\zzr^d} f(x)\ZD x-\int_{\zzr^d}f_{+;\,(R,N)}(x)\ZD x <\ZEP/2\,. 
 \]
 Then we use absolute continuity which holds when the integrand is the bounded function $f_{+;\,(R,N)}$: we fix $\ZDE>0$ such that
 \[
\zl(\mathcal{O})<\ \ZDE\ \implies 0\leq \int _{\mathcal{O}}f_{+;\,(R,N)}(x)\ZD x<\ZEP/2\,.
 \]
Then we have
 
\begin{multline*}
 \int_{\mathcal{O}}  f(x) \ZD x\leq \underbrace{\int_{\mathcal{O}} [f(x)-f_{+;\,(R,N)}(x)]\ZD x}_{
 \tiny \leq \int_{\tiny \zzr^d}[f(x)-f_{+;\,(R,N)}(x)]{\ZD x\, \leq\, \ZEP/2}
 }\\
+
 \int_{\mathcal{O}}  f_{+;\,(R,N)}(x) \ZD x<\ZEP\,.\zdiaform
\end{multline*}
 \section{\ZLA{ch3SecDimoConve}The Limit of  Sequences and the Lebesgue Integral} 
 In this section we examine the theorems concerning limits and integrals. The statements and the proofs  are the same as those in Chap.~\ref{Cha1bis:INTElebeUNAvariaBIS}, provided that we use  the 
 Egorov-Severini Theorem~\ref{P2CH3TeoEgoSEVER} and absolute continuity of the integral, i.e. Theorem~\ref{Teo:CH3ABSOluContiPERiSoliAperti}. So, we confine ourselves to state the results.
 
 The first theorem concerns bounded sequences on bounded sets:

 \begin{Theorem}\ZLA{teo:ch3:LebeCAROPartIntervLIMIT}
 Let $\{f_n\}$ be a bounded sequence of quasicontinuous functions on a bounded closed rectangle $R\subseteq\zzr^d$. If 
 \[
\limn f_n(x)=f(x)\qquad {\rm a.e. \ on}\  R 
 \]
 then we have also
 \[
\limn \int_R f_n(x)\ZD x=\int_R f(x)\ZD x\,. 
 \]
 \end{Theorem}
 
  As a second step we examine sequences of nonegative functions and, as in Sect.~\ref{Cha1bisSecLIMINTposiSEQUE}, we derive  
  Fatou Lemma  and Beppo Levi Theorem: 
\begin{Lemma}[\sc Fatou Lemma]\index{Fatou Lemma}\index{Lemma of Fatou}
If $\{f_n\}$ is a sequence of nonnegatrive functions which are integrable on a set $A$ and if $f_n(x)\to f(x)$ a.e. on $A$ then we have
\begin{equation}\ZLA{eq:DiseqFATOUch3}
\int _A f(x)\ZD x\leq  \liminf_{n\to+\ZIN} \int_{A} f_{n}(x)\ZD x \,.
\end{equation}
\end{Lemma}

\begin{Theorem}[{\sc Beppo Levi} or {\sc monotone convergence}]\ZLA{TeoCH3Beppo}\index{Beppo Levi Theorem}\index{Theorem!Beppo Levi} 
 \index{monotone convergence theorem}\index{Theorem!monotone convergence}
Let $\{f_n\}$ be a sequence of integrable nonnegative functions on $A$ and let
\[
0\leq f_n(x)\leq f_{n+1}(x)\quad {\rm a.e.}\ x\in A\\,\qquad \limn f_n(x)=f(x)\ {\rm a.e}\  x\in A\,.
\]
Then we have
\begin{equation}\ZLA{eq:TeoCH3Beppo}
\limn \int_A f_n(x)\ZD x=\int_A f(x)\ZD x\,.
\end{equation}
\end{Theorem}

We recall Remark~\ref{Chap1BisRemaLemmaBeppoeFatou} of Sect.~\ref{Cha1bisSecLIMINTposiSEQUE}: the statement of Beppo Levi Theorem does not hold for decreasing sequences, and the inequality in Fatou Lemma in general is strict.

Finally we consider the general case and we state:

\begin{Theorem}[{\sc Lebesgue} or {\sc    dominated convergence}]\index{Theorem!Lebesgue!}\index{Theorem!dominated convergence}\index{dominated convergence}\ZLA{teo:CH3TeoLebeRlimitato}
Let $\{f_n \}$ be a sequence of summable functions a.e. defined on    $A\subseteq \zzr^d$ and let $f_n\to f  $ a.e. on $A$. If there exists a summable nonnegative function $g $ such that
\[
|f_n ( x)|\leq g(x)\qquad a.e.\ x\in A
\]
then $f(x)$ is summable and
\begin{equation}
\ZLA{eq:CH3LebeDOMIconve0}
\lim \int_A f_n ( x)\ZD x=\int_A f(x)\ZD x\,.
\end{equation}
\end{Theorem}

\section{\ZLA{P2appeCH3EGOSEVE}\textbf{Appendix:}  Egorov-Severini: Preliminaries   in Several Variables}
 
In this appendix we prove Theorem~\ref{P2CH1TeoEgoSEVERante}. The proof 
  follows closely   the proof we have seen in the case of functions of one variable (seen  in the Appendix~\ref{AppeCha1bisEgSEV}). So, we repeat the statement and we give the details of significant  proofs or when there are significant differences.

\subsection{\ZLA{P2S1:appeCH3EGOSEVE}Preliminary Observations}
We need the following quite obvious  observations on the measure of multirectangle sets. These observations are similar to those in Sect.~\ref{P1Ch2AppeSEctSulMultiint}, just a bit more elaborate since the multirectangle we must consider when $d>1$ are not disjoint but almost disjoint.

We recall that a set $A$ is a {\sc    multirectangle set}\index{set!multirectangle set}\index{multirectangle!set} when
\[
A=\mathcal{I}_{\Delta}\qquad \mbox{and $\Delta$ is almost disjoint }\,.
 \]
 \emph{Although  the results below holds without this condition, in this appendix
we explicitly assume that $\Delta $ is composed of rectangles of positive measure since this is the case we shall need.}

We recall the definition
 \[
\zl(A)=L(\Delta)\,. 
 \]

\subparagraph{Observation~1.} 
 The {\sc monotonicity   of the measure}\index{measure!monotonicity} i.e. Lemma~\ref{CH3:lemmaMONOTONmisura} statement~\ref{I1CH3:lemmaMONOTONmisura}:
if  $A$ and $B$ are multirectangle sets and if $A\subseteq B$ then $\zl(A)\leq \zl(B)$.
 
 \subparagraph{Observation~2.} The {\sc additivity of the measure}\index{measure!additivity}
 i.e. Lemma~\ref{CH3:lemmaMONOTONmisura} statement~\ref{I2CH3:lemmaMONOTONmisura}: 
  let $\Delta_1=\{R_{1,n}\}$ and $\Delta_2=\{R_{2,n}\}$ be almost disjoint multirectangles. We have
 \begin{align*}
 R_{1,n}\cap R_{2,j}=\emptyset \ \ \forall n\,,\ j\quad \implies\quad 
 \zl(\mathcal{I}_{{\Delta_1\cup\Delta_2}})=L(\Delta_1\cup\Delta_2)&=L(\Delta_1)+L(\Delta_2)\\
&=\zl(\mathcal{I}_{\Delta_1})+\zl(\mathcal{I}_{\Delta_2})
\,.
 \end{align*} 
 An equivalent formulation is as follows: if $A$ and $B$ are multirectangle sets and if $A\cap B=\emptyset$ then $\zl(A\cup B)=\zl(A)+\zl(B)$.
\subparagraph{Observation~3.}  
 
Let $S$ and $\hat S$ be rectangles. If 
nonempty, both $S\cap \hat S$ and $S\setminus  S$ are finite unions of rectangles \emph{which are pairwise almost disjoint.}
If $\Delta=\{R_n\}_{n\geq 1}$ is an almost disjoint multirectangle we can represent any $S\cap R_n$ which is nonvoid as  just described:
\[
S\cap R_n=\bigcup _{\nu=1}^{N(n)} R_{n,\nu}\,.
\]
\emph{In order to simplify the notations, if $
 S\cap R_n=\emptyset$ we put $N(n)=0$ and we intend that 
\[
\bigcup _{\nu=1}^{0} R_{n,\nu}=\emptyset\,,\qquad L\left (\{  R_{n,\nu}\}_{\nu=1}^0\right )=0\,.
\]
}
This way we get an almost disjoint  multirectangle
    which we shortly denote $S\cap \Delta$:
\[
S\cap \Delta=\{S\cap R_n\}=\{R_{n,\nu}\}\quad \mbox{(indexed by $(n,\nu)$)} 
\]

and
\[
L(S\cap \Delta)=\sum _{n>0}\sum _{\nu=1}^{N(n)}R_{n,\nu}\leq \min\{L(S)\,,\ L(\Delta)\}\,.
\]

 \subparagraph{Observation~4.}  
 Let 
$\Delta=\{R_n\}$ and $\tilde\Delta=\{S_n\}$ be two almost disjoint multirectangles and let
\[
R_n\cap S_k=\bigcup_{\nu=1}^{N(n,k)}  R_{n,k;\nu }\quad   \mbox{(finite union of almost disjoint  rectangles)}\,.
\] 
\emph{As in the {\bf Observation~3,} if  $ R_n\cap S_k=\emptyset$ we put formally $N(n,k)=0$ and we intend that the sequence indexed by $\nu $ from $1$ to $0$ does not exist.}

 The sequence (indexed by the three indices $n$, $k$  and $\nu $)
 $\{R_{n,k;\nu }\} $ is an almost disjoint multirectangle and
 \[
L\left ( \{R_{n,k;\nu }\}\right )=\sum _{n,k}\sum _{\nu =1}^{N(n,k) }L(R_{n,k;\nu }) \leq\min\{L(\Delta)\,,\ L(\tilde\Delta)\}\,.
 \]

 If it happens that
\[
\mathcal{I}_{\Delta}\subseteq \mathcal{I}_{\tilde\Delta}
\]
then
\[\mathcal{I}_\Delta
=\bigcup_{n\,,\ k}\left[ \bigcup_{\nu =1}^{N(n,k)}  R_{n,k;\nu } \right]=\bigcup _{k\geq 1} \left [ \bigcup _{n\geq 1}  \bigcup _{\nu =1}^{N(n,k)}R_{n,k;\nu } \right ]
=\bigcup _{k\geq 1} S_k\cap\mathcal{I}_\Delta =\mathcal{I}_{\{S_k\cap\Delta\}_{k\geq 1}}
\]
and
\begin{equation}\ZLA{P2eq:CH4FormSommaObser4}
L(\Delta)=\sum _{n,k}\sum _{\nu =1}^{N(n,k) }L(R_{n,k;\nu } ) =\sum _{k\geq 1}  L(S_k\cap\Delta)\,.
\end{equation}
  \subparagraph{Observation~5.} 
If it happens that $\mathcal{I}_\Delta\subseteq R$ and $S\subseteq R$ ($R$ and $S$ both rectangles) then we have\footnote{$R\setminus S=\{R_1\,,\dots\,, R_k\}$
and the rectangles are almost disjoint. Then $\Delta_2=\cup _{i=1}^k R_i\cap\Delta$.}
\[
\Delta=\underbrace{  [S \cap \Delta]
}_{\Delta_1}\bigcup \underbrace{[(R\setminus S)\cap \Delta]}_{\Delta_2}\,.
\] 
We represent $R\setminus S$ as the almost disjoint union of finitely many rectangles.

When $\Delta$ is almost disjoint 
both $\Delta_1$ and $\Delta_2$ are almost disjoint   and no rectangle of $\Delta_1$ intersects a rectangle of $\Delta_2$. The additivity property in {\bf Observation~2} gives
\begin{equation}\ZLA{P2eq:3APPeSec1Punto5}
\zl(\mathcal{ I}_\Delta)=
L(\Delta)=L(\Delta_1)+L(\Delta_2)
=\zl(\mathcal {I}_{\Delta_1})+\zl(\mathcal{I}_{\Delta_2})
\,.
\end{equation}

 The previous considerations permits to extend the arguments in the Appendix~\ref{AppeCha1bisEgSEV} from the case $d=1$ to the case $d>1$.
 
 We state\footnote{we leave to the reader the reformulation of the lemma in terms of open sets, as   explicitly done in the Lemma~\ref{LemmaiAppeCh3SUOmultirettANNID} of the Appendix~\ref{AppeCha1bisEgSEV}.}:

\begin{Lemma}\ZLA{P2LemmaiAppeCh3SUOmultirettANNID}
Let  
  $\{\Delta_n\}$ be a sequence of     multirectangles in $\zzr^d$. We assume:
\begin{enumerate}
\item\ZLA{I1LemmaiAppeCh3P2} the existence of a bounded rectangle $R$ such that $\mathcal{I}_{\Delta_n}\subseteq R$ for every $n$;
\item\ZLA{P2I0LemmaiAppeCh3}    every multirectangle  $\Delta_n$ is almost disjoint;
 
\item\ZLA{P2I2LemmaiAppeCh3} 
$\mathcal{I}_{\Delta_{n+1}}\subseteq \mathcal{I}_{\Delta_n}$ for every $n$;
 \item\ZLA{P2I3LemmaiAppeCh3}   there exists     $l >0$  such that $L(\Delta_n)>l $ for every $n$.
  
\end{enumerate}
Under these conditions, there exists $x_0\in\zzr^d$ which is an interior points of every~$\mathcal{I}_{\Delta_n}$.

\end{Lemma}

The proof is world by world equal to that of Lemma~\ref{LemmaiAppeCh3SUOmultirettANNID}, provided that ``interval'' is replaced with ``rectangle'' and ``disjoint'' with ``almost disjoint''. So, the proof is not repeated.

\begin{Remark}
{\rm
We recall that the boundedness assumption in the statement~\ref{I1LemmaiAppeCh3P2}   of the Lemma  cannot be removed, see Remark~\ref{Rema:AppeP1Ch2RemaRSullaNeceLimit}.\zdia
}
\end{Remark}

We shall also use Lemma~\ref{LemmaAppe3ProprVnmWm}. As noted in Remark~\ref{LemmaAppe3ProprVnmWm}, this lemma, proved in the Appendix~\ref{AppeCha1bisEgSEV}, holds for functions of any number of variables. 

We repeat  the definitions of the functions $v_{n,m}$ and $w_m$ and the statement of the lemma.

\begin{Lemma}\ZLA{P2LemmaAppe3ProprVnmWm}
Let $\{f_n\}$ be a sequence of functions defined on a rectangle $R$ and let $v_{n,m}(x)$, $w_m(x)$ be the functions 
\begin{equation}\ZLA{P2Ch3AppeDifivNM}
\begin{array}{l}\displaystyle
 v_{n,m} (x)=\max\{
  |f_{m+r}(x)-f_{m+s}(x)| \quad 1\leq r<n\,,\qquad 1\leq s<n\}\,,\\[2mm]
\displaystyle
  w_m(x)=\sup\{v_{n,m}(x)\quad n\geq 1\}\leq +\ZIN\,.
  \end{array}
\end{equation}
 We have:

\begin{enumerate}
\item\ZLA{I1P2LemmaAppe3ProprVnmWm} monotonicity properties:
\begin{enumerate}\ZLA{I1VdipNP2LemmaAppe3ProprVnmWm}
\item
for every $x\in R$ and every $m$, the sequence $n\mapsto  v_{n,m}(x)$ is increasing.

\item \ZLA{I1WdipMP2LemmaAppe3ProprVnmWm}  the sequence $m\mapsto w_m(x)$ is decreasing.
\end{enumerate}
(monotonicity of the two sequences needs not be strict).
 
\item \ZLA{I2P2LemmaAppe3ProprVnmWm} the convergence of the sequence $n\mapsto  f_n(x) $ for a fixed value of $x$:

\begin{enumerate}
\item\ZLA{I2DIREP2LemmaAppe3ProprVnmWm} 
let   $n\mapsto f_n(x)$ converge. Then the sequence $n\mapsto v_{n,m}(x)$ is bounded so that  
\begin{equation}\ZLA{EqAI2DIREP2LemmaAppe3ProprVnmWm} 
w_m(x)=\sup _{n>0}v_{n,m}(x)=\lim _{n\to+\ZIN}v_{n,m}(x)   \in\zzr
\end{equation}
and we have also
\begin{equation}\ZLA{EqABI2DIREP2LemmaAppe3ProprVnmWm}
\lim_{m\to+\ZIN  } w_m(x)=0\,.
\end{equation}

 \item \ZLA{I2VICEVP2LemmaAppe3ProprVnmWm} 
let $\lim _{m\to +\ZIN} w_m(x)=0$. Then the sequence $n \mapsto f_n(x)$ converges.
 \item\ZLA{I3VICEVP2LemmaAppe3ProprVnmWm}let $S$ be a subset of $R$. The sequence $\{f_n\}$ converges uniformly on $S$ if and only if the sequence $\{w_m\}$ converges to $0$ uniformly on $S$.
\end{enumerate}

\item\ZLA{I3P2LemmaAppe3ProprVnmWm} 
 let  $S\subseteq R$. Let us assume  that each     function $\left (f_n\right )_{|_{S}}$ be continuous on $S$ and let $\ZEP>0$.  
 The set
\[
A_{m,\ZEP}=\{x\in S\,, w_m(x)>\ZEP\}\,.
\]
  is relatively open in $S$.
\end{enumerate}
\end{Lemma}
 
\ZREF{  
We stress again the fact that the statements~\ref{I2DIREP2LemmaAppe3ProprVnmWm},~\ref{I2VICEVP2LemmaAppe3ProprVnmWm}  and~\ref{I3VICEVP2LemmaAppe3ProprVnmWm} of Lemma~\ref{P2LemmaAppe3ProprVnmWm} recast pointwise convergence of the sequence $\{f_n(x)\}$ in terms of the convergence to $0$ of the sequence $\{w_m(x)\}$ and uniform convergence of $\{f_n\}$ in terms of uniform convergence to zero of $\{w_m\}$.

}

 \subsection{The Proof of Theorem~\ref{P2CH1TeoEgoSEVERante}} 
  As in Section~\ref{P1Ch2AppeSEctSulMultiint} of the Appendix~\ref{AppeCha1bisEgSEV}, first we state and prove a lemma which is the core of Egorov-Severini Theorem.  The proof is similar to that of Lemma~\ref{LemmaiAppeCh3PRELIsucceFUNcont} of the Appendix~\ref{AppeCha1bisEgSEV} but we find convenient to repeat this proof.

\begin{Lemma}\ZLA{P2LemmaiAppeCh3PRELIsucceFUNcont}
  Let $R$ be a bounded closed rectangle and let $\{f_n\}$ be a sequence of \emph{continuous} functions   defined on $R$. We assume that the sequence $\{ f_n(x)\}$ converges for every $x\in R$. 
 
We prove that for every pair of positive numbers $\zg >0$ and $\eta >0$ there exists
\begin{description}
\item[\bf 1:] an almost disjoint c-multirectangle $\Delta_{\zg ,\eta } $ which is a
multirectangle of an open set $A_{\zg,\eta}=\mathcal{I}_{\Delta _{\zg,\eta}}$ and such that
\[
L(\Delta_{\zg ,\eta })=\zl(A_{\zg,\eta})<\zg\,.
\]
\item[\bf 2:] a number $M_{\zg ,\eta }$  such that 
  
\[
\left\{\begin{array}{l}
 x\in R\setminus\mathcal{I}_{ \Delta_{\zg ,\eta }}\\
 m> M_{\zg ,\eta } \\ 
 r>0\,,\ s>0
 \end{array}\right.\quad \ \implies \ |f_{m+r}(x)-f_{m+s}(x)|<\eta  \,.
\]
\end{description}
\end{Lemma}
 \zProof    
We recast the thesis of the lemma   in terms of the functions  $w_m(x)$ defined in~(\ref{P2Ch3AppeDifivNM}) as follows:
  \emph{ for every $\zg >0$ and $\eta >0$ there exists an almost disjoint c-multirectangle $\Delta_{\zg ,\eta } $ such that
\[
L(\Delta_{\zg ,\eta })<\zg \quad\mbox{and}\quad
\left\{\begin{array}{l}
\displaystyle x\in R\setminus\mathcal{I}_{\Delta_{\zg ,\eta }}\ \implies w_m(x)<\eta \\
\displaystyle
m>M_{\zg,\eta}
\end{array}\right.
\] 
 }
This  we prove now.

Let $\eta >0$ and
\[
A_{m,\eta }=\{x\in {\rm int}\, R\,, w_m(x)>\eta/2 \}\,.
\]
Statement~\ref{I3P2LemmaAppe3ProprVnmWm}  of Lemma~\ref{P2LemmaAppe3ProprVnmWm} shows that the set $A_{m,\eta }$ is   open.    So, 
  there exists an almost disjoint c-multirectangle $\Delta_{m,\eta }$ such that
$A_{m,\eta }=\mathcal{I}_{\Delta_{m,\eta }}$. Moreover, the component rectangles of $\Delta_{m,\eta }$ have nonempty interior (see Theorem~\ref{TeoCH3STRUttAPERTI}).

If $x\in R\setminus \mathcal{I}_{\Delta_{m,\eta }}= 
R\setminus A_{  m,\eta } 
$ we have
\[
  w_m(x)\leq \eta/2  \,.
\]
We   note
\begin{equation}\ZLA{eq2P2LemmaiAppeCh3PRELIsucceFUNcont}
\limm L( \Delta_{m,\eta } )=0\,.
\end{equation}
In fact, the sequence $m\mapsto  w_m(x)$ is decreasing so that
\[
\mathcal{I}_{\Delta _{{m+1},\eta }}=A_{m+1,\eta }\subseteq A_{m,\eta }\subseteq \mathcal{I}_{\Delta_{m,\eta }}\,.
\]
It follows that $\{ L(\Delta_{m,\eta })\}=\{\zl(A _{m,\eta })\}$ is decreasing too
and $\lim _{m\to+\ZIN}L(\Delta_{m,\eta })$ exists.

If the limit is positive then there exists $l >0$ such that
\[
L(\Delta _{m,\eta })>l>0 \qquad \forall m
\]
and, from Lemma~\ref{P2LemmaiAppeCh3SUOmultirettANNID} there exists $x_0\in  \mathcal{I}_{\Delta_m,\eta }$ for every $m$. So, for every $m$ we have
$w_m(x_0)>\eta $.
Statement~\ref{I2DIREP2LemmaAppe3ProprVnmWm} of Lemma~\ref{P2LemmaAppe3ProprVnmWm} shows that the sequence $\{ f_n(x_0)\}$ is not convergent, in contrast with the assumption. So it must be
\[
\limm L(\Delta_{m,\eta })=\limm\zl(A_{m,\eta })=0\,.
\]
It follows that there exists $M_{\zg ,\eta } $ such that when $m>M_{\zg ,\eta }$
then we have
\[
\begin{array}{l}
\displaystyle
L(\Delta _{m,\eta })<\zg \,,\\
x\in R\setminus \mathcal{I}_{\Delta(m,\eta )}\  \implies\ w_m(x)\leq \eta/2<\eta \,.\zdiaform
\end{array}
\]

It is worth repeating the observation in Remark~\ref{ChP1Ch2RenaImpoNOunic}:

\begin{Remark}[Important observation]\ZLA{ChP2Ch2RenaImpoNOunic}{\rm
The sets $A_{m,\eta}$ can be chosen with different laws, for example by replacing $\eta/2$ in~(\ref{ChP1Ch2eqRenaImpoNOunic}) with $\eta/3$, or by replacing  $A_{m,\eta}$ in~(\ref{ChP1Ch2eqRenaImpoNOunic}) with larger open sets, provided that~(\ref{ChP1Ch2eqBISRenaImpoNOunic}) holds. So, the number $M_{\zg,\eta }$  does depend also the chosen set  $A_{m,\eta}$.\zdia

}
\end{Remark}

 Now we examine Theorem~\ref{P2CH1TeoEgoSEVERante}. As in Sect.~\ref{AppeP2CH2FiNalStep} of the Appendix~\ref{AppeCha1bisEgSEV},
in the next box we report the statement of Lemma~\ref{P2LemmaiAppeCh3PRELIsucceFUNcont} and that of Theorem~\ref{P2CH1TeoEgoSEVERante} (recasted
 in terms of the functions $w_m$ defined 
  in~(\ref{P2Ch3AppeDifivNM})).

   {\footnotesize
  \begin{center}
  \ZREF{
	The functions $w_m$ are defined on a bounded closed interval $R$.
	\begin{center}
  \fbox{\parbox{2.3in}{ 
\begin{center}{\bf Under the assumptions of  Lemma~\ref{P2LemmaiAppeCh3PRELIsucceFUNcont} we proved}\end{center}

For every $\zg>0$ and $\eta>0$ there exist an \emph{open} set $A_{\zg,\eta}$ and a
number $M_{\zg,\eta}$ such that 
\[
\left\{
\begin{array}{l}
\zl(A_{\zg,\eta})<\zg\\
\left\{\begin{array}{l}
\mbox{if $m>M_{\zg,\eta}$ and $x\notin A_{\zg,\eta}$}\\
\mbox{then $w_m(x)<\eta$}\,.
\end{array}\right.
\end{array}\right.
\]
 }} 
  \fbox{\parbox{2.3in}{
\begin{center}{\bf We must prove}\end{center}
 
 \vskip -.1cm
  \[\begin{array}{l}
\mbox{$\forall\ZEP>0\ \exists\  \mathcal{O}_\ZEP$   such that} \\
\mbox{$\mathcal{O}_\ZEP$ is open and  $\zl(\mathcal{O}_\ZEP)<\ZEP$ and} \\
\mbox{$\forall\ZSI>0\ \exists M>0$ such that}\\
\left\{\begin{array}{l}
\mbox{if $x\in R\setminus \mathcal{O}_\ZEP$, $m>M$}\\
\mbox{then  $  \ w_m(x)<\ZSI$ }\,.
\end{array}\right.
\end{array} 
  \]
  The number $M$ depends
on the previously chosen and fixed set $\mathcal{O}_\ZEP$ and  
   on~$\ZSI$. 
   }} 
   \end{center}
    \emph{The important fact to be proved is that $\mathcal{O}_\ZEP$ does not depend on $\ZSI$.}
  }
 {~}
  \end{center}
 }
 \smallskip

  By looking at this table, we see that it is precisely equal to that in   Sect.~\ref{AppeP2CH2FiNalStep} of the Appendix~\ref{AppeCha1bisEgSEV} and so the proof is concluded with precisely the same steps as in the case of the functions of one variable.
  
  \begin{Remark}
 {\rm 
We repeat that the proof of Egorov-Severini Theorem uses the assumption that $R$ is bounded, hidden in the    use 
of  Lemma~\ref{P2LemmaiAppeCh3SUOmultirettANNID}.\zdia
 
 }
 \end{Remark}
 

 \chapter{\ZLA{ch4CHAPFUBINI}Reduction of Multiple Integrals}
The calculation of Riemann integrals of functiuons of several variables is simplified if it can be reduced to a chain of computations of integrals of functions of one variable. Whether this can be done depends on the properties of the domain of integration and, in the context of Riemann integral, it is quite difficult to derive general conditions under which reduction can be achieved. In this chapter we see that similar reduction formulas exist  also for the Lebesgue integral with the further bonus that in the case of Lebesgue integration general conditions can be given.

 \section{\ZLA{ch4SectFUBINI}Discussion and the Reduction Formulas}
 We defined the Lebesgue measure and integrals in $\zzr^d$ for every dimension $d$. There is an obvous relation among the measures in $\zzr^d$ for different values of $d$. For example if $R=[a,b]\times [c,d]\subseteq\zzr^2$ is a rectangle  then its measure in $\zzr^2$ is the product of the measures of the intervals $I=[a,b]$ and $J=[c,d]$ both computed in $\zzr^1$: $\zl(R)=\zl(I)\times \zl(J)$ and this equality has a corresponding formulation in terms of
 the (Riemann) integral:
 \[
\zl(R)=\int_{[a,b]\times[c,d]} 1\ZD x\,\ZD y= \int_a^b\left [\int_c^d 1\ZD y\right ]\ZD x \,.
 \]
This equality holds also if the constant function $1$ is replaced by any continuous function $f(x,y)$.
 
 In this section we study conditions  under which   the following {\sc reduction formula}\index{reduction formula}  holds for the Lebesgue integral:
 \begin{equation}\ZLA{eq:CH4RedUForMUPreli}
 \int_{\zzr^d} f(x,y) \ZD x\,\ZD y=\int_{\zzr^{d_1}}\left [\int_{\zzr^{d_2} } f (x,y) \ZD y\right ]\ZD x\,,\qquad d=d_1+d_2\,.
 \end{equation}

We note a difficulty with the equality~(\ref{eq:CH4RedUForMUPreli}): let $R=[-1,1]\times[-1,1]$ and let $f(x,y)=0$ when $y\neq 0$ while $f(x,0)$ is a 
on measurable function\footnote{the existence of such  kind of  functions will be seen in Remark~\ref{RemaFINAppeCH4InseNONLmisur} of the   Appendix~\ref{AppeCH4InseNONLmisur}.}. Then $f$ is a.e. equal zero on the rectangle $R$, hence   quasicontinuous,    since a segment is a null set, but the restrictions of $f$ to segments need not be quasicontinuous on the segments. So, the inner integral on the right side cannot be always computed while the integral on the left exists and it is equal zero.

   The results we are going to prove are stated  in Sect.~\ref{ch4SECnotatEstatementFUbini} while the proofs are in Sect.~\ref{sec:CH4SecDimoFunbini}. Preliminary results are in Sect.~\ref{Ch4PreliFubToneNuovo}.

 \subsection{\ZLA{ch4SECnotatEstatementFUbini}Notations and the Statements of the Theorems}
 We use the following notation:
\[\mbox{ 
   $\zl_n$ denotes the Lebesgue measure in $\zzr^n$}\,. 
   \]

 We work in dimension $n=d+1$,   $\zzr^n=\zzr^{d+1}$. We represent  $\zzr^{d+1}$ as the product space
 \[
\zzr^{d+1}=\zzr\times \zzr^d\,. 
 \]
 The elements of $\zzr^{d+1}$ are represented as $(x,\zb  y)$ where $\zb y =(y_1\,,\, \dots\,,\, y_d)\in\zzr^d$.
 
 An integral on $\zzr$, $\zzr^d$ or $\zzr^{d+1}$ is denoted\footnote{in Corollary~\ref{Coro:Ch4FubinProdRETTANG} we use   an obvious extension of  these notations.}
 \[
 \int \cdots\ZD x\,,\qquad \int\cdots \ZD\zb y\,,\qquad
\int\cdots  \ZD(x,\zb y)  \,.
 \]

  With these notations we can state the following Fubini and Tonelli Theorems.  
   \begin{subequations}
  \begin{Theorem}[{\sc Fubini}]\index{Theorem!Fubini}\index{Fubini Theorem}\ZLA{Theo:Ch4FUBINI}
  Let $(a,b)\times R\subseteq \zzr^{d+1}$ be a (possibly unbounded) rectangle and let $f$ be a \emph{summable} function on $(a,b)\times R$. We have:
\begin{enumerate}
\item  \ZLA{I0Theo:Ch4FUBINI}
 the function $\zb y\mapsto f(x,\zb y)$ is summable for a.e. $x\in (a,b)$.
\item  
  \ZLA{I1Theo:Ch4FUBINI}
  the function
  \[
x\mapsto \int _R f(x,\zb y)\ZD\zb y  
  \]
  is a.e. defined and   summable on $(a,b)$.
  
  \item\ZLA{I2Theo:Ch4FUBINI} the following reduction formula holds:
  
  \begin{equation}\ZLA{eq:LI2Theo:Ch4FUBINI}
\int_{(a,b)\times R} f(x,\zb y)\ZD(x,\zb y)  =\int_a^b \left [  \int_R f(x,\zb y)\ZD \zb y\right ]\ZD x\,.
  \end{equation}
  \end{enumerate}
  \end{Theorem}
  
  \medskip
  {~}
  \smallskip
  We note:
  \begin{itemize}
  \item The formulation of the Fubini Theorem looks restrictive since it concerns integration on rectangles. In fact, after the introduction of suitable notations   we can see that the formulation is general. This observation is in Sect.~\ref{sec:CH4SecDimoFunbini}, before the proof of the theorem.
  
  \item
 The rectangle $ (a,b)\times R$ is
 \[
  (a,b)\times \left (\prod_{k=1}^d(a_k,b_k)\right )\,. 
 \]
 We singled out the ``first  component'' $(a,b)$ in the statement of the theorem but permutations of the order of the   components  do not change neither the rectangle $(a,b)\times R$ nor  the   value of the integral. 
  \item
 Equality~(\ref{eq:LI2Theo:Ch4FUBINI}) can be iteratively applied to the inner integral   on the right side of~(\ref{eq:LI2Theo:Ch4FUBINI}).  
\end{itemize}
By taking these observations into account we deduce:
\[
\int_{(a,b)\times R} f(x,\zb y)\ZD(x,\zb y)=\int_a^b\left [      
\int _{a_1}^{b_1} \left [
\int _{a_2}^{b_2} 
\dots f(x,y_1,\dots,y_d)\dots
 \ZD y_2
\right ]\ZD y_1
\right ]\ZD x
\]
The order of the integrals can be arbitrary changed and by suitably collecting them we get the following equality. Let $(x,\zb y)=(\zb{x_1},\zb{x_2} )$,   $\zb{x_i}\in\zzr^{d_i}$  with $d_1+d_2=d+1$, and let $(a,b)\times R=R_1\times R_2$ be the corresponding decomposition. We have:
\begin{Corollary}\ZLA{Coro:Ch4FubinProdRETTANG}Let the assumptions of Theorem~\ref{Theo:Ch4FUBINI} hold. We have:
\begin{enumerate}
\item The functions $\zb {x_2}\mapsto f(\zb {x_1},\zb{x_2})$ and 
$\zb {x_1}\mapsto f(\zb {x_1},\zb{x_2})$ are summable respectively for a.e. $\zb{x_1}$ and for a.e. $\zb{x_2}$.
\item  
the functions 
\[
 \zb{x_1}\mapsto \int _{R_2}f(\zb{x_1},\zb{x_2})\ZD\zb{x_2}\,,\qquad 
 \zb{x_2}\mapsto \int _{R_1}f(\zb{x_1},\zb{x_2})\ZD\zb{x_1}
\]
 are a.e. defined and summable;
\item
The following equality holds:
\begin{multline}\ZLA{eq:LI2Theo:Ch4FUBINIblocchi}
\int_{R_1\times R_2} f(\zb{x_1},\zb{x_2})\ZD(\zb{x_1},\zb{x_2})=\int _{R_1}\left [\int _{R_2}f(\zb{x_1},\zb{x_2})\ZD\zb{x_2}\right ]\ZD\zb{x_1} \\
=\int _{R_2}\left [\int _{R_1}f(\zb{x_1},\zb{x_2})\ZD\zb{x_1}\right ]\ZD\zb{x_2}\,;
\end{multline}
 
\item
In particular we  have:
\begin{equation}\ZLA{eq:LI2Theo:Ch4FUBINIscambiato}
\int_{(a,b)\times R} f(x,\zb y)\ZD(x,\zb y)  =\int_R \left [  \int_a^b f(x,\zb y)\ZD x\right ]\ZD \zb y\,.
\end{equation}
\end{enumerate}
\end{Corollary}
\end{subequations}

  \begin{Remark}
   {\rm \emph{This remark is of interest to readers who have been introduced to integrability respect to arbitrary measures, a topic we did not touch.}
  
   Needless to say, Fubini Theorem holds because the measure in a space of higher dimension is related to that of its subspaces. If this relation does not hold then the reduction formula does not hold too, as the following example show. Let $\ZDE$ be the Dirac measure in $\zzr^2$:
   \[
\ZDE(A)=\left\{\begin{array}{ll}
1&\mbox{ if $(0,0)\in A$}\\
0 &\mbox{otherwise}\,.
\end{array}   \right.
   \]
Every   function which is continuous in a neighborhood of  $(0,0)$ is summable   respect to $\ZDE$ and its integral on a set $A$ is
  \[
 \left\{\begin{array}{ll}
f(0,0)&\mbox{ if $(0,0)\in A$}\\
0 &\mbox{otherwise}\,.
\end{array}  \right.
  \]
  It is then easily checked that the reduction formula does not hold 
for the function $f(x,y)=1/[(1+x)(1+y)]$  
  on the rectangle $[-1/2,1/2]\times[-1/2,1/2]$ with the Dirac measure, if we use the Lebesgue measure on the segments.
  
 In fact the   abstract treatment     of Fubini theorem first considers  two   spaces, $X$ and $Y$, each one with its measure. The next step is the constructs the ``product measure'' in $X\times Y$. After that Fubini Theorem, respect to the product measure, is stated and proved.\zdia 
  
   }
   \end{Remark}

The second    result we shall prove is:
   \begin{Theorem}[{\sc Tonelli}]\index{Theorem!Tonelli}\index{Tonelli Theorem}\ZLA{Theo:Ch4Tonelli}
  Let $f(x,\zb y)$ be quasicontinuous and nonnegative  on $(a,b)\times R\subseteq \zzr^{d+1}$ and suppose that $\zb y\mapsto f(x,\zb y)$ is summable for a.e. $x\in (a,b)$.   The function $(x,\zb y)\mapsto f(x,\zb y)$ is summable on $R$ (and so 
the reduction formula~(\ref{eq:LI2Theo:Ch4FUBINI}) holds) if  
\[
x\mapsto  \left [  \int_R f(x,\zb y)\ZD y\right ]
\]
is summable.
  \end{Theorem}
  
  \begin{Remark}
  {\rm In order to appreciate Theorem~\ref{Theo:Ch4Tonelli} it is convenient to keep in mind the following fact: there exists   positive \emph{but not summable} functions $f(x,\zb y)$ such that $\zb y\mapsto f(x,\zb y)$ is summable for every $x$. An example in $ (-\ZIN,+\ZIN)\times(-\ZIN,+\ZIN)= \zzr^2$ is the function
  \[
f(x,y)=e  ^{-|y|}\,.
  \]
This function $f$ is not summable in $\zzr^2$. In spite of this, for every $x $ we have
  \[
  \int _{-\ZIN}^{+\ZIN} f(x,y)\ZD y=2 
  \]
  but note that the function
  \[
x\mapsto    \int _{-\ZIN}^{+\ZIN} f(x,y)\ZD y\quad\mbox{is not summable}\,.\zdiaform 
  \]
  }
  \end{Remark}
 
 \subsection{\ZLA{Ch4PreliFubToneNuovo}Preliminary Results}  
 
First we present   observations on null sets and on sequences of Tietze extensions  and then we examine the sections with planes and lines of multirectangle sets.
 \subsubsection{\ZLA{ch4SectPreliOBSEaFubini}Observation on Null Sets and Quasicontinuous Functions}
 
We use statement~\ref{coroCH2copriOPENNULLset}
of Corollary~\ref{Ch3CoriFrontiRectNULLO}: a set $N\subseteq \zzr^n$ is a null set when  there exists a sequence $\{\Delta_k\}$ of open multirectangles\footnote{according to the definition in the table~\ref{tableCH3ADDITIONALshortNOTATIONS},  a multirectangle $\{R_k\}$ is an open multirectangle when each $R_k$ is open.}   such that 
 
\begin{equation}\ZLA{Eq;CH4SequePerNullNONrego} 
 \mbox{$  N\subseteq \mathcal{I}_{\Delta_k}$ and $\lim _{k\to+\ZIN}L(\Delta_k)=0$ } 
 \,.
  \end{equation}
 It is possible to choose a new sequence $\{\hat \Delta_k\}$ of open multirectangles with the  properties   
 
  \begin{equation}
 \ZLA{Eq;CH4SequePerNullNONdecreasing}  
  N\subseteq \mathcal{I}_{\hat \Delta_k} \  \mbox{and} \ 
  \left\{\begin{array}{l} \lim_{k\to+\ZIN}L(\hat \Delta_k)=0\,,\\  
  \mathcal{I}_{\hat\Delta_{k+1}}\subseteq  \mathcal{I}_{\hat\Delta_{k } } \,.
  \end{array}\right.
 \end{equation}
  The multirectangles $\hat\Delta _k$ are constructed   as follows from the sequence $\{\Delta_k\}$ in~(\ref{Eq;CH4SequePerNullNONrego}): by passing to subsequences, we can assume that the sequence $\{\Delta_k\}$ has the following property:
  \[
L(\Delta_k)<\dfrac{1}{k}\,.  
  \]
Then we define
 \[
\hat\Delta _k=\bigcup _{r\geq 1} \Delta_{2^{k+r}}=\bigcup _{\nu\geq k+1}\Delta_{2^\nu}\quad\mbox{so that $L(\hat\Delta_k)<1/2^k\to 0$}\,.
 \]
 
\begin{Remark} \ZLA{REMA;CH4SequePerNullNONdecreasingBIS}  
{\rm
We combine  the following observations:
\begin{enumerate}
\item the sets $ \mathcal{I}_{\hat \Delta_k} $ are open sets.
\item any open set is a multirectangle set and it   can be ``approximated'' from outside by   multirectangle sets composed by open rectangles, see~(\ref{ch3DefiLebeMISaperti}).
\end{enumerate} 
By using these observations we can reformulate the property of being a null set as follows: a set $N$ is a null set when there exists a   sequence $\{\mathcal{O}_k\}$ of open sets such that
  \begin{equation}
 \ZLA{Eq;CH4SequePerNullNONdecreasingBIS}  
  N\subseteq \mathcal{O}_k \  \mbox{and} \ 
  \left\{\begin{array}{l} \lim_{k\to+\ZIN}\zl(\mathcal{O}_k)=0\,,\\  
  \mathcal{O}_{ _{k+1}}\subseteq  \mathcal{O}_{{k } } \,.\zdiaform
  \end{array}\right. 
 \end{equation}
}
\end{Remark} 
 
 A similar idea we apply to the associated Tietze extensions of a quasicontinuous function. 
 
 Let $f$ be quasicontinuous on a rectangle $R$ and let $\{\Delta_k\}$ be a sequence of open multirectangles such that
$L(\Delta_k)\to 0$ and 
  $f_k=f_{|_{R\setminus\Delta_k}}$ is continuous. We denote $f_{k,e} $ a Tietze extension of $f_k$.
 
 The sequence $\{f_{k,e}\}$ in general does not converge to $f$ not even a.e. on $R$ but it is possible to select a suitable sequence of Tietze extensions   which a.e. converges  to $f$. As above, first we choose $\{\Delta_k\}$ such that $L(\Delta_k)<1/k$ and then we denote
\[
 \hat\Delta _k=\bigcup _{\nu\geq k+ 1} \Delta_{2^\nu}\quad\mbox{so that $\lim_{k\to+\ZIN}L(\hat \Delta_k)=0$ and $\hat\Delta_{k+1}\subseteq\hat\Delta_k$}\,.
  \]
   Then we choose the sequence $\{\hat f_{k,e}\}$
 \[
 \hat f_{k,e} =f_{2^{k},e}\quad\mbox{so that \ \ $\hat f_{k,e}(x)=f(x)$ if $x\in R\setminus \Delta_{2^k}$} \,.
 \]
 
  Let
 \[
 N=
\underbrace{ \bigcap_{k\geq 1}\hat\Delta_k}_{\tiny  
 {
 \begin{array}{l}     \subseteq \hat\Delta_\nu\ \forall \nu \\
\mbox{hence a null set} \end{array}
 } 
}\,.
\]
  If $x\notin N$
there exists $K=K_x$ such that
$x\notin  \hat\Delta _{K_x}$. So we have also $x\notin \hat\Delta _\nu$ for every $\nu>K_x$ since the sequence $\{\hat\Delta_\nu\}$ is decreasing. Hence
\[
\nu> K_x \ \implies\ 
\hat f_{\nu,e}(x)=f_{2^\nu,e}(x)=f(x) 
\quad\mbox{so that}\quad 
\lim _{\nu\to+\ZIN} \hat f_{\nu,e}(x)=f(x)
\]  
as wanted.

\subsubsection{\ZLA{ch4SECprojeMesSETS}Sections of Multirectangle Sets}
Let $A\subseteq \zzr^{d+1} $. We define its sections 
 \begin{equation}\ZLA{eq:Ch4DEFIsezione}
  \begin{array}{l}\displaystyle 
A_x=\{ \zb y  \ \mbox{such that}\ (x,\zb y)\in A\}\subseteq \zzr^d \,,\\
\displaystyle
 A_{\zb y} 
=\{x  \ \mbox{such that}\ (x,\zb y)\in A\} \subseteq\zzr\,.
\end{array}
  \end{equation}

  The first result of this section, Lemma~\ref{Ch4LemmaRiduzMIsurINSIEM}, gives information  on the sections of multirectangle sets. Note that
  Lemma~\ref{Ch4LemmaRiduzMIsurINSIEM}  is the statement of the reduction formula in a special case.

  \begin{Lemma}
  \ZLA{Ch4LemmaRiduzMIsurINSIEM}
  Let $A\subseteq \zzr^{d+1}$ be a bounded multirectangle set, $A\subseteq [a,b]\times R$ where $R$ is a bounded rectangle of $\zzr^{d }$. 
We assume $A=\cup _{n\geq 1}R_n$ where $R_n$ are rectangles such that $L(R_n)>0$.  
  
  For every $x$ and $\zb y$, let $A_x$ and $ A_{\zb y} $ be its sections defined in~(\ref{eq:Ch4DEFIsezione}). 
Then:
  \begin{enumerate}
  
  \item \ZLA{I1Ch4LemmaRiduzMIsurINSIEM}
 the section $A_x$ is a multirectangle set of $\zzr^d$ \ a.e. $x\in [a,b]$  and $A_{\zb y}$ is a multiinterval a.e. $\zb y\in R$.
    \item \ZLA{I2Ch4LemmaRiduzMIsurINSIEM} the functions 
  \[
x\mapsto \int   _R \charfun _{A_x}(\zb y)\ZD\zb y=\zl_d (A_x)\,,\qquad 
\zb y\mapsto \int   _a^b \charfun _{A_{\zb y}}(x)\ZD x=\zl_1(A_{\zb y})
  \]
 are bounded quasicontinuous.
  \item \ZLA{I4Ch4LemmaRiduzMIsurINSIEM} we have
  \begin{subequations}
\begin{align} \nonumber
  \zl_{d+1}(A)& = \int_{[a,b]\times R} \charfun_A(x,\zb y) \ZD(x,\zb y)  
   =\int _a^b \left [ \int_ R \charfun _A(x,\zb y)\ZD\zb y   \right ]\ZD x\\
\ZLA{eq:Ch4Ridu1FUNZcar}   &=
  \int _a^b \zl_d(A_x)\ZD x\,, \\
  \nonumber
   \zl_{d+1}(A)&= \int_{[a,b]\times R} \charfun_A(x,\zb y) \ZD(x,\zb y)   
   =\int _R \left [ \int  _a^b \charfun _A(x,\zb y)\ZD  x   \right ]\ZD \zb y\\
\ZLA{eq:Ch4Ridu2FUNZcar}   =&
  \int _R \zl_1(A_{\zb y})\ZD \zb y  
  \,.
  \end{align}
  \end{subequations}
  \end{enumerate}
  \end{Lemma}
  \zProof 
  
We prove
the statement which concern the sections $A_x$ and
  the equality~(\ref{eq:Ch4Ridu1FUNZcar}). 
In an analogous way we can prove the statements which concern $A_{\zb y}$
  and  equality~(\ref{eq:Ch4Ridu2FUNZcar}).
  
 We note:
\begin{align*}
&\mbox{the assumption}\\
& \mbox{ 
  $A=\mathcal{I}_\Delta $ where $\Delta=\{R_n\}$ is almost disjoint and $L(R_n)>0$ } \\
  &\mbox{implies}
\\
&A_x=\{ \zb y  \ \mbox{such that}\ (x,\zb y)\in A\} =
\bigcup _{n>1}\{ \zb y  \ \mbox{such that}\ (x,\zb y)\in R_n\} =\bigcup _{n\geq 1} (R_n)_x 
\end{align*}
  and $(R_n)_x$ is a rectangle of $\zzr^d$. 
  It follows that $A_x\in\mathcal{L}(\zzr^d)$ and so its characteristic function   is quasicontinuous (see theorem~\ref{TeoCH4InsieLebeUgualeFUNcharQUASICONT}).
  
 \emph{ This proves that the integral of $\charfun_{A_x}$  in the statement~\ref{I2Ch4LemmaRiduzMIsurINSIEM} exists.  }

 We prove that $A_x$ is a multirectangle set for a.e. $x$ so that its measure can be computed by adding the measures of the rectangles.  
  We divide the proof in the following steps.
  \paragraph{Step~1: we prove the statement~\ref{I1Ch4LemmaRiduzMIsurINSIEM} of the lemma}
  We must prove that the rectangles $(R_n)_x$, considered as rectangles on the affine iperplane of the section, are almost disjoint for a.e. $x$. We need a notation.  
An open ball in $\zzr^{d+1}$ of radius $r$ and center    $(x,\zb y)$ is denoted $B((x,\zb y),r)$. Let $A\subseteq \zzr^{d+1}$. We fix $x$ and we consider the points $\zb y\in A_x$ which have the following property: there exists  $B((x,\zb y),r)$ such that its section is contained in $A_x$:
\[
[B((x,\zb y),r)]_{x}\subseteq A_x \,.
\] 
The set of the point $\zb y$ with this property is denoted ``${\rm r.int}\, A_x$''.
 
 Let $\zb {y_0}\in (R_k)_x\cap {\rm r.int }\, (R_n)_x$. Then $(x,\zb{y_0})\in R_k\cap R_n$. These rectangles  are quasidisjoint. So, the value of $x$ must correspond to the coordinate of an ``upper'' or ``lower''  face of $R_n$. The rectangles and so their faces are a numerable set. So, \emph{the family $\{(R_n)_x\}$ can be not almost disjoint only for a numerable set of values of $x$; and numerable sets are null sets.}
 
  \paragraph{Step~2: we prove the statement~\ref{I2Ch4LemmaRiduzMIsurINSIEM} of the lemma}
Egorov-Severini Theorem implies that $x\mapsto \zl_d(A_x)$ is quasicontinuous since
 \[
\zl_d(A_x)=\lim_{k\to+\ZIN}\left [\sum _{n=1} ^k \zl_d((R_n)_x)\right ]\qquad  {\rm a.e.}\ x\in [a,b]\,,
 \]
the limit of a sequence of quasicontinuous functions.

Boundedness follows from the assumption that $A$ is contained in a bounded rectangle of $\zzr^{d+1}$.

\emph{We recapitulate: these observations prove  the properties of $A_x$ in the statements~\ref{I1Ch4LemmaRiduzMIsurINSIEM} and~\ref{I2Ch4LemmaRiduzMIsurINSIEM} of the lemma.}
 
 \paragraph{Step~3: we prove the statement~\ref{I4Ch4LemmaRiduzMIsurINSIEM} of the lemma} In particular, we prove~(\ref{eq:Ch4Ridu1FUNZcar}).
  We note the following facts:
  \begin{enumerate}
  \item For every $n$ we have $\zl_{d+1}(R_n)=L(R_n)=L({\rm int}\,R_n)=\zl_{d+1}({\rm int}\,R_n)$;
  \item the set $\cup\partial(R_n)$ is a null set, so that integrals computed on $A$ and on $A\setminus\left [ \cup\partial(R_n)  \right ]$ have the same value; in particular, $\zl_{d+1}(A)=\zl_{d+1}\left (A\setminus\left [ \cup\partial(R_n)  \right ]\right )$.
  \item So, by removing the boundaries of the rectangles $R_n$ we can assume that $R_k\cap R_j=\emptyset$ for every $k$ and $j$ and that the set $A$ is a numerable union of disjoint rectangles.
  \end{enumerate}
  
 Now, 
equality~(\ref{eq:Ch4Ridu1FUNZcar}) follows from the following chain of equalities, which is justified below:
  \begin{multline*}
  \zl_{d+1}(A)=
  \ZSUno \zl_{d+1}(R_n)=
  \ZSUno \underbrace{ \int _{R_n} 1\ZD(x,\zb y)}
_{\tiny \begin{tabular}{l}
    Riemann    integral  \end{tabular}}
  = 
 \ZSUno \underbrace{\int _a^b \left [ \int_ R \charfun _{R_n}(x,\zb y)\ZD\zb y   \right ]\ZD x}
_{\tiny \begin{tabular}{l}
both    Riemann    integrals  \end{tabular}}
=\\
 \underbrace{\int _a^b \left [\ZSUno\overbrace{\int_ R \charfun _{R_n}(x,\zb y)\ZD\zb y}^{\tiny \begin{tabular}{l}
    Riemann  \\    integral  \end{tabular}}   \right ]\ZD x}_{\tiny \begin{tabular}{l}
    Lebesgue      integral  \end{tabular}}
=
  \underbrace{\int _a^b \left [\ZSUno\zl_d((R_n)_x)\right ]\ZD x}_{\tiny \begin{tabular}{l}
    Lebesgue       integral  \end{tabular}}= 
\underbrace{   \int _a^b \zl_d(A_x)\ZD x }_{\tiny \begin{tabular}{l}
    Lebesgue      integral  \end{tabular}} \,.  
  \end{multline*}
These equalities are justified by the following observations:
\begin{enumerate}
\item the integrals denoted ``Riemann integrals'' are integrals of piecewise continuous functions. Hence they are bona fide Riemann integrals and so also Lebesgue integrals.
\item the sequence   \[
k\mapsto   \sum_{n=1}^k \int_ R \charfun _{R_n}(x,\zb y)\ZD\zb y = \sum_{n=1}^k 
\zl_d((R_n)_x)
  \]
  is a bounded increasing sequence of nonnegative functions. 
It is bounded by $\zl_d(R)$ since we reduced ourselves to the case that the rectangles $R_n$ are pairwise disjoint. 
  \emph{ Its limit is $\zl_d(A_x)$ and   $x\mapsto \zl_d(A_x) $ is quasicontinuous since it is the limit of a sequence of quasicontinuous functions.}
  \item  The exchange of the series and the integral is justified by Beppo Levi Theorem. Note that once the series and  the integral has been exchanged, the exterior integral is an integral in the sense of Lebesgue.
  \item the equality   
  \[
  \ZSUno\zl_d((R_n)_x)=\zl_d(A_x) 
  \]
   follows (a.e. $x\in [a,b] $) from the fact that the sets $(R_n)_x$ are pairwise disjoint.
\end{enumerate}  

These observations complete the proof.\zdia
  
The second result of this section concerns sections of null sets.  
   
  \begin{Lemma}\ZLA{teoCh4:SectNULLediInsiNulli}
  Let $N\subseteq [a,b]\times R\subseteq \zzr^{d+1}$ be a bounded null set.
  Let $N_x$ and $N_{\zb y} $ be its sections  defined as in~(\ref{eq:Ch4DEFIsezione}). Then $N_x$ is a null set a.e. $x\in [a,b]$
and $N_{\zb y}$ is a null set a.e. $\zb y\in R$.
  \end{Lemma}
  \zProof  
We use Remark~\ref{REMA;CH4SequePerNullNONdecreasingBIS} and we prove that $N_x$ is a null set. The proof uses~(\ref{eq:Ch4Ridu1FUNZcar}). A similar argument, based on~(\ref{eq:Ch4Ridu2FUNZcar}),  shows that $N_{\zb y}$ is a null set.
We use~(\ref{Eq;CH4SequePerNullNONdecreasingBIS}). Because $N$ is a null set, there exists a sequence $\{\mathcal{O}_k\}$  of open sets with the following properies:
  \[
N\subseteq \mathcal{O}_{k+1 }  \subseteq \mathcal{O}_k\quad \mbox{for all $k$ and}\quad \zl_{d+1}(\mathcal{O}_k)=     \int _a^b \zl_d\left ((\mathcal{O}_k)_x  \right )\ZD x\to 0\,.
  \]
Note that~(\ref{eq:Ch4Ridu1FUNZcar}) can be used since nonempty open sets are multirectangle sets, union of rectangles $R_n$ such that $L(R_n)>0$ (see Theorem~\ref{TeoCH3STRUttAPERTI}). 
  
The inclusion
  \[
N_x  \subseteq (\mathcal{O}_{k+1})_x
 \subseteq ( \mathcal{O}_{k })_x 
  \]
shows that $k\mapsto  \zl_d\left ((\mathcal{O}_k)_x  \right )$ is bounded decreasing, hence convergent  for every~$x$, 
  \[
\lim _{k\to +\ZIN}  \zl_d\left ((\mathcal{O}_k)_x  \right )=f(x)\geq 0\,.  
  \]
  Theorem~\ref{teo:ch3:LebeCAROPartIntervLIMIT} implies
  \[
 \zl_{d+1}(\mathcal{O}_k)=    \int _a^b \zl_d\left ((\mathcal{O}_k)_x  \right )\ZD x\to \int_a^b f(x)\ZD x\quad\mbox{and so}\quad
   \int_a^b f(x)\ZD x  
     =0\,.
   \]
Theorem~\ref{CH4FunzCONinteNulloNULLA} shows that $f(x)=0$ a.e. $x\in (a,b)$ and so 
   
   \[
\zl_d\left (( \mathcal{O}_k)_x  \right )\to 0 \quad {\rm a.e.}\ x\in (a,b)\,.   
   \]
 So   we have a.e. on $[a,b]$:
   \[
N_x\subseteq    (\mathcal{O}_k)_x  \,,\qquad \mbox{ $(\mathcal{O}_k)_x $ is open in $\zzr^d$ and $\zl( (\mathcal{O}_k)_x )\to 0$ }\,.
   \]
 
   It follows that $N_x$ is a null set a.e. $x\in (a,b)$.\zdia
   
 The previous lemma can be lifted from bounded to unbounded null sets, by intersecting the set with an increasing sequence of rectangles, but we don't need this observation.

   \section{\ZLA{sec:CH4SecDimoFunbini}Fubini and Tonelli Theorems: the Proofs}

   Before proving the theorems, we note that the statements, expressed in terms of rectangles, are not restrictive. We recall that we defined
   \[
\int_{A} f(x,\zb y)\ZD (x,\zb y  ) 
   \]   
   when the set $A$ has the property that $\charfun_A$ is quasicontinuous and when\footnote{we recall that strictly speaking the notation  $f\charfun_A$  makes sense if  the domain of $f$ contains $A$. If not, $f$ is replaced by its extension to $\zzr^{d+1}$ with $0$.} $f\charfun_A$  is integrable.  By definition:
  \[
\int_{A} f(x,\zb y)\ZD (x,\zb y  ) =\int_{\zzr^{d+1}} f(x,\zb y)\charfun_A(x,\zb y) \ZD (x,\zb y  )\,.
   \]
 Fubini Theorem can be applied when $f \charfun_A$ is summable. Under this condition formula~(\ref{eq:LI2Theo:Ch4FUBINIblocchi}) takes the following form\footnote{$R_1\times R_2=\zzr^{d_1}\times \zzr^{d_2}$ can be replaced by a bounded rectangle $R_1\times R_2$   if $A$ is bounded.}:
 \begin{multline*}
 \int_{A} f(x,\zb y)\ZD (x,\zb y  ) =\int_{\zzr^{d_1}\times \zzr^{d_2}} 
 f(\zb{x_1},\zb{x_2})\charfun_A (\zb{x_1},\zb{x_2})
 \ZD(\zb{x_1},\zb{x_2})\\
 =\int _{\zzr^{d_1}}\left [\int _{\zzr^{d_2}}  f(\zb{x_1},\zb{x_2})\charfun_A (\zb{x_1},\zb{x_2})\ZD {\zb {x_2}}\right ] \ZD\zb {x_1}\\=\int _{\zzr^{d_1}}\left [ \int _{A_{\zb {x_2}}} 
  f(\zb{x_1},\zb{x_2})\charfun_A (\zb{x_1},\zb{x_2})\ZD {\zb {x_1}}\right ]\ZD \zb {x_2}  
  \\
  =\int _{A_{\zb {x_1}}} \left [ \int _{A_{\zb {x_2}}} 
  f(\zb{x_1},\zb{x_2})\ZD {\zb {x_1}}\right ]\ZD \zb {x_2}
  \,.
 \end{multline*}
 
After this observation we prove   Theorem~\ref{Theo:Ch4FUBINI}.
Linearity of the integral shows that we can prove the theorem separately for the functions
\[
f_+(x,\zb y)=\max\{f(x,\zb y)\,,\ 0\}\,,\qquad f_-(x,\zb y)=\min\{f(x,\zb y)\,,\ 0\}
\]   
i.e. we can prove the theorem when the function has constant sign, say when

\[
f\geq 0\,.
\]
   
First we prove the theorem when $f\geq 0$ is a bounded functions on a bounded rectangle $(a,b)\times R$. Then we extend to the general (unbounded) case.
 
   \subsubsection{Fubini Theorem: $f$   Bounded on a Bounded Rectangle}
   We use the following facts:
   \begin{itemize}
   \item The Lebesgue and the Riemann integrals coincide for continuous functions;
   \item the reduction formula on a rectangle holds for the Riemann, hence  the Lebesgue, integral of continuous functions;
   \item if $f$ is quasicontinuous on a bounded rectangle $(a,b)\times R\subseteq \zzr^{d+1}$ then\footnote{as noted in Sect.~\ref{ch4SectPreliOBSEaFubini}.}  there exists a
decresing sequence $\{\mathcal{O}_k\}$ of open sets in $(a,b)\times R$, such that $\zl_{d+1}(\mathcal{O}_k)\to 0$, and a   
    bounded sequence $\{\hat f_k\}$ of \emph{continuous} functions, such that
   \[
   \begin{array}{l}
   f(x,\zb y)=\hat f_k(x,\zb y) \quad {\rm if}\quad (x,\zb y)\notin \mathcal{O}_k\,,\\
f(x,\zb y)=\lim _{k\to+\ZIN} \hat f_k(x,\zb y)\quad x\in   \left [(a,b)\times R\right ]\setminus N
\end{array}     
   \]
   where $N$ is a null set.
We recall, from Lemma~\ref{teoCh4:SectNULLediInsiNulli}, that the section $N_x$ is a null set in $R$ a.e. $x\in (a,b)$ and $N_{\zb y}$ is a null set in $(a,b)$ a.e. $\zb y\in R$.
    \end{itemize}
 
 Now we proceed with the following steps:
 \begin{description}
 \item[\bf Step~1)] \emph{$\zb y\mapsto f(x,\zb y)$ is a.e. quasicontinuous  on $(a,b)$.}
 We fix $x_0\notin N_{\zb y}$ and we consider the function $\zb y\mapsto f(x_0,\zb y)$. We have
 \[
 f(x_0,\zb y)=\lim _{k\to+\ZIN} \hat f_k(x_0,\zb y)\qquad {\rm a.e.}\ \zb y\in R\,.
 \]
 The function  $\zb y\mapsto\hat f_k(x_0,\zb y)$ is 
continuous since $\hat f_k(x,\zb y)$ is continuous. Hence, its a.e. limit $f(x_0,\zb y)$ is quasicontinuous.

 \item[\bf Step~2)] \emph{the inner integral in the right hand side of the reduction formula~(\ref{eq:LI2Theo:Ch4FUBINI}) is quasicontinuous.}
 In fact,  
 \begin{multline*}
\underbrace{\int_R f(x,\zb y)\ZD \zb y }_{
\tiny \begin{tabular}{l}
    Lebesgue      integral  \end{tabular}
    }
=\underbrace{\int_R\left [\limn  \hat f_n(x,\zb y)\right ] \ZD \zb y }_{\tiny 
\begin{tabular}{l}
    Lebesgue      integral  \end{tabular}
    }
    \\
= \limn\underbrace{ \overbrace{\int_R \hat f_n(x,\zb y)\ZD \zb y}^{\tiny 
\begin{tabular}{l}
  continuous function of $x$  \end{tabular}
    }  }_{\tiny 
\begin{tabular}{l}
    Riemann      integral  \end{tabular}
    } \quad \mbox{ a.e. on $(a,b)$}
 \end{multline*}
 and the limit of a sequence of continuous function is quasicontinuous.
 
 Note that we use boundedness of the domain and of the sequence to exchange the limit and the Lebesgue integral.
 \item[\bf Step 3)] the reduction formula~(\ref{eq:LI2Theo:Ch4FUBINI}) holds under the stated boundedness assumptions. In fact,
    boundedness of the domain and boundedness of the sequence imply
 {\footnotesize
    \begin{multline*}
 \underbrace{
  \int _{(a,b)\times R} f(x,\zb y)\ZD (x,\zb y) 
  }_{\tiny \begin{tabular}{l}
     Lebesgue      integral  \end{tabular}}
 =
  \underbrace{
  \int _{(a,b)\times R}\left [\lim _{k\to+\ZIN} \hat f_k(x,\zb y)\right ]\ZD (x,\zb y) 
  }_{\tiny \begin{tabular}{l}
     Lebesgue      integral  \end{tabular}}\\
    =\lim _{k\to+\ZIN}   
 \underbrace{
\int _{(a,b)\times R} \hat f_k(x,\zb y)\ZD (x,\zb y) 
}_{\tiny \begin{tabular}{l}
    Riemann      integral  \end{tabular}}
=\lim _{k\to+\ZIN}   
 \underbrace{
\int _a^b  \left [
\underbrace{\int_R \hat f_k(x,\zb y)\ZD\zb y}_{\tiny \begin{tabular}{l}
  Riemann      integral  \end{tabular}}
   \right ]
  \ZD x 
}_{\tiny \begin{tabular}{l}
  Riemann      integral  \end{tabular}}\\
=
 \underbrace{
\int _a^b\left [ \lim _{k\to+\ZIN}  
 \underbrace{
\int_R \hat f_k(x,\zb y)\ZD\zb y 
}_{\tiny 
\begin{tabular}{l}
   Riemann      integral  \end{tabular}
   }
     \right ]\ZD x
      }_{\tiny 
     \begin{tabular}{l}
    Lebesgue      integral  \end{tabular}
    }
    %
 =
 \underbrace{
\int _a^b    \left [\underbrace{\int_R \left (\lim _{k\to+\ZIN} \hat f_k(x,\zb y)\right )\ZD\zb y
 }_{\tiny \begin{tabular}{l}
    Lebesgue      integral  \end{tabular}}\right ] \ZD x  }_{\tiny \begin{tabular}{l}
    Lebesgue      integral  \end{tabular}}
\\
= 
 \underbrace{
\int _a^b  
 \left [
    \underbrace{  
\int_R f (x,\zb y) \ZD\zb y   
 }_{\tiny \begin{tabular}{l}
    Lebesgue      integral  \end{tabular}}
 \right ]
   \ZD x
     }_{\tiny \begin{tabular}{l}
    Lebesgue      integral  \end{tabular}}
   \end{multline*}
   }
   as wanted.
 
   \subsubsection{\ZLA{susbsubFubiGene}Fubini Theorem under General Assumptions}
 We still consider $f\geq 0$ but now it can be unbounded and it is  summable on a rectangle which can be unbounded too.  
 By definition
 \[
\int_{(a,b)\times R} f(x,\zb y)\ZD (x,\zb y) 
 \]
 is define in two steps: first we compute the integral on the bounded rectangle $(a_\nu,b_\nu)\times R_\nu\subseteq \zzr^{d+1}$
 where
 \[
(a_ \nu,b_\nu)= (a,b)\cap (-\nu, \nu)\,,\qquad R_ \nu=R\bigcap\left (\prod _{k=1}^d (-\nu,\nu)\right )\,.
 \]
 The function $f$ can be unbounded on this domain and so first we consider the integral of $f^N$:
 \[
f^N(x,\zb y)=\min\{ f(x,\zb y)\,,\ N\} \,.
 \]
 This way we reduce ourselves to the bounded case already studied and we proved
 the reduction formula
 \[
 \int _{(a_\nu,b_\nu)\times R_\nu} f^N(x,\zb y)\ZD (x,\zb y)=
 \int _{a_\nu}^{b_\nu}    \left [\int_{R_\nu}  f^N (x,\zb y) \ZD\zb y\right ] \ZD x\,.
 \]
 The sequence $\{f^N\}$ is increasing and a.e. convergent to $f$ so that we can compute the limit for $N\to+\ZIN$ and exchange the limit with the integrals (twice on the right hand side) thanks to Beppo Levi Theorem. We get
 \begin{equation}\ZLA{eq:Ch4FubiSOILefILLI}
  \int _{(a_\nu,b_\nu)\times R_\nu} f (x,\zb y)\ZD (x,\zb y)=
 \int _{a_\nu}^{b_\nu}    \left [\int_{R_\nu}   f  (x,\zb y) \ZD\zb y\right ] \ZD x\,.
 \end{equation}
 
Then we compute the limit for $\nu\to +\ZIN$. We note that the integrals in~(\ref{eq:Ch4FubiSOILefILLI}) are the integrals on $(a,b)\times R$ (on the left side) and on $R$ and on $(a,b)$ (on the right side) of the function
\begin{equation}\ZLA{eq:CH4NellaDimpFubCnFunzCara}
f(x,\zb y)\charfun_{(a_\nu,b_\nu)\times R_\nu}(x,\zb y)
\end{equation}
and the sequence $\{f(x,\zb y)\charfun_{(a_\nu,b_\nu)\times R_\nu}(x,\zb y)\}$ is increasing and a.e. convergent to $f$. So we can  use   Beppo Levi Theorem  again  and finally get the reduction 
formula~(\ref{eq:LI2Theo:Ch4FUBINI}).

  \subsubsection{Tonelli Theorem: the Proof}

  The proof consists in the observation that under the assumption of Tonelli Theorem the function $f$ is summable, hence the conditions of Fubini Theorem are satisfied.
  
  Due to the fact that $f$ is nonnegative, it is sufficient to show that the integral on the left side of~(\ref{eq:LI2Theo:Ch4FUBINI})  is not $+\ZIN$.
  
 We use~(\ref{eq:CH4NellaDimpFubCnFunzCara}). We see that 
 \[
\int_{(a,b)\times R}f(x,\zb y)\ZD(x,\zb y) =\lim_{\stackrel{N\to+\ZIN}{\nu \to+\ZIN}} 
\int_a^b \left [\int_R f_N(x,\zb y)\charfun_{(a_\nu,b_\nu)\times R_\nu}(x,\zb y)\ZD\zb y\right ]\ZD x\,.
 \]
 The right hand side is an increasing sequence of $(N,\nu)$ which, under the assumption of Tonelli Theorem,  is bounded so that
 \[
\int_a^b \left [\int_R f (x,\zb y) (x,\zb y)\ZD\zb y\right ]\ZD x <+\ZIN\,,
 \]
 as wanted.
 \end{description}
  \part{\ZLA{PART3TheMeasure}Recovering Lebesgue Measure}

\chapter{\ZLA{Chap:4Measure}Borel and Lebesgue Measure}
 
On purpose, measure theory is not used in the Tonelli approach to Lebesgue integration.   But, once the Lebesgue integral has been defined, the Lebesgue measure of sets  can be recovered. This is the  goal of this chapter.

\section{\ZLA{ch4SECDEFIprimPromeasurSETS}Open Sets and Lebesgue Measurable Sets}

The following result is well known and easily proved
\begin{Theorem}\ZLA{CH4:TeoPRELIcontFUN}
Let $K\subseteq \zzr $ and let    $f$:  $\zzr^d\supseteq K\mapsto    \zzr^{m} $. The following properties are equivalent:
 
\begin{enumerate}
\item\ZLA{I1CH4:TeoPRELIcontFUN}   the function $f$  is continuous on $K$;
 
\item\ZLA{I2CH4:TeoPRELIcontFUN}  the set  $f^{-1}(A)$ is relatively open in $K$ for every open set $A\subseteq \zzr^m$.

 \item\ZLA{I3CH4:TeoPRELIcontFUN}  the set  $f^{-1}(C)$  is relatively closed in $K$ for every closed set  $C \subseteq \zzr^m$. 
\end{enumerate}
 
 Furthermore, when $m=1$:
 \begin{itemize}
 \item   in order to check the properties in the statement~\ref{I2CH4:TeoPRELIcontFUN} it is sufficient to consider the case that $A$ is any open interval or even solely the case that $A$ is any open half line, both $(a,+\ZIN)$ and $(-\ZIN,b)$. 
 
 \item
in order to check the properties in the statement~\ref{I3CH4:TeoPRELIcontFUN} it is sufficient to consider the case that $A$ is any closed interval or even solely the case that $A$ is any closed half line, both $[a,+\ZIN)$ and $(-\ZIN,b]$. 
 \end{itemize}
\end{Theorem}

This theorem shows a relation between the lattice structure of the sets and that of the functions. In fact, the following observation (which we used already) holds:   if $f$ and $g$ are continuous   then
\[
\phi(x)=\max\{f(x)\,,\ g(x)\}\,,\qquad \psi(x)=\min\{f(x)\,,\ g(x)\} 
\]
are continuous. Continuity is easily seen from Theorem~\ref{CH4:TeoPRELIcontFUN} since
\begin{align*}
&\left\{\begin{array}{l}
\{x\,:\ \phi(x)>a\} =\{x\,:\ f(x)>a\} \cup \{x\,:\ g(x)>a\} \\
\{x\,:\ \phi(x)<a\} =\{x\,:\ f(x)<a\} \cap \{x\,:\ g(x)<a\} \,, 
\end{array}\right.\\
&\left\{\begin{array}{l}
\{x\,:\ \psi(x)>a\} =\{x\,:\ f(x)>a\} \cap \{x\,:\ g(x)>a\}\\
 \{x\,:\ \psi(x)<a\} =\{x\,:\ f(x)<a\} \cup \{x\,:\ g(x)<a\}\,.
\end{array}\right.
\end{align*}
and finite unions  and intersections of open (or closed) sets are open (or closed) sets.

Of course, $\phi$ and $\psi$ are continuous if they are the maximum or the minimum of any \emph{finite} set of functions. Instead,  nothing can be said of the  functions
 \[
\phi(x)=\sup\{f_n(x)\,,\  1\leq n<+\ZIN\}\,,\qquad \psi(x)=\inf\{f_n(x)\,,\ 1\leq n<+\ZIN\} \,.
\]
In fact, let
\[
f_n(x)=\left\{\begin{array}{lll}
n|x|&{\rm if}& |x|\leq 1/n\\
1&{\rm if}& |x|> 1/n\,,
\end{array}\right. \qquad f(x)=
\left\{\begin{array}{lll}
0&{\rm if}& x=0\\
1& & \mbox{otherwise}\,.
\end{array}\right. 
\]
 Then we have
 \begin{align*}
& \phi(x)=\sup\{f_n(x)\,,\  1\leq n<+\ZIN\}=f(x)\,,\\
 & \psi(x)=\inf\{-f_n(x)\,,\ 1\leq n<+\ZIN\} =-f(x)\,,
 \end{align*}
both discontinuous.

In contrast with this, we proved   
  that \emph{quasicontinuity is preserved   when computing both the supremum and the infimum of bounded sequences of functions} (see    Corollary\ref{Coro:Ch3SUPinfPERch4} and, when $d=1$, Corollary~\ref{Coro:Cha1bisSUPinfPERch4}).  This observation suggest that we study the sets
$
 f^{-1}(J) 
$
 when $J$ is an open interval  and $f$ is quasicontinuous. 
We define:
\begin{Definition} \ZLA{che4DEFIlebeMisu}
 {\rm 
  {\sc Lebesgue measurable sets}\index{set!measurable!Lebesgue}\index{measurable!set (Lebesgue)}\index{Lebesgue!measurable!set}
   are
    the sets $A\subseteq \zzr ^d $ which have the following property: there exists a quasicontinuous function $f$: $\zzr^d\mapsto \zzr $  and an open interval $J$ such that
  \[
 A =f^{-1}(J ) \,.
  \]
   
  The family of the  subsets of $\zzr ^d $ which are  Lebesgue measurable  is  denoted $\mathcal{L}(\zzr  ^d)$ (or, simply, $\mathcal{L}$)\index{$\mathcal{L}$}\index{$\mathcal{L}(\zzr ^d)$}.\zdia
 }
 \end{Definition} 
 We note:
 \begin{Theorem}\ZLA{TeoCH4NULLsetsareLebe}
 The following properties hold
 \begin{enumerate}
 
  \item\ZLA{I1TeoCH4NULLsetsareLebe}
 $\zzr ^d \in \mathcal{L}(\zzr ^d )$ and $\emptyset\in \mathcal{L}(\zzr ^d )$ since 
 \[
\zzr ^d =\heaviside_{\zzr ^d }^{-1}((0,2)) \,,\qquad \emptyset =
\heaviside_{\zzr  ^d}^{-1}((3,4)) \,.
 \]
  
    \item\ZLA{I3TeoCH4NULLsetsareLebe} If $N$ is a null set then $\charfun_N$ is quasicontinuous and $N=\charfun_N^{-1}((1/2,3/2))$. So, \emph{any null set is Lebesgue measurable.}
  
\item\ZLA{I2TeoCH4NULLsetsareLebe}  The interval $J$ in the definition of Lebesgue measurable sets can be fixed at will, for example $J=(0,1)$ or  $J=(0,+\ZIN)$.
  \end{enumerate}
  \end{Theorem}
 \zProof The statements~\ref{I1TeoCH4NULLsetsareLebe} and~\ref{I3TeoCH4NULLsetsareLebe} are self explanatory. The statement~\ref{I2TeoCH4NULLsetsareLebe} follows from the following facts:
 \begin{enumerate}
 \item
 two intervals   are homeomorphic: they are transformed the first over the second by a continuous and continuously invertible function $g$; 
\item 
   the function $x\mapsto g(f(x))$ is quasicontinuous if and only if $f$ is quasicontinuous
and $g$ is continuous,   
    see the statement~\ref{I3TheoCH2LinearQC} of Theorem~\ref{TheoCH2LinearQC} (and, when $d=1$, the statement~\ref{I3BisTheoCH1PropELEMqcfunc} of Theorem~\ref{TheoCH1PropELEMqcfunc}). 
    \item  So, when $g$ is continuous and with continuous inverse, the function $f$ is quasicontinuous if and only if the composition $x\mapsto g(f(x))$ is quasicontinuous.\zdia
  \end{enumerate}
 We consider a sequence
$\{A_n\}$ of Lebesgue measurable sets. From the property~\ref{I2TeoCH4NULLsetsareLebe} of Theorem~\ref{TeoCH4NULLsetsareLebe}, there exist quasicontinuous functions $f_n$ such that
\[
A_n=f_n^{-1}((0,+\ZIN))\,.
\] 
  
Let
\[
\phi(x)=\sup\{ f_n(x)\,,\quad n>0\}
\,.
\]
The  function $\phi$ is quasicontinuous (see 
Corollary~\ref{Coro:Ch3SUPinfPERch4}).
 It is clear that
 \[
\phi^{-1}((0,+\ZIN))=\bigcup _{n\geq 1} A_n
\,.
 \]
 So we have
 \begin{Theorem}\ZLA{teoCH4SIGmaLebe}
The union of a sequence of Lebesgue measurable sets is a  Lebesgue measurable set.
 \end{Theorem}
 
A convenient characterization of   Lebesgue measurable sets is in the next theorem:
 \begin{Theorem}\ZLA{TeoCH4InsieLebeUgualeFUNcharQUASICONT}
 The set $A$ is Lebesgue measurable if and only if $\heaviside_A$ is quasicontinuous.
 \end{Theorem} 
 \zProof 
If $\heaviside _A$ is quasicontinuous then $A\in \mathcal{L}(\zzr ^d )$ since
\[
A=\heaviside_A^{-1}((1/2,3/2))\,.
\] 
Conversely we prove that if $A\in\mathcal{L}(\zzr ^d  )$ then $\heaviside_A$ is quasicontinuous.
 Let $f$ be a quasicontinuous function such that
 \[
A=f^{-1}((0,+\ZIN))\,. 
 \]
 By replacing $f(x)$ with $\max\{0,f(x)\}$ we can assume $f\geq 0$.
 
 We consider the functions
  \[
f_n(x)=\dfrac{f(x)}{f(x)+1/n}\quad \mbox{so that}\quad \limn f_n(x)=\left\{\begin{array}{lll}
0&{\rm if}& f(x)=0\\
1&{\rm if}& f(x)>0\,.
\end{array}\right.
\]
So,
\[
\limn f_n=\heaviside_A \,.
\]
The function $f$ is quasicontinuous and so the functions $f_n$ are quasicontinuous too\footnote{see   
Theorem~\ref{TheoCH2LinearQC}.}.

  Egorov-Severini Theorem shows that
  $\heaviside_A$ is quasicontinuous.\zdia
  
The Assumption~\ref{ASSUch2QuascontINSa} (in the case of functions of one variable, the Assumption~\ref{ASSUch1QuascontINSa})   implies:
 
 \begin{Corollary}
Let $f$ be defined on a set $A$. If the function is integrable then  $A$ is Lebesgue measurable.
 \end{Corollary}
 
We use $\tilde A$ to denote the complement of $A$:
\[
\tilde A=\zzr  ^d \setminus A 
\]
and we note that 
\[
\charfun_{\tilde A} (x)=1-\charfun_A(x)
\]
is quasicontinuous if and only if $\charfun_A$ is quasicontinuous. So we have:
 
\begin{Theorem}\ZLA{teoCH4SIGmaLebeCOMPLeMenTo}
The complement of every Lebesgue measurable set  is Lebesgue measurable.
\end{Theorem}

We use the following formula which holds for every sequence of sets:
\[
\bigcap _{n\geq 1} A_n=\widetilde{\ \left (\cup _{n\geq 1} \tilde A_n  \right ) \ }\,.
\]
Let the set $A_n$ be measurable. Theorems~\ref{teoCH4SIGmaLebe} and~\ref{teoCH4SIGmaLebeCOMPLeMenTo} give:
\begin{Theorem}\ZLA{Teo:Ch4InTeRnuMeLebEareLebe}
The intersection of a sequence of Lebesgue measurable sets is a  Lebesgue measurable set.
\end{Theorem}

We collect  the properties already stated and few more which can be deduced from them and from the elementary properties of the operations among sets:
\begin{Theorem}\ZLA{TeoCH4CollectingPROP}
The family of sets $\mathcal{L}(\zzr  ^d )$ has the following properties:
\begin{description}
\item[complement of sets:] if $A\in \mathcal{L}(\zzr  ^d )$ then $\tilde A\in \mathcal{L}(\zzr  ^d )$.
\item[union of sets:] if $A\in \mathcal{L}(\zzr ^d  )$ and $B\in \mathcal{L}(\zzr ^d  )$
then $A\cup B\in \mathcal{L}(\zzr ^d  )$ {\rm (since $\charfun_{A\cup B}(x)=\max\{\charfun_A(x),\charfun_B(x)\}$)}.
\item[intersection of sets:] if $A\in \mathcal{L}(\zzr  ^d )$ and $B\in \mathcal{L}(\zzr  ^d )$
then $A\cap B\in \mathcal{L}(\zzr ^d  )$  {\rm (since $\charfun_{A\cap B}(x)=\min\{\charfun_A(x),\charfun_B(x)\}$)}.
\item[difference of sets:] if $A\in \mathcal{L}(\zzr ^d  )$ and $B\in \mathcal{L}(\zzr ^d  )$
then $A\setminus B\in \mathcal{L}(\zzr  ^d )$
 {\rm (since $\charfun_{A\setminus B}(x)=\max\{0,\charfun_A(x)-\charfun_B(x)\}$)}.
\item[symmetric difference of sets:] if $A\,,\ B\in \mathcal{L}(\zzr  ^d )$ 
then $A\triangle B\in \mathcal{L}(\zzr  ^d )$ {\rm (since $A\triangle B=(A\cup B)\setminus(A\cap B$)}.
\item[sequences of sets:] if $A_n\in\mathcal{L}(\zzr ^d  )$ then $\cup _{n\geq 1}A_n$ and 
 $\cap _{n\geq 1}A_n$  both belong to $\mathcal{L}(\zzr^d)$ 
  {\rm (since $\charfun_{\cup A_n}=\sup\{\charfun_{A_n}\}$ and  
$\charfun_{\cap A_n}=\inf\{\charfun_{A_n}\}$. See also Theorems~\ref{teoCH4SIGmaLebe} and~\ref{Teo:Ch4InTeRnuMeLebEareLebe})}.
 \item[furthermore we recall:] $\zzr  ^d  $,  $\emptyset  $ both belong to $\mathcal{L}(\zzr ^d  )$ and null sets are Lebesgue measurable sets too {\rm  (Theorem~\ref{TeoCH4NULLsetsareLebe})}.
\end{description}
\end{Theorem}
  
Any subset of a null set is a null set too. So, the last statement of Theorem~\ref{TeoCH4CollectingPROP} has the following consequence:
\begin{Corollary}
\ZLA{CoroCH4LebeMISUcompleta} Any subset of a null set is Lebesgue measurable.
\end{Corollary}  
  
  We recall that an {\sc algebra}\index{algebra} is a ring with unity. The previous property show that $\mathcal{L}(\zzr  ^d )$ is an {\sc algebra of sets}\index{algebra!of sets} respect to the operations
  \[
+=\triangle\,,\qquad \cdot =\triangle
  \]
  (the unity is $\zzr  ^d $) For this reason we say that $\mathcal{L}(\zzr  ^d )$ is an {\sc algebra of sets} and, more precisely, we say that it is a {\sc $\ZSI$-algebra}\index{algebra!$\ZSI$-}\index{$\ZSI$-algebra} to indicate that it is closed under  numerable unions (and intersections) of its elements, a fact that we prove in Theorem~\ref{TeoCH4ProprieLebeMeasu} below.
  \begin{Remark}
{\rm
General or abstract study of measure theory   can be found in many books. We refer the reader to~\cite{DOOBbookMEASURE94,Halmos50MEASURE,MalliavinBOKKintegr95,Royden66LibroANALreale,Tao11Measure}.\zdia 
}
\end{Remark}

\subsection{\ZLA{sectCH4MEASUREgeneral}The Measure of Lebesgue Measurable Sets}
 
Let $A\in\mathcal{L}(\zzr ^d  )$. The function $\heaviside _A$ is quasicontinuous. We put ${B_N}=\{x\,:\ \|x\|<N\}$ and we see that
 
\begin{equation}\ZLA{eq:ch4PRELIdefiMISA}
\int_A1 \ZD x=\int _{\zzr  ^d } \heaviside_A(x)\ZD x=\lim _{N\to+\ZIN}\int _{B_N} \heaviside_A(x)\ZD x \,.
\end{equation}
    The limit~(\ref{eq:ch4PRELIdefiMISA}) exists, either a number or $+\ZIN$.
So we define:
\begin{Definition}
Let $A$ be Lebesgue measurable.
We put
\[
\zl(A)=\int_A1 \ZD x= \int _{\zzr ^d  } \heaviside_A(x)\ZD x \in [0,+\ZIN]\,.
\]
We call $\zl(A)$ the {\sc Lebesgue measure of the (Lebesgue measurable) set $A$.}\index{measure!Lebesgue}\index{set!Lebesgue measure}\index{Lebesgue!measure}
 
In the contest of the Lebesgue measure, a function which is quasicontinuous is called a {\sc (Lebesgue) measurable function.}\index{Lebesgue!measurable!function}\index{function!measurable!Lebesgue}\index{measurable!function}
\end{Definition}

From the known properties of the Lebesgue integral we get:
\begin{Theorem}\ZLA{TeoCH4ProprieLebeMeasu}
The following properties hold:
\begin{enumerate}
\item\ZLA{I1TeoCH4ProprieLebeMeasu}for any Lebesgue measurable set $A$ we have $\zl(A)\geq 0$.
\item\ZLA{I0TeoCH4ProprieLebeMeasu}if $\mathcal{O}$ is an open set and
\[
\mathcal{O}=\bigcup R_n\quad \mbox{($R_n$ pairwise almost disjoint closed rectangles)}
\]
then we have\footnote{i.e. the definition of the measure of an open set given here agree with that given in the definition~\ref{DefiCH1DeFiMiSuOpen} of Chap.~\ref{Ch1:INTElebeUNAvaria}.}
\[
\zl(\mathcal{O})=\sum_{n\geq 1}\zl(R_n)\,.
\]
\item\ZLA{I2TeoCH4ProprieLebeMeasu}monotonicity: $A\subseteq B$ (both Lebesgue measurable) we have $\zl(A)\leq\zl(B)$.
\item\ZLA{I3TeoCH4ProprieLebeMeasu}additivity:
if $A\cap B=\emptyset$ and both  $A$ and $B$ are Lebesgue measurable then
\[
\zl(A\cup B)=\zl(A)+\zl(B)\,.
\]
\item\ZLA{I4TeoCH4ProprieLebeMeasu}$\ZSI$-additivity: if $\{A_n\}$ is a sequence of \emph{pairwise disjoint} Lebesgue measurable sets then
\[
\zl(\cup A_n)=\sum\zl(A_n)\,.
\]
\item\ZLA{I5TeoCH4ProprieLebeMeasu}let $\{A_n\}$   be Lebesgue measurable sets and let us assume
\[
A_n\subseteq A_{n+1}\,.
\]
Let $A=\cup _{n\geq 1}A_n$. We have:
\[
\zl(A)=\limn \zl(A_n)\,.
\]
\item\ZLA{I6TeoCH4ProprieLebeMeasu}let $\{A_n\}$   be Lebesgue measurable sets and let us assume:

\[
\left\{\begin{array}{l}
\displaystyle \mbox{there exists a \emph{bounded} rectangle $R$ such that $A_1\subseteq R$;}
\\
 A_{n+1}\subseteq A_{n}  \,.
\end{array}\right.
\]
Let $A=\cap _{n\geq 1}A_n$. We have
\[
\zl(A)=\limn \zl(A_n)\,.
\]
\end{enumerate}
 \end{Theorem}  
 \zProof The proof is obvious. Only the following observations can be useful.
Property~\ref{I0TeoCH4ProprieLebeMeasu} is clear if the union is finite. Otherwise it follows from Beppo Levi Theorem since
\[
\charfun_{\mathcal{O}}=\lim _{N\to+\ZIN}\sum_{n=1}^N \charfun _{R_n}\qquad {\rm a.e.}\,\  x\in\zzr^d
\] 
and the sequence
\[
n\mapsto \sum_{n=1}^N \charfun _{R_n}
\] 
is increasing.

 As regard to the properties~\ref{I5TeoCH4ProprieLebeMeasu} and~\ref{I6TeoCH4ProprieLebeMeasu}: 
 in both the cases  we know $A\in\mathcal{L}(\zzr^d)$ from Theorem~\ref{TeoCH4CollectingPROP}. 
We consider the statement~\ref{I5TeoCH4ProprieLebeMeasu}. 
  The sequence of the functions $\charfun_{A_n}$ is increasing,  $\charfun_{A_n} (x)\geq 0$ for every $x$ and for every $x$
 \[
\limn  \charfun_{A_n} (x)=\charfun_A(x)\,.
 \]
Beppo Levi Theorem implies
\[
\zl(A)=\int _{\zzr  ^d } \charfun_A(x)\ZD x=\limn 
\int _{\zzr ^d  } \charfun_{A_n}(x)\ZD x=\limn \zl(A_n)\,. 
\]

We consider the statement~\ref{I6TeoCH4ProprieLebeMeasu}.  
The sequence of the functions $\charfun_{A_n}$ is decreasing,  $\charfun_{A_n} (x)\geq 0$ for every $x$ and for every $x$
 \[
\limn  \charfun_{A_n} (x)=\charfun_A(x)\,.
 \]
 Furthermore, the sequence $\{\charfun_{A_n}\}$ is \emph{bounded on a bounded rectangle.} So, from Theorem~\ref{teo:Cha1bis:LebeCAROPartIntervLIMIT} we have
 \[
\zl(A)=\int _{\zzr ^d  } \charfun_A(x)\ZD x=\limn 
\int _{\zzr ^d  } \charfun_{A_n}(x)\ZD x=\limn \zl(A_n)\,. 
\]
The proof is completed.\zdia
 
\begin{Remark}
{\rm
The boundedness assumption in the statement~\ref{I6TeoCH4ProprieLebeMeasu} is crucial. It can be slightly relaxed and we can assume that the sets $A_n$ are contained in a possibly unbounded set, provided that its measure is finite. But it cannot be completely removed as the following example shows: if $A_n=[n,+\ZIN)\subseteq \zzr$ then $\cap A_n=\emptyset $ but $\zl(A_n)=+\ZIN$ for every $n$.\zdia
}
\end{Remark}
 
Now we observe that if $A$ and $B$ are Lebesgue measurable then $A\cup B=A\cup\left (B\setminus A\right )$, the disjoint union of two Lebesgue measurable sets. 
 Properties~\ref{I2TeoCH4ProprieLebeMeasu} and~\ref{I3TeoCH4ProprieLebeMeasu} of Theorem~\ref{TeoCH4ProprieLebeMeasu}
 give: 
 \begin{Corollary}  \label{COROTeoCH4ProprieLebeMeasu}
 Let $A$ and $B$ be Lebesgue measurable. We have
\[
\zl(A\cup B)=\zl(A)+\zl(B\setminus A)\leq \zl(A)+\zl(B)\,. 
\] 
 \end{Corollary}

 Let us consider a bounded set $A\in\mathcal{L}(\zzr  ^d )$. Its characteristic function is quasicontinuous. For every $\ZEP>0$ there exists a  multirectangle $\Delta_\ZEP$ such that   
 \[
 \begin{array}{l}
 \mbox{$\mathcal{I}_{\Delta_\ZEP}$ is open and $ 
L(\Delta_\ZEP) <\ZEP$}
\\  \mbox{$(\heaviside_A)_{\zzr ^d  \setminus \mathcal{I}_{\Delta_\ZEP}}$ is continuous} \,.
\end{array}
 \]
 We denote \[
 \mathcal{O}_\ZEP=\mathcal{I}_{\Delta_\ZEP}\,.
 \]
 The set $ \mathcal{O}_\ZEP$ is open, $\zl( \mathcal{O}_\ZEP)=L(\Delta_\ZEP)<\ZEP$ and
 the set $A\cup\mathcal{O}_\ZEP $ is open. In fact, let  $x\in A\cup\mathcal{O}_\ZEP $.  If $x\in A$ is not an interior point then  $x\in \partial A$ and it is a point where $\heaviside_A$ is discontinuous so that $x\in \mathcal{O}_\ZEP $, hence it is an interior point of $A\cup\mathcal{O}_\ZEP $. 
 
 Monotonicity of the measure and Corollary~\ref{COROTeoCH4ProprieLebeMeasu} gives
 \[
\zl(A)\leq \zl(\underbrace{A\cup\mathcal{O}_\ZEP}_{\mbox{\tiny open}} )\leq \zl(A)+\zl(\mathcal{O}_\ZEP )\leq  \zl(A)+\ZEP\,. 
 \]
 At the same conclusion we arrive when $A$ is not bounded   by considering the sequence of sets $A\cap \{x\,\ \|x\|<n\}$.
 So we have the following theorem:
 \begin{Theorem}\ZLA{TEO:eq:CH4RappreMeasMISESTER-A}
 Let $A\in\mathcal{L}(\zzr ^d  )$. We have:
 \begin{equation}\ZLA{eq:eq:CH4RappreMeasMISESTER-A}
\zl(A)= \inf\{\zl(\mathcal{O})\,,\quad \mbox{$\mathcal{O}\supseteq A $ and open}\}\,. 
 \end{equation}
 \end{Theorem}
 The equality~(\ref{eq:eq:CH4RappreMeasMISESTER-A}) was already stated in
 Lemma~\ref{ch3LemmaPrimaDefiLunivDisjquasiMult}  when $A$ is a multirectangle set (Theorem~\ref{Teo:Ch1P1MeasMultirec} in dimension $1$).

Translation invariance of the integral has a reformulation in terms of the measure. Let $A\in \mathcal{L}(\zzr ^d  )$ and let $x_0\in\zzr ^d  $. We consider the translation of $A$:
 \[
x_0+A=\{x+x_0\,,\quad x\in A\}\,. 
 \]
 Then:
 \[
\charfun_{x_0+A}(x)=\charfun_A(x-x_0) 
 \]
 and we have
 \[
\zl(A+x_0)=\int _{\zzr ^d  } \charfun_{x_0+A}(x) \ZD x=\int _{\zzr  ^d }\charfun_A(x-x_0) \ZD x=\int_{\zzr ^d  }\charfun_A(x)\ZD x= \zl(A)\,.
 \]
This property is the {\sc translation invariance}\index{transaltion invariance!of the Lebesgue measure}    of the Lebesgue measure.

\begin{Remark}{\rm
Any set function
\begin{equation}\ZLA{EqCh4DefiMISUconDensi}
A\mapsto \int_A f(x)\ZD x
\end{equation}
with $f\geq 0$ has the property listed in Theorem~\ref{TeoCH4ProprieLebeMeasu} and can be considered a measure ``weighted'' by the function $f$, for example the quantity of material when $f$ denotes a density. We note the existence of set functions with the properties in theorem~\ref{TeoCH4ProprieLebeMeasu}  which cannot be represented as in~(\ref{EqCh4DefiMISUconDensi}). An example is the {\sc Dirac  measure}\index{measure!Dirac}\index{Dirac measure} 
\[
\delta(A)=\left\{\begin{array}{lll}
1 &{\rm if}& 0\in A\\
0&{\rm if}& 0\notin A\,.
\end{array}\right.
\] 
Any set function which enjoys the properties~\ref{I1TeoCH4ProprieLebeMeasu}-\ref{I3TeoCH4ProprieLebeMeasu} in Theorem~\ref{TeoCH4ProprieLebeMeasu} is a {\sc (positive) measure}\index{measure!positive} and it is a {\sc $\ZSI$-additive measure}\index{measure!$\ZSI$-additive} when also property~\ref{I4TeoCH4ProprieLebeMeasu} holds. So, the Dirac measure is indeed a $\ZSI$-additive measure.

The Dirac measure is not translation invariant.\zdia
}
\end{Remark}

 \subsection{Absolute Continuity of the Integral: the General Case}
We note:   
\begin{Theorem}\ZLA{CH4FunzCONinteNulloNULLA}
Let $f\geq 0$ be quasicontinuous  and let $A$ be a measurable set. Let
\[
\int_A f(x)\ZD x=\int_{\zzr ^d  } \charfun_A(x)f(x)\ZD x=0\,.
\]
The function $f$ is a.e. zero on $A$, i.e. $\charfun_A f $ is a.e. zero on $\zzr  ^d $.
\end{Theorem} 
\zProof  We proceed by contradiction and we prove that if $f$ is positive on a set of positive measure then its integral is positive.   Let
\[
A_+=\{x\,:\  f(x)>0\}=\bigcup _{n\geq 1} A_{+,n}\quad \mbox{where}\quad A_{+,n}=
\{x\,:\ \ f(x)>1/n\}\,.
\]
If $\zl(A_+)>0$ then\footnote{see statement~\ref{I5TeoCH4ProprieLebeMeasu} of Theorem~\ref{TeoCH4ProprieLebeMeasu}.} there exists $n_0$ such that $\zl(A_{+,n_0})>0$. Monotonicity of the integral gives
\[
\int_{A} f(x)\ZD x\geq \int_{A_{+,n_0}} f(x)\ZD x\geq \dfrac{1}{n_0}\zl(A_{+,n_0})\,.\zdiaform
\]
 
 We recall that the integral of a function whose support is a null set is equal zero (property~\ref{I3TEO:Ch2InteNULLOfunzQOnulla} of Theorem~\ref{TEO:Ch2InteNULLOfunzQOnulla}).   We combine this fact with Theorem~\ref{CH4FunzCONinteNulloNULLA}
   and we get the following result which justifies the term ``null set'':
 \begin{Theorem}
 The set $N$ is a null set if and only if $\zl(N)=0$\,.
 \end{Theorem}
 

A related and important result is a  
  consequence of Theorem~\ref{TEO:eq:CH4RappreMeasMISESTER-A}. This result extends {\sc absolute continuity of the integral,}\index{set function!absolute continuity}\index{absolute continuity}  already proved for open sets in Theorem~\ref{Teo:CH3ABSOluContiPERiSoliAperti} (and in Theorem~\ref{Teo:Cha1bisABSOluContiPERiSoliAperti} when $d=1$), to  
 Lebesgue measurable sets   $A$:  
 
 \begin{Theorem}  \ZLA{TeoCH4ABsolContINgene}
  Let $f $ be summable on $\zzr  ^d $. For every $\ZEP>0$ there exists $\ZDE>0$ such that 
 \[ 
A\in \mathcal{L}(\zzr  ^d )\,,\quad \zl(A)\leq\ZDE\ \implies\ 
 \int_A |f(x)|\ZD x<\ZEP\,.
 \]
 \end{Theorem}  
 \zProof 
  Theorem~\ref{Teo:CH3ABSOluContiPERiSoliAperti} 
states that for every $\ZEP>0$ there exists $\ZSI>0$ such that

\[
\int_{\mathcal{O}} |f(x)|\ZD x<\ZEP \quad {\rm if}\quad \left\{
\begin{array}{l}
\mbox{$\mathcal{O}$ is  open}\\
\zl(\mathcal{O})<\ZSI\,.
\end{array}\right.
\]  
 
Let $A\in\mathcal{L}(\zzr  ^d )$ satisfy
\[
\zl(A)<\ZDE=\ZSI/2\,.
\]
Theorem~\ref{TEO:eq:CH4RappreMeasMISESTER-A} shows the existence  of an open set $\mathcal{O}$ such that
\[
A\subseteq\mathcal{O}\,,\qquad \zl(\mathcal{O})<\ZSI
\]
so that
\[
\int_A |f(x)|\ZD x\leq  \int_{\mathcal{O}} |f(x)|\ZD x<\ZEP\,.\zdiaform
\]

\section{\ZLA{Ch4SECTbor3elLEBE}Borel Sets and Lebesgue Sets}

We begun this chapter with Theorem~\ref{CH4:TeoPRELIcontFUN} which shows the relations between   continuous functions and  open or closed sets. The family of open and closed sets is  not an algebra of sets.  In fact in general the difference of two (open or closed) sets is neither open nor closed. We can construct a  $\ZSI$-algebra by taking any  numerable union or intersection of open and closed sets. This way we obtain, among all the  $\ZSI$-algebras of subsets of $\zzr  ^d $, the smallest one which contains all the open sets and all the closed sets.
This  $\ZSI$-algebra is denoted 
$\mathcal{B}(\zzr ^d )$\index{$\mathcal{B}(\zzr^d)$} (or simply  
  $\mathcal{B}$\index{$\mathcal{B}$})   and it is the  $\ZSI$-algebra of the {\sc Borel   sets.}\index{set!Borel}\index{algebra!Borel} A Borel set is also called a {\sc Borel measurable set.}\index{set!measurable!Borel}\index{measurable!set (Borel)}\index{Borel!measurable set}

It is clear that 
\[
\mathcal{B}(\zzr  ^d )\subseteq\mathcal{L}(\zzr  ^d ) \,.
\]
because any open or closed set  is Lebesgue measurable and $\mathcal{B}$ is the \emph{smallest} $\ZSI$-algebra which contains   open and closed sets.

We observe\footnote{we recall that ~$\tilde{}$~ denotes the complement.}:
\begin{equation}\ZLA{Eq:ch4ProInsieFunzINV}
\left\{\begin{array}{l}
 f^{-1}(\cup A_n) =\cup f^{-1}(A_n)\,,\qquad f^{-1}(\cap A_n) =\cap f^{-1}(A_n)\\
  A=f^{-1}(B)\ \implies \tilde A=f^{-1}(\tilde B)\quad\mbox{(provided that ${\rm dom}\, f=\zzr ^d  $)} \,.
\end{array}\right.
\end{equation}
 
 We combine the equalities~(\ref{Eq:ch4ProInsieFunzINV}) and Theorem~\ref{CH4:TeoPRELIcontFUN}. We get:

\begin{Theorem}\ZLA{TeoCH4ContrimmBOrelianiBIS}
Let $f$: $\zzr ^d\mapsto \zzr ^m $ be continuous. If $A\in\mathcal{B}(\zzr^m)$     then $f^{-1}(A)\in\mathcal{B}(\zzr^d)$.
\end{Theorem}

Even more:
\begin{Theorem}\ZLA{TeoCH4ContrimmBOrelianiTER}
We assume:
\begin{enumerate}
\item   $B\in\mathcal{B}(\zzr^d)$;
\item  $K$ is a closed set of $\zzr^m$;
\item   $f$: $B\mapsto \zzr^m$ is continuous.
\end{enumerate}
Under these conditions, the set 
\[
f^{-1}(K)=\{x\in B\,:\ f(x)\in K\}
\]
is a Borel set.
\end{Theorem}
\zProof We know from Theorem~\ref{CH4:TeoPRELIcontFUN} that $f^{-1}(K)$ is  relatively closed in $B$, i.e.  we know that $f^{-1}(K)=B\cap C$ where $C$ is a closed set. The set $B\cap C$ is the intersection of two Borel set, and it is a Borel set.\zdia
 
It is a fact, to be seen in Appendix~\ref{APPEch4NullNONborel}, that \emph{there exists Lebesgue measurable sets which are not Borel sets, i.e. the inclusion $\mathcal{B}\subseteq \mathcal{L}$ is strict.} In spite of this, the two $\ZSI$-algebras $\mathcal{B}$ and $\mathcal{L}$ are closely related: 
 
\begin{Theorem}
Let $A\in\mathcal{L}= \mathcal{L}(\zzr ^d )$. there exist a null set $ N $ and a set $B   \in\mathcal{B}(\zzr ^d  )$ such that
\[
A=B\cup N\,.
\]
\end{Theorem}
\zProof  
We can confine ourselves to prove the theorem in the case that $A$ is bounded.

The proof is by iteration. So, it is convenient to rename $A_0$ the set $A$ and $\heaviside_0$ its characteristic function.
The assumption is that $A_0\in\mathcal{L}(\zzr ^d  )$ so that $\heaviside_0$ is quasicontinuous on $\zzr ^d  $: there exists an open set  $\mathcal{O}_1$ such that
\[
\zl(\mathcal{O}_1)<1\,,\qquad( \heaviside_0)_{|_{\zzr ^d  \setminus \mathcal{O}_1}}
\quad \mbox{is continuous}\,.
\]
 We have:
 \[
A_0=\underbrace{\left (A_0\setminus\mathcal{O}_1\right )}_{=B_0\in\mathcal{B}}\cup \underbrace{\left (A_0\cap\mathcal{O}_1\right )}_{\stackrel{=A_1\in\mathcal{L}}{\zl(A_1)\leq\zl(\mathcal{O}_1) <1}   }\,.
 \] 
 Note that $B_0\in\mathcal{B}$ from Theorem~\ref{TeoCH4ContrimmBOrelianiTER} since 
 \[
 B_0= \left (
  {
  (
  \charfun_0)_{
  |_{\zzr ^d  \setminus\mathcal{O}_1
  }
  }
}   
\right )
^{-1}  
  (1)  
 \]
 and $ (\charfun_0)_{
  |_{\zzr ^d  \setminus\mathcal{O}_1
  }}$ is a continuous function defined on a closed, hence a Borel, set. 
  
We repeat this construction for the set $A_1$ buth with an open set $\mathcal{O}_2$ such that $\zl(\mathcal{O}_2)<1/2$:
\[
A_1=\underbrace{\left (A_1\setminus\mathcal{O}_2\right )}_{=B_1\in\mathcal{B}}\cup \underbrace{\left (A_1\cap\mathcal{O}_2\right )}_{\stackrel{=A_2\in\mathcal{L}}{\zl(A_2)\leq\zl(\mathcal{O}_2)<1/2}   }\,.
 \]
So we have
\[
 A_0=\underbrace{B_0\cup B_1
}_{\in\mathcal{B}}\cup  
\underbrace{\left (A_0\cap\left (\mathcal{O}_1\cap \mathcal{O}_2 \right ) \right )}_{\stackrel{=A_2\in\mathcal{L}}{\zl(A_2)\leq\zl(\mathcal{O}_2)<1/2}   }\,.
\]
We iterate and we get 
\[
A=A_0=\underbrace{\left [ \bigcup_{k\geq 1} B_k\right ]
}_{\in\mathcal{B}}
\bigcup  
\underbrace{\left (A_0\cap \left ( \cap_{k\geq 1} \mathcal{O}_k\right )  \right )}_{\stackrel{=N\in\mathcal{L}}{\zl(N)=0}   }\,.
\]
The fact that    $N$ is a null set follows from Theorem~\ref{TEO:eq:CH4RappreMeasMISESTER-A} since  for every $k$ we have
\[
\mbox{
$N\subseteq \mathcal{O}_k$,  \ \
  $ \mathcal{O}_k$ open and $\zl (\mathcal{O}_k)<1/k$ for every $k$}\,.\zdiaform
\]

 In conclusion, any Lebesgue  set is ``almost'' a  Borel set: the difference is a null set. 

Of course, to every Borel set we can associate its Lebesgue measure: the restriction of the Lebesgue measure to $\mathcal{B}$ is the {\sc Borel measure.}\index{measure!Borel}\index{set!Borel measure}\index{Borel!measure}
\emph{We note an important difference between the $\ZSI$-algebras of Lebesgue and Borel: the $\ZSI$-algebra of Lebesgue contains any subset of a null set (see Corollary~\ref{CoroCH4LebeMISUcompleta})
while there exists \emph{null sets which are Borel sets} and which contains subsets which do not belong to $\mathcal{B}$. This is seen in Appendix~\ref{APPEch4NullNONborel}.} 

By definition, a $\ZSI$-additive measure  is  {\sc complete}\index{measure!complete} 
 when every subset of a null set is measurable.  So, $\mathcal{L}(\zzr^d)$ is complete while $\mathcal{B}(\zzr^d)$ is not.

  \subsection{\ZLA{secCAP4:multirElebeMis}Multirectangle Sets and Lebesgue Measurable Sets}
  In the table~\ref{tableCH3ADDITIONALshortNOTATIONS} of sect~\ref{CH3MoreonMULTIREC} we defined the multirectangle sets: a multirectangle set is a set $\mathcal{J}_\Delta$ where $\Delta$ is almost disjoint. So, a multirectangle set is   a Borel set, hence a Lebesgue set. Its measure
  can be computed as in~(\ref{eq:CH3RappreMeasMISESTER-A})
and this formula extends to every Lebesgue measurable set, see~(\ref{eq:eq:CH4RappreMeasMISESTER-A}).  The interpretation is that  \emph{ a Lebesgue measurable set can be  approximated \emph{from outside} by open  multirectangles.}

In the special case of the multirectangle sets: 
\emph{any multirectangle set can be approximate \emph{from inside} by    disjoint  multirectangles. } This is the interpretation of Remark~\ref{P2Ch2RemaModiDefiMultRwSET},
 Theorem~\ref{TEOeq:CH3RappreMeasMISESTER-C} and formula~(\ref{eq:CH3RappreMeasMISESTER-C}).  

 The goal of this section is to show the existence of Lebesgue measurable sets (even of Borel sets) which are not multirectangle sets and whose measure \emph{cannot} be computed from inside by using   formula~(\ref{eq:CH3RappreMeasMISESTER-C})  \emph{not even we allow $ \Delta_{\rm ins}$ to be composed by infinitely many rectangles.} 
This goal is achieved by constructing an example in dimension~$1$. The construction is similar to that of the standard Cantor set and we present the two constructions in parallel. Moreover, \emph{we shall see that the standard Cantor set is a null set which is not numerable.
This observation confirms the statement  in item~\ref{I3Exe:Ch0Qnullset1} of Remark~\ref{REMA:Ch0Qnullset}.}

It is convenient  to state first the following {\bf preliminary observations:}

  \begin{enumerate}
  \item\ZLA{I2secCAP4:multirElebeMis}
  Representation of the numbers of the interval $[0,1]$ in the base $k$. A number $x\in [0,1]$ is represented by a sequence  $\{c_n\}$ of nonnegative integers, where $c_n$ are the numerators which appear in the equality
  \[
x=\sum _{n=1}^{+\ZIN}   \frac{c_n}{k^n}\,.
  \]
  
  Any element $x\in [0,1]$ has a unique representation of this form unless $x=c/k^j$ with $j$ and  $c$ in $ \mathbb{N}$ (and $c\leq k^j$) at least for large enough $j$. In this case the number is represented by two sequences 
whose terms are  (for    large enough $j$)\footnote{here $\ZDE_{n,j}$ is the Kronecker delta: $\ZDE_{j,j}=1$ and $\ZDE_{n,j}=0$ if $n\neq j$.}:
  \[
\mbox{the sequence $\left \{c \ZDE_{n,j}  \right \}$ and the stationary sequence   $\{c(k-1) \}$}\,.
  \]
 In order to describe the sets we are going to construct, it is convenient to discard the first representation and to keep the second one. 
\emph{ 
 This way any element of $[0,1]$ admits a unique representation.}  
 For example the number $1$ is represented by the stationary sequence $\{ k-1\}$ while the number zero is represented by the stationary sequence $0$.
  
    \item\ZLA{I1secCAP4:multirElebeMis} Let $k$, $k_1$ and $k_2$  be   fixed positive integers. We represent the elemets of $[0,1]$ in base $k$. Any open interval $(a,b)\subseteq [0,1]$ contains infinitely many numbers of the form $ x=c/k^n$ where $c$ and $n$ are positive integers and $c\neq k_1$, $c\neq k_2$. The reason is that   $c/k^n-c/k^{n+1}<(b-a)/2$ if $n$ is sufficiently large.
  \end{enumerate}
  
 Now we construct the two sets. The first one is the usual ternary Cantor set. The elements of this set are best represented if we choose $3$ as the base. The second set is a modification of the Cantor set first constructed by  K.~Smith in~\cite{Smith1874} and later by V.~Volterra in~\cite{VolterraJmatBatta1881} (both these interesting papers are available on line). Correspondingly, we denote $K_C$ and $K_{SV} $ the two sets.
 
 We present the two constructions   in the following tables. The left column is for the set $K_C$ and the right one for $K_{SV}$.

We proceed with the following steps:
\bigskip

\begin{center}
\begin{tabular}{||l||c||}    
\hline\hline  
&
\\              
{\bf Step~0}& \parbox{4.3in}{
 \begin{tabular}{l|l}
  \parbox{2in}{
   \begin{center}\fbox{The construction of $K_C$}   \end{center}
The numbers are represented in base $3$. The interval $[0,1]$ is
divided in $3$ equal intervals.} 
& 
  \parbox{2in}{
   \begin{center}\fbox{The construction of $K_{SV}$}   \end{center}  
  The numbers are represented in base $8$. The interval $[0,1]$ is
divided in $8$ equal intervals.}
 \end{tabular} 
 }   \\
 & \\
 \hline
\hline
 \end{tabular}

\medskip

\begin{tabular}{||l||c||}    
\hline\hline   
&
\\              
{\bf Step~$\bf 1$}& \parbox{4.3in}{
 \begin{tabular}{l|l}
  \parbox{2in}{
$C_1$ is the open middle interval of $[0,1]$  of length $1/3$
and we put   $D_1=[0,1]\setminus C_1$. We have  $\zl(C_1)=1/3$ and  $D_1$ is the union of two closed intervals. 

 {\bf The representation of the ele\-ments of $D_1$ in base $3$ is $\sum d_n/3^n$, $d_1\neq 1$.}
} 
& 
  \parbox{2in}{
$C_1$ is   the open middle interval  of $[0,1]$, of length $1/4$
and we put   $D_1=[0,1]\setminus C_1$ We have  $\zl(C_1)=1/4$ 
  and $D_1=[0,1]\setminus C_1$ is the union of two closed intervals. 
 
   {\bf The representation of the e\-lements of $D_1$ in base $8$ is $\sum d_n/8^n$, $d_1\notin\{ 3,\, 4\}$.}
  }
 
 \end{tabular} 
 }   \\
 & \\
 \hline
\hline
 \end{tabular}
\begin{tabular}{||l||c||}    
\hline\hline   
&
\\              
{\bf Step~$\bf 2$}& \parbox{4.3in}{
 \begin{tabular}{l|l}
  \parbox{2in}{
$C_2$ is the union of the open middle intervals of diameter $1/3^2$ of the two  which compose $D_1$ and we put   $D_2=D_1\setminus C_2$ ($D_2$ is the union of $2^2$ closed intervals). We have  $\zl(C_2)=2/3^2$. 
 
  {\bf The representation of the elements of $D_2$ in base $3$ is $\sum d_n/3^n$, $d_n\neq 1$ if $n\leq 2$.}
} 
& 
  \parbox{2in}{
$C_2$ is the union of the open middle intervals
of diameter $1/4^2$ of the two which compose
  $D_1$ 
and we put   $D_2=D_1\setminus C_2$ ($D_2$ is the union of $2^2$ closed intervals). We have  $\zl(C_2)=2/4^2$

  {\bf  The representation of the elements of $D_2$ in base $8$ is $\sum d_n/8^n$, $d_n\notin\{ 3,\, 4\}$ if $n\leq 2$.}
}
 \end{tabular} 
 }   \\
 & \\
 \hline
\hline
 \end{tabular}

\medskip

\begin{tabular}{||l||c||}    
\hline\hline   
&
\\              
{\bf Step~$\bf 3$}& \parbox{4.3in}{
 \begin{tabular}{l|l}
  \parbox{2in}{
$C_ 3$ is the union of the open middle intervals
 of length $1/3^3$ of those which compose
 of $D_2$ 
and we put   $D_3=D_2\setminus C_3$. We have  $\zl(C_3)=2^2/3^3$ 
and $D_3$ is the union of $2^3$ closed intervals.

  {\bf The representation of the elements of $D_3$ in base $3$ is $\sum d_n/3^n$, $d_n\neq 1$ if $n\leq 3$.} 
 } 
& 
  \parbox{2in}{
$C_3$ is the union of the open middle intervals
of length $1/4^3$ of  of those which compose
  $D_2$ 
 and we put   $D_3=D_2\setminus C_3$. We have  $\zl(C_3)=2^2/4^3$ 
 and $D_3$ is the union of $2^3$ closed intervals.

  {\bf The representation of the elements of $D_3$ in base $8$ is $\sum d_n/8^n$, $d_n\notin\{ 3,\, 4\}$ if $n\leq 3$.}
}
 \end{tabular} 
 }   \\
 & \\
 \hline
\hline
 \end{tabular}

 \medskip

\begin{tabular}{||l||c||}    
\hline\hline   
&
\\              
{\bf Step~$\bf j$}& \parbox{4.3in}{
 \begin{tabular}{l|l}
  \parbox{2in}{
$C_ {j}$ is the union of the open middle intervals
of diameter $1/3^j$ of each one which compose 
 of $D_ {j-1}$ 
and we put  $D_j=D_{j-1}\setminus C_{j}$. We have: $\zl(C_{j})=2^{j-1}/3^j$ 
and $D_j$ is the union of $2^j$ closed intervals.

  {\bf The representation of the elements of $D_j$ in base $3$ is $\sum d_n/3^n$, $d_n\neq 1$ if $n\leq j$.} 
 } 
& 
  \parbox{2in}{
$C_ {j}$ is the union of the open middle intervals
of diameter $1/4^j$ of each one which compose 
 of $D_ {j-1}$ 
and we put   $D_j=D_{j-1}\setminus C_{j}$. We have: $\zl(C_{j})=2^{j-1}/4^j$ 
and $D_j$ is the union of $2^j$ closed intervals.

  {\bf the representation of the elements of $D_j$ in base $8$ is $\sum d_n/8^n$, $d_n\notin\{ 3,\, 4\}$ if $n\leq j$.}
}
 \end{tabular} 
 }   \\
 & \\
 \hline
\hline
 \end{tabular}
\end{center}

\bigskip
 
We denote $C_C$, respectively $C_{SV}$, the open sets $\cup C_n$ of the two constructions and
\[
K_C=[0,1]\setminus C_C\,,\qquad K_{SV}=[0,1]\setminus C_{SV}\,.
\]

These sets $C_C$, $C_{SV}$, $K_C$ and $K_{SV}$  are Borel sets and $K_C$ is the standard Cantor set.

The set $K_C$ is the set of those numbers of $[0,1]$ whose representation in the base $3$ does not contains $1$ while the elements of $K_{SV}$ are characterized by the fact that their representation in base $8$ does not contains neither $3$ nor $4$ (and are contained in $[0,1]$).

The sets $C_C$ and $C_{SV} $ are numerable unions of non overlapping intervals, so that their measure can be computed by adding the measure of the single intervals. We have:
\begin{align*}
&\zl(C_C)=1\quad \mbox{so that}\quad \zl(K_C)=0\\
&\zl(C_{SV})=\dfrac{1}{2}\quad \mbox{so that}\quad \zl(K_{SV})=\dfrac{1}{2}\,.
\end{align*}
The set $K_{SV}$ is not numerable since its measure is positive. Also $ K_C$
is not numerable\footnote{a well known fact, easily seen thanks to
the fact that any of its element is represented by a sequence $\{c_n\}$ with $c_n$ equal either to $0$ or to $2$. So, the function
  $x=\{c_n\}\mapsto (c_n/2)$: $K_C\mapsto [0,1]$, with the elements of $[0,1]$ represented in base $2$, is surjective.}. 
\emph{So, $K_C$ is an example of a null set which is not numerable.}

Thanks to the statement~\ref{I1secCAP4:multirElebeMis} above, neither $K_C$ nor $K_{SV}$ contains nondegenerate intervals. In particular, the multiintervals $\Delta$ such that $J_\Delta\subseteq K_{SV}$ have the form
$\Delta=\{[q_k,q_k]\}$ and $L(\Delta)=0$. Formula~(\ref{eq:CH3RappreMeasMISESTER-C}) does not hold for the Borel set $K_{SV} $: in general, \emph{a Lebesgue (or Borel) measurable set \emph{cannot} be approximated from inside---in the sense of the measure---by almost disjoint multirectangles (or, in dimension $d=1$, by   disjoint multiintervals).}
 
A result on the approximation of measurable sets from inside is Theorem~\ref{TeoFINpa3Ch1ApproInte} below.
 
 \begin{Remark}\ZLA{RemaCH4InsOpeFrontPOSI}
 {\rm
  The statement in the {\bf Preliminary observation}~\ref{I1secCAP4:multirElebeMis} implies also that $K_C=\partial C_C$ and $K_{SV}=\partial C_{SV} $. \emph{This last equality shows that  the boundary of the open set $C_{SV}$ has positive measure. }

This observation explains also the error in Remark~\ref{Ch1RemeCONerrore}: the points of discontinuity of the function $\charfun_{\mathcal{O}}$ in this remark are not only the points $a_n$  
  and $b_n$. The function $\charfun_{\mathcal{O}}$ is discontinuous on the boundary of the open set ${\mathcal{O}}$ and in general the boundary of an open set is not a null set.\zdia 
 }\end{Remark}

\section{Littlewood's Three Principles}
We recapitulate the three main properties we have seen:
\begin{enumerate}
\item
Theorem~\ref{TEO:eq:CH4RappreMeasMISESTER-A} shows that \emph{any Lebesgue measurable set is ``close''---from the point of view of the measure---to an open set} in the sense that Lebesgue measurable sets can be ``approximated'' (from   outside) by open sets.
\item
Egorov-Severini Theorem asserts that any pointwise convergent sequence of  quasicontinuous (i.e. Lebesgue measurable) functions is ``close'' to being uniformly convergent. 
\item
From the very beginning of our treatment, quasicontinuous functions are ``close'' to being continuous.
\end{enumerate}

These three informal statements are the {\sc ``Littlewood's three  principles''.}\index{Littlewood's three principles} It is useful to keep    them in mind when working with Lebesgue integral.

A recent analysis of the Littlewood's three  principles is in~\cite{16LittlewPrincMAGNANINI}.

We noted that the usual route to the integral is the converse way around. First the class  of Lebesgue measurable sets   is defined and studied.
The Lebesgue measurable functions are then defined as those functions 
with the property that $f^{-1}(I)$ is Lebesgue measurable for every open set $I\subseteq \zzr$. 

Finally the integral and its properties is studied.

When going this way, the relation of Lebesgue measurable functions and continuous function has to be separately proved: it has to be proved that any Lebesgue measurable functions, defined in this way, is
 quasicontinuous. This statement is {\sc Lusin Theorem}\index{Theorem!Lusin} first proved in~\cite{Lusin1912}.

 \section{Lebesgue Definition of Measurable Sets}
 In order to conclude this presentation, we investigate whether Lebesgue measurable sets can be ``approximated'' not only from outside, as stated by Theorem~\ref{TEO:eq:CH4RappreMeasMISESTER-A}, but also ``from inside''. 
 The response is positive and gives a characterization of Lebesgue measurable sets which is precisely the way measurable sets where originally defined by Lebesgue.
 
 We note that 
closed sets are measurable sets and now we prove that closed sets can be used to approximate a measurable set from inside:
\begin{Theorem}\ZLA{TeoFINpa3Ch1ApproInte}
\ZLA{CH4TEOapproMeasDAinteCHIUSI}
Let $A$ be a bounded Lebesgue measurable set.
We have:
\begin{equation} 
\ZLA{EQCH4TEOapproMeasDAinteCHIUSI}
\zl(A)=\sup\{\zl(K)\,,\ \mbox{$K$ compact subset of $A$}\}\,.
\end{equation}
\end{Theorem}
\zProof
Monotonicity of the measure implies
\begin{equation}
\ZLA{EQCH4TEOapproMeasDAinteCHIUSIP0}
\zl(A)\geq \sup\{\zl(K)\,,\ \mbox{$K$ compact subset of $A$}\}\,.
\end{equation}
We prove that the equality cannot be strict by finding a sequence $\{K_n\}$ of 
compact subsets of $A$ such that 
\[
\limn \zl(K_n)=\zl(A)\,.
\]

 Let $R$ be a bounded closed rectangle such that $A\subseteq R$. The   set $R\setminus A$  is measurable as the difference of two measurable 
 sets and
 \begin{equation}
 \ZLA{EQCH4TEOapproMeasDAinteCHIUSIP1}
 \zl(R)=\zl(R\setminus A)+\zl(A)\,.
 \end{equation}
We use~(\ref{eq:eq:CH4RappreMeasMISESTER-A}): there exists a sequence $\mathcal{O}_n$ of open sets such that
 \[
R\setminus A\subseteq \mathcal{O}_n\,,\quad \zl( \mathcal{O}_n)\to \zl( R\setminus A)\,.
 \] 
 Note that $K_n=R\setminus \mathcal{O }_n\subseteq A$ and that
 $K_n $ is a closed set. We have:  
 \[R=
  \mathcal{O }_n\bigcup K_n\quad \mbox{(disjoint union of measurable sets)}
 \]
 so that
 \[
\zl(R)=\zl(   \mathcal{O }_n)+\zl(K_n)\,.
 \]
The limit for $n\to+\ZIN$ gives
 \[
 \zl(R)=  \zl( R\setminus A)+\lim_{n\to+\ZIN }\zl(K_n)\,.
 \]
 We compare with~(\ref{EQCH4TEOapproMeasDAinteCHIUSIP1}) and we see
 \begin{equation}\ZLA{EQCH4TEOapproMeasDAinteCHIUSIP2}
\zl(A)= \lim_{n\to+\ZIN }\zl(K_n) 
 \end{equation}
 as we wished to achieve.\zdia
 
 This result on the ``approximation from inside'' does not have an intuitive appeal since Cantor set shows that closed sets have a complex structure. But, it suggests the following characterization of Lebesgue measurable sets:
 \begin{Theorem}
 Let $A$ be a bounded subset of $\zzr^d$, $A\subseteq R$ where $R$ is a bounded rectangle. The set $A $ is Lebesgue measurable if and only if
 \begin{equation}\ZLA{Eq:P3CH1ApproSIAinteCHEeste}
 \sup\{\zl(K)\,,\ \mbox{$K$ compact subset of $A$}\}=\inf\{\zl(\mathcal{O})\,,\quad \mbox{$\mathcal{O}\supseteq A $ and open}\}\,. 
 \end{equation}
 \end{Theorem}
 \zProof If $A$ is Lebesgue measurable then the equality~(\ref{Eq:P3CH1ApproSIAinteCHEeste}) follows from~(\ref{eq:eq:CH4RappreMeasMISESTER-A}) and~(\ref{EQCH4TEOapproMeasDAinteCHIUSI}). Conversely, we prove that~(\ref{Eq:P3CH1ApproSIAinteCHEeste}) implies that $\charfun_A$ is quasicontinuous and so that $A$ is Lebesgue measurable.
 
 We observe:

\medskip 
\begin{center} 
  \fbox{\parbox{2.5in}{
  Let $\{K_n\}$ be a sequence of compact subsets of $A$ such that
  \begin{multline*}
  \lim _{k\to+\ZIN}\zl(K_n) = \sup\{\zl(K)\,,\\ \mbox{$K$ compact subset of $A$}\}
  \end{multline*}
 By replacing $K_n$ with $\cup _{r=1}^n K_r$ we can assume that 
 \[
K_{n}\subseteq K_{n+1}  
 \]
 and so
 \begin{enumerate}
 \item 
  the numerical sequence  $\{\zl(K_n) \}$ is increasing.
  \item the sequence of the characteristic functions $\{ \charfun _{K_n} \}$ is increasing. Hence the following limit exists for every $x$:
   \[
  f(x)=\limn \charfun _{K_n}(x) \,.
   \]
    \end{enumerate}
 }}  
  \fbox{\parbox{2.5in}{
 Let $\{\mathcal{O}_n\}$ be a sequence of open sets which contain  $A$ and such that
 \begin{multline*}
 \lim _{k\to+\ZIN}\zl(\mathcal{O}_n) =\inf\{\zl(\mathcal{O})\,,\\
  \mbox{$\mathcal{O}\supseteq A $ and open}\}\,. 
 \end{multline*}
 By replacing $\mathcal{O}_n$ with $\cap _{r=1}^n \mathcal{O}_r$ we can assume that 
 \[
\mathcal{O}_{n+1} \subseteq \mathcal{O}_n
 \]
 and so:
\begin{enumerate}
\item 
  the numerical sequence $\{\zl(\mathcal{O}_n)\}$ is decreasing.
  \item
  the sequence of the characteristic functions $\{\charfun _{\mathcal{O}_n}\}$ is decreasing. Hence the following limit exists for every $x$:
  \[
 g(x)=\limn\charfun _{\mathcal{O}_n}(x) \,.
  \]
  \end{enumerate}
 }}  
  \end{center}
 
 \bigskip
 
 It is clear that 
 \[
f(x)\leq\charfun_A(x)\leq g(x) \,.
 \]

Measurability of open and closed sets and Egorov-Severini Theorem imply that 
the functions $f$ and $g$ are quasicontinuous since they are the limit of sequences of quasicontinuous functions.

Measurability of $A$ follows since now we prove $\charfun_A=f=g$ a.e. $x\in R$, so that the characteristic function $\charfun_A$ is quasicontinuous too.

We use Theorem~\ref{teo:ch3:LebeCAROPartIntervLIMIT} and we exchange  the limit and the integral:
 \[
\limn\left [\zl(\mathcal{O}_n)-\zl(K_n)\right ]=\limn \int_R \left [\charfun _{\mathcal{O}_n}(x)- 
\charfun _{K_n}(x)\right ]\ZD x  =\int_R [g(x)-f(x)]\ZD x \,. 
 \]
 But,
 \[
\limn\left [ \zl(\mathcal{O}_n)-\zl(K_n) \right ]=0
 \]
 so that the integral of the \emph{nonnegative} function  $g-f$ is zero:
 \[
 \int_R [g(x)-f(x)]\ZD x=0\quad\mbox{so that  $f(x)=g(x)$ a.e. $x\in R$ (see Theorem~\ref{CH4FunzCONinteNulloNULLA})} 
 \] 
 and 
 
 \[
\charfun_A(x)=  f(x)=g(x) \ \mbox{a.e. $ x\in R $; i.e. $\charfun_A$ is quasicontinuous}
 \]
 as we wanted.\zdia

 Now we can explain the original definition of the measure as given by Lebesgue in its thesis: in its essence, it is the characterization~(\ref{Eq:P3CH1ApproSIAinteCHEeste}) taken as a definition. Lebesgue proceeds with the following steps to define the measure of a bounded set $E$. We recast Lebesgue terminology in the form we have used up to now. In particular we note that Lebesgue prefers to consider as a basic ``bricks'' of its construction  not the rectangles but the triangles.
 \begin{enumerate}
 \item  we fix any bounded rectangle\footnote{in fact,  Lebesgue uses a triangle.}  $R\supseteq E$.  $\zl(R)$ is the standard ``volume'' (i.e. length, area, 
 volume,\dots as we defined in the Chapt.s~\ref{Ch1:INTElebeUNAvaria} and~\ref{Chap:2quasicontPIUvar}) of the rectangle. The number $\zl(R)$ does not depend on the topological properties of $R$ and so we can assume that $R$ is closed.
 \item in this second step Lebesgue defines the {\sc exterior measure}\index{measure!exterior} of the set $E$ as follows
 \[
m_e(E)=\inf\left\{ L(\Delta)\,,\quad  
\mbox{$\Delta$ almost disjoint and }\quad  
E\subseteq \mathcal{I}_\Delta
 \right\}
\] 
 By using the characterization of the open sets in Theorem~\ref{TeoCH3STRUttAPERTI}, we can recast this definition as follows:
 \begin{multline}
 \ZLA{eq:P3Ch1DEFImisUSE}
 m_e(E)=\inf
  \{m_e(\mathcal{O})\,,\quad
    \mbox{
  $\mathcal{O}$
    open and  $
     E\subseteq
    \mathcal{O}
    $
   }
  \}\\
 =\inf
  \{\zl(\mathcal{O})\,,\quad
    \mbox{
  $\mathcal{O}$
    open and  $
     E\subseteq
    \mathcal{O}
    $
   }
  \}\,.
 \end{multline}
 \item In the third step, the {\sc interior measure}\index{measure!interior} is defined as follows:
 \[
m_i(E)=\zl(R)-m_e(\tilde E) 
 \]
 where $\tilde E$ denotes the complement of $E$ respect to $R$,
 \[
\tilde E=R\setminus E\,. 
 \]
We elaborate on this definition:
\begin{multline*}
m_i(E)=\zl(R)-
\inf \{
\zl(\mathcal{O})\,,\quad \,
 \mbox{$\mathcal{O}$ open and $\tilde E\subseteq\mathcal{O}$}
\}\\
=
 \zl(R)+\sup  \{-\zl(\mathcal{O})\,,\quad \,\mbox{$\mathcal{O}$ open and $\tilde E\subseteq\mathcal{O}$}\}\\
 =
 \sup 
  \{
 \zl(R)-\zl(\mathcal{O})  
\,,\ \  \mbox{$\mathcal{O}$ open and $\tilde E\subseteq\mathcal{O}$}\}\\
=
\sup\{\zl(R\setminus\mathcal{O})\,,\ \mathcal{O}\ \mbox{open and $\tilde E\subseteq \mathcal{O}$}\}
\,.
\end{multline*}
Note that $K=R\setminus \mathcal{O}\subseteq E$ is compact.
So, this chain of equalities suggests   putting\footnote{consequence of the additivity of the measure, but in the process of the definition of the measure used by Lebesgue additivity of the measure is not yet proved at this stage, and this equality is taken as  the definition of the measure of a compact set.}
\[
\zl(K)=\zl(R\setminus \tilde K)=\zl(R)-\zl(\tilde K) 
\]
for every compact subset  $K$ of $R$.
 So,
 \[
m_i(E)=\sup\{ \zl (K)\quad \mbox{ K compact and $K\subseteq E$              }  \} \,.
 \]
 
 \item the final step is the definition of the (Lebesgue) measureble sets: \emph{the set $E$ is measurable when $m_i(E)=m_e(E)$ and Lebesgue definition of the measure is}
 \[
\zl(E)=m_i(E)=m_e(E)\,. 
 \]
 \end{enumerate}

\section{\ZLA{AppeCH4InseNONLmisur}\textbf{Appendix:} A Set Which is not Lebesgue Measurable}
Vitali in~1905 constructed the following example of a subset $E\in(0,1/2)$ which is not Lebesgue measurable\footnote{even more: $E$ is not measurable respect to any $\ZSI$-additive measure which is translation invariant and such that the measure of a segment is its length.}. The example is in~\cite{VITALI1905NONmisur}.

We introduce the following equivalence relation in $\zzr$:
\[
x\sim y\qquad \mbox{if $x-y\in\zzq$}\,.
\]
We denote with Greek letters the equivalence classes: $\zaa=x+\zzq$ is an equivalent class. We write $\zaa_x$ if we want to stress that $\zaa$ is the equivalence class which contains $x$.

Note that  
\[
\mbox{$\zaa$ is a numerable set and $\zzr=\cup\zaa$}
\]
so that \emph{the family of the equivalence classes is not numerable.}

Let us fix any $x\in\zzr$. There exist numbers $q\in\zzq$   such that
\[
x+q\in(0,1/2)
\] 
and then $\zaa_x\cap(0,1/2)\neq\emptyset$ for every $x\in\zzr$. 
So,  every equivalence class intersects $(0,1/2)$.

\emph{For every $\zaa$ we choose one element $x_\zaa\in\zaa\cap(0,1/2)$.
The set $E$ is the set of these elements $x_\zaa$. }

The set $E$ is not numerable since $\{\zaa\}$ (the set of the equivalence classes) is not numerable.

For every $q\in\zzq$ we denote $E_q$ the translation of $E$:
\[
E_q=E+q=\{x+q\,, \ x\in E\}\,.
\]
We note:
\begin{enumerate}
\item if $q_1$ and $q_2$ are different rational  numbers then $E_{q_1}\cap E_{q_2}=\emptyset$. In fact, if $x\in E_{q_1}\cap E_{q_2}$ then
\[
x=y_1+q_1=y_2+q_2\,,\qquad y_1\in E\,,\quad y_2\in E\,.
\]
This is not possible since the equality implies $y_1\sim y_2$ while different elements of $E$ are taken from different equivalence classes;
\item we have
\begin{equation}\ZLA{eq:CH4AppeArgoVITALI}
R=\bigcup _{q\in \zzq} E_q 
\end{equation}
since every $x\in\zzr$ belongs to its equivalence class $\zaa_x$, so it is of the form $y+q$ with $y\in\zaa_x\cap(0,1/2)$.

\item 
 Equality~(\ref{eq:CH4AppeArgoVITALI}) shows that $\zzr$ is a numerable union of sets and so at least one of them either is not Lebesgue measurable or it is not a null set.
 
 Now we use translation invariance of the measure: 
 If $E\in\mathcal{L}(\zzr)$ then we have also $E_{q }\in\mathcal{L}(\zzr) $ and 
 \[
\zl(E_{q })=\zl(E)\,. 
 \]  
 \emph{It follows from~(\ref{eq:CH4AppeArgoVITALI}) that $E$, if measurable, is not a null set} since $\zzr$ is not a null set.
\end{enumerate}

\emph{The proof that $E$  is \emph{not   measurable}  consists in devising a different argument which shows that $E$,  if measurable,    must be a null set.}  The argument is as follows: We recall that $E\subseteq (0,1/2)$ so that
\[
E+1/n\in (0,2)\qquad \forall n\in\zzn
\]
and we saw already that
\[
\left (E+\dfrac{1}{n}\right )\bigcup 
\left (E+\dfrac{1}{m}\right )=\emptyset \quad\mbox{(if $n\neq m$)}\,.
\]
So
\[
\bigcup _{n\geq 1}\left (E+\dfrac{1}{n}\right )\subseteq (0,2)\quad\mbox{(the sets are pairwise disjoint)}\,.
\]
  We use monotonicity and $\ZSI$-additivity of the measure and we use again translation invariance.  We find:
\[
2\geq \sum _{n\geq 1}\zl(E+1/n)=\lim _{k\to+\ZIN}\sum_{n=1}^k \zl(E+1/n)=\lim _{k\to+\ZIN} k\zl(E)\ \implies\ \zl(E)=0\,.
\]

\emph{The contradiction shows that the set $E$ is not Lebesgue measurable.}

\begin{Remark}\ZLA{RemaFINAppeCH4InseNONLmisur}
{\rm
Note a consequence of this example: the function which is $1$ on the set $E$ and $-1$ on its complement is not quasicontinuous, but its absolute value is constant, hence it is quasicontinuous.\zdia
}
\end{Remark}
 
Vitali construction of the set $E$  uses the {\sc axiom of the choice,}\index{axiom of the choice} i.e. the fact that we can arbitrarily choose one element  from each one of infinitely many sets. It is a fact that up to now an example of a non measurable set constructed without using this axiom   is not known. 


\section{\ZLA{APPEch4NullNONborel}A Null Set Which is not Borel Measurable}
We prove the existence of Lebesgue measurable sets in $\zzr^2$ which are not Borel set.   
\begin{enumerate}
\item let $N\subseteq \zzr^2$. If for every $\ZEP>0$ there exists a sequence $R_n$ of rectangles   whose areas sum to a number less the $\ZEP$, then $N$ is a null set.

\item Borel sets of $\zzr^2$ are countable unions and intersections of open or closed sets.
\item if $B\subseteq \zzr^2$ is a Borel set and if $f$: $\zzr\mapsto \zzr^2$ is continuous then $f^{-1}(B)$ is a Borel set of $\zzr$ (see Theorem~\ref{TeoCH4ContrimmBOrelianiBIS}).
\end{enumerate}

We use these observations and we prove the existence of subsets af $\zzr^2$ which are Lebesgue but not Borel measurable.  More precisely we show a set    $N\subseteq \zzr^2$ which is not Borel measurable but such that $\zl(N)=0$.

Let  $E\subseteq \zzr$  be a set which is not Lebesgue measurable, for example the Vitali set constructed  in  Appendix~\ref{AppeCH4InseNONLmisur}. Let
\[
G=\{(x,0)\,,\ x\in E\}\subseteq \zzr^2\,.
\]
The set $G$ is a null set in $\zzr^2$, hence it is Lebesgue measurable, \emph{but it is not Borel measurable.} In fact, the function
\[
x\mapsto f(x)=(x,0)\,,\qquad \zzr^1\mapsto \zzr^2
\]
is continuous. Hence, if $G\in\mathcal{B}(\zzr^2)$ then   $f^{-1}(G)\in\mathcal{B}(\zzr)\subseteq\mathcal{L}(\zzr)$ while  $f^{-1}(G)=E\notin\mathcal{L}(\zzr)$.

\begin{Remark}
{\rm
We note the following consequences of this example:
\begin{enumerate}\ZLA{rema:AppeCap4NONcompos}
 
\item 
the set   $\{(x,0)\,,\ x\in\zzr\}\in \mathcal{B}(\zzr^2)$  and it is a null set. 
It contains the set $G\notin\mathcal{B}(\zzr^2)$
So, \emph{$\mathcal{B}(\zzr^2)$ is not a complete $\ZSI$-algebra. }
\item the function $f$ is continuous and $G$ is Lebesgue measurable. But, $f^{-1}(G)$ \emph{is not} Lebesgue measurable. The inverse image  of a Lebesgue measurable set under continuous function (and, a fortiori, under Lebesgue measurable functions) in general is not Lebesgue measurable.
  A consequence is that in general the composition   of Lebesgue measurable functions  is not Lebesgue measurable.   
    \item In a similar way it is possible to construct examples in any dimension $d\geq 2$. The previous arguments cannot be adapted when $d=1$. Nevertheless, examples exists also in dimension~$1$. We refer to~\cite[p.~56]{Benedetto76BOOK}.\zdia
\end{enumerate}

}
\end{Remark}


%

%
	 
\backmatter

\printindex
 

\end{document}